\documentclass[11pt]{amsart}

\usepackage{amssymb, amsmath}
\allowdisplaybreaks
\usepackage[all]{xy}
\usepackage[inline]{enumitem}
\usepackage{hyphenat}
\usepackage[colorlinks,linkcolor=blue,citecolor=blue,urlcolor=teal]{hyperref}
\usepackage{xcolor}
\usepackage{mathtools}
\usepackage{tikz-cd}
\usepackage[labelformat=empty]{caption}
\usepackage{cleveref} 

\usepackage[ left=3cm, right=3cm, top = 2.8cm, bottom = 2.8cm ]{geometry}

\newcommand{\eq}[2]{\begin{equation}\label{#1}#2 \end{equation}}
\newcommand{\ml}[2]{\begin{multline}\label{#1}#2 \end{multline}}
\newcommand{\mlnl}[1]{\begin{multline*}#1 \end{multline*}}

\newcommand{\xr}[1] {\xrightarrow{#1}}

\newcommand{\ra}{\rightarrow} 
\newcommand{\lra}{\longrightarrow}
\newcommand{\inj}{\hookrightarrow}
\newcommand{\surj}{\twoheadrightarrow}

\newcommand{\A}{\mathbb{A}}
\newcommand{\D}{\mathbb{D}}
\newcommand{\F}{\mathbb{F}}
\newcommand{\G}{\mathbb{G}}

\renewcommand{\P}{\mathbb{P}}
\newcommand{\N}{\mathbb{N}}
\newcommand{\Q}{\mathbb{Q}}
\newcommand{\Z}{\mathbb{Z}}

\newcommand{\sB}{\mathcal{B}}
\newcommand{\sE}{\mathcal{E}}
\newcommand{\sF}{\mathcal{F}}

\newcommand{\sH}{\mathcal{H}}

\newcommand{\sO}{\mathcal{O}}
\newcommand{\sU}{\mathcal{U}}
\newcommand{\sV}{\mathcal{V}}

\newcommand{\sX}{\mathcal{X}}

\newcommand{\sZ}{\mathcal{Z}}

\newcommand{\sD}{\mathcal{D}}
\newcommand{\sL}{\mathcal{L}}

\newcommand{\sP}{\mathcal{P}}
\newcommand{\sQ}{\mathcal{Q}}
\newcommand{\sR}{\mathcal{R}} 
\newcommand{\cE}{\mathcal{E}} 
\newcommand{\cH}{\mathcal{H}} 
\newcommand{\cO}{\mathcal{O}} 
\newcommand{\cD}{\mathcal{D}} 
\newcommand{\cR}{\mathcal{R}} 

\newcommand{\tF}{\widetilde{F}}

\newcommand{\fm}{\mathfrak{m}}

\newcommand{\bF}{\bar{F}}

\newcommand{\Hom}{\operatorname{Hom}}
\newcommand{\cHom}{\operatorname{\cH om}}
\newcommand{\cExt}{\operatorname{\cE xt}}

\newcommand{\End}{\operatorname{End}}
\newcommand{\Ker}{\operatorname{Ker}}
\newcommand{\Coker}{\operatorname{Coker}}
\newcommand{\coker}{\operatorname{Coker}}
\renewcommand{\Im}{\operatorname{Im}}

\newcommand{\id}{\operatorname{id}}

\newcommand{\Tor}{\operatorname{Tor}}

\newcommand{\uHom}{\operatorname{\underline{Hom}}}

\newcommand{\Frac}{\operatorname{Frac}}

\newcommand{\CH}{{\operatorname{CH}}}

\newcommand{\xra}{\xrightarrow}
\newcommand{\colim}{{\operatorname{colim}}}

\newcommand{\ol}{\overline}

\newcommand{\cone}{\operatorname{cone}}

\newcommand{\Tr}{\operatorname{Tr}}
\newcommand{\tr}{{\operatorname{tr}}}
\newcommand{\Nm}{\operatorname{Nm}}
\newcommand{\Pic}{\operatorname{Pic}}

\newcommand{\Br}{\operatorname{Br}}
\newcommand{\Spec}{\operatorname{Spec}}

\newcommand{\Div}{\operatorname{Div}}
\renewcommand{\div}{\operatorname{div}}

\newcommand{\Res}{\operatorname{Res}}

\newcommand{\red}{{\operatorname{red}}}

\newcommand{\dlog}{\operatorname{dlog}}
\newcommand{\AS}{\mathrm{AS}} 
\newcommand{\Ztr}{\mathbb{Z}_{\rm tr}}
\newcommand{\hra}{\hookrightarrow}
\renewcommand{\hat}{\widehat} 

\newcommand{\Zar}{{\operatorname{Zar}}}
\newcommand{\Nis}{{\operatorname{Nis}}}
\newcommand{\et}{{\operatorname{\acute{e}t}}}

\newcommand{\FRP}{\operatorname{FRP}}

\newcommand{\fillog}{\operatorname{fil}^{\rm log}}
\newcommand{\fil}{{\operatorname{fil}}}
\newcommand{\Fil}{{\operatorname{Fil}}}

\newcommand{\gr}{{\operatorname{gr}}}
\newcommand{\e}{{\epsilon}}

\newcommand{\ul}{\underline}

\newcommand{\Sm}{\operatorname{\mathbf{Sm}}}
\newcommand{\Cor}{\operatorname{\mathbf{Cor}}}
\newcommand{\PST}{{\operatorname{\mathbf{PST}}}}

\newcommand{\NST}{\operatorname{\mathbf{NST}}}

\newcommand{\uMCor}{\operatorname{\mathbf{{\underline{M}}Cor}}}

\newcommand{\bcube}{{\ol{\square}}}

\newcommand{\uomega}{{\ul{\omega}}}

\newcommand{\uMPST}{\operatorname{\mathbf{{\underline{M}}PST}}}

\newcommand{\CI}{{\operatorname{\mathbf{CI}}}}

\newcommand{\RSC}{{\operatorname{\mathbf{RSC}}}}

\newcommand{\idrw}{{\operatorname{inv-dRW}}}
\newcommand{\ddrw}{{\operatorname{dir-dRW}}}

\newcommand{\uMNST}{\operatorname{\mathbf{{\underline{M}}NST}}}

\newcommand{\ulMDM}{\operatorname{\mathbf{\underline{M}DM}}}
\newcommand{\uMDM}{\operatorname{\mathbf{\underline{M}DM}}}

\newcommand{\pdd}{$p$-divisibility decomposition}

\theoremstyle{plain}
\newtheorem{introthm}{Theorem}
\newtheorem{introcor}{Corollary}
\newtheorem{prop}{Proposition}[section]
\newtheorem{lem}[prop]{Lemma}
\newtheorem{cor}[prop]{Corollary}
\newtheorem{thm}[prop]{Theorem}
\newtheorem{claim}{Claim}[prop]

\theoremstyle{definition}
\newtheorem{defn}[prop]{Definition}
\newtheorem{notation}[prop]{Notation}
\newtheorem{example}[prop]{Example}

\newtheorem{para}[prop]{}

\newtheorem{rmk}[prop]{Remark}
\newtheorem{remark}[prop]{Remark}

\newcommand{\beq}{\begin{equation}}
\newcommand{\eeq}{\end{equation}}

\numberwithin{equation}{prop}

\title[Duality for Hodge-Witt cohomology with modulus]{Duality for Hodge-Witt cohomology with modulus}

\author{Fei Ren}
\address{Bergische Universit\"at Wuppertal\\ Gau\ss str. 20, D-42119 Wuppertal, Germany}
\email{renfei@uni-wuppertal.de}

\author{Kay R\"ulling}
\address{Bergische Universit\"at Wuppertal\\ Gau\ss str. 20, D-42119 Wuppertal, Germany}
\email{ruelling@uni-wuppertal.de}

\thanks{ }
\begin{document}
\begin{abstract}
Given an effective Cartier divisor $D$ with simple normal crossing support on
a  smooth and proper scheme $X$ over a perfect field of positive characteristic $p$,
there is a natural notion of de Rham-Witt sheaves on $X$ with zeros along $D$.
We show that these sheaves correspond under Grothendieck duality for coherent sheaves to
de Rham-Witt sheaves on $X$ with modulus $(X,D)$, 
as defined in the theory of cube invariant modulus sheaves with transfers developed by Kahn-Miyazaki-Saito-Yamazaki.
From this we deduce refined versions of Ekedahl -  and Poincaré duality for crystalline cohomology generalizing results
of Mokrane and Nakkajima for reduced $D$, and a modulus version of Milne-Kato duality 
for étale motivic cohomology with $p$-primary torsion coefficients, which refines  a result of Jannsen-Saito-Zhao.
We furthermore get new integral models for rigid cohomology  with compact supports on the complement of $D$ 
and a modulus version of Milne's perfect Brauer group pairing for  smooth projective surfaces over  finite fields.
\end{abstract}
\maketitle

\tableofcontents

\section{Introduction}

Let $k$ be a perfect field of positive characteristic $p>0$. 
There is a wealth of duality results for  smooth and proper $k$-schemes with $p$-torsion - or $p$-adic coefficients, such as
Serre duality for differential forms, Poincar\'e duality for crystalline cohomology \cite{Berthelot-Cris}, or
Milne  duality \cite{Milne76}, see also \cite{Kato-dualityI}, \cite{Kato-dualityII},  which by \cite{GL} 
can be interpreted as a duality for \'etale motivic cohomology with mod $p^n$ coefficients.
For surfaces over a finite field the latter induces a pairing for the Brauer group \cite{Milne-ConjArtin}.
These parings can all be understood via Grothendieck-Ekedahl duality \cite{Ekedahl}, \cite{Gros} which is a duality statement for
the de Rham-Witt complex of Bloch-Deligne-Illusie \cite{Bloch-DRW}, \cite{IlDRW}.

A generalization of these duality statements to open subschemes (say of a smooth and proper $k$-scheme)
requires the use  of cohomology with compact support
and relies on considering pairings between certain pro-groups (on the compact support side) and certain ind-groups.
For example there is the Deligne-Hartshorne duality for compactly supported cohomology of coherent sheaves \cite{Ha72},
Berthelot's duality for rigid cohomology \cite{Berthelot-PD-KunnethRig}, and the generalization of 
Milne  duality by Jannsen-Saito-Zhao \cite{JSZ} and Gupta-Krishna \cite{Gupta-Krishna}, which generalizes 
the $p$-primary torsion part  of geometric global class field theory.

As the (compact supported) cohomology groups of an open in the above situations are very often not 
finitely generated over the base ring at hand, e.g., $k$, $W_n(k)$, $W(k)$, or $\Z/p^n$, 
it is desirable to have more precise pairings which hold before taking the colimit or the limit 
to  go all the way to compactly supported cohomology. 
This requires to have a sturdy notion of poles and pole orders along the complement of an open immersion $U\inj X$
between smooth $k$-schemes and also a dual notion of zeros and vanishing order along the same closed subset. 
Candidates appear at various places in the literature: in \cite{Brylinski}, \cite{Kato-Swan}, \cite{Matsuda},
Brylinski, Kato, and Matsuda introduce  pole order filtrations on Witt vectors which under the Artin-Schreier-Witt sequence
measure the ramification of $\Z/p^n$-Galois coverings of $U$ along $X\setminus U$. In \cite{KeS-Lefschetz} these filtrations
are used to prove a Lefschetz-type result for the abelian fundamental group with modulus. In \cite{Tanaka} Tanaka introduces 
Witt divisorial sheaves to obtain Witt versions of classical vanishing theorems, similar sheaves appear also in 
\cite{JSZ}, \cite{Gupta-Krishna}, \cite{KSY}.  Hodge-Witt sheaves with vanishing order along
a (non-reduced) closed subscheme also play an important role  in the work \cite{Morrow}.

In all these cases we would like to have a precise duality pairing  which matches Witt differential forms with 
certain poles along $X\setminus U$ with corresponding zeros along the same closed subset. 
In this paper we are going to construct such pairings in the case $X\setminus U$ is the support of a simple normal crossing divisor.

Let $X$ be a smooth $k$-scheme (not necessarily proper) of pure dimension $N$ and let $D$ be an effective Cartier divisor 
such that the underlying reduced divisor  $D_{\red}$ has simple normal crossings and denote the complement by $U=X\setminus D$. 
The question is now how to define  appropriate notions of poles bounded by $D$ and vanishing along $D$
for the Hodge-Witt sheaves $W_n\Omega^q_X$. One natural way is to consider  
the multiplicative lift $\Pic(X)\to \Pic(W_n X)$.
Then $D$ gives rise to an invertible $W_n\sO_X$-module $W_n\sO(D)$ and we can consider 
$W_n\Omega^q_X\otimes_{W_n\sO_X} W_n\sO(\pm D)$ for poles ($+$) and  zeros ($-$).
This approach is taken for example in \cite{Tanaka} and at least for the pole side also in \cite{JSZ}, \cite{Gupta-Krishna}, 
\cite[Appendix]{KSY}. The advantage of this choice is that by Ekedahl duality multiplication induces immediately an isomorphism 
\[W_n\Omega^q_X\otimes_{W_n\sO_X} W_n\sO(D)\cong R\cHom_{W_n\sO_X}(W_n\Omega^{N-q}_X\otimes_{W_n\sO_X} W_n\sO(-D), W_n\Omega^N_X),\]
and this clearly works for any Cartier divisor $D$. The drawback of this choice is that the differential $d$, 
the Verschiebung $V$, and the Frobenius $F$ which are defined on the de Rham-Witt complex do not extend to
endomorphisms of the graded pro-object $W_\bullet\Omega^*_X\otimes_{W_\bullet\sO_X} W_\bullet \sO(D)$. 
Hence we cannot extend this to a pairing for crystalline cohomology or to a Milne type pairing.
Furthermore, this makes it hard to understand the duality in the limit over $n$, see, e.g, \cite{Lemcke}, where the 
product $\prod_r  W_\bullet(\sO(p^r D))$ is considered to accomplish this for $q=0$.  
Also  considering $D$ as a (non-reduced) scheme the complex $W_n\Omega^*_D$  is defined and 
studying it requires to consider $\Ker(W_n\Omega^*_X\to W_n\Omega^*_D)$, which is not the same as 
$W_\bullet\Omega^*_X\otimes_{W_\bullet\sO_X} W_\bullet\sO(-D)$.\footnote{Though the pro-systems over $n$ are isomorphic.}
In case $D$ is {\em reduced} the log de Rham-Witt complex \cite{HK}, \cite{Matsuue}
has all the desired properties. Furthermore, \cite{Mokrane} constructed a de Rham-Witt complex with zeros along the reduced
divisor $D$, which by \cite{Hyodo} and \cite{Nakkajima} is dual to the log de Rham-Witt. But at least for the pole side 
it was not clear how to define a notion of poles bounded by a non-reduced divisor. 

The notion of pole orders we choose in this paper relies on the theory of motives with modulus developed 
by Kahn-Miyazaki-Saito-Yamazaki in \cite{KMSYI}, \cite{KMSYII}, \cite{KMSYIII}. 
We will explain the relation more precisely a little later.
For now we make the following a bit more direct definition which comes from the work \cite{RS-AS} and is
inspired by  \cite{AbbesSaito11}:
Let $F$ be a Nisnevich sheaf of abelian groups on all finite type and separated $k$-schemes. Let $X$ and $D$ be as above.
We denote by $F_X$  the Nisnevich sheaf on the small site $X_{\Nis}$.

\medskip

{\bf Poles along $D$.}
Denote by ${\rm Bl}_D(X\times X)$ the blow-up in $D$ diagonally embedded into $X\times X$ 
and denote by $P^{(D)}_X$ the complement of the union of the strict transforms of $X\times D$ and $D\times X$ in ${\rm Bl}_{D}(X\times X)$.
The open embedding $U\times U\inj X\times X$ extends to an embedding $U\times U\inj P^{(D)}_X$.
We set
\[F(X,D):= \text{Equalizer}\left(\xymatrix{F(U)\ar@<1ex>[r]^-{p_1^*}\ar@<-1ex>[r]_-{p_2^*} & F(U\times U)/F(P^{(D)}_X)}\right)\]
and denote by $F_{(X,D)}$\footnote{In the body of the text these sheaves will be denoted by $F^{\rm AS}_{(X,D)}$, but we prefer to have
this simpler notation in the introduction.} the Nisnevich sheaf $(V\to X)\mapsto F(V, D_{|V})$, see \ref{para:AS} for details.

\medskip

{\bf Zeros along $D$.}
Let $D=\sum_i D_i$ with $D_{i,\red}$ smooth.
We set
\[F_{(X,-D)}:= \Ker \left(F_X\to \bigoplus_i F_{D_i}\right).\]

Note that both definitions are functorial in $F$.
One of the main results of this paper is now:

\begin{introthm}[Theorem \ref{thm:duality-mod}]\label{introthm1}
Multiplication induces isomorphisms, for all $q\ge 1$ and $n$,
\[W_n\Omega^{q}_{(X,D)}\xr{\simeq} R\cHom_{W_n\sO_X}(W_n\Omega^{N-q}_{(X,-D)}, W_n\Omega^N_X)\]
and 
\[W_n\Omega^{N-q}_{(X,-D)}\xr{\simeq}R\cHom_{W_n\sO_X}(W_n\Omega^{q}_{(X,D)}, W_n\Omega^N_X).\]
\end{introthm}

From Ekedahl's trace isomorphism
\eq{eq:intro1}{\pi_n^! W_n(k)\cong W_n\Omega^N[N],}
where $\pi_n: W_n(X)\to\Spec W_n(k)$ is the structure map, together with Grothendieck duality we immediately get:
\begin{introcor}\label{introcor1}
Assume additionally that $X$ is proper over $k$. Then there is a canoncial isomorphism of finite $W_n(k)$-modules 
\[H^{N-i}(X, W_n\Omega^{N-q}_{(X,-D)})\cong H^{i}(X, W_n\Omega^{q}_{(X,D)})^\vee,\]
where $(-)^\vee$ denotes the dual $W_n(k)$-module.
\end{introcor}

We make some comments on the above result:
\begin{itemize}
    \item It is not obvious from the definition given above that $W_n\Omega^q_{(X,D)}$ is a
    $W_n\sO_X$-module. And in fact this is not true for $q=0$, which is shown in \cite{RS-AS}.
    Using a certain conductor as defined in \cite{RS21}, we can also
    define the correct $W_n\Omega^0_{(X,D)}$, so that the isomorphisms in Theorem \ref{introthm1} hold for $q=0$ as well.
    \item We have $Rj_*W_n\Omega^q_{U}=\colim_r W_n\Omega^q_{(X,rD)}$, where $j:U\inj X$ is the open immersion,
    hence taking the limit over $\{rD\}_r$ in
    Corollary \ref{introcor1} yields 
    \[H^{N-i}_c(U, W_n\Omega^{N-q}_U)\cong H^{i}(U, W_n\Omega^{q}_U)^\vee, \]
    where the left hand side is the compactly supported cohomology from \cite{Ha72}.
    \item By Ekedahl's  isomorphism \eqref{eq:intro1} the $W_n\sO_X$-module $W_n\Omega^N_X$ 
    is a dualizing complex in the sense of Grothendieck. 
    Hence the two isomorphisms in Theorem \ref{introthm1} imply each other.
    \item For $D=\emptyset$, Theorem \ref{introthm1} is due to Ekedahl. 
    We have $W_n\Omega^q_{(X, D_{\red})}=W_n\Omega^q_X(\log D_{\red})$ and hence for $D$ reduced Theorem \ref{introthm1} recovers 
     the duality from   \cite[Theorem 5.3(1)]{Nakkajima}, \cite[(3.3.1)]{Hyodo}.
    \item Part of the statement of Theorem \ref{introthm1} is also the vanishing
    \[\cExt^i(W_n\Omega^q_{(X,\pm D)}, W_n\Omega^N_X)=0\quad \text{for } i \ge 1.\]
    Note that this is not automatic as $W_n\Omega^q_{(X,\pm D)}$ is not a locally free $W_n\sO_X$-module. 
    The idea to prove this vanishing is similar as in \cite{Ekedahl} but the actual computations are  more involved,
    see Theorem \ref{introthm7} below.
    \item It is not hard to guess from the existing literature, e.g.,
    \cite{JSZ}, \cite{Morrow}, \cite{Gupta-Krishna},\cite{RS-AS}, that there is an isomorphism of pro-systems
    over $\{rD\}_r$ as in  the second isomorphism of  Theorem \ref{introthm1}, but that these two abstract notions of
    poles and zeros introduced above for general sheaves correspond precisely under multiplication came as a surprise
    to the authors and indicates that there might be some more abstract yet to discover duality in the background. 
    We want to stress however, that our proof essentially proceeds by computing the pole and the zero side
    completely and then show by hand that they match.    
\end{itemize}

Before we  give more applications we consider the pole part in greater detail.
Recall the notion of cube invariant sheaves with transfers $\CI_{\Nis}^{\tau,{\rm sp}}$ 
from \cite{KSY-RSCII}, \cite{Saito-Purity}, see also \ref{para:RSC} and \ref{para:MNST-gamma} 
for some more details and references. In particular we have a functor 
\[\ul{\omega}^{\CI}:\NST\to \CI^{\tau,{\rm sp}}_{\Nis}\]
from the category of Nisnevich sheaves with transfers to (certain) cube invariant sheaves with transfers 
which in case $X$ is  proper is given by 
\[\uomega^{\CI}(F)(X,D)=\Hom_{\NST}(h_0(X,D), F).\]
Here $h_0(X,D)$ is a certain Nisnevich sheaf with transfers which on a function field $K/k$ is equal to 
$\CH_0(X_K,D_K)$ the Chow group with modulus as introduced in \cite{KeS}, 
where the index $K$ indicates the base change over $K$. 
Moreover, if for each  henselian discrete valuation field $L$ of geometric type over $k$,  
we are given an increasing filtration $\{\Fil_r F(L)\}_{r\ge 0}$ which defines a conductor 
$c=\{c_L: F(L)\to \N\}_L$ in the sense of \cite{RS21}, then there is a cube invariant sheaf with transfers
$F_c$ which for $X$ proper is given by
\[F_c(X,D)=\{a\in F(U)\mid \rho^*(a)\in \Fil_{v_L(\rho^*D)} F(L)\quad \text{for all } \rho: \Spec \sO_L\to X\},\]
see \ref{para:cond} for details and references. By \cite[Theorem 4.15]{RS21} and \cite[Theorem 2.6]{RS-AS}
we have 
\eq{eq:intro2}{F_{c, (X,D)}\subset \uomega^{\CI}F_{(X,D)}\subset F_{(X,D)},}
where the right hand side is {\em poles along D} as defined  above. 
Furthermore by \cite[Theorem 2.10]{RS-AS} the second inclusion is an equality in case $X$ is projective and smooth. 

For $L$ as above we define for $r,q\ge 0$
\[\fil_r W_n\Omega^q_L:= \fillog_{r-1}W_n(L)\cdot \dlog(K^M_q(L))
+V^{n-m}\left(\fillog_rW_n(L)\cdot \dlog K^M_q(\sO_L)\right),\]
where $\fillog_{r}W_n(L)$ is the Brylinski-Kato filtration, 
$m=\min\{n, v_p(r)\}$, and $K^M_q$ denotes the $q$th Milnor $K$-theory, and 
\[\Fil^p_rW_n\Omega^q_L= \sum_{s\ge 0} p^s \left(\fil_{p^sr}W_n\Omega^q_L+ d(\fil_{p^sr}W_n\Omega^q_L)\right),\]
see section \ref{sec:p-sat-fil} for details.
The second main result of this paper, which is also essential for the proof of Theorem \ref{introthm1}, is:
\begin{introthm}[Theorem \ref{thm:fil-cond} and Theorem \ref{thm:HW-modulus}]\label{introthm2}
The filtration $\{\Fil^p_rW_n\Omega^q_L\}_{r,L}$ defines a conductor $c$. Moreover, for
$q\ge 1$ and $X$ smooth (not necessarily proper)
\[W_n\Omega^q_{c, (X,D)}= \uomega^{\CI}W_n\Omega^q_{(X,D)}=  W_n\Omega^q_{(X,D)}.\]
\end{introthm}
For $q=0$, the correct definition for the pole side from the point of view of Theorem \ref{introthm1} is 
given by $W_n\Omega^0_{c,(X,D)}$. 
In view of the first part of this theorem the proof of the second part is essentially by \eqref{eq:intro2} 
reduced to show that the Nisnevich stalk
$(W_n\Omega^q_{(X,D)})_\eta^h$ in a generic point $\eta$ of $D$  is contained in $\Fil^p_rW_n\Omega^q_L$,
where $L=\Frac(\sO_{X,\eta}^h)$. The proof of this result takes the sections \ref{sec:p-sat-fil} - \ref{sec:HWM}.
It relies crucially on the construction of some ad hoc characteristic form (see Lemma \ref{lem:char}).
We remark that for $n=1$ the above result is proven in \cite[Theorem 6.6]{RS-AS} and our proof makes use of this result.
Note that the $p$-saturation is not visible in the case $n=1$.

Though the definition of the filtration $\Fil^p_rW_n\Omega^q_L$ 
looks a bit unmotivated it has several universal interpretations, e.g.,
Theorem \ref{introthm2} states ($q\ge 1$)
\[\Fil^p_rW_n\Omega^q_L=\uomega^{\CI}W_n\Omega^q(\sO_L, \fm_L^{-r}),\]
which also implies that it defines the motivic conductor in the sense of \cite{RS21}.
Furthermore Theorem \ref{introthm1} gives the isomorphism (which also holds for $q=0$)
\eq{eq:intro3}{\Fil^p_rW_n\Omega^q_L \xr{\simeq}\Hom_{W_n\sO_L}(W_n\Omega^{N-q}_{(\sO_L,\fm^r_L)}, W_n\Omega^N_L), }
where $W_n\Omega^{N-q}(\sO_L,\fm_L^r)=\Ker(W_n\Omega^{N-q}_{\sO_L}\to W_n\Omega^{N-q}_{\sO_L/\fm_L^r})$ and 
$N={\rm trdeg}(L/k)$.  If $f_n: \Spec W_n(\sO_L/\fm_L^r)\to \Spec W_n(k)$ is the structure map,
then we get an isomorphism 
\[f_n^! W_n(k)\cong \frac{\Fil^p_rW_n\Omega^N_L}{W_n\Omega^N_{\sO_L}}[N-1],\]
see Corollary \ref{cor:twisted-inverse-image-Fil}.
Also note that this filtration is related to Kato's ramification filtration 
for \'etale motivic cohomology with $p$-torsion coefficients from \cite{Kato-Swan}, see Lemma \ref{lem:relKato} 
for a precise statement. However we remark,  that there are several ways to lift $W_n\Omega^q$ to a cube invariant sheaf
with transfers (corresponding to different conductors), and that this particular lift we are using here
is the maximal one and turns out to be suitable for duality, but has the disadvantage, 
as already remarked in \cite[6.9]{RS-AS} (there for $n=1$), that the functor 
\[(X,D)\mapsto R\Gamma(X, \uomega^{\CI}W_n\Omega^q_{(X,D)})\]
does not factor through the triangulated category of motives with modulus $\ulMDM^{\rm eff}$ defined in \cite{KMSYIII}.
A sufficient criterion for a cube invariant sheaf to have this property is given by Koizumi in \cite{Koizumi}.
There it is also shown that the filtration $\{\fillog_r W_n(L)\}_{r}$ gives rise to a structure of cube invariant
Nisnevich sheaf with transfers on $W_n\sO$ which defines a realization of Witt vector cohomology from $\uMDM^{\rm eff}$.

In the following,  we write $W_n\Omega^q_{(X,D)}$ for $W_n\Omega^q_{c, (X,D)}$ and allow also $q=0$. 
In particular  we get that $W_n\Omega^q_{(X,D)}$ is a $W_n\sO_X$-module,
that $\colim_r W_n\Omega^q_{(X,rD)}= j_*W_n\Omega^q_U$, and that the equality
\[W_n\Omega^q_{(X,D_{\red})}=W_n\Omega^q_X(\log D_{\red})\]
holds. By Theorem \ref{introthm2} this sheaf  has transfers, 
which gives a new proof of a result by Merici \cite{Merici}, see Remark \ref{rmk:HW-modulus}. 
We also recall the formula for $n=1$, which was proven in \cite{RS-AS},
\[W_1\Omega^q_{(X,D)}=\Omega^q_{(X,D)}= \Omega^q_X(\log D_0)\otimes_{\sO_X}\sO_X(D_0- D_{0,\red}+pD_1),\]
where we write $D=D_0+pD_1$ such that $D_0$ and $D_1$ have no irreducible components in common 
and  $D_0$ has none of its multiplicities  divisible by $p$.

Using the explicit description of the filtration we  can check that multiplication induces
the map in \eqref{eq:intro3} which also induces the natural maps in Theorem \ref{introcor1}. 
To prove that the pairing is perfect and that the higher ext-groups vanish requires some further analysis of the structure
of $W_n\Omega^q_{(X,\pm D)}$, see Theorem \ref{introthm7} below.  We first give some applications.

\medskip

In Theorem \ref{lem:HWM-gamma} we use Theorem \ref{introthm7} plus a result from \cite{BRS} to show
\[\uHom_{\uMNST}(K^M_r, \uomega^{\CI}W_n\Omega^q)_{(X,D)}\cong W_n\Omega^{q-r}_{(X,D)}.\]
For $D=\emptyset$, this was proven in  \cite{BRS}. 
Hence we can apply the results of \cite{BRS} to obtain
a projective bundle formula and a blow-up formula for $W_n\Omega^q_{(X,D)}$, see \ref{para:pbf-buf-gt},
and we can apply Theorem  \ref{introthm1} to get these formulas as well for $D$ replaced by $-D$, 
see Corollary \ref{cor:pbf-buf-zeros}. 
These generalize results by Gros \cite{Gros} for the case $D=\emptyset$. 
We furthermore get a Gysin triangle (for poles):
Let $Z\inj X$ be a closed immersion of pure codimension $r$ between smooth $k$-schemes, which intersects 
$D$ transversally ($D$ may be non-reduced). Denote  by $\rho:\tilde{X}\to X$ the blow-up of $X$ in $Z$ and by 
$E$ the exceptional divisor. Then there is a canonical distinguished triangle in $D(W\sO_X)$
\[i_*W_n\Omega^{q-r}_{(Z, D_{|Z})}[-r]\xr{g}  
W_n\Omega^q_{(X,D)}\xr{\rho^*} R\rho_*W_n\Omega^q_{(\tilde{X}, \rho^*D+E)}\xr{\partial} 
i_*W_n\Omega^{q-r}_{(Z, D_{|Z})}[-r+1]. \]
As an application we obtain the following Lefschetz-type statement for the cohomology of the top Hodge-Witt forms:
\begin{introthm}[Theorem \ref{thm:top-Lef}]\label{introthm3}
Assume  that $X$ is additionally projective and let $H\subset X$ 
be a smooth hypersurface section which intersects $D$ transversally and satisfies
\[H^j(X, \Omega^N_{(X,D)}\otimes_{\sO_X}\sO_X(H))=0, \quad \text{for all } j\ge 1.\]
Then  the Gysin map
\[g:H^{j-1}(H, W_n\Omega^{N-1}_{(H, D_{|H})})\lra H^j(X, W_n\Omega^N_{(X,D)})\]
is an isomorphism for $j\ge 2$ and is surjective for $j=1$.
\end{introthm}
We will not spell it out in the following, but the reader should keep in mind that the projective bundle formula,
the blow-up formula and the Gysin triangle are compatible with Frobenius, Verschiebung, differential, etc. 
and  that these formulas continue to hold in a suitably modified sense in the crystalline cohomology  and 
\'etale motivic cohomology with poles and zeros which we will discuss next.

\medskip

The two isomorphisms in Theorem \ref{introthm1} behave quite differently when taking 
the limit over $n$. Indeed,
by Proposition \ref{prop:modulus-limit},
\[R\lim_n W_n\Omega^q_{(X, D)}=W\Omega^q_X(\log D_{\red}).\]
Hence taking limits over $n$, the first isomorphism in Theorem \ref{introthm1} together with Ekedahl's duality formalism
yields in case $X$ is proper the Poincar\'e duality, see \ref{para:conseq-coh}\ref{para:conseq-coh2},
\eq{eq:intro4}{R\Gamma(X, W\Omega^\bullet_{(X, -D_{\red})})\cong R\Hom_W(R\Gamma_{\text{log-crys}}((X,D_{\red})/W), W)[-2N],}
where we use that the de Rham-Witt complex with log poles along $D_{\red}$ computes log-crystalline cohomology
of the smooth log-scheme $(X,D_{\red})$. On finite level this statement is also proven in 
\cite[Proposition (3.3)]{Hyodo} and \cite[Theorem 5.3]{Nakkajima}. 
From the above isomorphism (or rather a version before taking the total complex)
and weight filtration arguments \cite{Mokrane}, \cite{Nakkajima}
we deduce that $R\Gamma(X, W\Omega^*_{(X, \pm D_{\red})})$
\footnote{If $C^*$ is a dga of sheaves, we use $R\Gamma(X,C^*)$ to denote the dga in complexes, and we use $R\Gamma(X,C^\bullet)$ to denote the total complex of the double complex $R\Gamma(X,C^*)$.
} 
are coherent complexes of modules
over the Cartier-Dieudonn\'e-Raynaud ring $\sR$, in the sense of  Illusie-Raynaud, Ekedahl. 
We do not know if $R\Gamma(X, W\Omega^*_{(X, -D)})$ is a coherent $\sR$-complex for $D$ non-reduced, in fact 
we do not know if it is a complete $\sR$-complex, see Remark \ref{rmk:not-coh}.
As the natural map 
\[R\Gamma(X, W\Omega^*_{(X, -D)})\to R\Gamma(X, W\Omega^*_{(X, -D_{\red})})\]
is an isomorphism up to 
bounded $p$-primary torsion (see \Cref{para:conseq-coh}\ref{para:conseq-coh3})
it follows from Shiho's comparison of log-crystalline cohomology with rigid cohomology, 
Berthelot's Poincar\'e duality for rigid cohomology, and the isomorphism
\eqref{eq:intro4} that $R\Gamma(X, W\Omega^\bullet_{(X, -D)})$ is an integral model 
for compactly supported rigid cohomology of $U$, for any $D$ with $U=X\setminus D_{\red}$.
There are now several corollaries one can draw from this concerning the degeneration of the slope spectral sequence and 
vanishing results and such, but which all work up to bounded torsion and hence only require to work with $D_\red$
in which case they can be found at least implicitly in the literature, see \ref{para:conseq-coh}
for more details. 

However the multiplicities of $D$ remain visible when we take the limit in the second isomorphism of 
Theorem \ref{introthm1}. To state the result set 
\[W_\infty\Omega^*_{(X,D)}:=\colim_{\ul{p}} W_n\Omega^q_{(X,D)}\quad \text{and}\quad
W\Omega^*_{(X,-D)}:=\Ker \left(W\Omega^*_X\to \bigoplus_i W\Omega^*_{D_i}\right),\]
where the colimit is over the map $\ul{p}$ which is "lift and multiply by $p$" and $D=\sum_i D_i$ with $D_{i,\red}$ smooth.
We define {\em co-crystalline cohomology with modulus $(X,D)$} by
\[R\Gamma_{\rm crys}((X,D)/W_\infty):= R\Gamma(X, W_\infty \Omega^\bullet_{(X,D)}),\]
which has an  operator $\sV$ on it induced by $p^{N-q}V$ in degree $q$. We furthermore define
{\em crystalline cohomology with zeros along $D$} by
\[R\Gamma_{\rm crys}((X,-D)/W):= R\Gamma(X, W\Omega^\bullet_{(X,-D)}),\]
which has a Frobenius $F$ on it induced by the absolute Frobenius on $X$.
Theorem \ref{introthm1} together with Ekedahl's duality formalism yields the following version of Poincar\'e duality.
\begin{introthm}[Corollary \ref{cor:limit-crys}]\label{introthm4}
Assume $X$ is proper. There is a canonical isomorphism in $D^b(W[F])$
\[R\Gamma_{\rm crys}((X,-D)/W)\cong R\Hom_W(R\Gamma_{\rm crys}((X,D)/W_\infty), K/W)[-2N],\]
where the action of $F$ on the right hand side is induced by $\sV^\vee$.
\end{introthm}
We  have  
$R\Gamma_{\rm crys}((X,D_{\red})/W_n)= R\Gamma(X, W\Omega^\bullet_{(X,D_{\red})})\otimes^L_W W_n$ 
(cf. Corollary \ref{thm:Dred-coh}) and from this one can deduce that the isomorphism 
in Theorem \ref{introthm4} for $D_\red$ is induced by \eqref{eq:intro4}. 
This is not the case  if $D$ is not reduced, as the left hand side of the isomorphism in Theorem \ref{introthm4}
depends on the multiplicities of $D$. See \ref{para:conseq-coh} for more comments and consequences 
such as vanishing statements. Finally we remark that it is an intriguing problem  to 
investigate the relation of $R\Gamma_{\rm crys}((X,-D)/W)$ or $R\Gamma_{\rm crys}((X,D)/W_\infty)$ with 
the edged crystalline cohomology from \cite{DAddezio}.

\medskip

Next we explain a modulus version of Milne duality which we get from Theorem \ref{introthm1}.
Consider the following two complexes of abelian sheaves on $X_{\et}$
\[\Z/p^n(q)_{(X, D)}:= \left(\uomega^{\CI}(FW_{n+1}\Omega^q)\xr{1-C} \uomega^{\CI}W_n\Omega^q\right)_{(X,D)}[-q]\]
and 
\[\Z/p^n(q)_{(X,-D)}:=\left(W_n\Omega^q_{(X,-D)}\xr{C^{-1}-1} 
\frac{W_n\Omega^q_{(X,-D)}}{dV^{n-1}\Omega^{q-1}_X\cap W_n\Omega^q_{(X,-D)}}\right)[-q].\]
Both complexes sit in degree $[q, q+1]$ and are for $D=\emptyset$ quasi-isomorphic 
to the \'etale motivic complex  with $\Z/p^n$-coefficients by \cite{GL}.
Complexes similar to $\Z/p^n(q)_{(X,-D)}$  are for example considered in \cite{JSZ}, \cite{Morrow}, and \cite{Gupta-Krishna}.
A complex similar to $\Z/p^n(0)_{(X,D)}$ is considered in \cite{KeS-Lefschetz}.
We remark that $\Z/p^n(q)_{(X,-D)}$ is 
quasi-isomorphic to a sheaf concentrated in degree $q$, which is denoted by $W_n\Omega^q_{(X,-D),\log}$.
In the {\em pole-case} this is also true if $D$ is reduced but not in general. 
The local sections of the \'etale sheaf 
$\sH^{q+1}(\Z/p^n(q)_{(X,D)})$ only vanish after a finite cover which ramifies along $D$, 
see Lemma \ref{lem:Zpnq-} and Remark \ref{rmk:Zpnq+}. 

In order to obtain a version of Kato's generalization of Milne's  duality 
we denote by  $\Z/p^n(q)^{\FRP}_{(X,\pm D)}$ the extension to the flat relatively perfect site $X_{\FRP}$.

\begin{introthm}[Theorem \ref{thm:Zpnqfrp-duality}]\label{introthm5}
There are isomorphisms 
$$\Z/p^n(q)^{\FRP}_{(X,D)}[q]
\xr{\simeq}
\D_{n,X}\left(\Z/p^n(N-q)_{(X,-D)}^{\FRP}[N-q]\right),$$
$$\Z/p^n(N-q)_{(X,-D)}^{\FRP}[N-q]
\xr{\simeq}
\D_{n,X}\left(\Z/p^n(q)^{\FRP}_{(X,D)}[q]\right),$$
where $\D_{X,n}(-)$ is Kato's dualizing functor, see \ref{para:frp}.
\end{introthm}

\begin{introcor}[Corollary \ref{cor:Zpnq-duality-finite}]\label{introcor2}
Assume additionally that $k$ is a finite field and $X$ is proper over $k$.
Then there is an isomorphism in $D(\Z/p^n)$
\[R\Gamma(X_{\et}, \Z/p^n(q)_{(X,D)})
\cong R\Hom_{\Z/p^n}\left(R\Gamma(X_{\et}, \Z/p^n(N-q)_{(X,-D)}), \Z/p^n\right)[-2N-1].\]
In particular we obtain isomorphisms of finite groups for all $i$
\[H^{i+q}(X_{\et}, \Z/p^n(q)_{(X,D)})\cong 
\Hom_{\Z/p^n}(H^{2N-i-q+1}(X_{\et}, \Z/p^n(N-q)_{(X,-D)}), \Z/p^n).\]
\end{introcor}
Taking the limit over $\{rD\}_r$ in the second isomorphism above yields \cite[Theorem 2]{JSZ},
see Remark \ref{rmk:fin-duality}. 
Denote by $\pi_1^{\rm}(X,D)^p$ the \'etale fundamental group which classifies abelian $p$-covers of $U$ with 
ramification bounded by $D$ from \cite{KeS}. Then the case $q=0$ and $i=1$ in the above corollary yields 
(assuming $k$ finite and $X$ proper)
\[\pi_1^{\rm ab}(X,D)^p=H^N(X_{\et}, W\Omega^N_{(X,-D),\log}).\]

As last application we mention a refinement of Milne's pairing for the Brauer group 
of a smooth projective  surface over a finite field to the ramified situation.
The Brauer group of $X$ with ramification bounded by $D$, is defined by
\[\Br(X,D):= H^0(U, R^2\e_* \Q/\Z(1)'_{U})\oplus H^0(X, R^2\e_*\Q_p/\Z_p(1)_{(X,D)}),\]
where $\e:X_\et\to X_{\Nis}$ is the change of sites map, and 
\[\Q/\Z(1)'_{U}= \colim_{n'}\, \mu_{n',U}, \quad  
\Q_p/\Z_p(1)_{(X,D)}=\colim_{\ul{p}}\, \Z/p^n(1)_{(X,D)}\]
with the colimit on the left over all $n'$ which are prime to $p$.
The Brauer group with zeros along $D$ is defined by
\[\Br(X,-D):= H^2(X_\et, \sO_{(X,-D)}^\times),\]
where $\sO_{(X,-D)}^\times=\Ker(\sO_X^\times\to \sO_D^\times)$. 
There are natural inclusions $\Br(X)\subset \Br(X,D)\subset \Br(U)$ and 
$\colim_r \Br(X,rD)=\Br(U)$. Furthermore there is  an exact sequence
\[\Pic(X)\to \Pic(D)\to \Br(X,-D)\to \Br(X)\to \Br(D),\]
where  $\Br(Z)=H^2(Z_{\et},\G_m)$.
For $D=\emptyset$ the following duality statement is the $p$-part of Milne's duality  \cite[Theorem 2.4]{Milne-ConjArtin}.

\begin{introthm}[Theorem \ref{thm:duality-Br-surf}]\label{introthm6}
Assume $k$  is a finite field and  $X$ is a smooth proper surface. 
Then there is a canonical isomorphism of profinite groups
\[ \frac{\Br(X,-D)[p^\infty]}{(\Br(X, -D)[p^\infty])_{\rm div}}\stackrel{\simeq}{\lra}
\Hom\left(\frac{\Br(X,D)[p^\infty]}{(\Br(X,D)[p^\infty])_{\rm div}}, \Q/\Z\right),\]
where the index ``div" refers  to the divisible part and $M[p^\infty]$ is the $p$-primary torsion in $M$.
\end{introthm}

For example, if $A$ is an abelian surface over a finite field, 
it follows from the finiteness of $\Br(A)$, proved by Tate and Milne,  
that $\Br(A,-D)[p^\infty]$ is finite as well and hence the above duality yields that
$\Br(A, D)[p^\infty]$ is the direct sum of a divisible group with a finite group, 
for all effective Cartier divisors $D$ on $A$ with $D_{\red}$ 
a simple normal crossing divisor, see Remark \ref{rmk:BrA}.

\medskip

Finally we state the  key results which are needed to complete the proof of Theorem \ref{introthm1}.

\begin{introthm}[Theorem \ref{thm:strHWM}, Theorem \ref{thm:HW-zeros}]\label{introthm7}
There are short exact sequences of $W_{n+1}\sO_X$-modules
\[0\lra B_n\Omega^{q+1}_{(X,D)}\lra \frac{W_{n+1}\Omega^q_{(X,D)}}{\ul{p}W_n\Omega^q_{(X,D)}}\xr{F^n} Z_n\Omega^q_{(X,D)}
\lra 0,\]
and
$$0\lra (\Omega/B)^q_{n, (X,-D)} \xra{V^n}
\Ker\left(W_{n+1}\Omega^q_{(X,-D)}\xra{R} W_n\Omega^q_{(X,-D)}\right)
\lra
(\Omega/Z)^{q-1}_{n, (X,-D)}
\lra 0,
$$
where in both cases the outer terms are locally free $\sO_X$-modules.
Moreover, multiplication induces $\sO_X$-linear isomorphisms 
\[
(\Omega/B)^q_{n, (X,-D)} \xra{\simeq} \sH om(Z_n\Omega^{N-q}_{(X,D)},\Omega_X^N)
,\quad 
(\Omega/Z)^{q-1}_{n, (X,-D)}
\xra{\simeq}
\sH om(B_n\Omega^{N-q+1}_{(X,D)},\Omega_X^N).\]
\end{introthm}

Recall that using Theorem \ref{introthm2} we can show that there are well-defined maps as in Theorem \ref{introthm1}.
With Theorem \ref{introthm7} at hand it is direct to check that the same arguments from Ekedahl's proof 
of Theorem \ref{introthm1} for $D=\emptyset$ work in the general case as well.

The outer terms of the two exact sequences in Theorem \ref{introthm7} are defined in 
\eqref{thm:strHWM1} and \eqref{para:Omega/BZ3}. 
For $D=\emptyset$ the exact sequences were proven by Illusie and Ekedahl and the fact that the outer terms
are locally free $\sO_X$-modules is proven by Cartier operator calculus. 
The  strategy for $D\neq \emptyset$ is the same, but it turns out to be more involved.
For example, by definition of $(\Omega/B)^q_{n, (X,-D)}$ there is a cartesian diagram of $W_{n+1}\sO_X$-modules
\[\xymatrix{
(\Omega/B)^q_{n, (X,-D)}\ar[r]\ar[d] & W_{n+1}\Omega^q_{(X,-D)}\ar[d]\\
\frac{F_{X*}^n\Omega^q_{X}}{B_n\Omega^q_{X}}\ar[r]^{V^n} & W_{n+1}\Omega^q_X.
}\]
But in general $(\Omega/B)^q_{n, (X,-D)}$  is not a quotient of 
\[\Omega^q_{(X,-D)}= \Omega^q_X(\log D_{0,\red})(-D),\]
where $D=D_0+pD_1$ and $D_0$ has none of its multiplicities divisible by $p$.
It rather depends on the whole $p$-divisibility decomposition of $D$ up to length $n+1$.
More precisely, write 
\[D=D_0+pD_1+\ldots+ p^{n+1}D_{n+1},\]
where $D_i$ and $D_j$ have no irreducible components in common for $i\neq j$, and 
$p$ does not divide any multiplicity of $D_i$, for $i=0,\ldots, n$. Then there is a surjection
\[\bigoplus_{j=0}^{n-j} 
F^{n-j}_{X*}\left(\Omega^q_X(\log \ul{D}_{j,\red})\otimes_{\sO_X}\sO_X(-\lceil \ul{D}_j/p^j\rceil- \ul{D}^{n+1-j})\right)
\xr{\oplus_j C^{-j}} (\Omega/B)^q_{n, (X,-D)},\]
where $\ul{D}_j=D_0+\ldots+p^jD_j$ and $\ul{D}^{n+1-j}=pD_{j+1}+\ldots+ p^{n+1-j}D_{n+1}$,
see Remark \ref{rmk:pHWM-VdV-zeros}. The proof of Theorem \ref{introthm7} takes almost all 
of the sections \ref{sec:strHWM} and \ref{sec:strHWZ}.

\subsection*{General conventions}
Throughout the whole paper,  {\em $k$ is a perfect field of characteristic $p>0$}.
We denote by $\Sm$ the category of separated schemes which are smooth and of finite type over $k$.
If $F$ is a Nisnevich sheaf on $\Sm$ and $X\in \Sm$, then we denote by $F_{X}$ the restriction of $F$ to the small
Nisnevich site $X_{\Nis}$. If $x\in X$ is a point, then we denote by $F_{X,x}^h$ the Nisnevich stalk at $x$.

\section{The {$p$}-saturated filtration}\label{sec:p-sat-fil}

We start by fixing some standard notation and recalling some results on the de Rham-Witt complex.
In \ref{para:fil} and Definition \ref{defn:p-sat} we introduce a filtration 
which will play an essential role  throughout the rest of the paper.

\begin{para}\label{para:DRW}
Let $X\in \Sm$.
We denote by $W_n\Omega_X^\bullet$ the de Rham-Witt complex of Bloch-Deligne-Illusie of length $n$  on $X$ 
(see \cite{IlDRW}). We denote by $R: W_{n+1}\Omega_X^\bullet\to W_n\Omega_X^\bullet$, 
$V: W_n\Omega_X^\bullet\to W_{n+1}\Omega_X^\bullet$,  $F:W_{n+1}\Omega_X^\bullet\to W_n\Omega_X^\bullet$, 
the restriction, the Verschiebung, and the Frobenius morphism, respectively,
 which are part of the structure of the de Rham-Witt complex. Furthemore, we have the map
 \eq{para:DRW0}{\ul{p}: W_n\Omega^\bullet_X\to W_{n+1}\Omega^\bullet_X}
 which is given by ``lifting to level $n+1$ and multiply by $p$", it is well-defined and injective by \cite[I, Proposition 3.4]{IlDRW}.
Recall that $W_n\Omega^\bullet_X$ is a differetial graded $W_n(k)$-algebra;
we denote by $d: W_n\Omega_X^\bullet\to W_n\Omega_X^{\bullet+1}$ the differential and by $W_n\Omega^q_X$
its degree $q$ part. Also recall that $W_n\Omega^0_X=W_n\sO_X$ is the sheaf of Witt vectors of length $n$ on $X$
and that we have the multiplicative Teichm\"uller lift $[-]=[-]_n: \sO_X\to W_n\sO_X$ at our disposal.
By \cite[Corollary 3.2.5(3)]{KSY-RSCII}  the functor
\[X\mapsto  H^0(X, W_n\Omega^q_X)\]
defines a Nisnevich sheaf with transfers for all $q\ge 0$, which has SC-reciprocity;
 we also denote it by $W_n\Omega^q\in \RSC_{\Nis}$.
See also \cite{CR12} for details on how to define the transfers structure and see \ref{para:RSC} for the definition of $\RSC_{\Nis}$.
If $f: X\to Y$ is a morphism in $\Sm$, then the morphism 
\[\Gamma_f^*=f^*:  W_n\Omega^q(Y)\to  W_n\Omega^q(X)\]
induced by its graph $\Gamma_f\in \Cor(X,Y)$, is the natural pullback morphism induced 
by the functoriality of the de Rham-Witt complex. If $f$ is finite and surjective, then the transpose of the graph
defines an element $\Gamma^t_f\in \Cor(Y,X)$ and  $\Gamma_f^{t*}=f_*$, where $f_*$ is the pushforward defined using 
duality theory. We list some properties of the transfers structure which will be used later,
we refer to \cite[Lemmas 7.7]{RS21} for  proofs of the first two and further references to the literature:
\begin{enumerate}[label=(\alph*)]
    \item\label{para:DRW1}  The restriction, Verschiebung, Frobenius, $\ul{p}$, and the differential define morphisms in $\RSC_{\Nis}$
\[R: W_{n+1}\Omega^q\to W_n\Omega^q, \quad V: W_n\Omega^q\to W_{n+1}\Omega^q,\quad F: W_{n+1}\Omega^q\to W_n\Omega^q,\]
\[\ul{p}: W_n\Omega^q\to W_{n+1}\Omega^{q+1},  \quad d: W_n\Omega^q\to W_n\Omega^{q+1}.\]
   \item\label{para:DRW2} The Nisnevich sheaf with transfers $W_n\Omega^0=W_n\sO$ coincides with the Nisnevich sheaf with transfers
defined by the algebraic group $W_n$ in  \cite[Proof of Lemma 3.2]{SpSz}.
\item\label{para:DRW4} Denote by $K^M_q$, $q\ge 0$, the restriction of the improved Milnor $K$-theory from \cite{Ke10} to $\Sm$.
It is an $\A^1$-invariant Nisnevich sheaf with transfers, in particular $K^M_q\in\RSC_{\Nis}$
and the map
\[\dlog: K^M_q\to W_n\Omega^q, \quad u=\{u_1,\ldots, u_q\}\mapsto \dlog u:=\dlog_n u:= \frac{d[u_1]_n}{[u_1]}\cdots \frac{d[u_q]_n}{[u_q]}, \quad u_i\in \sO^\times,\]
defines a morphism in $\RSC_{\Nis}$, e.g. \cite[11.1(4)]{BRS}. 
In particular if $f: Y\to X$ is a finite  and surjective morphism in $\Sm$, then 
\[f_*\dlog u= \dlog( \Nm(u)) \quad \text{in } W_n\Omega^1(X),\quad u\in \sO_Y^\times(Y),\]
where $\Nm: f_*\sO^\times_Y\to \sO_X^\times$ denotes the usual norm.
\end{enumerate}
Another property of the de Rham-Witt forms which will be important in the following is that the natural map of Nisnevich sheaves (without transfers)
\eq{para:DRW5}{(W_n\sO\otimes_{\Z} K^M_q) \oplus (W_n\sO\otimes_{\Z} K^M_{q-1})\to W_n\Omega^q,\quad (a\otimes u, b\otimes v)\mapsto a\dlog u+db\dlog v,}
is surjective. This follows easily by induction over $n$ from the corresponding fact for $\Omega^q$, the exact sequence of Nisnevich sheaves on $\Sm$
(see \cite[I, Proposition 3.2]{IlDRW})
\[\Omega^q\oplus \Omega^{q-1}\xr{V^{n-1}+dV^{n-1}} W_n\Omega^q\xr{R} W_{n-1}\Omega^q\to 0,\]
and the formulas  $V(a\dlog u)=V(a)\dlog u$ and $Vd=pdV$, cf. also  \cite[Proposition (4.6)]{HK}.
\end{para}

The following notation will be used throughout.
\begin{notation}\label{nota:hdvf}
A {\em henselian discrete valuation field of geometric type over $k$} is a field 
\[L=\Frac(\sO_{U,z})^h,\]
where $U\in \Sm$, 
$z\in U^{(1)}$ is a point of codimension 1, and $(\sO_{U,z})^h$  denotes the henselization of the local ring $\sO_{U,z}$.
We denote the set of all  such $L$ by $\Phi$.
For $L\in \Phi$ we denote by $\sO_L$, $v_L$, $\fm_L$, $\kappa_L$ the ring of integers, the normalized discrete valuation, the maximal ideal, 
and the residue field, respectively.  In case there is no ambiguity we also write $\fm$ instead of $\fm_L$ and $\kappa$ instead of $\kappa_L$. 
For each $L\in \Phi$, we pick  a local parameter $z=z_L\in \sO_L$.
\end{notation}

\begin{para}\label{para:fillog}
Let $L\in \Phi$ (see \ref{nota:hdvf}) and pick a local parameter $z\in \fm_L$. 
Recall that the {\em Brylinski-Kato filtration} of $W_n(L)$ is defined as follows (see \cite{Brylinski}):
\begin{align*}
\fil^{\log}_r W_n(L)& :=\{(a_0,\ldots, a_{n-1})\mid p^{n-1-i} v_L(a_i)\ge -r\}\\
                          &= \{a\in W_n(L)\mid [z]^r \bF^{n-1}(a)\in W_n \sO_L\},\quad r\ge 0.
\end{align*}
We have
\[\fil^{\log}_0 W_n(L)=W_n\sO_L, \quad \fil^{\log}_r L=\fil^{\log}_r W_1(L)= \fm_L^{-r}.\]
Furthermore, the maps $V$, $R$, $F$ restrict to 
\eq{para:fillog1}{\xymatrix{ \fil^{\log}_r W_n(L)\ar[r]^{V} &  \fil^{\log}_r W_{n+1}(L)\ar@/^1pc/[l]^F\ar@/_1.5pc/[l]_R    }}
and more precisely, 
\eq{para:fillog1.5}{
R(\fil^{\log}_r W_{n+1}(L))\subset \fil^{\log}_{\lfloor r/p\rfloor}W_n(L) \subset \fil^{\log}_r W_n(L).}
It follows from the formula $V^i([a])\cdot V^j([b])=V^{i+j}([a^{p^j}][b^{p^i}])$ that we have 
\eq{para:fillog2}{\fil^{\log}_rW_n(L)\cdot \fil^{\log}_sW_n(L)\subset \fil^{\log}_{r+s} W_n(L), \quad r,s\ge 0.}
In particular, $\fil^{\log}_r W_n(L)$ is a $W_n \sO_L$-submodule of $W_n(L)$. 
Moreover the injective map $\ul{p}: W_n(L)\to W_{n+1}(L)$ induces an isomorphism
\eq{para:fillog3}{\ul{p}: \fil^{\log}_r W_n(L)\xr{\simeq} p\cdot \fil^{\log}_{rp}W_{n+1}(L).}
\end{para}

\begin{defn}\label{defn:fil}
Let $L\in \Phi$. For $q$, $r\ge 0$ we define
\[\fil^{\log}_r W_n\Omega^q_L:=\Im\left( \fil^{\log}_rW_n(L)\otimes_{\Z}K^M_q(L)\xr{\id\otimes\dlog} W_n\Omega^q_L\right),\]
\[\fil^{\log'}_r W_n\Omega^q_L:=\Im\left( \fil^{\log}_rW_n(L)\otimes_{\Z}K^M_q(\sO_L)\xr{\id\otimes\dlog} W_n\Omega^q_L\right),\]
and  
\[\fil_r W_n\Omega^q_L:=\fil_{r-1}^{\log}(W_n\Omega^q_L)+ V^{n-m}(\fil^{\log'}_rW_m\Omega^q_L),\]
where $m=\min\{v_p(r), n\}$  with $v_p$ the normalized $p$-adic valuation. We have $\fil_0 W_n\Omega^q_L=\fil_0^{\log'}W_n\Omega^q_L$ (by convention).
Finally
\[\Fil_r W_n\Omega^q_L:= \fil_r W_n\Omega^q_L+ d(\fil_r W_n\Omega^{q-1}_L).\]
\end{defn}

\begin{para}\label{para:fil}
We make some comments and list some easy properties of the above defined filtrations:
\begin{enumerate}[label=(\arabic*)]
\item\label{para:fil1} The filtration $\{\Fil_r W_n\Omega^0_L\}_{r\ge 0}=\{\fil_r W_n(L)\}_{r\ge 0}$ coincides with 
    the Kato-Matsuda filtration \cite{Matsuda} (with the notational conventions from \cite[2.1]{KeS}).
    \item\label{para:fil2}  The family $\{\Fil_r W_n\Omega^q_L\}$ defines an  increasing and exhaustive filtration of $W_n\Omega^q_L$
    \[\Fil_0 W_n\Omega^q_L\subset \Fil_1 W_n\Omega^q_L\subset \ldots \subset \Fil_r W_n\Omega^q_L\subset\ldots \subset W_n\Omega^q_L.\]
    (The filtration is exhaustive by the surjectivity of \eqref{para:DRW5}.)
    \item\label{para:fil3} The surjectivity of \eqref{para:DRW5} yields
    \[\Fil_0 W_n\Omega^q_L=W_n\Omega^q_{\sO_L}, \quad \Fil_1W_n\Omega^q_L= W_n\Omega^q_{\sO_L}(\log),\]
    where 
    \[W_n\Omega^q_{\sO_L}(\log)= W_n\Omega^q_{\sO_L}+ W_n\Omega^{q-1}_{\sO_L}\dlog z\subset W_n\Omega^q_L,\]
    which is also equal to  the $q$th de Rham-Witt forms of the log ring $(\sO_L, \lambda)$  over $k$ with trivial log structure, 
    where $\lambda$ is the  logarithmic structure  associated to $\N\to \sO_L$, $1\mapsto z$, see e.g. \cite{Matsuue}. 
    \item\label{para:fil4} We have 
    \[\fil_r W_n\Omega^q_L=\begin{cases}
      \fil^{\log}_{r-1}W_n\Omega^q_L & \text{if } v_p(r)=0\\
       \fil_{r-1}^{\log}W_n\Omega^q_L +\fil_r^{\log'}W_n\Omega^q_L & \text{if } v_p(r)\ge n.
      \end{cases}\]
    \item\label{para:fil5} For $j\ge 0$ we have a natural map 
    \[\Fil_r W_n\Omega^{j}_L\otimes_{\Z} K^M_{q-j}(\sO_L)\to \Fil_r W_n\Omega^q_L, \quad a\otimes u\mapsto a\dlog u,\]
    which is surjective for $j\ge 1$. (Use the formulas $V(a\dlog u)=V(a)\dlog u$ and $d \dlog u=0$ and the surjectivity of
    $L^\times\otimes_{\Z} K^M_{q-1}(\sO_L)\surj K^M_q(L)$, $x\otimes u\mapsto \{x,u\}$.)
    \end{enumerate}
\end{para}

\begin{defn}\label{defn:p-sat}
Let $\{G_r\}_{r\ge 0}$ be a family of subgroups of $W_n\Omega^q_L$ with $G_r\subset G_{r+1}$.
The {\em $p$-saturation of $G_r$} is then defined by 
\[G_r^p:=\sum_{s= 0}^{n-1} p^s G_{rp^s}\subset W_n\Omega^q_L.\]
In particular we have 
\[\Fil^p_r W_n\Omega^q_L= \sum_{s\ge 0} p^s\Fil_{rp^s}W_n\Omega^q_L= \fil_r^p W_n\Omega^q_L+ d(\fil^p_r W_n\Omega^{q-1}_L).\]
\end{defn}

We immediately get from \ref{para:fil} that $\{\Fil^p_rW_n\Omega^q_L\}_{r\ge 0}$ is an increasing and exhaustive filtration of $W_n\Omega^q_L$, which 
comes with a natural map
\eq{defn:p-sat1}{\Fil^p_rW_n\Omega^j_L\otimes_{\Z} K^M_{q-j}(\sO_L)\to \Fil^p_{r}W_n\Omega^q_L,}
which is surjective for $j\ge 1$. 

\begin{lem}\label{lem:filp-round}
For all $0\le s\le n-1$ and all $r\ge 1$ we have 
\[p^s\fil^{\log}_{(r-1)p^s}W_n\Omega^q_L= p^s\fil^{\log}_{rp^s-1}W_n\Omega^q_L,\quad \text{and}\quad
p^s\fil^{\log'}_{(r-1)p^s}W_n\Omega^q_L= p^s\fil^{\log'}_{rp^s-1}W_n\Omega^q_L.\]
\begin{proof}
Clearly we have this $\subset$ inclusion in both cases. For the other direction it suffices to consider $q=0$.
Let $y\in L$ with $p^{n-j-1}v_L(y)\ge -rp^s+1$. If $s\ge n-j$, then $p^s V^{j}([y])=0$; 
if $s\le n-j-1$, then  $p^{n-j-s-1}v_L(y)$ is an integer $\ge -r+1/p^s$ and hence it is  $\ge -r+1$. 
This yields this $\supset$ inclusion in both cases. 
\end{proof}

\end{lem}

\begin{lem}\label{lem:Filp}
We have
\[\Fil^p_0W_n\Omega^q_L=W_n\Omega^q_{\sO_L},\quad \Fil^p_1W_n\Omega^q_L=W_n\Omega^q_{\sO_L}(\log).\]
\end{lem}
\begin{proof}
By \ref{para:fil}\ref{para:fil3} it remains to show the inclusion $p^s\fil_{p^s}W_n\Omega^q_L\subset W_n\Omega^q_{\sO_L}(\log)$, for $s\in [1, n-1]$.
Since $p^s V^{n-s}(W_s\Omega^q_L)=0$, it suffices to show 
\[p^s \fil^{\log}_{p^s-1}W_n\Omega^q_L\subset W_n\Omega^q_{\sO_L}(\log),\]
which follows from Lemma \ref{lem:filp-round}.
\end{proof}

\begin{lem}\label{lem:fil-FVR}
Set $w_n:=W_n\Omega^q_L$. We have 
\eq{lem:fil-FVR1}{V(\fil^{\log}_r w_n)\subset \fil^{\log}_r w_{n+1}, \quad V(\fil^{\log'}_r w_n)\subset \fil^{\log'}_r w_{n+1},
\quad V(\fil_r w_n)\subset \fil_r w_{n+1},}
\eq{lem:fil-FVR1.5}{\ul{p}(\fil^{\log}_r w_n)= p\cdot \fil^{\log}_{rp} w_{n+1}, \quad \ul{p}(\fil^{\log'}_r w_n)= p\cdot \fil^{\log'}_{rp} w_{n+1},
\quad \ul{p}(\fil_r w_n)= p\cdot\fil_{rp} w_{n+1},}
and for $\varphi\in \{R,F\}$ 
\eq{lem:fil-FVR2}{\varphi(\fil^{\log}_r w_{n+1})\subset \fil^{\log}_r w_n, \quad 
\varphi(\fil^{\log'}_r w_{n+1})\subset \fil^{\log'}_r w_n,
\quad \varphi(\fil_r w_{n+1})\subset \fil_r w_n.}
Furthermore,
\eq{lem:fil-FVR3}{Fd(\fil^{\log}_{r}W_{n+1}\Omega^{q-1}_L)\subset \fil^{\log}_r W_n\Omega^q_L+ d(\fil^{\log}_r W_n\Omega^{q-1}_L)}
and if $v_p(r)\ge n+1$, then 
\eq{lem:fil-FVR4}{Fd(\fil^{\log'}_r W_{n+1}\Omega^{q-1}_L)\subset p\cdot \fil^{\log}_r W_n\Omega^q_L+\Fil_r W_n\Omega^q_L.}
\end{lem}
\begin{proof}
The inclusions \eqref{lem:fil-FVR1} and \eqref{lem:fil-FVR2} follow immediately from
\eqref{para:fillog1} and the formula $\varphi(a \dlog x)=\varphi(a)\dlog x$, for $\varphi\in \{F,R,V\}$,
the first two equalities in \eqref{lem:fil-FVR1.5} follow from \eqref{para:fillog3}, the last equality from this and Lemma \ref{lem:filp-round}.
Note that any element in $\fil^{\log}_r W_{n+1}\Omega^{q-1}_L$ is a sum of elements 
$\omega= a\dlog x$, with  $a\in \fil^{\log}_r W_{n+1}(L)$ and $x\in K^M_{q-1}(L)$. Write $a=[a_0]+V(b)$, so that 
$[a_0]\in \fil^{\log}_{r}W_{n+1}(L)$ and $b\in \fil^{\log}_{r}W_n(L)$. Thus  \eqref{lem:fil-FVR3} follows from
\[Fd(\omega)= F([a_0])\dlog\{a_0,x\}+ d(b\dlog x) \in \fil^{\log}_r W_n\Omega^q_L+ d(\fil^{\log}_r W_n\Omega^{q-1}_L).\]
Finally, \eqref{lem:fil-FVR4}. We assume $v_p(r)\ge n+1$. Let $\omega$ be as above but this time $x\in K^M_{q-1}(\sO_L)$, 
and write $a_0= z^e u$, with $u\in \sO_L^\times$ and $ep^n \ge -r$. We get
\[Fd(\omega)= \underbrace{eF([a_0])\dlog\{z,x\}}_{=\omega_1} + \underbrace{F([a_0])\dlog\{u,x\}+ d(b\dlog x)}_{\in \Fil_r W_n\Omega^q_L}.\]
If $ep^n  > -r$, then $[a_0]\in \fil^{\log}_{r-1}W_{n+1}(L)$ and hence $\omega_1\in \fil^{\log}_{r-1}W_n\Omega^q_L$.
If $ep^n=-r$, then $p|e$ and hence $\omega_1\in p\cdot\fil^{\log}_{r}W_n\Omega^q_L$. This yields the statement. 
\end{proof}

\begin{cor}\label{cor:Filp-FVR}
The maps $F$, $R$, $V$, $\ul{p}$, $d$ on $W_\bullet\Omega^*_L$ induce well-defined maps
\[\xymatrix{ \Fil^p_r W_n\Omega^q_L\ar@<1ex>[r]^-{V,\, \ul{p}} &  \Fil^p_r W_{n+1}\Omega^q_L\ar@<1ex>[l]^-{F,\, R}}, \quad 
             d: \Fil^p_r W_n\Omega^q_L \to \Fil^p_r W_n\Omega^{q+1}_L.\]
Furthermore,
\[ \ul{p}(\Fil^p_r W_{n-1}\Omega^q_L)= p\Fil^p_{pr}W_n\Omega^q_L.\]
\end{cor}
\begin{proof}
This holds for $d$ by definition. The well-definedness of $\ul{p}$ and the final statement follow directly from \eqref{lem:fil-FVR1.5}.
For $V$ and $R$ the statement follows from \eqref{lem:fil-FVR1} and \eqref{lem:fil-FVR2}
and the formulas $Vd=pdV$ and $Rd=dR$. For $F$ it follows from \eqref{lem:fil-FVR2} - \eqref{lem:fil-FVR4}
and the observation
\[\fil^{\log}_{rp^s} W_n\Omega_L^q\subset \fil^{\log}_{rp^{s+1}-1} W_n\Omega_L^q \subset  \fil_{rp^{s+1}}W_n\Omega^q_L.\]
\end{proof}

It will be useful to have a different presentation of $\Fil^p_rW_n\Omega^q_L$ which is provided in the following lemma.
\begin{lem}\label{lem:pres-filp}
Let $r\ge 1$ and let $z\in \sO_L$ be a local parameter.
\begin{enumerate}[label=(\arabic*)]
\item\label{lem:pres-filpI} For  $0\le j\le n-1$, $i\in \Z$,  $\alpha\in W_{n-j}\Omega^q_{\sO_L}$, and $\beta\in W_{n-j}\Omega^{q-1}_{\sO_L}$ we have 
\eq{lem:pres-filp1}{ V^j([z]^i \alpha)\in \fil^{\log}_{r-1}W_n\Omega^q_L+ d(\fil^{\log'}_{r-1}W_n\Omega^{q-1}_L),
\qquad \text{if}\quad ip^{n-1-j}\ge -r+1,}
\eq{lem:pres-filp2}{ V^j([z]^i\dlog z\cdot \beta)\in \fil^{\log}_{r-1}W_n\Omega^q_L + d(\fil^{\log}_{r-1}W_n\Omega^{q-1}_L),
\qquad \text{if} \quad  ip^{n-1-j}\ge -r+1,}
Furthermore,  if  $n-m\le j\le n-1$, where $m=\min\{v_p(r), n\}\ge 1$,  then 
\eq{lem:pres-filp3}{V^j([z]^i\cdot \alpha)\in 
p \cdot V^{n-m}(\fil^{\log}_r W_m\Omega^q_L) + 
dV^{n-m}(\fil^{\log'}_r W_m\Omega^{q-1}_L),
\qquad \text{if}\quad ip^{n-j-1}=-r.
}
\item\label{lem:pres-filpII} Let $H^q_r\subset W_n\Omega^q_L$ be the subgroup generated by the elements \eqref{lem:pres-filp1} - \eqref{lem:pres-filp3}.
Then 
\[\fil_r W_n\Omega^q_L\subset H^q_r\subset \Fil_r W_n\Omega^q_L+ p\cdot \Fil_{rp} W_n\Omega^q_L.\]
In particular
\[\Fil^p_r W_n\Omega^q_L= \sum_{s=0}^{n-1}\ p^s (H^q_{rp^s}+d(H^{q-1}_{rp^s})).\]
\end{enumerate}
\end{lem}
\begin{proof}
\ref{lem:pres-filpI}. 
If $\alpha$ (resp. $\beta$) is equal to  $a\dlog u$ with $a\in W_{n-j}(\sO_L)$ and  $u\in K^M_{q}(\sO_L)$ (resp. $u\in K^M_{q-1}(\sO_L)$),
then \eqref{lem:pres-filp1} - \eqref{lem:pres-filp3} follow directly from \ref{para:fillog} and Definition \ref{defn:fil}.
By the surjectivity of \eqref{para:DRW5} and \ref{para:fil}\ref{para:fil5} 
it remains to consider the case where $\alpha$ (resp. $\beta$) is equal to 
$db$ with $b\in W_{n-j}(\sO_L)$ and $q=1$ (resp. $q=2$). 
We compute
\[V^j([z]^i db)= p^j d V^j(b[z]^i)- V^j(ib[z]^i)\dlog z,\]
where the equality follows from the Leibniz rule and the formula $Vd=pdV$.
This yields \eqref{lem:pres-filp1}. 
Observe that in the situation of \eqref{lem:pres-filp3} the integer $i$ is divisible by $p$, which yields \eqref{lem:pres-filp3}.
Similarly we get
\[V^j([z]^i\dlog(z)db)=-p^jd V^j(b[z]^i)\dlog z,\]
whence \eqref{lem:pres-filp2}. Concerning the chain of inclusions in \ref{lem:pres-filpII} 
observe that the left inclusion holds by definition and that the right inclusion follows 
from \ref{lem:pres-filpI}, where for the elements of type \eqref{lem:pres-filp3} 
we observe that 
\[p\cdot\fil_r^{\log}W_n\Omega^q_L\subset p\cdot\fil_{r+1}W_n\Omega^q_L
\subset p\cdot\fil_{rp}W_n\Omega^q_L.\]
This completes the proof.
\end{proof}

\begin{cor}\label{cor:fil-WO-mod}
The group $\Fil^p_rW_n\Omega^q_L$ is a finitely generated $W_n\sO_L$-submodule of $W_n\Omega^q_L$, for all $r, q\ge 0$ and $n\ge 0$.
Moreover  we have inclusions  of $W_n\sO_L$-modules
\eq{cor:fil-WO-mod1}{W_n\Omega^q_{\sO_L}\subset \Fil^p_r W_n\Omega^q_L\subset W_n\Omega^q_{\sO_L}\cdot \frac{1}{[z]^{rp^{n-1}}},}
where $z\in \sO_L$ denotes a local parameter.
\end{cor}
\begin{proof}
Let $H^q_r$ be as in Lemma \ref{lem:pres-filp}. It is a $W_n\sO_L$-submodule of $W_n\Omega^q_L$, by the formula $a V^j(x)=V^j(F^j(a)x)$.
Using additionally the Leibniz rule we find for $\gamma\in W_{n-j}\Omega^{q-1}_L$ and $a\in W_n\sO_L$ 
\[adV^j([z]^i\gamma)= dV^j([z]^i F^j(a)\gamma)-V^j([z]^i F^j(da)\gamma ).\]
Hence $H^q_r+d(H^q_r)$ is a $W_n\sO_L$-submodule of $W_n\Omega^q_L$, for all $r\ge 0$, and hence so is 
$\Fil^p_r W_n\Omega^q_L$, by Lemma \ref{lem:pres-filp}. 
The first inclusion in \eqref{cor:fil-WO-mod1} holds by Lemma \ref{lem:Filp}, the second  inclusion 
follows from $[z]^{t}\cdot \fil^{\log}_t W_n(L)\subset W_n\sO_L$, for all $t\ge 0$.
Since $W_n\sO_L$ is noetherian (e.g. \cite[Prop A.4]{LZ}) and $W_n\Omega_{\sO_L}^q$ is a finite $W_n\sO_L$-module
(where we use that $L$ is of geometric type), we obtain that $\Fil^p_r W_n\Omega^q_L$ is a finite 
$W_n\sO_L$-module as well.
\end{proof}

\section{The {$p$}-saturated filtration defines a conductor}
In this section we show that the $p$-saturated filtration from section \ref{sec:p-sat-fil} defines a conductor in the sense of
\cite{RS21}.

\begin{para}\label{para:mod-pairs}
Following \cite{KMSYI} we call a pair $(X,D)$ a {\em modulus pair}, if $X$ is separated and of finite type over $k$ 
and $D$ is an effective Cartier divisor on $X$, such that $X\setminus D$ is smooth. The modulus pair $(X,D)$ is called \textit{proper} if 
$X$ is a proper $k$-scheme.  
A {\em compactification of a modulus pair} $(X,D)$ is a proper modulus pair $(\ol{X}, D'+ X_\infty)$ such that
$\ol{X}\setminus X_{\infty}=X$ and $D'_{|X}=D$.  All compactifications of a fixed  modulus pair form a cofiltered set, see
\cite[Lemma 1.8.2]{KMSYI}.
\end{para}

\begin{para}\label{para:cond}
Recall from \cite[Definition 4.3]{RS21} that a {\em conductor} on  a presheaf with transfers $F$ on $\Sm$ is
a collection of set maps 
\[c=\{c_{L}: F(L)\to \N_0\}_{L\in\Phi},\]
where $\Phi$ denotes the set of henselian discrete valuation fields of geometric type over $k$, see Notation \ref{nota:hdvf},
satisfying the following properties for all $L\in\Phi$ and all $X\in \Sm$:
\begin{enumerate}[label = {(c\arabic*)}]
\item\label{c1} $c_{L}(a)=0 \,\Rightarrow\, a\in \Im(F(\sO_L)\to F(L))$.
\item\label{c2} $c_L(a+b)\le \max\{c_L(a), c_L(b)\}$.
\item\label{c3} For any finite extension  $L'/ L$ of ramification degree $e$ and any $a\in F(L')$ we have 
       \[c_L(\Tr_{L'/L}a)\le \left\lceil\frac{c_{L'}(a)}{e}\right\rceil,\]
        where $\Tr_{L'/L}$ is the trace which is given by the pullback along the transpose of the graph of $\Spec L'\to \Spec L$ 
        (viewed as a finite correspondence) and $\lceil -\rceil$ is the round up.
\item\label{c4}
      Let $a\in F(\A^1_X)$. Then 
      \[c_{k(x)(t)_\infty}(\rho_x^*a)\le 1, \quad \forall \, x\in X \,\Longrightarrow  a\in F(X),\]
      where $k(x)(t)_\infty:=\Frac(\sO_{\P^1_x,\infty}^h)$, $\rho_x :\Spec k(x)(t)_\infty \to \A^1_X$ is the natural map,
      and  we identify $F(X)$ with its image in $F(\A^1_X)$ under the  pullback along  $\A^1_X\to X$.
\item\label{c5}
        For any $a\in F(X)$ there exists a proper modulus pair $(\ol{X}, D)$
       with $X=\ol{X}\setminus D$, such that for all  $k$-morphisms $\rho:\Spec L \to X$ 
         we have  
            \[c_L(\rho^*a)\le v_L(D_L),\]
            where $D_L$ denotes the pullback of $D$ under the unique extension $\Spec \sO_L\to \ol{X}$ of $\rho$ and 
            $v_L(D_L)$ denotes its multiplicity.
\end{enumerate}
We say a conductor has {\em level $n$} if in \ref{c4} it suffices to consider points 
\[x\in X_{(\le n-1)}= \{\text{points of dimension } \le n-1\}.\footnote{
In fact in \cite{RS21} a conductor of level $n$ is only required to be defined for $L\in\Phi$ of transcendence degree over $k$ at most $n$. 
The restriction of conductor of level $n$ in the above sense to those $L$ is a conductor of level $n$ in the sense of \cite{RS21}.}\]
We say $c$ is  {\em semi-continuous} if it satisfies the following condition:
\begin{enumerate}[label = {(c6)}]
\item\label{c6} Let $X\in \Sm$ and let $Z\subset X$ be a smooth prime divisor with generic point $z$  
and $K=\Frac(\sO_{X,z}^h)$.  Then for any $a\in F(X\setminus Z)$ with $c_K(a_K)\le r$ there exists 
a Nisnevich neighborhood $u:U\to X$ of $z$ and a compactification $(Y, E)$ of the modulus pair $(U, r\cdot u^*Z)$ such that 
\[c_{L}(\rho^*a_{U})\le v_L(E_L), \quad \text{for all } L\in \Phi \text{ and all } \rho:\Spec L\to U,\]
 where $a_U$ (resp. $a_K$) denotes the restriction of $a$ to $U$ (resp. $K$).
\end{enumerate}
\end{para}

\begin{thm}\label{thm:fil-cond}
The collection  $c=\{c_L:W_n\Omega^q_L\to \N_0\}_{L\in \Phi}$ with
\eq{eq:dRW-cond1}{c_L(a)=\min\left\{r\ge 0\mid a\in \Fil^p_r W_n\Omega^q_L\right\}}
defines a semi-continuous conductor of level $q+1$. (See Definition \ref{defn:p-sat}  for $\Fil^p_r W_n\Omega^q_L$.)
\end{thm}

\begin{proof}
\ref{c1} holds by Lemma \ref{lem:Filp} and \ref{c2} is clear.
The proof of \ref{c3} requires some more lemmas and will be given at the end of the next subsection.
\ref{c4} Let $X\in \Sm$ and $a\in W_n\Omega^q(\A^1_X)$. For  $x\in X$  
let $\rho_x: \Spec k(x)(t)_{\infty}\to \A^1_X$ be as in \ref{c4}. 
Assume $c_{k(x)(t)_\infty}(\rho_x^*a)\le 1$, for all points $x\in X_{(\le q)}$.
Since $W_n\Omega^q_{\P^1_x}(\log \infty)$ is a Nisnevich sheaf we find by Lemma \ref{lem:Filp} 
\[a_{|\A^1_x}\in \Ker\left(W_n\Omega^q(\A^1_x)\to \frac{W_n\Omega^q_{k(x)(t)_\infty}(\log \infty)}{(W_n\Omega^q_{\P^1_x,\infty}(\log \infty))^h}\right)
= W_n\Omega^q(\log \infty)(\P^1_x).\]
We have  an exact sequence 
\[0\to W_n\Omega^q(\P^1_x)\to W_n\Omega^q(\log \infty)(\P^1_x)\xr{\Res} W_n\Omega^{q-1}(x) \xr{\delta} H^1(\P^1_x, W_n\Omega^q_{\P^1_x}),\]
see \cite[p. 68]{Matsuue}. By the projective bundle formula (see \cite[I, Corollaire 4.2.13]{Gros})
the map $\delta$ is a split injection and we find
\[a_{|\A^1_x}\in W_n\Omega^q(\log \infty)(\P^1_x)=W_n\Omega^q(\P^1_x)=W_n\Omega^q(x),\]
i.e., $a_{|\A^1_x}$ is pulled-back from $W_n\Omega^q(x)$.
Let $s:X\inj \A^1_X$ be a section of the projection $\pi: \A^1_X\to X$ and consider $b:= a-\pi^*s^*a\in W_n\Omega^q(\A^1_X)$.
By the above, 
\eq{thm:fil-cond2}{b_{|\A^1_x}=0, \quad \text{for all }x\in X_{(\le q)}.} 
We claim that \eqref{thm:fil-cond2}  implies $b=0$. To prove this latter claim we may assume $X=\Spec A$,  with $A$ a smooth integral $k$-algebra
and $b\in W_n\Omega^q_{A[t]}$. 
By \cite[Theorem B]{HeMa04} (for $p$ odd) and \cite[Theorem 4.3]{Costeanu} (for $p=2$) 
the {\em abelian group}  $W_n\Omega^q_{A[t]}$ is isomorphic to a direct sum of the groups $W_{m}\Omega^q_A$ and $W_m\Omega^{q-1}_A$,
for various $m\in [0, n]$. Furthermore this isomorphism is natural in $A$. The claim thus follows from Lemma \ref{lem:levelq} below.

\ref{c5}. Let $X\in \Sm$ and $a\in W_n\Omega^q(X)$. Let $(\ol{X}, D)$ be a proper modulus pair with $X=\ol{X}\setminus D$.
By the surjectivity of \eqref{para:DRW5} we find an open cover $\ol{X}=\cup_i V_i$ inducing an open cover $X=\cup_i U_i$, with $U_i=V_i\setminus D$,
and finitely many elements  $a_{i,j}\in W_n\sO(U_i)$, $b_{i,j}\in W_n\sO(U_i)$, $u_{i,j}\in K^M_q(U_i)$, and $v_{i,j}\in K^M_{q-1}(U_i)$
such that 
\[a_{U_i}= \sum_j a_{i,j} \dlog u_{i,j} + \sum_j db_{i,j} \dlog v_{i,j},\]
where $a_{U_i}$ denotes the restriction of $a$ to $U_i$. For $N\ge 1$ big enough we have 
\[a_{i,j}, \,b_{i,j}\in H^0(V_i, W_n \sO_{V_i}(N \cdot D_{|V_i})),\quad \text{for all } i,\, j.\]
Here  for a divisor $E$ on a finite type $k$-scheme $Y$, we denote by $W_n\sO_Y(E)$ the invertible $W_n\sO_Y$-modulue,
which on an open $V\subset Y$ with $E_{|V}=\Div(e)$ is equal to $W_n\sO_V\cdot \frac{1}{[e]}$.
By the proof of \cite[Claim 7.5.1]{RS21}  we obtain for $r>p^{n-1}N$ 
\[\rho^*(a_{i,j}), \, \rho^*(b_{i,j})\in \fil^{\log}_{rm_{i,L}-1}W_n(L), \quad 
\text{for all } \rho: \Spec \sO_L\to V_i, \text{with } L\in\Phi, \text{ all }i,\, j,\]
where $m_{i,L}=v_L(\rho^* D_{|V_i})$. 
By definition of $\Fil^p_r$ we find for all $\rho: \Spec \sO_L\to V_i$
\eq{thm:fil-cond3}{c_L(\rho^* a_{U_i})\le v_L(r D_{|V_i}).}
This yields \ref{c5} as any  $\Spec L\to X$ extends uniquely to $\Spec \sO_L\to \ol{X}$ and factors via some $V_i\inj \ol{X}$.

Finally \ref{c6}.
Let $Z\subset X$, $z\in Z$ and $K=\Frac(\sO_{X,z}^h)$ be as in \ref{c6}.
Let $a\in W_n\Omega^q(X\setminus Z)$ with $c_K(a_K)\le r$.
By Definition \ref{defn:p-sat} and \eqref{lem:fil-FVR1.5}
\[a_K=\sum_{s=0}^{n-1}\ul{p}^s (\alpha_{s} +  d\beta_{s})\]
with $\alpha_{s}\in \fil_{r} W_{n-s}\Omega^q_K$ and $\beta_{s}\in \fil_{r} W_{n-s}\Omega^{q-1}_K$. 
Replacing $X$ by a Nisnevich neighborhood of $z$, 
we can assume that $\alpha_{s}$ and $\beta_{s}$ are the restriction to $K$ 
of de Rham Witt forms on $X\setminus Z$.
As the set of compcatifications of $(U,rZ)$ is cofiltered, for $U$ a Nisnevich neighborhood of $z$,
see \ref{para:mod-pairs}, it suffices to prove the following, for all $q\ge 0$ and $n\ge 1$:
Assume $a$ as above satisfies $a_K\in \fil_{r} W_n\Omega^q_K$.
Then there exists a Nisnevich neighborhood
$U\to X$ and a compactification $(Y, E)$ of $(U,rZ_U)$ such that 
\eq{thm:fil-cond4}{\rho^*a\in \fil_{ v_L(E_L)} W_n\Omega^q_L, \quad \text{for all }\rho:\Spec L\to U.}

Let $m=\min\{v_p(r), n\}$. We find an open neighborhood  $U\subset X$ of $z$ such that 
\[a_U=\sum_j a_j \dlog(\lambda_j) + V^{n-m}(b_j \dlog(u_j)),\]
with $a_j\in W_n(U\setminus Z_U)$, $b_j\in W_m(U\setminus Z_U)$, $\lambda_j \in K^M_q(U\setminus Z_U)$, and 
$u_j\in K^M_q(U)$, and
\[a_{j, K}\in \fil^{\log}_{r-1}W_n(K), \quad \text{and}\quad b_{j,K}\in \fil^{log}_rW_m(K), 
\quad \text{all }j.\]
If $r=0$, then we can take $a_j=0=\lambda_j$. 
Let $(Y, Z'+ Y_\infty)$ be a compactification of $(U,Z_U)$ with $Y$ normal.
By the proof of (c6) in \cite[Proof of Proposition 7.5]{RS21} we have for any $N\gg 0$ with $p^n|N$,
and for any map $\rho: \Spec L\to U$, with $L\in\Phi$,
\[\rho^*(a_j)\in \fil^{\log}_{v_L((r-1)Z'+ (N-1)Y_\infty)} W_n(L) \subset \fil^{\log}_{v_L(rZ'+ N Y_\infty)-1} W_n(L)\]
and 
\[V^{n-m}(\rho^* b_j)\in V^{n-m_0}(\fil^{\log}_{v_L(rZ'+ N Y_\infty)} W_{m_0}(L)),\]
where $m_0=\min\{v_p(v_L(rZ'+ N Y_\infty)), n\}$.
Hence condition \eqref{thm:fil-cond4} is satisfied for $(Y, E)$ with 
$E= r Z'+ NY_\infty$,  $N$ as above, which is a modulus compactification
of $(U, rZ_U)$. 

It remains to prove \ref{c3}, which is done at the end of the next subsection.
\end{proof}

\begin{lem}\label{lem:levelq}
Let $X$ be a smooth $k$-scheme and $a\in W_n\Omega^q(X)$, for some $n\ge 1$ and $q\ge 0$.
If $a$ is nonzero, then there exists a point $x\in X_{(\le q)}$, such that the image of $a$ in $W_n\Omega^q(x)$ is nonzero.
\end{lem}
\begin{proof}
First note that the restriction $W_n\Omega^q(X)\to W_n\Omega^q(U)$ along a dense open subset $U\subset X$ is injective, 
which follows from the fact that $W_n\Omega^q_X$ is a successive extension of locally free $\sO_X$-modules, see \cite[I, Corollaire 3.9]{IlDRW}.
Since $W_n\Omega^q_X$ is moreover an \'etale sheaf we may assume $k$ is algebraically closed and $X=\Spec A$ with $A$ smooth
and \'etale over $k[t_1,\ldots, t_e]$ with $e=\dim X >q$. 

For $n=1$, the differential forms $dt_{i_1}\cdots dt_{i_q}$, $1\le i_1<\ldots< i_q\le e$, form a basis
of the $A$-module $\Omega^q_A$. Assume for $a\in \Omega^q_A$ we have $b:=(dt_{i_1}\cdots dt_{i_q})^{\vee}(a)\neq 0$ for some  sequence
$(i_1,\ldots, i_q)$ as above. As $k$ is algebraically closed we find elements $\lambda_j\in k$ such that $b\neq 0$ in 
$B:=A/(t_j-\lambda_j| j\in \{1,\ldots, e\}\setminus \{i_1,\ldots, i_q\})$. Then $a(x)\neq 0$ in $\Omega^q_{k(x)}$, 
for some generic point $x$ of $\Spec B$, which is a point of dimension $q$ in $X$. 

By induction over $n$ we assume the statement is proven for $W_m\Omega^q_Y$, for $m\le n-1$ and all $Y\in \Sm$.
For an \'etale sheaf of abelian groups $F$ on $\Sm$ consider the property:
\begin{enumerate}[label=$(*)$]
\item\label{Pq} For each non-zero element $a\in F(X)$, there exists  a morphism $f:Z\to X$  in $\Sm$ such that $\dim Z\le e-1$ and 
$f^*a\neq 0$ in $F(Z)$. 
\end{enumerate}
By the above, property \ref{Pq} holds for $W_m\Omega^q$, for $m\le n-1$, and it suffices  to show that 
$F=W_n\Omega^q$ has property \ref{Pq}.
To this end observe that if we have an exact sequence of sheaves 
\[0\to F'\to F\to F'',\]
and $F$ satisfies \ref{Pq}, then so does $F'$; moreover if $F''$ and $F'$ satisfy \ref{Pq}, then so does $F$.
By  \cite[I, Corollaire 3.9]{IlDRW} we are reduced to show that $\Omega^{q-1}/Z_n\Omega^{q-1}$ and  
$\Omega^{q}/B_n\Omega^{q}$ satisfy \ref{Pq}, for all $n$ (with the notation from {\em loc. cit.}).
Since we have an injection 
$\Omega^{q-1}/Z_1\Omega^{q-1}\inj  \Omega^q$ (via $d$) and  isomorphisms induced by the  Cartier operator
\eq{lem:levelq0}{Z_{n-1}\Omega^{q-1}/Z_n\Omega^{q-1}\cong Z_{n-2}\Omega^{q-1}/Z_{n-1}\Omega^{q-1}\cong \ldots \cong \Omega^{q-1}/Z_1\Omega^{q-1},}
the quotient $\Omega^{q-1}/Z_n\Omega^{q-1}$ satisfies \ref{Pq}.
Next we consider $\Omega^{q}/B_n\Omega^{q}$. Let $t:=t_1$ and set $A_0:= A/tA$. Up to replacing $X$ by a Nisnevich neighborhood of 
$\Spec A_0\subset X$, we can assume that we have a $k$-algebra morphism $\sigma: A_0\to A$ 
which composed with the quotient map $A\to A_0$ is the identity on $A_0$, e.g. \cite[Lemma 7.14]{BRS}. 
Note that $A_0[t]\to A$ is \'etale and  induces an identification $\hat{A}:=\varprojlim_n A/t^nA\cong A_0[[t]]$.
As   $B_n\Omega^q_X$ has the structure of a coherent locally free $\sO_X$-submodule of $F_{X*}^{n}\Omega^q_X$
we have 
\[(\Omega^{q}/B_n\Omega^{q})(X)= \Omega^q_A/B_n\Omega^q_A=:M_A.\] 
Let $a\in M_A$  be a nonzero element. We have to show there exists  a map $\varphi: A\to C$  with $C$ 
a smooth $k$-algebra of dimension $\le e-1$ such that $\varphi(a)\neq 0$ in  $M_C$. 
Since $M_A$ is a projective $A$-module of finite rank (with module structure induced by the $p^n$-power map)
we find that the natural map 
\[M_A\to \lim_{n} M/t^n M \cong M_{\hat{A}}\]
is injective, e.g., \cite[Theroem 8.9]{Matsumura}. Denote by $\hat{a}$ the image of $a$ in $M_{\hat{A}}$, which is nonzero.
It suffices to show:
\begin{claim}\label{lem:levelq-claim}
Either the image of $\hat{a}$ in $M_{\hat{A}/(t-\lambda)}$  is nonzero, for some $\lambda\in k$,
or  there exists a map $\varphi_0:A_0\to C_0$, with $C_0$ smooth of dimension $\dim C_0\le e-2$, such that the image of $\hat{a}$ in 
$M_{C_0[[t]]}$ is nonzero.
\end{claim}
 Indeed, in this case $C:= A/(t-\lambda)$ (resp. $C:=A\otimes_{A_0} C_0$) is smooth of dimension $\le e-1$ 
and the $t$-adic  completion of the natural map $\varphi: A\to  C$ is the morphism $\hat{A}\to \hat{A}/(t-\lambda)$ 
(resp. $A_0[[t]]\to C_0[[t]]$), hence $\varphi(a)\neq 0$ in $M_C$.   

We prove the claim. Note that we can write any element $\Omega^q_{A_0[[t]]}$ uniquely as a sum
\eq{lem:levelq1}{\sum_{l=0}^\infty b_l t^l+ \sum_{l=1}^\infty c_l t^{l} d\log t, \quad b_l\in \Omega^q_{A_0},\, c_l\in \Omega^{q-1}_{A_0}.}
We characterize those elements which are in $B_n\Omega^q_{A_0[[t]]}$.
To this end, note that $B_n\Omega^q_{A_0[[t]]}= F^{n-1}dW_n\Omega^{q-1}_{A_0[[t]]}$, by \cite[I, Proposition 3.11]{IlDRW}.
By \cite[Theorem B]{GH06} each element in $W_n\Omega^{q-1}_{A_0[[t]]}$ can uniquely be written as an infinite sum of elements of the form
\[b_0 [t]^i, \quad c_0 [t]^i d\log [t], \,(i\ge 1), \quad  V^s(b_s [t]^j), \quad dV^s(c_s [t]^j), \]
where $i\ge 0$, $j\ge 1$ with $(j,p)=1$, $s\in \{1, \ldots, n-1\}$, and $b_r\in W_{n-r}\Omega^q_{A_0}$, $c_r\in W_{n-r}\Omega^{q-1}_{A_0}$. 
It follows that an element \eqref{lem:levelq1} lies in $B_n\Omega^q_{A_0[[t]]}$ if and only if the following conditions are satisfied:
\begin{enumerate}
    \item[1st case:]\label{lem-levelq-1st} $l= p^r l_0$ with $r\in [0, n-2]$ and $(l_0,p)=1$. Then
               \[\exists\, \gamma\in W_{r+1}\Omega^{q-1}_{A_0} \text{ such that }  b_l= F^rd(\gamma) \text{ and } c_l= l_0 F^r(\gamma).\]
    \item[2nd case:]\label{lem-levelq-2nd} $l=p^{n-1} l_0$ with  $(l_0,p)=1$. Then 
     \[\exists\, \gamma\in W_n\Omega^q_{A} \text{ and } \delta\in B_{n-1}\Omega^{q-1}_{A_0} \text{ such that } b_l= F^{n-1}d(\gamma) \text{ and }
     c_l= l_0 F^{n-1}(\gamma)+\delta.\]
     \item[3rd case:]\label{lem-levelq-3rd} $p^n| l$  (including the case $l=0$). Then 
     \[b_l\in B_n\Omega^q_{A_0} \text{ and } c_l\in B_n\Omega^{q-1}_{A_0} \quad (\text{where } c_l=0, \text{ if } l=0).\]
\end{enumerate}
Let $\tilde{a}\in \Omega^q_{A_0[[t]]}$ be a representative of  $\hat{a}\in M_{A_0[[t]]}$ and write $\tilde{a}$ in the form \eqref{lem:levelq1}.
As $\hat{a}\neq 0$, we find  $l_1\ge 0$  minimal with the property 
\[b_{l_1} t^{l_1}+ c_{l_1} t^{l_1} dlog(t)\not\in B_n\Omega^q_{A_0[[t]]}.\] 
First assume $p^n|l_1$. Then either $b_{l_1}\not \in B_n\Omega^q_{A_0}$ or $c_{l_1}\not\in B_n\Omega^{q-1}_{A_0}$.
If $b_{l_1}\not \in B_n\Omega^q_{A_0}$, then taking $\lambda=1$ in Claim \ref{lem:levelq-claim} will work;
if $c_{l_1}\not\in B_n\Omega^{q-1}_{A_0}$ we find by induction over $q$ a map $\varphi_0: A_0\to C_0$ as in the claim such that 
$\varphi(c_{l_1})\not\in B_n\Omega^{q-1}_{C_0}$, whence $\hat{a}$ does not map to $0$ in $M_{C_0[[t]]}$.

Now assume $l_1= p^rl_0$ with $r\in [0, n-2]$ and $(l_0,p)=1$. If $b_{l_1}\not\in B_{r+1}\Omega^q_{A_0}$, again taking
$\lambda=1$ in the claim will work. Else there exists a $\gamma\in W_{r+1}\Omega^{q-1}_{A_0}$, such that
$b_{l_1}=F^rd(\gamma)$. If $\gamma'\in W_{r+1}\Omega^{q-1}_{A_0}$ has the same property, then $\gamma-\gamma'\in FW_{r+2}\Omega^{q-1}_{A_0}$,
by \cite[I, (3.21.1.1)]{IlDRW}. Since $Z_{r+1}\Omega^{q-1}= F^{r+1} W_{r+2}\Omega^{q-1}_{A_0}$, 
see \cite[I, (3.11.3)]{IlDRW}, we deduce from
the 1st case above that we have 
\[c_{l_1}-l_0 F^r(\gamma)\neq 0 \quad \text{in } \Omega^{q-1}_{A_0}/ Z_{r+1}\Omega^{q-1}_{A_0}.\]
By \ref{Pq} for $\Omega^{q-1}/Z_{r+1}\Omega^{q-1}$, proven above,
we find a map $\varphi_0 : A_0\to C_0$ as in the claim such that,
$\varphi_0(c_{l_1}-l_0 F^r(\gamma))\neq 0$ in $\Omega^{q-1}_{C_0}/ Z_{r+1}\Omega^{q-1}_{C_0}$.
Thus the image of $\hat{a}$ in $M_{C_0[[t]]}$ is nonzero.

Finally assume $l_1=p^{n-1} l_0$ with  $(l_0,p)=1$. In the case $b_{l_1}\not\in B_n\Omega^q_{A_0}$ we take $\lambda=1$ as above.
Otherwise  we can argue as in the 2nd case, as $B_{n-1}\Omega^{q-1}_{A_0}\subset Z_n\Omega^{q-1}_{A_0}$.
This completes the proof of Claim \ref{lem:levelq-claim} and hence of the lemma.
\end{proof}
\begin{rmk}\label{rmk:levelq}
As the above proof shows, the statement of  Lemma \ref{lem:levelq} also holds for $\Omega^q/B_n\Omega^q$ instead of $W_n\Omega^q$, for $n\ge 1$.
But note that the same statement does not hold for $\Omega^q/Z_n\Omega^q$, as the latter sheaf vanishes on all smooth $k$-schemes of dimension $\le q$.
\end{rmk}

\subsection{The proof of \ref{c3}}
In the following, we let $L\in \Phi$ and denote by $\sO_L$ its ring of integers with maximal ideal $\fm_L$ and 
we let $z\in \fm_L$ be a fixed local parameter.

\begin{lem}\label{lem:fillog-explicit}
Let $K\to \sO_L$ be a coefficient field. Then modulo $W_n\sO_L$ 
 any element of $\fil^{\log}_r W_n L$, $r\ge 1$, 
is a sum of elements
\[ V^j([\mu z^i])\quad \text{with}\quad \mu\in K,\, 0\le j\le n-1,\, -r\le ip^{n-j-1}<0.\]
\end{lem}
\begin{proof}
This is immediate for $n=1$. For $n\ge 2$ we can write
$a=[a_0]+V(b)$, where $a_0\in L$ with $p^{n-1}v_p(a_0)\ge -r$ and $b\in \fil^{\log}_rW_{n-1}L$.
We find $\mu_ i\in K$ and $c\in \sO_L$ such that  $a_0=\sum_{0>i p^{n-1}\ge -r} \mu_i z^i+c$.
Note that 
\[[a_0]- \sum_{0>i p^{n-1}\ge -r} [\mu_i z^i] -[c]\in 
(\fil^{\log}_r W_n\sO_L\cap VW_{n-1}L) = V\fil^{\log}_rW_{n-1}L.\]
Hence 
\[a= \sum_{0>i p^{n-1}\ge -r} [\mu_i z^i] \quad \text{mod } V\fil^{\log}_{r}W_{n-1}(L)+ W_n\sO_L.\]
The statement follows by induction over $n$.
\end{proof}

\begin{lem}\label{lem:fil-V1K}
For $q\ge 1$ denote by $\sV^1K^M_q(\sO_L)$ the image under the  map 
\[(1+\fm_L)\otimes_{\Z} K^M_{q-1}(\sO_L)\to K^M_q(\sO_L),\quad (1+b)\otimes u\mapsto\{1+b,u\}.\]
The multiplication \ref{para:fil}\ref{para:fil5} induces  maps 
\eq{lem:fil-V1K1}{\fil^{\log}_{r} W_n(L)\otimes_{\Z} \sV^1K^M_q(\sO_L)\lra 
\fil_{r-1}^{\log} W_n\Omega^q_L+ d(\fil_{r-1}^{\log} W_n\Omega^{q-1}_L)}
and 
\ml{lem:fil-V1K2}{\fil^p_r W_n(L)\otimes_{\Z} \sV^1K^M_q(\sO_L)\\
\lra \Fil^p_{r-1}W_n\Omega^q_L
           + \sum_{s=0}^{n-1} p^s \left(V^{n-s-1}(\fil^{\log'}_{(r-1)p^s}W_{s+1}\Omega^q_L)+
           dV^{n-s-1}(\fil^{\log'}_{(r-1)p^s}W_{s+1}\Omega^{q-1}_L)\right).}
\end{lem}
\begin{proof}
It suffices to consider the case $q=1$, cf. \ref{para:fil}\ref{para:fil5}.
For \eqref{lem:fil-V1K1} we have to show
\[V^j([x])\dlog(1-b)\in \left(\fil_{r-1}^{\log} W_n\Omega^1_L+ d(\fil_{r-1}^{\log} W_n(L))\right), \]
for $x\in L$ with  $p^{n-1-j}v_L(x)\ge -r $, and $b\in \fm_L$. 
In view of $V^j(y)\dlog z=V^j(y\dlog z)$ and  \eqref{lem:fil-FVR1} it suffices to consider the case $j=0$.
By \cite[Lemma 7.13]{RS21} we have in $W_n\Omega^q_{\hat{\sO}_L}$, where $\hat{\sO}_L=\varprojlim_s \sO_L/\fm_L^s$,
\[\dlog(1-b)=-\sum_{i\ge 0} [b]^i d[b]- \sum_{j=1}^{n-1} \sum_{(i,p)=1} \frac{1}{i} dV^j([b]^i).\]
Thus it remains to show for  $x\in L$ with $p^{n-1}v_L(x)\ge -r$ and $0\neq b\in\fm_L$
\begin{enumerate}[label=(\alph*)]
\item\label{lem:fil-V1Ka} $[x] [b]^id[b]\in \fil^{\log}_{r-1}W_n\Omega^1_L$, for all $i\ge 0$;
\item\label{lem:fil-V1Kb} $[x]dV^j([b]^i)\in \fil^{\log}_{r-1}W_n\Omega^1_L+d(\fil^{\log}_{r-1}W_n(L))$, 
for all $0\le j\le n-1$, and all $i\ge 1$ with $(i,p)=1$.
\end{enumerate}
For \ref{lem:fil-V1Ka} it suffices to consider $i=0$, in which case the claim is immediate, as
$[x]d[b]=[xb]\dlog b$ and $p^{n-1}v_L(xb)\ge -r+ p^{n-1}\ge -r+1$.
For \ref{lem:fil-V1Kb} we compute (using the Leibniz rule)
\[[x]dV^j([b]^i)= dV^j([x]^{p^j}[b]^i)- V^j([b]^i[x]^{p^j})\dlog x.\]
Since $p^{n-j-1}(p^jv_L(x)+i)\ge -r+1$ it lies in the right hand side of \eqref{lem:fil-V1K1}. 

We show \eqref{lem:fil-V1K2}. 
By Definition \ref{defn:p-sat} and Lemma \ref{lem:filp-round}
\[\fil^p_r W_n(L)=\sum_{s=0}^{n-1}p^s\left(\fil^{\log}_{(r-1)p^s}W_n(L) +V^{n-m_s}(\fil^{\log}_{r p^s}W_{m_s}(L))\right),\]
where $m_s=\min\{v_p(rp^s), n\}=\min\{v_p(r), n-s\}+s$.
By \eqref{lem:fil-V1K1} it suffices to consider $p^sV^{n-m_s}(\fil^{\log}_{p^s r}W_{m_s}(L))$ and 
by \eqref{lem:fil-FVR2} we  may assume $s=0$. We can furthermore assume  $m:=m_0\ge 1$. It remains to show
that $V^j([y])\dlog(1-b)$ lies in the right hand side of \eqref{lem:fil-V1K2}
for $n-m\le j\le n-1$, $y\in L$ with $p^{n-j-1}v_L(y)=-r$, and $b\in \fm_L$.
By  Corollary \ref{cor:Filp-FVR} we may as before assume $j=0$ (and hence $m=n$) and consider the elements in the  cases 
\ref{lem:fil-V1Ka} and \ref{lem:fil-V1Kb} from above with $x$ replaced by $y$ 
with $p^{n-1}v_L(y)=-r$. Note $v_p(r)\ge n$ hence $p|v_L(y)=:e$. 
In case \ref{lem:fil-V1Ka} (with $i=0$) write 
$b= z^c u$ and $y=z^e v$, with $u,v\in \sO_L^\times$ and $c\ge 1$.
We obtain 
\[[y]d[b]= c[uv] [z]^{e+c}\dlog z+ [uv] [z]^{e+c}\dlog u.\]
It follows that, if $n\ge 2$ or $c\ge 2$, then 
\[[y]d[b]\in \fil^{\log}_{r-2}W_n\Omega^1_L\subset \Fil^p_{r-1}W_n\Omega^1_L.\]
On the other hand if  $n=1=c$, then 
\[[y]d[b]= ydb= \tfrac{1}{e+1} d(uv z^{e+1})-  uvz^{e+1}(\tfrac{1}{e+1}\dlog(uv)-\dlog u),\]
which lies in $\fil^{\log'}_{r-1}W_{1}\Omega^1_L+ d(\fil^{\log'}_{r-1}W_{1}(L))$, 
which lies in the right hand side of \eqref{lem:fil-V1K2}.

Finally we consider \ref{lem:fil-V1Kb} for $y\in L$ as above.
We have 
\[[y]dV^j([b]^i)= dV^j([y]^{p^j}[b]^i)-  V^j([b]^i[y]^{p^j})\dlog y.\]
If $j< n-1$ or $i\ge 2$, then it follows
\[ [y]dV^j([b]^i)\in \fil^{\log}_{r-2}W_n\Omega^1_L+ d(\fil^{\log}_{r-2}W_n(L))\subset \Fil^p_{r-1}W_n\Omega^1_L.\]
If $j=n-1$ and $i=1$, then we write $y=z^ev$ with $v\in \sO_L^\times$ and obtain 
\[[y]dV^{n-1}([b]^i)= dV^{n-1}([y]^{p^{n-1}}[b])- V^{n-1}([b][y]^{p^{n-1}})\dlog v,\]
where we use that $p$ divides $e$ and hence $e\cdot V^{n-1}(W_{1}\sO_L)=0$. 
Thus $[y]dV^{n-1}([b])$ lies in the right hand side of \eqref{lem:fil-V1K2}, in this case as well.
\end{proof}

\begin{lem}\label{lem:fil-vs-G}
Let $K\inj \sO_L$ be a coefficient field.
For $r\ge 1$ set
\[E:=\fil^{\log'}_r W_n\Omega^1_L+ W_n\Omega^1_{\sO_L}(\log).\]
The following two subgroups of $W_n\Omega^1_L$ are equal:
\begin{enumerate}[label=(\arabic*)]
    \item $\fil^{\log}_r W_n\Omega^1_L+d(\fil^{\log}_r W_n(L))$.
    \item The subgroup $G$  generated by $E$  and the elements
    \[\mu F^e(da),\quad a\in \fil^{\log}_r W_{n+e}(L), \, \mu \in W_n(K), \, e\ge 0.\]
\end{enumerate}
\end{lem}
\begin{proof}
Clearly, $d(\fil^{\log}_r W_n(L))\subset G$. We show $\fil^{\log}_r W_n\Omega^1_L\subset G$.
Since $W_n\Omega^1_{\sO_L}(\log)$ is contained in both groups it suffices by Lemma \ref{lem:fillog-explicit} to show 
\[V^j([\nu z^i])\dlog b\in G, \quad \text{for}\quad 0\le j\le n-1,\, b\in L^\times,\, \nu\in K^\times, \, 0>ip^{n-j-1}\ge -r.\]
Write $b=z^c u$ with $u\in \sO_L^\times$ and $c\ge 0$ and $i=p^e i_0$ with $(i_0,p)=1$ and $e\ge 0$.
We obtain
\[V^j([\nu z^i])\dlog b= \tfrac{c}{i_0}V^j([\nu]) F^e d(V^j([z]^{i_0}))+ V^j([\nu z^i])\dlog u\in G,\]
where we use 
\[[z]^i \dlog z^c= F^e([z]^{i_0}\dlog z^c)= \tfrac{c}{i_0} F^e d([z]^{i_0})= \tfrac{c}{i_0} F^{e+j}dV^j([z]^{i_0}).\]

We show the other inclusion $G\subset \fil^{\log}_r W_n\Omega^1_L+d(\fil^{\log}_r W_n(L))$. 
By \ref{para:fil}\ref{para:fil5} and Lemma \ref{lem:fillog-explicit} it suffices to show
\eq{lem:fil-vs-G1}{\mu F^e dV^j([\nu z^i])\in \fil^{\log}_r W_n\Omega^1_L+d(\fil^{\log}_r W_n(L)), }
for $\mu\in W_n(K)$,  $e\ge 0$, $0\le j\le n+e-1$, $\nu\in K^\times$, and    $ip^{n+e-j-1}\ge -r$.
There are two cases. First assume $e<j$. Then
\[\mu F^e dV^j([\nu z^i])=\mu dV^{j-e}([\nu z^i])= dV^{j-e}(F^{j-e}(\mu)[\nu z^i]) - V^{j-e}([\nu z^i] F^{j-e}d\mu).\]
Hence \eqref{lem:fil-vs-G1} holds in this case by \eqref{lem:pres-filp1}. 
Now assume $e\ge j$. We obtain
\[\mu F^e dV^j([\nu z^i])= \mu [\nu z^i]^{p^{e-j}} \dlog \nu + i\mu [\nu z^i]^{p^{e-j}} \dlog z\]
which satisfies \eqref{lem:fil-vs-G1} by Lemma \ref{lem:pres-filp}.
\end{proof}

\begin{lem}\label{lem:gr}
Let $K\inj \sO_L$ be a coefficient field and  $r\ge 2$. The following elements are contained in $\Fil^p_r W_n\Omega^q_L$:
\eq{lem:gr1}{a\dlog u,\qquad a\in \fil_r^p W_n(L),\, u \in K^M_q(K),}
\eq{lem:gr2}{da\dlog v,\qquad a\in \fil_r^p W_n(L),\, v \in K^M_{q-1}(K), }
\eq{lem:gr3}{(F^e db)\cdot \gamma, \qquad b\in \fil^p_r W_{n+e}(L), \, \gamma \in W_n\Omega^{q-1}_K, \, e\ge 0.}
Furthermore  the abelian group
\[\frac{\Fil^p_r W_n\Omega^q_L}{\Fil^p_{r-1} W_n\Omega^q_L}\]
is generated by the elements \eqref{lem:gr1} and \eqref{lem:gr2} and those elements \eqref{lem:gr3}, which have 
$b\in \fil^{\log, p}_{r-1}W_{n+e}(L)$.
\end{lem}
\begin{proof}
As we have by definition $p^s \Fil^p_{p^sr}W_n\Omega^q_L\subset \Fil^p_r W_n\Omega^q_L$, it suffices to prove the first statement 
with $\fil^p_r$ in \eqref{lem:gr1} - \eqref{lem:gr3} replaced by $\fil_r$. In this case we find
$\eqref{lem:gr1}, \eqref{lem:gr2}\in \Fil_r W_n\Omega^q_L$. It remains to consider \eqref{lem:gr3} for which 
it suffices by the surjectivity of \eqref{para:DRW5} and by \ref{para:fil}\ref{para:fil5}  to show 
\begin{enumerate}[label=(\alph*)]
    \item\label{lem:gra} $F^e(db)d\mu\in \Fil^p_r W_n\Omega^2_L$,
    \item\label{lem:grb} $\mu F^e(db)\in \Fil^p_r W_n\Omega^1_L$,
\end{enumerate}
where  $b\in \fil_r W_{n+e}(L)$, $\mu \in W_n(K)$, and $e\ge 0$.
In case \ref{lem:gra} we have 
\[F^e(db)d\mu= F^ed(bdV^e(\mu)).\]
As $bdV^e(\mu)\in \Fil^p_rW_{n+e}\Omega^1_L$, by Lemma \ref{lem:pres-filp}, the compatibility of $\Fil^p_r$ with $F$ and $d$, 
see Corollary \ref{cor:Filp-FVR}, yields the claim in this case. Similarly we have $F^e(db)\in \Fil^p_r W_n\Omega^1_L$, since 
$\Fil^p_r W_n\Omega^1_L$ is an $W_n\sO_L$-module by Corollary \ref{cor:fil-WO-mod}, we get \ref{lem:grb}.

It remains to show the second statement on the generators of the quotient. We assume $q\ge 1$ as there is nothing to show for $q=0$.
We define the following subgroups of $W_n\Omega^q_L$
\[\sV^1\fil^p_rW_n\Omega^q_L:=\Im\left( \fil^p_rW_n(L)\otimes \sV^1K^M_q(\sO_L)\to W_n\Omega^q_L\right),\]
 see Lemma \ref{lem:fil-V1K}  for the definition of $\sV^1 K^M_q(\sO_L)$, and 
\[ \sV^1\Fil^p_rW_n\Omega^q_L:= \sV^1\fil^p_rW_n\Omega^q_L+ d(\sV^1\fil^p_rW_n\Omega^{q-1}_L).\]
Furthermore denote by $G_0^q\subset W_n\Omega^q_L$ the subgroup generated by all the elements \eqref{lem:gr1} and \eqref{lem:gr2}.
We claim
\eq{lem:gr-claim}{\sV^1\Fil^p_rW_n\Omega^q_L\subset G_0^q+\Fil^p_{r-1}W_n\Omega^q_L.}
By Lemma \ref{lem:fil-V1K} and the fact that 
$d(G_0^{q-1}+\Fil^p_{r-1}W_n\Omega^{q-1}_L)\subset G_0^q+\Fil^p_{r-1}W_n\Omega^q_L$ 
it suffices  to show 
\[p^s V^{n-s-1}(\fil^{\log'}_{(r-1)p^s} W_{s+1}\Omega^q_L)\subset G^q_0+\Fil^p_{r-1}W_n\Omega^q_L.\]
By \eqref{lem:fil-FVR1.5} we can assume $s=0$. 
By the decomposition  
\eq{lem:gr4}{\sO_L^\times=K^\times\cdot(1+\fm_L)}
every element in the left hand side (for $s=0$) can be written as a sum of elements
\begin{enumerate}[label=(\alph*)]
    \item\label{lem:gri} $V^{n-1}(a\dlog u)=V^{n-1}(a)\dlog u$, 
    \item\label{lem:grii} $V^{n-1}(a\dlog \beta)=V^{n-1}(a)\dlog\beta$,
\end{enumerate}
where $a\in \fil^{\log}_{r-1}W_n(L)$, $u\in K^M_q(K)$, and $\beta\in \sV^1K^M_{q}(\sO_L)$.
The elements \ref{lem:gri} are clearly in $G_0^q$, and the elements \ref{lem:grii} 
are in $\fil^{\log}_{r-2}W_n\Omega^q_L+d(\fil^{\log}_{r-2}W_n\Omega^q_L)\subset \Fil^p_{r-1}W_n\Omega^q_L$, by  \eqref{lem:fil-V1K1}.
This shows Claim \eqref{lem:gr-claim}.

As any element in $K^M_q(L)/\sV^1K^M_q(\sO_L)$ can be represented by a sum of elements
$u\in K^M_q(K)$ and $\{z,v\}$ with $v\in K^M_{q-1}(K)$ we find by \eqref{lem:gr-claim} and the definition of $G_0^q$ 
a surjection
\[H^q_0+ d(H^{q-1}_0)\surj\frac{\Fil^p_r W_n\Omega^q_L}{G_0^q+\Fil^p_{r-1}W_n\Omega^q_L},\]
where $H^q_0\subset W_n\Omega^q_L$ is the sugbroup generated  by elements
\[p^s a \dlog \{z,u\}\quad \text{for}\quad 0\le s\le n-1,\, a\in \fil^{\log}_{(r-1)p^s}W_n(L),\, u\in K^M_{q-1}(K).\]
Denote by $H^q_1$  the subgroup of $W_n\Omega^q_L$ generated by the elements \eqref{lem:gr3}, which have 
$b\in \fil^{\log, p}_{r-1}W_{n+e}(L)$. It remains to show 
\[H_0^q+ dH_0^{q-1}\subset R^q=:H^q_1+G_0^q+\Fil^p_{r-1}W_n\Omega^q_L.\]
As $d(R^{q-1})\subset R^q$ (for $q\ge 2$)
it suffices to show $H_0^q\subset R^q$, for all  $q\ge 1$.
Since the multiplication $W_n\Omega^1_L\otimes K^M_{q-1}(K)\to W_n\Omega^q_L$ induces a surjection
$H_0^1\otimes K^M_{q-1}(K)\surj H_0^q$ and a well-defined map $R^1\otimes K^M_{q-1}(K)\to R^q$ 
it suffices to consider the case $q=1$.
By \eqref{lem:fil-V1K1},  the decomposition \eqref{lem:gr4}, \eqref{lem:fil-FVR1.5}, and \ref{lem:Filp} we have 
\[p^s\fil^{\log'}_{(r-1)p^s}W_n\Omega^1_L + W_n\Omega^1_{\sO_L}(\log)\subset R^1,\quad \text{for all } s.\]
Hence the inclusion $H_0^1\subset R^1$ follows from Lemma \ref{lem:fil-vs-G}.
\end{proof}

\begin{lem}\label{lem:c3fil}
Let  $L'/ L$ be a finite extension with ramification degree $e$.  
Denote by $\Tr: W_n(L')\to W_n(L)$ the trace.
Then 
\[\Tr(\fil_rW_n (L'))\subset \fil_{\lceil\frac{r}{e}\rceil} W_n(L) \quad \text{and} \quad 
\Tr(\fil^p_rW_n (L'))\subset \fil^p_{\lceil\frac{r}{e}\rceil} W_n(L).\]
\end{lem}
\begin{proof}
By \eqref{lem:fil-FVR1.5} we have $\fil^p_{r} W_n(L')= \sum_{s=0}^{n-1} \ul{p}^s \fil_{r}W_{n-s}(L')$. Hence it suffices to prove the first
statement. This is proven in \cite[(7.5.1)]{RS21} for the $F$-saturated filtration $\fil^{F}_rW_n (L)$. 
We check that it also works without $F$-saturation.
Set $s=\lceil\frac{r}{e}\rceil$. By the same argument as in {\em loc. cit.} (below (7.5.1)) we find 
$\Tr(a)\in \fil^{\log}_s W_n(L)$.   Set $m:=\min\{v_p(s), n\}$. If $m=n$, then $\fil^{\log}_s W_n(L)=\fil_s W_n(L)$ and we are done.
If $m\le n-1$, then it follows directly from  the definition that we have an isomorphism
\eq{lem:c3fil1}{K\xr{\simeq} \frac{\fil^{\log}_s W_n(L)}{\fil_s W_n(L)}, \quad x\mapsto V^{n-m-1}([\tilde{x} z^{-s_0}]),}
where $s_0= s/p^m$ (which is prime to $p$) and $\tilde{x}\in \sO_L^\times$ is any lift of $x$.
The map
\[F^{n-1}d: W_n (L)\to \Omega^1_L, \quad (a_0, \ldots, a_{n-1})\mapsto \sum_{j=0}^{n-1} a_j^{p^{n-j-1}}\dlog a_j, \]
clearly induces maps
\[\fil^{\log}_s W_n(L)\to \fm_L^{-s}\Omega^1_{\sO_L}(\log) \quad \text{and} \quad \fil_s W_n(L)\to \fm_L^{-s}\Omega^1_{\sO_L}.\]
By \eqref{lem:c3fil1} the induced map  on the quotient is injective,
\[\frac{\fil^{\log}_s W_n(L)}{\fil_s W_n(L)}\longrightarrow \frac{\fm_L^{-s}\Omega^1_{\sO_L}(\log)}{\fm_L^{-s}\Omega^1_{\sO_L}}, \quad 
V^{n-m-1}(\tilde{x} z^{-s_0})\mapsto -s_0\tilde{x}^{p^m} z^{-s}\dlog z.\]
Thus the statement follows from $\fm^{s}_L\cdot F^{n-1}d(\Tr(a))\subset \Omega^1_{\sO_L}$, see \cite[(7.5.2)]{RS21}.
\end{proof}

Now we can complete the proof of Theorem \ref{thm:fil-cond}.

\begin{proof}[Proof of \ref{c3}.]
Let $L'/L$ be a finite extension with ramification index $e$ and denote by $\Tr: W_n\Omega^q_{L'}\to W_n\Omega^q_{L}$ the trace.
We have to show 
\eq{pfc31}{\Tr(\Fil^p_r W_n\Omega^q_{L'})\subset \Fil^p_{\lceil\frac{r}{e}\rceil} W_n\Omega^q_L.} 
As $\Tr$ restricts to $\Tr: W_n\Omega^q_{\sO_{L'}}\to W_n\Omega^q_{\sO_{L}}$ we can assume $r\ge 1$. 
We consider several cases.
In the following we will use that $\sO_L$ is excellent and that henceforth for every finite extension $L'/L$  we have $[L':L]= e(L'/L) f(L'/L)$.

{\em 0th case: We have finite field extensions $L'\supset E \supset L$,
and \eqref{pfc31} holds for $L'/E$ and $E/L$. Then it  holds for $L'/L$.} 
This follows from the transitivity of the trace and the formula $\lceil\lceil r/e_1\rceil /e_2\rceil=\lceil r/e_1e_2\rceil$.

{\em 1st case: $e=1$.}  In this case the local parameter $z\in \sO_L$ is local parameter of $\sO_{L'}$ as well 
and  hence \eqref{pfc31} follows from Lemma \ref{lem:pres-filp}\ref{lem:pres-filpII} and the fact that $\Tr$ commutes with  $V$ and $d$ 
(see \ref{para:DRW}) and satisfies a projection formula.

{\em 2nd case: $e=[L':L]$.} Let $K\inj \sO_L$ be a coefficient field, by assumption the composition
$K\inj \sO_L\inj \sO_{L'}$ is a coefficient field of $\sO_{L'}$. Let $t\in \sO_{L'}$ be a local parameter.
Any element in $W_n\Omega^q_{\sO_{L'}}(\log)/W_n\Omega^q_{\sO_{L'}}$ admits a representative
which is a sum of elements $\alpha \dlog t$ with $\alpha\in W_n\Omega^{q-1}_{K}$.
We have $\Tr(\alpha \dlog t)= \alpha \dlog (\Nm (t))$, see \ref{para:DRW}\ref{para:DRW4}. 
Thus in this case \eqref{pfc31} holds for $r=0,1$, by Lemma \ref{lem:Filp}.
For $r\ge 2$ condition \eqref{pfc31} follows  by induction over $r$, from the Lemmas \ref{lem:gr} and \ref{lem:c3fil}, and the fact that 
$\Tr$ commutes with $F$, $d$, and satisfies a projection formula.

{\em 3rd case: $L'/L$ purely inseparable.} In this case we can refine $L'/L$ into a tower of subextensions 
of degree $p$. Thus the statement follows from the above cases.

{\em 4th case: $L'/L$ separable.} In this case we use the  $p$-extension trick from \cite[5.9]{BaTa}.
Let  $H$ a $p$-Sylow subgroup in the Galois group of some finite Galois extension $M/L$ containing $L'$.
Set $E:= M^{H}$. Then every finite extension of $E$  inside $M$ has  $p$-power degree and and $([E:L],p)=1$.
Since $L'/L$ is separable we obtain decomposition $L'\otimes_L E= \oplus_i E_i$ with $E_i/E$ a finite field extension of $p$-power degree. 
We obtain a commutative diagram 
\eq{prop:fil-sc2}{\xymatrix{
W_n\Omega^q_{L'}\ar[d]_{\Tr_{L'/L}}\ar[r]^-{\oplus_i \varphi_i} & 
\bigoplus_i W_n\Omega^q_{E_i}\ar[d]^{\sum_{i} \Tr_{E_i/E}}\\
W_n\Omega^q_{L}\ar[r]^{\varphi} & W_n\Omega^q_E. 
}}
By the cases 0 - 2, \eqref{pfc31} holds for  $E/L$ since it is tamely ramified;
similarly \eqref{pfc31} holds for $E_i/E$ for all $i$, since a $p$-power extension can be refined
to a tower of degree $p$-extensions. 
Now let $a\in \Fil_{r}^p W_n\Omega^q_{L'}$,
by definition, or  \ref{c6}, we have
\[c_{E_i}(\varphi_{i}(a))\le e(E_i/L')\cdot r,\]
by \ref{c3} for $E_i/E$ and $e(L'/L) e(E_i/L')= e(E/L) e(E_i/E)$ we have 
\[c_E(\Tr_{E_i/E}(\varphi_i(a)))\le \bigg\lceil  \frac{e(E/L) \cdot r}{e(L'/L)}\bigg\rceil,\]
by \ref{c3} for $E/L$ we have 
\eq{prop:fil-sc3}{c_L(\Tr_{E/L}\Tr_{E_i/E}(\varphi_i(a)))\le \bigg\lceil \frac{e(E/L)\cdot r}{e(L'/L) e(E/L)}\bigg\rceil 
=\bigg\lceil \frac{r}{e(L'/L)}\bigg\rceil.}
Thus altogether
\begin{align*}
 c_L(\Tr_{L'/L}(a))              &=c_L\left([E:L]\cdot \Tr_{L'/L}(a)\right), && \text{since } [E:L]\in \Z_p^\times,\\
                                      & = c_L\left(\Tr_{E/L} \varphi(\Tr_{L'/L}(a))\right), && \text{by }\Tr_{E/L}\circ\varphi=[E:L],\\
                                      & = c_L\left(\sum_i \Tr_{E/L} \Tr_{E_i/E} \varphi_i(a)\right), && \text{by } \eqref{prop:fil-sc2}\\
                                      & \le \max \left\{c_L\left(\Tr_{E/L} \Tr_{E_i/E} \varphi_i(a)\right)\right\}, && \text{by \ref{c2}}\\
                                      &\le \bigg\lceil \frac{r}{e(L'/L) }\bigg\rceil, && \text{by }\eqref{prop:fil-sc3}.
\end{align*}

{\em 5th case:} For an arbitrary finite field extension $L'/L$ property \ref{c3} follows from the cases 0, 3, and  4 above. 
This completes the proof of Theorem \ref{thm:fil-cond}.
\end{proof}

\section{Some more properties of the {$p$}-saturated filtration}\label{sec:propFil}

For later use we analyze in this section the $p$-saturated filtration further.

\begin{lem}\label{lem:char}
Let $L\in \Phi$  and $r\ge 2$. The  map 
\eq{lem:char1}{F^{n-1}\oplus F^{n-1}d: \frac{\Fil_r^p W_n\Omega^q_L}{\Fil^p_{r-1}W_n\Omega^q_L+\ul{p}(\Fil_r^pW_{n-1}\Omega^q_L)}\lra 
\frac{\Omega^q_L}{\Fil^p_{r-1}\Omega^q_L} \oplus \frac{\Omega^{q+1}_L}{\Fil^p_{r-1}\Omega^{q+1}_L}}
is injective.
\end{lem}
\begin{proof}
First note that the map in the statement is well-defined by Corollary \ref{cor:Filp-FVR}. Let $K\inj \sO_L$ be a coefficient field.
Let $\phi$ be an element in the source of \eqref{lem:char1} and assume
\eq{lem:char2}{(F^{n-1}(\phi), F^{n-1}d\phi)=0 \quad \text{in} \quad 
\frac{\Omega^q_L}{\Fil^p_{r-1}\Omega^q_L} \oplus \frac{\Omega^{q+1}_L}{\Fil^p_{r-1}\Omega^{q+1}_L}.}
By Lemma \ref{lem:gr} we can assume that $\phi$ is a sum
\eq{lem:char2.5}{\phi=\sum_i a_i\dlog u_i +\sum_j da'_j\dlog v_j +\sum_l (F^{c_l}db_l)\cdot \delta_l,}
with
$$a_i, a'_j\in \gr_rW_nL:=\frac{\fil_rW_nL}{\fil_{r-1}W_nL}, 
\qquad 
b_l\in \frac{\fil^{\log}_{r-1}W_{n+c_l}L}{\fil_{r-1}W_{n+c_l}L \cap \fil^{\log}_{r-1}W_{n+c_l}L},$$ 
$$u_i\in K^M_q(K), \quad 
v_j\in K^M_{q-1}(K), \quad 
\delta_l\in W_n\Omega^{q-1}_K, \quad 
c_l\ge 0.$$
Write 
\[e=v_p(r),\quad  e_1=v_p(r-1),\quad r=r_0p^e, \quad r-1=r_1p^{e_1}.\]
We consider four cases. 

\medskip

\noindent {\em 1st case:} $e=0$, $e_1\in [0,n-1]$.
By \cite[7.18(1)]{RS21} we can assume 
\[a_i=V^{n-1-e_1}([\alpha_i z^{-r_1}]),\quad a'_j=V^{n-1-e_1}([\alpha'_j z^{-r_1}]),\quad b_l=V^{n+c_l-1-e_1}([\beta_l z^{-r_1}]),\] 
for some $\alpha_i,\alpha'_j,\beta_l\in K$. Hence
\[\phi=V^{n-1-e_1}([z]^{-r_1}A)+ dV^{n-1-e_1}([z]^{-r_1}B),\]
where 
\[A:=\sum_i [\alpha_i]\dlog u_i- \sum_l [\beta_l] F^{n-1-e_1}d\delta_l\,\,\in W_{e_1+1}\Omega^q_K,\]
\[B:=\sum_j[\alpha_j']\dlog v_j +\sum_l [\beta_l]F^{n-1-e_1}(\delta_l)\,\,\in W_{e_1+1}\Omega^{q-1}_K.\]
Therefore
\begin{align*}
F^{n-1}d\phi
&=
z^{-(r-1)} F^{e_1}dA +(-r_1)z^{-(r-1)}\dlog z \cdot F^{e_1}(A),
\\
F^{n-1}(\phi)
&=
z^{-(r-1)}\left(p^{n-1-{e_1}} F^{e_1}(A)+ F^{e_1}d B\right) +(-r_1)z^{-(r-1)}\dlog z \cdot F^{e_1}(B).
\end{align*}
If $e_1=0$, then
$$\Fil^p_{r-1}\Omega^q_L=\frac{1}{z^{r-2}}\Omega^q_{\sO_L}(\log z)$$
and the vanishing \eqref{lem:char2} yields $A=0$ and $B=0$, whence $\phi=0$.
If $e_1\in [1,n-1]$, then 
$$\Fil^p_{r-1}\Omega^q_L=\frac{1}{z^{r-1}}\Omega^q_{\sO_L}$$
and the vanishing \eqref{lem:char2} yields $F^{e_1}(A)=0$ and $F^{e_1}(B)=0$.
By \cite[I, (3.11.3)]{IlDRW} we find $A'\in W_{e_1}\Omega^{q}_K$ and $B'\in W_{e_1}\Omega^{q-1}_K$ with 
$V(A')=A$ and $V(B')=B$. 
Therefore
$$\phi=V^{n-e_1}([z^{-r_1p}] A')+dV^{n-e_1}([z^{-r_1p}] B').$$
Thus $\phi\in \Fil^p_{r-1}W_n\Omega^q_L$, by Lemma \ref{lem:pres-filp}.

\medskip

\noindent {\em 2nd case:} $e=0$, $e_1\ge n$. 
In this case $\gr_rW_nL=0 =\gr_rW_{n+c_l}L$, for all $c_l$ with $e_1\ge n+c_l$, by \cite[Cor. 7.18(1)]{RS21}.
Hence we can assume  
$\phi=\sum_l (F^{c_l} db_l)\cdot \delta_l$,  with $e_1<n+c_l$ for all $l$, 
and $b_l=V^{n+c_l-1-e_1}([\beta_l z^{-r_1}])$, with $\beta_l\in K$.
Hence
\[\phi= [z]^{-\frac{r-1}{p^{n-1}}} A+ [z]^{-\frac{r-1}{p^{n-1}}} \dlog z\cdot B,\]
where
\[A= \sum_l F^{e_1-(n-1)}d[\beta_l]\cdot \delta_l\,\, \in W_n\Omega^q_K\quad \text{and} \quad 
B=\sum_l F^{e_1-(n-1)}([\beta_l])\delta_l\,\,\in W_n\Omega^{q-1}_K.\]
By Lemma \ref{lem:pres-filp} we have 
\[[z]^{-\frac{r-1}{p^{n-1}}} A\in p \cdot\fil^{\log}_{r-1}W_n\Omega^q_L + d(\fil^{\log'}_r W_n\Omega^{q-1}_L)\subset \Fil^p_{r-1} W_n\Omega^q_L.\]
Thus we can assume $A=0$ and hence
\[F^{n-1}d\phi= - z^{-(r-1)}\dlog z\cdot F^{n-1}dB \quad \text{and} \quad F^{n-1}\phi= z^{-(r-1)}\dlog z \cdot F^{n-1}(B).\]
As in the  case under consideration
$$\Fil^p_{r-1}\Omega^q_L=\frac{1}{z^{r-1}}\Omega^q_{\sO_L}$$
the vanishing \eqref{lem:char2} yields
\eq{lem:char3}{B\in \Ker( F^{n-1})\cap \Ker(F^{n-1}d)= pW_n\Omega^{q-1}_K,}
see e.g. \cite[(0.6.3)]{Ekedahl} for the equality\footnote{It is easily deduced from \cite[I, Proposition 3.11]{IlDRW}}.
By Lemma \ref{lem:pres-filp}  and Corollary \ref{cor:Filp-FVR} we find
\[\phi \in p\cdot \Fil^p_rW_n\Omega^q_L\subset p\Fil_{pr}^pW_n\Omega^q_L=\ul{p}(\Fil^p_rW_{n-1}\Omega^q_L).\]

\medskip

\noindent {\em 3rd case:} $e\in [1,n-1]$, $e_1=0$.
By \cite[Cor. 7.18(2) and Cor 7.17(1)]{RS21}, we can assume
\[a_i=V^{n-1}(\alpha_i z^{-(r-1)})+V^{n-e}(\beta_i[z]^{-r_0p}), \quad 
a'_j=V^{n-1}(\alpha'_j z^{-(r-1)})+V^{n-e}(\beta'_j[z]^{-r_0p}),\] 
and 
\[b_l=V^{n+c_l-1}(\gamma_l z^{-(r-1)}),\] 
for some $\alpha_i,\alpha'_j,\gamma_l\in K$, $\beta_i, \beta'_j\in W_eK$.
Hence
\[\phi= V^{n-1}(z^{-(r-1)} A)+ dV^{n-1}(z^{-(r-1)}B) +V^{n-e}([z]^{-r_0p} C)+ dV^{n-e}([z]^{-r_0p} D),\]
where
\[A:=\sum_i \alpha_i\dlog u_i-\sum_l\gamma_l F^{n-1}d\delta_l\,\, \in \Omega^q_K,\]
\[B:= \sum_j\alpha_j'\dlog v_j + \sum_l \gamma_l F^{n-1}(\delta_l)\,\,\in \Omega^{q-1}_K,\]
\[C:= \sum_i\beta_i\dlog u_i\,\,\in W_e\Omega^q_K,\qquad D:=\sum_j \beta_j'\dlog v_j\,\,\in W_e\Omega^{q-1}_K.\]
Therefore
\[F^{n-1}d\phi= z^{-(r-1)}dA -(r-1) z^{-(r-1)}\dlog z\cdot A + z^{-r} F^{e-1}dC,\]
\[F^{n-1}\phi= z^{-(r-1)} dB -(r-1) z^{-(r-1)}\dlog z \cdot B + z^{-r} F^{e-1}dD.\]
As in the  case under consideration
$$\Fil^p_{r-1}\Omega^q_L=\frac{1}{z^{r-2}}\Omega^q_{\sO_L}(\log z)$$
the vanishing \eqref{lem:char2} yields
\[A=0,\quad B=0,\quad  F^{e-1}d C=0,\quad F^{e-1}d D=0.\]
By \cite[(3.11.4)]{IlDRW} we  find $C'\in W_{e+1}\Omega^q_K$ and $D'\in W_{e+1}\Omega^{q-1}_K$ such that 
\[C= F(C')\quad \text{and}\quad D=F(D').\]
Altogether 
\[\phi= p V^{n-e-1}([z]^{-r_0} C') +p dV^{n-e-1}([z]^{-r_0} D').\]
Hence $\phi\in p \Fil^p_{r+1}W_n\Omega^q_L$, by Lemma \ref{lem:pres-filp}.
As $r+1\le pr$, Corollary \ref{cor:Filp-FVR} yields
\[\phi\in p\Fil^p_{pr}W_n\Omega^q_L= \ul{p} (\Fil^p_r W_{n-1}\Omega^q_L).\] 

\medskip

\noindent {\em 4th case:} $e\ge n, e_1=0$.    
By \cite[Cor. 7.18(3) and Cor. 7.17(1)]{RS21}, we can assume
\[a_i=V^{n-1}(\alpha_iz^{-(r-1)})+\beta_i[z]^{-\frac{r}{p^{n-1}}},\quad
a'_j=V^{n-1}(\alpha'_jz^{-(r-1)})+\beta'_j[z]^{-\frac{r}{p^{n-1}}},\]
and 
\[\beta_l=V^{n+c_l-1}(\gamma_l z^{-(r-1)}),\] 
for some $\alpha_i$, $\alpha'_j$, $\gamma_l\in K$, $\beta_i$, $\beta_j\in W_nK$.
Similarly as in the 3rd case we find
\[\phi= V^{n-1}(z^{-(r-1)} A)+ dV^{n-1}(z^{-(r-1)}B)+ [z]^{-\frac{r}{p^{n-1}}} C + p [z]^{-\frac{r}{p^{n-1}}} \dlog z\cdot D,\]
where $A\in \Omega^{q}_K$, $B\in \Omega^{q-1}_K$, $C\in W_n\Omega^q_K$, and $D\in W_n\Omega^{q-1}_K$ and 
where we use that $r/p^{n-1}$ is divisible by $p$ to get the $p$ in front of the last summand.
Hence
\[ F^{n-1} d\phi= -(r-1)z^{-(r-1)}\dlog z\cdot A + z^{-{(r-1)}}dA + z^{-r} F^{n-1}(dC),\]
\[ F^{n-1} \phi =  -(r-1)z^{-(r-1)}\dlog z\cdot B+ z^{-(r-1)}dB + z^{-r} F^{n-1}(C).\]
As in the  case under consideration
$$\Fil^p_{r-1}\Omega^q_L=\frac{1}{z^{r-2}}\Omega^q_{\sO_L}(\log z)$$
the vanishing \eqref{lem:char2} yields
\[ A=0,\quad B=0,\quad C\in \Ker(F^{n-1})\cap \Ker(F^{n-1}d).\]
Thus 
\[\phi\in p\Fil^p_{r+1}W_n\Omega^q_L\subset \ul{p}(\Fil^p_rW_{n-1}\Omega^q_L),\]
by \eqref{lem:char3}, Lemma \ref{lem:pres-filp}, and Corollary \ref{cor:Filp-FVR}.
This completes the proof.
\end{proof}

\begin{lem}\label{lem:gr-ex}
Let $L\in \Phi$ and $r\ge 2$. The sequences    
\eq{lem:gr-ex1}{F(W_{n+1}\Omega^q_L)\cap\Fil_r^p W_n\Omega^q_L\lra 
\frac{\Fil_r^p W_n\Omega^q_L}{\Fil^p_{r-1}W_n\Omega^q_L+\ul{p}(\Fil_r^pW_{n-1}\Omega^q_L)} \stackrel{F^{n-1}d}{\lra}
\frac{\Omega^{q+1}_L}{\Fil^p_{r-1}\Omega^{q+1}_L}}
and
\eq{lem:gr-ex2}{V(W_{n-1}\Omega^q_L)\cap\Fil_r^p W_n\Omega^q_L\lra 
\frac{\Fil_r^p W_n\Omega^q_L}{\Fil^p_{r-1}W_n\Omega^q_L+\ul{p}(\Fil_r^pW_{n-1}\Omega^q_L)} \stackrel{F^{n-1}}{\lra}
\frac{\Omega^{q}_L}{\Fil^p_{r-1}\Omega^{q}_L}}
are exact, the maps on the left  being induced by the quotient maps.
\end{lem}
\begin{proof}
We use the notation from the proof of Lemma \ref{lem:char}.
Clearly the sequences from the statement form a complex. 
Write 
\[e=v_p(r),\quad  e_1=v_p(r-1),\quad r=r_0p^e, \quad r-1=r_1p^{e_1}.\]

We show the exactness of \eqref{lem:gr-ex1}.
Let $\phi\in \Fil^p_rW_n\Omega^q_L$ be as in \eqref{lem:char2.5} and assume 
\[F^{n-1}d(\phi)=0 \quad \text{in } \frac{\Omega^{q+1}_L}{\Fil^p_{r-1}\Omega^{q+1}_L}.\]
We  consider the same  four cases as in Lemma \ref{lem:char} to show that the image of $\phi$ in the quotient in the middle of \eqref{lem:gr-ex1}
lies in the image of $F(W_{n+1}\Omega^q_L)$. 

\medskip

\noindent {\em 1st case:} $e=0$, $e_1\in [0,n-1]$.
With the notation from the corresponding case in Lemma \ref{lem:char} we find: 
if $e_1=0$, then $A=0$ and hence
\[\phi= dV^{n-1}([z]^{-r_1}B)= FdV^{n}([z]^{-r_1}B);\]
if $e_1\in[1, n-1]$, then 
\[A\in \Ker(F^{e_1})\cap \Ker(dF^{e_1})\subset p W_n\Omega^q_K,\]
see  \eqref{lem:char3} for the inclusion, hence writing $A=pA'=FV(A')$ yields 
\[\phi=F\Big(V^{n-e_1}([z]^{-r_1}A')+dV^{n-e_1}([z]^{-r_1}B)\Big).\]

\medskip

\noindent {\em 2nd case:} $e=0$, $e_1\ge n$. 
The argument from the second case of the proof of Lemma \ref{lem:char} yields 
\[\phi=[z]^{-\frac{r-1}{p^{n-1}}}\dlog z\cdot B, \quad \text{with } F^{n-1}d(B)=0.\]
By \cite[I, (3.11.4)]{IlDRW} we have $B=F(B')$ and hence 
\[\phi=F([z]^{-\frac{r-1}{p^n}}\dlog z\cdot B').\]

\medskip

\noindent {\em 3rd case:} $e\in [1,n-1]$, $e_1=0$.
As in the third case of  the proof of Lemma \ref{lem:char} we find $A=0$ and $C=F(C')$.
Thus
\[\phi= F\Big(dV^n(z^{-(r-1)}B)+ V^{n-e}([z]^{-r_0}C')+dV^{n-e}([z]^{-r_0}V(D))\Big).\]

\medskip

\noindent {\em 4th case:} $e\ge n, e_1=0$.    
Similarly as in the fourth case of Lemma \ref{lem:char}, we get $A=0$ and $C=F(C')$.
Thus
\[\phi= F\Big(dV^n(z^{-(r-1)}B)+ [z]^{-\frac{r}{p^n}}C'+ [z]^{-\frac{r}{p^n}}\dlog z \cdot V(D)\Big).\]
This completes the proof of the exactness of \eqref{lem:gr-ex1}.

We show the exactness of \eqref{lem:gr-ex2}. We assume $n\ge 2$ as the statement is trivial for $n=1$.
Let $\phi\in \Fil^p_rW_n\Omega^q_L$ be as in \eqref{lem:char2.5} and assume 
\[F^{n-1}(\phi)=0 \quad \text{in } \frac{\Omega^{q}_L}{\Fil^p_{r-1}\Omega^{q}_L}.\]
We consider the same  four cases as above to show that the image of $\phi$ in the quotient 
lies in the image of $V(W_{n-1}\Omega^q_L)$ under the quotient map. 

\medskip

\noindent {\em 1st case:} $e=0$, $e_1\in [0,n-1]$.
With the notation from the corresponding case in Lemma \ref{lem:char} we find: 
if $e_1=0$, then $B=0$ and hence
\[\phi= V^{n-1}([z]^{-r_1}A);\]
if $e_1\in[1, n-1]$, then  $F^{e_1}(B)=0$, hence  $B=V(B')$, by \cite[I, (3.11.3)]{IlDRW}.
As
\[dV^{n-e_1-1}([z^{-r_1}B])= dV^{n-e_1}([z]^{-r_1 p}B')\in \Fil^p_{r-1}W_n\Omega^q_L,\]
by Lemma \ref{lem:pres-filp}, we obtain that the image of $\phi$ in the quotient is equal to
\[\phi=V^{n-1-e_1}([z]^{-r_1}A).\]

\medskip

\noindent {\em 2nd case:} $e=0$, $e_1\ge n$. 
The argument from the second case of the proof of Lemma \ref{lem:char} yields 
\[\phi=[z]^{-\frac{r-1}{p^{n-1}}}\dlog z\cdot B, \quad \text{with } F^{n-1}(B)=0.\]
By \cite[I, (3.11.3)]{IlDRW} we have $B=V(B')$ and hence 
\[\phi=V([z]^{-\frac{(r-1)p}{p^n}}\dlog z\cdot B').\]

\medskip

\noindent {\em 3rd case:} $e\in [1,n-1]$, $e_1=0$.
As in the third case of  the proof of Lemma \ref{lem:char} we find $B=0$ and $F^{e-1}dD=0$.
Thus $D=F(D')$ and
\[\phi= V^{n-1}(z^{-(r-1)}A) +V^{n-e}([z]^{-r_0p}C)+ V(FdV^{n-e-1}([z]^{-r_0}D')).\]

\medskip

\noindent {\em 4th case:} $e\ge n, e_1=0$.    
Similarly as in the fourth case of Lemma \ref{lem:char}, we get $B=0$ and $C=V(C')$.
Thus
\[\phi= V^{n-1}(z^{-(r-1)}A) + V([z]^{-\frac{pr}{p^{n-1}}} C')+ V F([z]^{-\frac{r}{p^{n-1}}}\dlog z\cdot D).\]
This completes the proof of the exactness of \eqref{lem:gr-ex2}.
\end{proof}

\section{Hodge-Witt sheaves with modulus}\label{sec:HWM}
We begin this section by recalling the notion of reciprocity sheaves and several a priori different ways to assign a modulus to a section
of certain presheaves, which appear in \cite{KSY-RSCII}, \cite{RS21}, and  \cite{RS-AS}. The main result of this section
is Theorem \ref{thm:HW-modulus}, which says that these different notions of modulus agree for $W_n\Omega^q$, $q\ge 1$.
Theorem \ref{thm:fil-cond} and Lemma \ref{lem:char} are crucial ingredients in the proof of that theorem.

\begin{para}\label{para:RSC}
We recall the definition of a {\em reciprocity sheaf} from \cite[Definition 2.2.4 and Section 2.4]{KSY-RSCII},
where these are called sheaves with SC-reciprocity.
Denote by $\PST$ the category of presheaves with transfers on $\Sm$ in the sense of Voevodsky and by 
$\NST$ its full subcategory of Nisnevich sheaves.
For $U\in \Sm$ we denote by $\Ztr(U)$ the presheaf with transfers represented by $U$. 
Let $F\in \PST$ and let $(X,D)$ be a proper modulus pair (see \ref{para:mod-pairs}) with $X\setminus D=U$.
A section $a\in F(U)$ {\em has  modulus $(X,D)$} if the map $a: \Ztr(U)\to F$, defined by $a$ via the Yoneda embedding,
factors in $\PST$ as
\[\xymatrix{\Ztr(U)\ar[r]\ar[d] & F\\  h_0(X,D)\ar@{.>}[ur]}.\]
Here $h_0(X,D)$ is the presheaf with transfers given on $S\in\Sm$ by
\[h_0(X,D)(S)= \Coker(\Cor((\P^1_S,\infty_S), (X,D))\xr{i_0^*-i_1^*} \Ztr(U)(S)),\]
where $i_\epsilon: \{\epsilon\}_S\inj \P^1_S$ denotes the closed immersion, for $\epsilon\in \{0,1\}$, 
and $\Cor((\P^1_S ,\infty_S), (X,D))$ denotes the free abelian group generated on integral closed subschemes $Z\subset \A^1_S\times_k U$
which are surjective and finite over a connected component of $U$ such that 
\[\nu^*\infty_{S}\ge \nu^* D,\]
where $\nu : \widetilde{Z}\to \P^1_S\times_k X$ is the normalization of the closure of $Z$.

For $(X,D)$ any modulus pair with $U=X\setminus D$ and $(\ol{X}, D'+ X_\infty)$ a compactification of $(X,D)$ (see \ref{para:mod-pairs}) we set
\[\ul{\omega}^{\CI}F(X,D)=\{a\in F(U)\mid a \text{ has modulus } (\ol{X}, D'+ N \cdot X_\infty)\text{ for some } N\gg 0\}.\]
In fact the assignment   $(X,D)\mapsto \ul{\omega}^{\CI}F(X,D)$  defines a cube invariant semipure modulus presheaf with M-reciprocity 
$\ul{\omega}^\CI F$ as defined in \cite{KSY-RSCII} and \cite{Saito-Purity} and with the notation from there
we have $\ul{\omega}^{\CI}F=\tau_! h_0^{\bcube}\omega^* F$. 

A presheaf with transfers $F$ is a {\em reciprocity presheaf} if for any $U\in \Sm$  any section $a\in F(U)$ has a modulus, 
i.e., if $F(U)=\bigcup \ul{\omega}^{\CI}F(X,D)$, where the union is over all proper modulus pairs $(X,D)$ 
with $X\setminus D=U$. 
A {\em reciprocity sheaf} is a Nisnevich sheaf on $\Sm$ which is also a reciprocity presheaf. 
The  category of reciprocity sheaves is denoted by $\RSC_{\Nis}$.

If $F\in \RSC_{\Nis}$ and $(X,D)$ is a modulus pair, then the assignment
\[(\text{\'etale $X$-schemes})\ni (u:V\to X)\mapsto \ul{\omega}^{\CI}F(V, u^*D)=:\Gamma(V, \ul{\omega}^{\CI}F_{(X,D)})\]
is a  sheaf on $X_{\Nis}$ by \cite[Corollary 4.16]{RS21} and is denoted by $\ul{\omega}^{\CI}F_{(X,D)}$.
\end{para}

\begin{para}\label{para:Fc}
Let $F\in \NST$ and let $c=\{c_L: F(L)\to \N_0\}_{L\in\Phi}$ be a semi-continuous conductor (see \ref{para:cond}).
Let $(X,D)$ be a {\em proper} modulus pair with $U=X\setminus D$. For $a\in F(U)$ we write
\[c_X(a)\le D\]
as a shorthand for 
\[c_L(\rho^*a)\le v_L(\rho^*D),\quad \text{for all } L\in \Phi \text{ and  all } \rho:\Spec L\to U.\]
Following \cite[4.8]{RS21} we define for any modulus pair $(X,D)$ with a compactification $(\ol{X}, D'+ X_\infty)$
\[\tF_c(X,D):=\{a\in F(U)\mid c_{\ol{X}}(a)\le D'+ N\cdot X_\infty, \text{ for some } N\gg 0\}.\]
By \cite[(4.8.1)]{RS21} the assignment
\[(\text{\'etale $X$-schemes})\ni (u:V\to X)\mapsto \tF_c(V, u^*D)=:\Gamma(V, \tF_{c,(X,D)})\]
defines a Nisnevich sheaf $\tF_{c, (X,D)}$ on $X_\Nis$.
By \cite[Theorem 4.15(4)]{RS21} we have
\eq{para:Fc1}{\tF_c\subset \ul{\omega}^{\CI}F, \qquad \text{for }F\in \RSC_{\Nis}.}
For $L\in \Phi$ with ring of integers $\sO_L$, maximal ideal $\fm_L$, and $r\ge 1$ set
\[\tF_c(\sO_L, \fm_L^{-r})= \varinjlim_{(U,Z)} \tF_{c}(U,rZ),\]
where the filtered colimit is over all modulus pairs $(U,Z)$ with $U$ and $Z$ smooth and connected 
such that $\sO_L= \sO_{U,z}^h$, with $z\in Z$ is the generic point. As $c$ is assumed to be semi-continuous we have 
\eq{para:Fc2}{\tF_c(\sO_L,\fm_L^{-r})=\{a\in F(L)\mid c_L(a)\le r\},}
see \cite[Lemma 4.23]{RS21}. Finally we recall that if the conductor has level $m\le \infty$ 
and $(X,D)$ is a modulus pair with $U=X\setminus D$, then 
\[\tF_c(X,D)=\{a\in F(U)\mid f^*a\in \tF_c(S, f^*D), \text{for all } f: S\to X \text{ with } S\in \Sm,\, \dim S\le m\},\]
see \cite[Corollary 4.18]{RS21}.
\end{para}

\begin{para}\label{para:AS}
The following notion of modulus using dilatations is motivated by work of Abbes and (Takeshi) Saito, who used this approach
to study the ramification of Galois torsors, see \cite{AbbesSaito11} and \cite{TakeshiSaito}. 
In the framework of reciprocity sheaves this notion was studied in \cite{RS-AS}.
Let $(X,D)$ be a modulus pair with $X\in \Sm$ and assume that the reduced divisor $D_{\red}$ underlying $D$ has only simple normal crossings,
i.e., $D_{\red}$ is an SNCD divisor. Denote by ${\rm Bl}_D(X\times X)$ the blow-up in $D$ diagonally embedded into $X\times X$ 
and denote by $\widetilde{X\times D}$ and $\widetilde{D\times X}$ the strict transforms of $X\times D$ and $D\times X$, respectively.
Set 
\[P^{(D)}_X= {\rm Bl}_D(X\times X)\setminus (\widetilde{X\times D}\cup \widetilde{D\times X}).\]
We note that $P^{(D)}_X$ is smooth and that the open embedding $U\times U\inj X\times X$, with $U=X\setminus D$, induces
an embedding $U\times U\inj P^{(D)}_X$, see e.g. \cite[Lemma 2.3]{RS-AS}.
We denote by $p_1$, $p_2: P^{(D)}_X\to X$ the maps induced by the projection $X\times X\to X$ to the first and second factor, respectively;
they restrict to the projection maps on $U\times U$. For $F\in \RSC_{\Nis}$ we define
\eq{para:AS0}{F^{\AS}(X,D)=\{a\in F(U)\mid p_1^*a-p_2^*a \in F(P^{(D)}_X)\},}
for this definition to make sense we use that $F(U\times U)$ is a subgroup of $F(P^{(D)}_X)$ 
by \cite[Theorem 6]{KSY} together with \cite[Corollary 3.2.3]{KSY-RSCII}.
As an \'etale map $u: V\to X$ induces a morphism $P_V^{(u^*D)}\to P_X^{(D)}$ we obtain 
a Nisnevich sheaf $F^{\AS}_{(X,D)}$ on $X_\Nis$ given by
\[(u:V\to X)\mapsto F^{\AS}(V, u^*D)=\Gamma(V, F^{\AS}_{(V, u^*D)}).\]
By \cite[Theorem 2.6]{RS-AS}
\eq{para:AS1}{\ul{\omega}^{\CI}F_{(X,D)}\subset F^{\AS}_{(X,D)},}
for $X\in \Sm$ and $D_{\red}$ SNCD. By Theorem 2.10 in {\em loc. cit.} this is an equality if $(X,D)$ additionally admits 
a projective SNC-compactification.
 \end{para}
 
\begin{thm}\label{thm:HW-modulus}
Let $c=\{c_L:W_n\Omega^q_L\to \N_0\}_{L\in \Phi}$ with
\eq{thm:HW-modulus1}{c_L(a)=\min\left\{r\ge 0\mid a\in \Fil^p_r W_n\Omega^q_L\right\}, \qquad a\in W_n\Omega^q_L,}
be the  semi-continuous conductor of level $q+1$ from Theorem \ref{thm:fil-cond}.
Let $(X,D)$ be any modulus pair and let $n$, $q\ge 1$. Then
\[\ul{\omega}^{\CI}W_n\Omega^q_{(X,D)}= \widetilde{W_n\Omega^q}_{c, (X,D)}.\]
In particular
\[\ul{\omega}^{\CI}W_n\Omega^q(\sO_L, \fm_L^{-r})= \Fil^p_r W_n\Omega^q_L,\]
for all $L\in \Phi$ and $r\ge 0$.
If furthermore $X$ is smooth and $D_{\red}$ is  an SNCD, then 
\[\ul{\omega}^{\CI}W_n\Omega^q_{(X,D)}= \widetilde{W_n\Omega^q}_{c, (X,D)}= (W_n\Omega^q)^{\AS}_{(X,D)}.\]
\end{thm}

The proof will be given after the next remark.

\begin{rmk}\label{rmk:HW-modulus}
\begin{enumerate}
\item Note that the statement is not true for $q=0$. As the absolute Frobenius on $X$ induces an endomorphism of 
$W_n\sO_X$ one has to consider the Frobenius saturation of $\fil_r W_n(L)$, see \cite[Theorem 7.20]{RS21} and \cite[Proposition 5.3]{RS-AS}.
\item The case $n=1$ and $q\ge 1$ follows from \cite[Theorem 6.6 and Corollary 6.8]{RS-AS}. 
In particular this says that if $X$ and $D$ are smooth and $r\ge 0$ then
\eq{rmk:HW-modulus1}{(\Omega^q)^{\AS}(X,rD)=
\begin{cases}
\Omega^q_X(\log D)((r-1)D) & \text{if } p\nmid r\\
\Omega^q_X(rD) & \text{if } p|r,
\end{cases}}
which at the generic point $\eta$ of $D$  coincides with $\Fil^p_rW_1\Omega^q_{L_\eta}$, where $L_\eta=\Frac(\sO_{X,\eta}^h)$.
\item In the case $D_\red$ is an SNCD,  Theorem \ref{thm:HW-modulus} and Lemma \ref{lem:Filp} imply
\eq{rmk:HW-modulus2}{\ul{\omega}^{\CI}W_n\Omega^q_{(X,D_\red)}=W_n\Omega^q_X(\log D),}
where the right hand side is the logarithmic de Rham-Witt sheaf, see, e.g., \cite[Proposition-Definition 3.10]{Matsuue}.
(This $\supset$ inclusion is immediate, this $\subset$ inclusion  can be, for example, deduced from the isomorphism 
\cite[I, (2.14.8)]{IlDRW} together with  \cite[Corollary 4.4]{Matsuue}.)
Note that the above equality is also true for $q=0$ as in this case both sides are equal to $W_n\sO_X$.
As a consequence we obtain the equality
\[\sL og(\ul{\omega}^{\CI}W_n\Omega^q)(X,D_{\red})=W_n\Omega^q_X(\log D),\]
where $\sL og: \uMNST_{\rm log}\to \textbf{Shv}^{\rm ltr}_{\rm dNis}$ is the functor defined in
\cite[(6.0.2)]{Saito-Log}. This gives a new proof of \cite[Theorem 4.4]{Merici}.
\end{enumerate}
\end{rmk}

\begin{proof}
First of all we note that by \cite[Theorem 4.15]{RS21} and \eqref{para:Fc2} the first statement of Theorem \ref{thm:HW-modulus} 
is equivalent to show 
\[\ul{\omega}^{\CI}W_n\Omega^q(\sO_L,\fm_L^{-r})= \Fil^p_r W_n\Omega^q_L, \quad \text{for all }L\in \Phi.\]
Furthermore, if $X$ is smooth and $D_{\red}$ is an SNCD, then $P^{(D)}_X$ (see \ref{para:AS}) is smooth  and 
$W_n\Omega^q_{P^{(D)}_X}$ is a successive  extension of locally free sheaves (see \cite[I, Corollaire 3.9]{IlDRW}) and hence the question
whether an element $a\in W_n\Omega^q(U)$ lies in $(W_n\Omega^q)^{\AS}(X,D)$ 
is Nisnevich local around the 1-codimensional points of $D$. Thus we are reduced to the following situation:

Let $X=\Spec A$, with $A$ a smooth $k$-algebra, let $D=\div(z)$ be a smooth connected divisor on $X$ with generic point $\eta\in D^{(0)}$,  
and set $\sO_L= \sO_{X,\eta}^h$ and $L=\Frac(\sO_{X,\eta}^h)$.
Then we have to show 
\[\Fil^p_r W_n\Omega^q_L= \ul{\omega}^{\CI}W_n\Omega^q(\sO_L,\fm_L^{-r})= \left((W_n\Omega^q)^{\AS}_{(X,rD)}\right)^h_{\eta},\quad 
\text{for all }r\ge 0,\]
where the right hand side denotes the Nisnevich stalk of $(W_n\Omega^q)^{\AS}_{(X,rD)}$ in $\eta$.
By \eqref{para:Fc1} and \eqref{para:AS1} it remains to show
\[\left((W_n\Omega^q)^{\AS}_{(X,rD)}\right)^h_{\eta}\subset  \Fil^p_r W_n\Omega^q_L.\]
By \eqref{rmk:HW-modulus1} this is true for $n=1$.
Thus we assume $n\ge 2$ in the following. As  $\{\Fil^p_r W_n\Omega^q_L\}_{r \ge 0}$ is an exhaustive filtration of
$W_n\Omega^q_L$ it remains to show 
\eq{pfThm:HW1}{\Fil^p_r W_n\Omega^q_L\cap (W_n\Omega^q)^{\AS}(X,(r-1)D) \subset \Fil^p_{r-1} W_n\Omega^q_L,}
where we identify $(W_n\Omega^q)^{\AS}(X,(r-1)D)$ with its image in $W_n\Omega^q_L$.
By Lemma \ref{lem:Filp} this is true for $r=1$ and we can assume $r\ge 2$. 
For $\phi$ in the left hand side of \eqref{pfThm:HW1} we claim
\eq{claimphir}{\phi=0 \quad \text{in}\quad  \frac{\Fil_r^p W_n\Omega^q_L}{\Fil^p_{r-1}W_n\Omega^q_L+\ul{p}(\Fil_r^pW_{n-1}\Omega^q_L)}.}
Indeed,  by \cite[Theorem 6.6]{RS-AS} (see also Remark \ref{rmk:HW-modulus}) we have 
\[(\Omega^q)^{\AS}(X, (r-1)D)\cap \Omega^q_L= \Fil^p_{r-1}\Omega^q_L,\]
and since the formation $E\mapsto E^{\AS}$ is functorial with respect to any morphism of sheaves of abelian groups we find
\[ F^{n-1}\left((W_n\Omega^q)^{\AS}(X,(r-1)D)\right)\subset \Fil^p_{r-1}\Omega^q_L\]
and
\[F^{n-1}d\left((W_n\Omega^q)^{\AS}(X,(r-1)D)\right)\subset \Fil^p_{r-1}\Omega^{q+1}_L.\]
Thus the claim \eqref{claimphir} follows from Lemma \ref{lem:char}.

Therefore  (after possibly shrinking $X$ around $\eta$) we can write  $\phi=\psi+\ul{p}(\phi')$  with
\[\psi\in \Fil^p_{r-1}W_n\Omega^q_L\cap W_n\Omega^q_{A[\frac{1}{z}]} \quad \text{and} 
\quad \phi'\in \Fil_r^p W_{n-1}\Omega^q_L\cap W_{n-1}\Omega^q_{A[\frac{1}{z}]}.\]
By Theorem \ref{thm:fil-cond}, \eqref{para:Fc1}, and \eqref{para:AS1} the form $\psi$ lies in  $(W_n\Omega^q)^{\AS}(X, (r-1)D)$
and hence so does $\ul{p}(\phi')$.
Since $P_X^{((r-1)D)}\setminus U\times U$ is supported on a smooth divisor (e.g. \cite[Lemma 2.3]{RS-AS}),
the equality \eqref{para:AS0} together with \Cref{upreg} below and the injectivity of $\ul{p}$ 
imply that $\phi'\in (W_{n-1}\Omega^q)^{\AS}(X,(r-1)D)$. 
By induction on $n$ we thus have $\phi'\in \Fil^p_{r-1}W_{n-1}\Omega^q_L$.
This proves the containment \eqref{pfThm:HW1} and hence the theorem.
\end{proof}

\begin{lem}[{\cite[Lemma 6.7]{Gupta-Krishna}}]\label{upreg}
Let $C$ be a smooth $k$-algebra and $t\in C$ an element such that $C/tC$ is smooth. Then (see \eqref{para:DRW0} for the notation $\ul{p}$)
\[\underline p(W_{n-1}\Omega^q_{C[\frac{1}{t}]}) \cap W_n\Omega_C^q=\underline p(W_{n-1}\Omega_C^q)\quad \text{in } W_n\Omega^q_{C[\frac{1}{t}]}.\]
\end{lem}

\section{The structure of Hodge-Witt sheaves with modulus}\label{sec:strHWM}
In this section we investigate the structure of Hodge-Witt sheaves with modulus. The main result is Theorem \ref{thm:strHWM},
whose proof will occupy all of this section and 
which will be essential for the proof of the main statement in section \ref{sec:duality}.

Throughout this section we assume $X\in \Sm$ and we let $D$ be an effective Cartier divisor on $X$ such that $D_{\red}$ is an SNCD.

\begin{para}\label{para:div-not}
If $E$ is any divisor on $X$ and $m\in \Z\setminus\{0\}$, then we write
\[m \mid E :\Longleftrightarrow  m \text{ divides the multiplicity of every irreducible component of } E\]
and 
\[m\nmid E :\Longleftrightarrow m \text{ does not divide the multiplicity of any irreducible component of } E.\]

Given an increasing sequence of natural numbers $1\le r_1<\ldots< r_s$ we say 
\[E= E' + p^{r_1}E_1+\ldots + p^{r_s}E_s\]
is a {\em $p$-divisibility decomposition  of $E$ (with respect to $r_1<\ldots< r_s$)},
if $p^{r_1}\nmid E'$ and $p^{r_i}\nmid p^{r_{i-1}}E_{i-1}$, for $i=2,\ldots, s$, 
and $E'_{\red}+\sum_{i=1}^s E_{i,\red}$ is a reduced divisor.
Note that $p\mid E_s$ is allowed and that a $p$-divisibility decomposition of $E$ always exists 
and is uniquely determined by the sequence $r_1<\ldots< r_s$. 
For example if $E=\sum_i n_i \sE_i$ with $\sE_i$ the irreducible components of $E$, $\sE_i\neq\sE_j$, for $i\neq j$,
then to say that 
\[E=E'+ p^rE_r\]
is a $p$-divisibility decomposition of $E$ means that 
\[E'=\sum_{i, \, p^r\nmid n_i} n_i \sE_i\quad \text{and}\quad E_r=\sum_{i,\, p^r\mid n_i} \frac{n_i}{p^r} \sE_i.\]
Finally for $E$ and $m$ as above we set
\[\lceil E/m\rceil:= \sum_i \lceil n_i/m\rceil \sE_i \quad \text{and} \quad \lfloor E/m\rfloor:= \sum_i \lfloor n_i/m\rfloor \sE_i,\]
where $\lceil - \rceil$ (resp. $\lfloor -\rfloor$) denotes the  round-up (resp. round-down).
\end{para}

\begin{para}\label{para:HW-not}
Under the assumptions of this section all three ways to define Hodge-Witt sheaves with modulus discussed in Section \ref{sec:HWM} coincide
by Theorem \ref{thm:HW-modulus}. We therefore set
\eq{para:HW-not1}{W_n\Omega^q_{(X,D)}:= \ul{\omega}^{\CI}W_n\Omega^q_{(X,D)}, \quad \text{for } q\ge 1,}
which coincides with $\widetilde{W_n\Omega^q}_{c, (X,D)}$ and $(W_n\Omega^q)^{\AS}_{(X,D)}$.
For $q=0$ we set 
\eq{para:HW-not2}{W_n\Omega^0_{(X,D)}:= W_n\sO_{(X,D)}:= (\widetilde{W_n})_{c, (X,D)}}
in the notation of \ref{para:Fc} and where $c$ is the conductor defined in \eqref{thm:HW-modulus1} for $q=0$,
i.e., it is defined by  the  $p$-saturated Kato-Matsuda filtration $\fil^p_* W_n(L)$, see \ref{para:fil} and Definition \ref{defn:p-sat}.
By \cite[Theorem 7.20]{RS21} we have 
\eq{para:HW-not3}{\ul{\omega}^{\CI}W_n\sO_{(X,D)}=W_n\sO_{(X,D)}^F,}
where the upper index $F$ denotes the saturation under the absolute Frobenius acting on $W_n\sO_U$, with $U=X\setminus D$.
Note that $\ul{\omega}^{\CI}W_n\sO_{(X,D)}$ is not a $W_n\sO_X$-module.

Let $D=D_0+pD_1$ be a \pdd, see \ref{para:div-not}. 
Recall from \cite[Theorem 6.6]{RS-AS} (see also Theorem \ref{thm:HW-modulus}) that with the above notation
\eq{para:HW-not5}{\Omega^q_{(X,D)}= 
\begin{cases} 
\Omega^q_X(\log D_0)(D_0-D_{0,\red}+pD_1), & q\ge 1;\\
\sO_X(D_0-D_{0,\red}+pD_1),& q=0. 
\end{cases}}
For us the notation $\Omega^q_X(\log E)$ for a possibly non-reduced effective Cartier divisor with 
simple normal crossing support $E$ will always mean $\Omega^q_X(\log E_{\red})$.
\end{para}

We record: 
\begin{prop}\label{prop:strHWM}
The sheaves $W_n\Omega^q_{(X,D)}$ are coherent $W_n\sO_X$-modules, for all $q\ge 0$.
Furthermore, let  $a\in W_n\Omega^q(U)$, with $U=X\setminus D$. 
Then 
\[a\in \Gamma(X, W_n\Omega^q_{(X,D)})\, \Longleftrightarrow\, a_\eta\in \Fil^p_{m_\eta}W_n\Omega^q_{L_\eta},\quad \text{for all } \eta\in |D|^{(0)},\]
where $L_\eta=\Frac(\sO_{X,\eta}^h)$, $a_\eta\in W_n \Omega^q_{L_\eta}$ denotes the image of $a$, and $m_\eta$ is the multiplicity of $\eta$ in $D$.
\end{prop}
\begin{proof}
Denote by $W_n\sO_X(D)$ the invertible  $W_n\sO_X$-submodule of $j_*W_n\sO_U$, with $j:U\inj X$ the open immersion, 
which is uniquely characterized by the equality
\[W_n\sO_X(D)_{|U_0}=W_n\sO_{U_0}\cdot \frac{1}{[f]},\]
for any open $U_0\subset X$ with $D_{|U_0}=\Div(f)$, $f\in \sO(U_0)$. By Corollary \ref{cor:fil-WO-mod} 
we have the following inclusions of $W_n\sO_X$-modules
\eq{para:HW-not4}{W_n\Omega^q_X\subset W_n\Omega^q_{(X,D)}\subset W_n\Omega^q_X\otimes_{W_n\sO_X} W_n\sO_X(p^{n-1}D),\qquad \text{for all }q\ge 0.}
As $W_nX=(|X|, W_n\sO_X)$ is a scheme  of finite type over $\Spec W_n(k)$ and as $W_n\Omega^q_X$ is a coherent $W_n\sO_X$-module,
we see that $W_n\Omega^q_{(X,D)}$ is a coherent $W_n\sO_X$-module for all $q\ge 0$. 

The second statement follows for $q\ge 1$ from the equality $W_n\Omega^q_{(X,D)}= (W_n\Omega^q)^{\AS}_{(X,D)}$, 
see the beginning of the proof of Theorem  \ref{thm:HW-modulus}. For $q=0$, it is direct to check from the definition.
\end{proof}

The aim of this section is to prove the following theorem:

\begin{thm}\label{thm:strHWM}
For $n$, $q\ge 0$ set
\eq{thm:strHWM1}{Z_n\Omega^q_{(X,D)}:= F^n(W_{n+1}\Omega^q_{(X,D)}) \quad \text{and}\quad  B_n\Omega^q_{(X,D)}:= 
\begin{cases}
F^{n-1}d(W_n\Omega^{q-1}_{(X,D)}), &n\ge 1;\\
0, & n=0.
\end{cases}
}
Then $Z_n\Omega^q_{(X,D)}$ and $B_n\Omega^q_{(X,D)}$ are locally free $\sO_X$-submodules of $F^n_{X*}\Omega^q_{(X,D)}$, 
where $F_X:X\to X$ denotes the absolute Frobenius. Furthermore, there is  a short exact sequence of coherent $W_{n+1}\sO_X$-modules
\[0\lra B_n\Omega^{q+1}_{(X,D)}\lra \frac{W_{n+1}\Omega^q_{(X,D)}}{\ul{p}W_n\Omega^q_{(X,D)}}\xr{F^n} Z_n\Omega^q_{(X,D)}\lra 0,\]
where the first map on the left is given by $F^{n-1}d(\beta)\mapsto V(\beta)$ and the $W_{n+1}\sO_X$-module structure on $Z_n$ and $B_n$
is induced by $F^n:W_{n+1}\sO_X\to \sO_X$.
\end{thm}

The above theorem follows from Proposition \ref{prop:BZnn} and Proposition \ref{prop:HWM-ex} proven below.
In Proposition \ref{prop:BZnn} and Lemma \ref{lem:BnZn} we also give alternative descriptions of $B_n\Omega^q_{(X,D)}$ and $Z_n\Omega^q_{(X,D)}$, 
in the former using certain twisted Cartier operators, in the latter using the functors $\ul{\omega}^{\CI}$ and $(-)^{\AS}$ 
introduced in section \ref{sec:HWM}. 

\begin{rmk}\label{rmk:strHWM}
For $D=\emptyset$, the first part of the above theorem  holds by \cite[I, Proposition 3.11]{IlDRW} and \cite[0, Proposition 2.2.8]{IlDRW},
and the second part, the exactness of the sequence, holds by \cite[(0.6.2)]{Ekedahl}. 
Note that in this case there is also a second short exact sequence which is dual (under Grothendieck duality) to the one above, 
it describes $\Ker(R: W_{n+1}\Omega^q_X\to W_n\Omega^q_X)$ as an extension of locally free $\sO_X$-modules,
see \cite[I, (3.9.1)]{IlDRW}. For $D\neq\emptyset$ this latter sequence  does not seem to work, cf. Remark \ref{rmk:BZ-loc-free}.
\end{rmk}

\begin{lem}\label{lem:twisted-qis}
Let $A$ be an SNCD on $X$. Let $B$ be a (not necessarily effective) Cartier divisor whose support is contained in $A$ 
and let $A_0$  be an irreducible component of $A$, which is not a component of $|B|$. 
Let $s\in \Z$ and $\rho\in \N$ with $0\le \rho\le p-1$. Then the natural map
\[\Omega^\bullet_{X/k}(\log A)(B+psA_0)\lra \Omega^\bullet_{X/k}(\log A)(B+(ps+\rho)A_0)\]
is a quasi-isomorphism. Here $\Omega^q_X(\log E)(D):= \Omega^q_X(\log E)\otimes_{\sO_X} \sO_X(D)$. 
\end{lem}
\begin{proof}
By \cite[Lemma 4.1]{MS} (applied to $A=D_1+D_2$, $A_0=D_1$ and $D=B+A_0,\ldots, B+ \rho A_0$) we have a quasi-isomorphism
as in the statement for $s=0$. Applying $F_{X*}$ we obtain a quasi-isomorphism of $\sO_X$-modules; twisting with $\sO_X(sA_0)$ 
therefore yields the statement.
\end{proof}

\begin{lem}\label{lem:twisted-Cartier}
Let $E$ be an effective Cartier divisor on $X$ such that $D_{\red}+E_\red$ is a reduced SNCD.
Then the inverse Cartier operator 
\eq{lem:twisted-Cartier1}{C^{-1}: \Omega^q_X(\log D)\xr{\simeq} \sH^q(F_{X*}\Omega^\bullet_X(\log D)),}
see, e.g., \cite[Theorem (7.2)]{Katz}, induces an isomorphism 
\eq{lem:twisted-Cartier2}{
C^{-1}: \Omega^q_X(\log D)(\lceil D/p\rceil -\lceil D/p\rceil_{\red} +E)\xr{\simeq} \sH^q\Big(F_{X*}\Omega^\bullet_X(\log D)(D-D_{\red}+pE)\Big).}
\end{lem}
\begin{proof}
Observe
\[D-D_{\red}= p(\lceil D/p\rceil- \lceil D/p\rceil_{\red}) + \sum_i \rho_i D_i,\]
for certain $\rho_i\in \N$ with $0\le \rho_i\le p-1$,  where the $D_i$ are the irreducible components of $D$.
Indeed, for  $r\ge 1$ we have
$r-1=p(\lceil r/p\rceil-1)+ \rho$, for some  $0\le \rho\le p-1$.
Thus tensoring the isomorphism of $\sO_X$-modules \eqref{lem:twisted-Cartier1} with $\sO_X(\lceil D/p\rceil -\lceil D/p\rceil_{\red})$
and applying Lemma \ref{lem:twisted-qis} yields the  isomorphism \eqref{lem:twisted-Cartier2} in the case $E=0$.
Tensoring further with $\sO_X(E)$ yields the general case.
\end{proof}

\begin{para}\label{para:Omega-n}
Let $E$ be an effective Cartier divisor on $X$ such that $D_{\red}+E_\red$ is a reduced SNCD.
We will use the following notation
\eq{para:Omega-n1}{
\Omega^q_n(D,E):=\Omega^q_X(\log D)(\lceil D/p^n\rceil- \lceil D/p^n\rceil_{\red}+E),\qquad q\ge 0,\, n\ge 0.}
This is a locally free coherent $\sO_X$-module. 
We observe:
\begin{enumerate}[label= (\alph*)]
\item\label{para:Omega-n2}
Let $n\ge 1$.
If $D= D'+ p^n D_n$ is a {\pdd} (see \ref{para:div-not}), then
\[\Omega^q_n(D,E)= \Omega^q_X(\log (D'+D_n))(\lceil D'/p^n\rceil-\lceil D'/p^n\rceil_{\red}
 +D_n-D_{n,\red}+E).\]
Hence we have an inclusion
\eq{para:Omega-n3}{\Omega^q_n(D,E)\subset \Omega^q_n(D', D_n+E).}
Note that this inclusion is strict if $D_n\neq 0$ and $q<\dim X$
and that the cokernel is annihilated by $\sO_X(-D_{n,\red})$. In particular it is not an inclusion of vector bundles.
\item\label{para:Omega-n4}
Let $n\ge 0$.
There is a well-defined map $d:\Omega^q_n(D, pE)\to \Omega^{q+1}_n(D,pE)$ induced by the usual differential of the de Rham complex.
(Note that we need to have $pE$ instead of $E$ for this.) 
We therefore get a subcomplex $\Omega^\bullet_n(D, pE)\subset j_*\Omega^\bullet_U$, 
where $j:U=X\setminus |D+E|\inj X$ denotes the open immersion.
Given a {\pdd} $D=D'+p^{n+1}D_{n+1}$ we obtain an inclusion of complexes
\eq{para:Omega-n5}{\Omega^\bullet_n(D, pE)\subset \Omega^\bullet_n(D', p(D_{n+1}+E)).}
\item 
Let $n\ge 1$.
Using the equality 
\[\lceil D/p^n \rceil= \lceil \lceil D/p^{n-1}\rceil /p \rceil\]
we find that the inverse Cartier operator \eqref{lem:twisted-Cartier2} induces an isomorphism of $\sO_X$-modules
\eq{para:Omega-n6}{ C^{-1}: \Omega^q_n(D,E)\xr{\simeq} \sH^q(F_{X*}\Omega^\bullet_{n-1}(D,pE)).}
Let $n\ge 1$.
If $D=D'+p^nD_n$ is a \pdd, then we obtain a commutative diagram of locally free $\sO_X$-modules
\eq{para:Omega-n7}{
\xymatrix{
\Omega^q_n(D,E)\ar[r]^-{C^{-1}}_-{\simeq}\ar@{^(->}[d]   & \sH^q(F_{X*}\Omega^\bullet_{n-1}(D, pE))\ar@{^(->}[d]\\
\Omega^q_n(D', D_n+E)\ar[r]^-{C^{-1}}_-{\simeq} & \sH^q(F_{X*}\Omega^\bullet_{n-1}(D', p(D_n+E))).
}
}
Note that the injectivity of the vertical map on the right follows from the injectivity of the vertical map on the left.
\end{enumerate}
\end{para}


\begin{para}
For a complex $C^{\bullet}$ of $\sO_X$-modules on a scheme $X$, we denote 
\[\sZ^q(C^\bullet)= \Ker(d: C^q\to C^{q+1}), \quad \text{and}\quad \sB^q(C^\bullet)=\Im(d: C^{q-1}\to C^q),\]
which we view as $\sO_X$-modules.
\end{para}

\begin{para}\label{para:C-exs-eq-poles}
Let  $X$, $D$, and $E$ be as in \ref{para:Omega-n} above.
The Cartier isomorphism defines the following exact sequence of locally free $\sO_X$-modules
\eq{para:C-ex-seq-poles1}{ 0\to\sB^q(F_{X*}\Omega^\bullet_{n-1}(D,pE))\to \sZ^q(F_{X*}\Omega^\bullet_{n-1}(D,pE))\xr{C}
\Omega^q_n(D,E)\to 0.}
The fact that these modules are locally free follows via descending induction on $q$ by considering as well the exact sequence
\[0\to \sZ^q\to F_{X*}\Omega^q_{n-1}(D,pE)\xr{d}\sB^{q+1}\to 0.\]
We obtain a class  in the Ext-group which we denote by the symbol
\[\e_n^q(D,E)\in {\rm Ext}^1(\Omega^q_n(D,E), \sB^q(F_{X*}\Omega^\bullet_{n-1}(D,pE))).\]
\end{para}

\begin{para}\label{para:twisted-BZ}
Let  $X$, $D$, and $E$ be as in \ref{para:Omega-n} above. Let $n\ge 0$ and let 
\[D=D_0+pD_1+\ldots+ p^n D_n\]
be a \pdd.  Set
\eq{para:twisted-BZ1}{\ul{D}_j:= D_0+\ldots +p^j D_j,\qquad j=0,\ldots, n,}
and 
\eq{para:twisted-BZ2}{\ul{D}^0:=0 \quad\text{and}\quad \ul{D}^j:= pD_{n-j+1}+ p^2 D_{n-j+2}+\ldots+  p^jD_n, \quad j=1,\ldots, n.}
Thus we have $p\mid \ul{D}^j$, for all $j$, and 
\[D=\ul{D}_{n-j}+ p^{n-j+1}\left(\frac{1}{p}\ul{D}^j\right)\]
is a \pdd in the sense of \ref{para:div-not}. By abuse of notation we will also say
that $D=\ul{D}_{n-j}+ p^{n-j}\ul{D}^j$ is a \pdd, keeping in mind that there is an extra $p$ hidden in $\ul{D}^j$.
For $j=1,\ldots, n$, we will use the following notations repeatedly
\[\Omega^q_{j,n}:= F^j_{X*}\Omega^q_{n-j}(\ul{D}_{n-j}, \ul{D}^j+p^jE), \quad \sB^q_{j,n}:=\sB^q(\Omega^\bullet_{j,n}), \quad 
\sZ^q_{j,n}:=\sZ^q(\Omega^\bullet_{j,n}).\]
We define $\cO_X$-submodules
\[B_{j,n}^q(D,E)\subset Z^q_{j,n}(D,E)\subset 
\Omega^q_{j,n}\]
by setting 
\[B_{0,n}^q(D,E):=0 \subset 
\Omega^q_{0,n}:= \Omega^q_n(\ul{D}_n, \ul{D}^0+E)=:Z_{0,n}^q(D,E),\]
and for $j\ge 1$ recursively by the condition that the squares on the right in the following two diagrams are pullback squares
\eq{para:twisted-BZ3}{\xymatrix{
0\ar[r] & \sB^q_{j,n}\ar@{=}[d]\ar[r] &
B^q_{j,n}(D,E)\ar[r]\ar[d] &
B^q_{j-1,n}(D,E)\ar@{^(->}[d]\ar[r] & 0\\
0\ar[r] & \sB^q_{j,n}\ar@{=}[d]\ar[r] &
Z^q_{j,n}(D,E)\ar[r]\ar[d] &
Z^q_{j-1,n}(D,E)\ar@{^(->}[d]\ar[r] & 0\\
0\ar[r] & \sB^q_{j,n}\ar[r]&
\sZ^q_{j,n}\ar[r]^-C&
F^{j-1}_{X*}\Omega^q_{n-j+1}(\ul{D}_{n-j}, D_{n-j+1}+\ul{D}^{j-1}+p^{j-1}E)\ar[r] & 0,
}}
where the vertical inclusion in the lower right corner  is induced by the natural inclusion
\[ \Omega^q_{j-1,n}
\inj F^{j-1}_{X*}\Omega^q_{n-j+1}(\ul{D}_{n-j}, D_{n-j+1}+\ul{D}^{j-1}+p^{j-1}E).\]
Phrased differently $B_{j,n}^q(D,E)$, $Z^q_{j,n}(D,E)$ are obtained by pulling back the exact sequence 
\[F^{j-1}_{X*}\e^q_{n-j+1}(\ul{D}_{n-j}, D_{n-j+1}+\ul{D}^{j-1}+p^{j-1}E),\]
see \eqref{para:C-ex-seq-poles1}, along the respective inclusions on the right in the diagram \eqref{para:twisted-BZ3}. 

Some consequences from this definition:
\begin{enumerate}[label= (\alph*)]
    \item\label{para:twisted-BZa} We have for $n\ge 1$
    \[B_{1,n}^q(D,E)= 
    \sB_{1,n}^q=
    \sB^q(F_{X*}\Omega^\bullet_{n-1}(\ul{D}_{n-1}, \ul{D}^1+pE))\]
    and 
    \[Z_{1,n}^q(D,E)= \sZ^q(F_{X*}\Omega^\bullet_{n-1}(D, pE)) + B_{1,n}^q(D,E).\]
    This follows directly from the definition and \eqref{para:Omega-n7}.
    \item
    \label{para:twisted-BZb} 
    The sheaves $B_{j,n}^q(D,E)$ and $Z_{j,n}^q(D,E)$ are locally free coherent $\sO_X$-modules, which we can identify
    with submodules of $\sigma_*F^j_{X*}\Omega^q_U$, where $\sigma: U=X\setminus D+E\inj X$ is the open immersion.
   Furthermore  the $j$-fold iteration of the inverse Cartier operator  induces an isomorphism
    \[\Omega^q_n(D,E)=\frac{Z_{0,n}^q(D,E)}{B_{0,n}^q(D,E)}  \stackrel{\simeq}{\lra}\frac{Z_{j,n}^q(D,E)}{B_{j,n}^q(D,E)}.\]
    \item
    \label{para:twisted-BZc} 
    For $j\ge 1$, there are  inclusions of $\sO_X$-submodules of $\Omega^q_{j,n}$
    \[F_{X*}B^q_{j-1,n-1}(\ul{D}_{n-1}, \ul{D}^1+pE)\subset B_{j,n}^q(D,E)\]
    and 
    \[Z^q_{j,n}(D,E)\subset F_{X*}Z^q_{j-1,n-1}(\ul{D}_{n-1}, \ul{D}^1+pE).\]    
    This follows from the inclusions in the case $j=1$ and the fact that for $j\ge 2$ the bottom  sequence in \eqref{para:twisted-BZ3}
    constructed from $D$, $E$, $n$, $j$ is equal to $F_{X*}$ applied to the same sequence constructed 
    from $\ul{D}_{n-1}$, $\ul{D}^1+pE$, $n-1$, $j-1$.
    Furthermore, the $(j-1)$-fold iteration of the inverse Cartier operator  induces an isomorphism
    \[B^q_{1,n}(D,E)= \frac{B^q_{1,n}(D,E)}{F_{X*}B^q_{0,n-1}(\ul{D}_{n-1},\ul{D}^1+pE)}\stackrel{\simeq}{\lra} 
    \frac{B^q_{j,n}(D,E)}{F_{X*}B^q_{j-1,n-1}(\ul{D}_{n-1},\ul{D}^1+pE)}\]
    as follows from the above and the snake lemma.
\end{enumerate}
\end{para}

\begin{rmk}\label{rmk:BZ-loc-free}
Classically one defines $B_j\Omega^q_X\subset Z_j\Omega^q_X\subset F^j_{X*}\Omega^q_X$    as above in the case $D=E=0$.
More precisely, we have $B_j\Omega^q_X= B^q_{j,n}(0,0)$, for all $n\ge j$, and similarly for $Z_j\Omega^q_X$.
In this case the quotients $F^j_{X*}\Omega^q_X/Z_j\Omega^q_X$ are locally free as well, e.g. \cite[I, Corollaire 3.9]{IlDRW}. 
We warn the reader that this is in general not the case for 
\eq{rmk:BZ-loc-free1}{
\Omega^q_{j,n}/ 
Z^q_{j,n}(D,E).}
For example, for $j=1$ it follows from \eqref{para:twisted-BZ3} that we have an exact sequence 
\[0\to \frac{\Omega^q_n(\ul{D}_{n-1},D_n+E)}{\Omega^q_n(D,E)}\to 
\frac{F_{X*}\Omega^q_{n-1}(\ul{D}_{n-1}, \ul{D}^1+pE)}{Z^q_{1,n}(D,E)}\to 
\frac{F_{X*}\Omega^q_{n-1}(\ul{D}_{n-1}, \ul{D}^1+pE)}{\sZ^q_{1,n}}\to 0.\]
It is easy to see that the quotient on the right is locally free, but the quotient on the left is 
non-zero if $D_n\neq 0$ and $q<\dim X$ and  is annihilated by
$\sO_X(-D_n)$ and hence is not locally free. Thus the quotient in the middle is not locally free as well.
\end{rmk}

\begin{lem}\label{lem:RFil}
Let $L\in \Phi$ (see Notation \ref{nota:hdvf}). Then for $q$, $r$, $n\ge 0$ 
\[R^n(\Fil^p_r W_{n+1}\Omega^q_L)=\begin{cases}
\Omega^q_{\sO_L}(\log)\cdot\fm_L^{-\lceil r/p^n\rceil +1}, &\text{if } v_p(r)\le n,\\
\Omega^q_{\sO_L}\cdot \fm_L^{-r/p^n}, &\text{if } v_p(r)\ge n+1,
    \end{cases}\]
where $\Fil^p_rW_{n+1}\Omega^q_L$ is the filtration from Definition \ref{defn:p-sat} and 
$R^n:W_{n+1}\Omega^q_L\to \Omega^q_L$ is the restriction.
\end{lem}
\begin{proof}
By Definition of $\Fil^p$ we have 
\[R^n(\Fil^p_rW_{n+1}\Omega^q_L)= R^n(\fil_rW_{n+1}\Omega^q_L)+dR^n(\fil_rW_{n+1}\Omega^{q-1}_L).\]
Assume $v_p(r)\le n$. In this case we have by Definition \ref{defn:fil}
\[R^n(\fil_r W_{n+1}\Omega^q_L)= R^n(\fil^{\log}_{r-1}W_{n+1}(L))\cdot \dlog(K^M_q(L))
= \fm_L^{\lceil - \frac{r-1}{p^n}\rceil}\cdot\Omega^q_{\sO_L}(\log).\]
It is direct to check that $\lceil - \frac{r-1}{p^n}\rceil=-\lceil \frac{r}{p^n}\rceil+1$, which yields the statement in this case.
Assume $v_p(r)\ge n+1$. In this case we have by Definition \ref{defn:fil} and the above 
\begin{align*}
R^n(\fil_r W_{n+1}\Omega^q_L) &= 
\fm_L^{\lceil - \frac{r-1}{p^n}\rceil}\cdot\Omega^q_{\sO_L}(\log)+
R^n(\fil^{\log}_{r}W_{n+1}(L))\cdot \dlog(K^M_q(\sO_L))\\
&= \fm_L^{ - \frac{r}{p^n}+1}\cdot\Omega^q_{\sO_L}(\log)+ \fm_L^{-\frac{r}{p^n}}\cdot \Omega^q_{\sO_L}=
\fm_L^{-\frac{r}{p^n}}\cdot \Omega^q_{\sO_L}.
\end{align*}
As $r/p^n$ is divisible by $p$, we have 
$d(\fm_L^{-r/p^n}\cdot \Omega^{q-1}_{\sO_L})\subset \fm_L^{-r/p^n}\cdot \Omega^q_{\sO_L}$.
Hence the statement.
\end{proof}

\begin{cor}\label{cor:RHWM}
Let $D=D'+ p^{n+1}D_{n+1}$ be a \pdd.
Then for $q$, $n\ge 0$ we have with the notation from \ref{para:Omega-n}
\[R^n(W_{n+1}\Omega^q_{(X,D)})=\Omega_n^q(D', pD_{n+1}).\]
\end{cor}
\begin{proof}
As $\Omega_n^q(D', pD_{n+1})$    is locally free it suffices to show this ``$\subset$'' inclusion for 
the  Nisnevich stalks in the generic points of $D$.
In this case the statement follows from Lemma \ref{lem:RFil} and the notation \eqref{para:Omega-n1},
see \ref{para:HW-not}, \ref{para:Fc}, and Theorem \ref{thm:HW-modulus}.
To show the other inclusion we can argue locally around  the {\em closed} points of $D$.
Let $x\in D$ be a point and set $A=\sO_{X,x}$. Then  we find an \'etale map $k[t_1,\ldots, t_d]\to A$, such that 
$t_1,\ldots,t_d$ form a regular sequence of parameters for $A$ and on $\Spec A$ we have
\[D'=\Div(t_1^{m_1}\cdots t_r^{m_r}) \quad\text{and}\quad D_{n+1}=\Div(t_{r+1}^{m_{r+1}}\cdots t_s^{m_s}),\]
with $0\le r\le s\le d$ and $p^{n+1}\nmid m_i$, for $i=1,\ldots,r$ (with the convention that $r=0$ means that $D'=\emptyset$ and $r=s$ means that 
$D_{n+1}=\emptyset$). Thus 
\[D=\Div(t_1^{m_1}\cdots t_r^{m_r}\cdot t_{r+1}^{p^{n+1}m_{r+1}}\cdots t_s^{p^{n+1}m_s}).\]
A basis of the free $A$-module $\Omega^q_n(D',pD_{n+1})_x$ (Zariski stalk at $x$) is given by 
\[e_{I,J} = t_1^{-\lceil\frac{m_1}{p^n}\rceil+1}\cdots t_r^{-\lceil\frac{m_r}{p^n}\rceil+1} t_{r+1}^{-pm_{r+1}}\cdots t_{s}^{-pm_{s}}\cdot
\dlog \ul t_I\dlog (\ul{1+t}_J),\]
where $I$ and $J$ run through the tuples $I=(1\le i_1<\ldots<i_{q_1}\le r)$ and $J=(r+1\le j_1<\ldots< j_{q_2}\le d)$ with $q_1+q_2=q$
and where we use the notations 
\[ \ul{t}_I= \{t_{i_1}, \ldots, t_{i_{q_1}}\}\in K^M_{q_1}(A[\tfrac{1}{t_{i_1}\cdots t_{i_{q_1}}}])
\quad \text{and}\quad \ul{1+ t}_J= \{ 1+t_{j_1},\ldots, 1+ t_{j_{q_2}}\}\in K^M_{q_2}(A).\]
As the  Zariski stalk $W_{n+1}\Omega^q_{(X,D),x}$ is a $W_{n+1}(A)$-module and $R^n(W_{n+1}(A))=A$, it suffices to show that
there exists an element $\tilde{e}_{I,J}\in W_{n+1}\Omega^q_{(X,D),x}$, for each $I$, $J$, 
with $R^n(\tilde{e}_{I,J})=e_{I,J}$. By the second part of Proposition \ref{prop:strHWM} it is direct to check  
from the definition of $\Fil^p_rW_n\Omega^q_L$ that the element
\[\tilde{e}_{I,J}:=[t_1]^{-\lceil\frac{m_1}{p^n}\rceil+1}\cdots [t_r]^{-\lceil\frac{m_r}{p^n}\rceil+1} 
[t_{r+1}]^{-pm_{r+1}}\cdots [t_{s}]^{-pm_{s}}\cdot
\dlog\ul{t}_I\dlog (\ul{1+t}_J),\]
which a priori is an element in $W_{n+1}\Omega^q_A[\frac{1}{t_1\cdots t_s}]$, has the  looked for properties.
\end{proof}

\begin{prop}\label{prop:BZnn}
Let $D=D'+p^{n+1}D_{n+1}$ be a \pdd. Then with the notation from \eqref{thm:strHWM1} and \ref{para:twisted-BZ} we have for $n\ge 0$
\[B_n\Omega^q_{(X,D)}= B^q_{n,n}(D', pD_{n+1}),\quad \text{and}\quad Z_n\Omega^q_{(X,D)}= Z^q_{n,n}(D', pD_{n+1}).\]
In particular, $B_n\Omega^q_{(X,D)}$ and $Z_n\Omega^q_{(X,D)}$ are locally free coherent $\sO_X$-modules.
\end{prop}
\begin{proof}
For $n=0$ the statement follows from the definitions and \eqref{para:HW-not5}.
Assume now $n\ge 1$. Let $D'= D''+p^nD_n$ be a \pdd.
By induction we have 
\eq{prop:BZnn1}{Z_{n-1}\Omega^q_{(X,D)}= Z^q_{n-1, n-1}(D'', pE),\quad \text{where }E=D_n+pD_{n+1}.}
and similarly with $Z$  replaced by $B$.
Consider the following diagram (in which we neglect the $\sO_X$-module structure)
\eq{prop:BZnn2}{\xymatrix{
W_n\Omega^{q-1}_{(X,D)}\ar@{->>}[r]^d\ar@{->>}[dd]^{R^{n-1}} & d(W_n\Omega^{q-1}_{(X,D)})\ar@{->>}[dd]^{R^{n-1}}\ar@{->>}[r]^-{F^{n-1}} & 
B_n\Omega^q_{(X,D)}\subset Z_{n-1}\Omega^q_{(X,D)}\ar@{->>}[dr] & \\
                         &              &       & \frac{Z_{n-1,n-1}^q(D'', pE)}{B^q_{n-1,n-1}(D'', pE)}.\\
\Omega^{q-1}_{n-1}(D'',pE)\ar@{->>}[r]^d & B^q_{1,n}(D', pD_{n+1})\ar[r]^{C^{-(n-1)}}_{\simeq} & 
                       \frac{B_{n,n}^q(D', pD_{n+1})}{B^q_{n-1,n-1}(D'', pE)}\ar@{^(->}[ur] &
}}
Here the map $d$ on the lower left is surjective by \ref{para:twisted-BZ}\ref{para:twisted-BZa} and the vertical map $R^{n-1}$ on the left
hand side is surjective by Corollary \ref{cor:RHWM}; clearly the square commutes and hence also the vertical $R^{n-1}$ in the middle is surjective.
The  lower horizontal map $C^{-{(n-1)}}$ is induced by the $(n-1)$-fold iteration of the inverse Cartier operator;  
it is an isomorphism, see \ref{para:twisted-BZ}\ref{para:twisted-BZc}.
The upwards diagonal inclusion on the right hand side is induced by the natural inclusions 
in \ref{para:twisted-BZ}\ref{para:twisted-BZc}, the downwards diagonal map is induced by
\eqref{prop:BZnn1}. The inclusion on the top right is given by
\eq{prop:BZnn2.5}{
B_n\Omega^q_{(X,D)}=F^{n-1}d(W_n\Omega^{q-1}_{(X,D)})=F^{n}dV(W_n\Omega^{q-1}_{(X,D)})\subset F^n (W_{n+1}\Omega^q_{(X,D)})=Z_n\Omega^q_{(X,D)}.}
To see that  the right part of the diagram commutes observe that by \ref{para:twisted-BZ}\ref{para:twisted-BZb}
\[\frac{Z_{n-1,n-1}^q(D'', pE)}{B^q_{n-1,n-1}(D'', pE)}\cong \Omega^q_{n-1}(D'', pE)\]
is a locally free $\sO_X$-module; hence it suffices to check the commutativity after restricting to $U=X\setminus D$
in which case it follows from the classical fact that the Frobenius on the de Rham-Witt complex is a lift of the inverse Cartier operator,
see \cite[I, Proposition (3.3)]{IlDRW}. Thus the whole diagram \eqref{prop:BZnn2} commutes and 
all the maps in the diagram are defined and are surjections, injections, or isomorphisms as indicated.
Furthermore we have the inclusions
\[B_{n-1,n-1}^q(D'', pE)= F^{n-2} d(W_{n-1}\Omega^{q-1}_{(X,D)})=F^{n-1} dV(W_{n-1}\Omega^{q-1}_{(X,D)})
\subset B_n\Omega^q_{(X,D)},\]
where the first equality holds by induction.
Therefore diagram \eqref{prop:BZnn2} yields 
\[\frac{B_{n,n}^q(D', pD_{n+1})}{B^q_{n-1,n-1}(D'', pE)}= \frac{B_n\Omega^q_{(X,D)}}{B^q_{n-1,n-1}(D'', pE)}\]
and hence we get the equality of subsheaves of $\Omega^q_{(X,D)}$
\eq{prop:BZnn2.6}{B_{n,n}^q(D', pD_{n+1})= B_n\Omega^q_{(X,D)}.}
Next we consider the diagram
\eq{prop:BZnn3}{\xymatrix{
W_{n+1}\Omega^q_{(X,D)}\ar@{->>}[r]^-{F^n}\ar@{->>}[dd]^{R^n} & Z_n\Omega^q_{(X,D)}\subset Z_{n-1}\Omega^q_{(X,D)}\ar@{->>}[dr] & \\
                                      &       & \frac{Z_{n-1,n-1}^q(D'', pE)}{B^q_{n,n}(D', pD_{n+1})}.\\
\Omega^q_n(D',pD_{n+1})\ar[r]^{C^{-n}}_{\simeq} &  \frac{Z_{n,n}^q(D', pD_{n+1})}{B^q_{n,n}(D', pD_{n+1})}\ar@{^(->}[ur] &
}}
Here $R^n$ on the left hand side is surjective by Corollary \ref{cor:RHWM}. The lower horizontal map $C^{-n}$
is an isomorphism by \ref{para:twisted-BZ}\ref{para:twisted-BZb}.
The upward diagonal injection is induced by \ref{para:twisted-BZ}\ref{para:twisted-BZc}, the downward diagonal map is induced by
\eqref{prop:BZnn1}. We claim this diagram  commutes. As above it suffices to show that the quotient on the right is locally free 
as the restriction of the diagram to $U=X\setminus D$ is known to commute. To this end   consider the isomorphisms 
\[\frac{Z_{n-1,n-1}^q(D'', pE)}{B^q_{n,n}(D', pD_{n+1})}\cong 
\frac{Z_{n-1,n-1}^q(D'', pE)/B_{n-1,n-1}(D'',pE)}{B^q_{n,n}(D', pD_{n+1})/B_{n-1,n-1}(D'',pE)}
\cong \frac{\Omega^q_{n-1}(D'',pE)}{\sB^q(\Omega^\bullet_{n-1}(D'',pE))},\]
where the isomorphisms hold by \ref{para:twisted-BZ}\ref{para:twisted-BZa}--\ref{para:twisted-BZc}.
Moreover we have the exact sequence
\[0\to\Omega^q_n(D'', E)\stackrel{C^{-1}}{\lra} \frac{F_{X*}\Omega^q_{n-1}(D'', pE)}{\sB^q(\Omega^\bullet_{n-1}(D'',pE))}
\stackrel{d}{\lra}  \sB^{q+1}(\Omega^\bullet_{n-1}(D'',pE))\to 0, \]
see \eqref{para:Omega-n6}. Since the two outer terms of this exact sequence are locally free so is the middle term.
Hence the quotient on the right in diagram \eqref{prop:BZnn3} is locally free and the diagram commutes.
Using the equality \eqref{prop:BZnn2.6} and the inclusion \eqref{prop:BZnn2.5}, the commutativity of \eqref{prop:BZnn3}
yields 
\[\frac{Z_n\Omega^q_{(X,D)}}{B_n\Omega^q_{(X,D)}}=\frac{Z_{n,n}(D',pD_{n+1})}{B_n\Omega^q_{(X,D)}}.\]
Hence we get the equality 
\[Z_n\Omega^q_{(X,D)}= Z_{n,n}(D',pD_{n+1})\]
of subsheaves of $\Omega^q_{(X,D)}$, which yields the statement.
\end{proof}

\begin{lem}\label{lem:BnZn}
Let $j:U=X\setminus D\inj X$ be the open immersion. With the notation from section \ref{sec:HWM} we have 
\eq{lem:BnZn1}{B_n\Omega^q_{(X,D)}= (\ul{\omega}^{\CI}B_n\Omega^q)_{(X,D)}= (B_n\Omega^q)^{\AS}_{(X,D)}= j_*(B_n\Omega^q_U)\cap \Omega^q_{(X,D)},\quad 
\text{for all } q\ge 0,}
and 
\eq{lem:BnZn2}{Z_n\Omega^q_{(X,D)}=(\ul{\omega}^{\CI}Z_n\Omega^q)_{(X,D)}= (Z_n\Omega^q)^{\AS}_{(X,D)}= j_*(Z_n\Omega^q_U)\cap \Omega^q_{(X,D)},\quad 
\text{for all } q\ge 1,}
and 
\eq{lem:BnZn3}{Z_n\Omega^0_{(X,D)}= j_*(Z_n\Omega^0_U)\cap \Omega^0_{(X,D)}= j_*F^n(W_{n+1}\sO_U)\cap \sO_{(X,D)}.}
\end{lem}
\begin{proof}
The second equality in \eqref{lem:BnZn3} holds by definition. The first equality reduces to 
\[F^n(\fil_r W_{n+1}\Omega^0_L) = F^n(W_{n+1}(L))\cap \fil_r \Omega^0_L,\quad \text{for }L\in \Phi, r\ge 0,\]
which follows directly from  Definition \ref{defn:fil}. As $B_n\Omega^0=0$ we assume  $q\ge 1$ in the following.
By the definition of $B_n\Omega^q_{(X,D)}$ and $Z_n\Omega^q_{(X,D)}$ in \eqref{thm:strHWM1}, by Theorem \ref{thm:HW-modulus}, 
and \eqref{para:AS1} we see that \eqref{lem:BnZn1} and \eqref{lem:BnZn2} hold with ``$=$" replaced by ``$\subset$".
Thus it remains to show 
\[j_*(B_n\Omega^q_U)\cap \Omega^q_{(X,D)}\subset B_n\Omega^q_{(X,D)} \quad \text{and}\quad 
j_*(Z_n\Omega^q_U)\cap \Omega^q_{(X,D)}\subset Z_n\Omega^q_{(X,D)}.\]
We first consider the $B_n$ case. Set  $M:=j_*(B_n\Omega^q_U)\cap \Omega^q_{(X,D)}$ and $N:=B_n\Omega^q_{(X,D)}$.
We note that $M$ and $N$ are $\sO_X$-submodules of $V:=j_*B_n\Omega^q_U$, with $N$ locally free of finite rank, 
by Proposition \ref{prop:BZnn}, and with $M$ coherent, as it is a submodule of the coherent module $F^n_*\Omega^q_{(X,D)}$. Moreover, 
$M_{|U}=N_{|U}=V_{|U}$ is a finite locally free $\sO_U$-module. By \cite[VII, \S 4, No. 2, Corollaire]{BourbakiCA} it suffices 
therefore to show the inclusion of stalks 
$M_{\eta}\subset N_\eta$, for all generic points $\eta$ of $D_\red$. As $N_{\eta}^h\cap V_\eta= N_\eta$ we may consider Nisnevich stalks instead 
and are thus reduced  to show the following claim for each henselian dvf $L\in \Phi$:
\eq{lem:BnZn-4}{F^{n-1}d(\Fil^p_sW_n\Omega^{q-1}_L)\cap \Fil^p_r\Omega^q_L\subset F^{n-1}d(\Fil^p_rW_n\Omega^{q-1}_L), 
\quad \text{for all }s\ge r\ge 1.}
Indeed, this holds trivially for $s=r$. Assume $s>r$ and let $\alpha\in \Fil^p_s W_n\Omega^{q-1}_L$ such that 
$F^{n-1}d(\alpha)\in \Fil^p_r\Omega^q_L\subset \Fil^p_{s-1}\Omega^q_L$.
By Lemma \ref{lem:gr-ex} we find elements
$\beta\in W_{n+1}\Omega^q_L$, $\alpha'\in \Fil^p_{s-1}W_n\Omega^q_L$, and $\alpha''\in \Fil^p_sW_{n-1}\Omega^q_L$
such that $\alpha= \alpha'+\ul{p}(\alpha'')+F(\beta)$.
Hence
\[F^{n-1}d(\alpha)=F^{n-1}d(\alpha')\in F^{n-1}d(\Fil^p_{s-1}W_n\Omega^q_L)\cap \Fil_r^p \Omega^q_L\]
and  claim \eqref{lem:BnZn-4} holds by induction on $s$. This completes the proof of \eqref{lem:BnZn1}. 

For $Z_n$ we can argue similarly as above to reduce to the following claim for each $L\in \Phi$:
\eq{lem:BnZn-claim5}{
F^n(\Fil^p_s W_{n+1}\Omega^q_L) \cap \Fil^p_r\Omega^q_L\subset F^n(\Fil^p_rW_{n+1}\Omega^q_L).}
Indeed, this holds trivially for $s=r$. We assume $s>r$. Let $\alpha\in \Fil^p_s W_{n+1}\Omega^q$ such that 
$F^n(\alpha)\in \Fil^p_r\Omega^{q}_L\subset \Fil^p_{s-1}\Omega^{q}_L$.
By Lemma \ref{lem:gr-ex} we find elements $\beta\in W_n\Omega^q_L$, $\alpha'\in \Fil^p_{s-1} W_{n+1}\Omega^q_L$ and 
$\alpha''\in \Fil^p_s W_n\Omega^q_L$ such that $\alpha= \alpha'+\ul{p}(\alpha'')+ V(\beta)$.
Hence
\[F^n(\alpha)= F^n(\alpha') \in F^n(\Fil^p_{s-1}W_{n+1}\Omega^q_L)\cap \Fil^p_r\Omega^q_L\]
and  claim \eqref{lem:BnZn-claim5} holds by induction. This completes the proof of \eqref{lem:BnZn2}. 
\end{proof}

\begin{lem}\label{lem:KerFnFnd-mod}
For the statement of this lemma denote by $F^n_{(X,D)}: W_{n+1}\Omega^q_{(X,D)}\to \Omega^q_{(X,D)}$    
the map induced by  the $n$th-power of the de Rham-Witt Frobenius $F^n$ (which otherwise is denoted by $F^n$ as well).
Then we have the following equality of subsheaves of $W_{n+1}\Omega^q_{(X,D)}$
\[\Ker (F^n_{(X,D)})\cap \Ker(F^n_{(X,D)}\circ d)=\ul{p}(j_*W_n\Omega^q_U)\cap W_{n+1}\Omega^q_{(X,D)}
=\ul{p}(W_n\Omega^q_{(X,D)}),\]
where $j:U=X\setminus D\inj X$ is the open immersion.
\end{lem}
\begin{proof}
Denote  by $F^n_U: j_*W_{n+1}\Omega^q_U\to j_*\Omega^q_U$ the Frobenius on $U$.  
As $\Omega^q_{(X,D)}\subset j_*\Omega^q_U$ we have 
 \[\Ker(F^n_{(X,D)})= \Ker(F^n_U)\cap W_{n+1}\Omega^q_{(X,D)}\quad\text{and}\quad 
 \Ker(F^n_{(X,D)}d)= \Ker(F^n_U d)\cap W_{n+1}\Omega^q_{(X,D)}.\]
As $\Ker(F^n_U)\cap \Ker(F^n_U d)=\ul{p}(j_*W_n\Omega^q_U)$ by \cite[(0.6.3)]{Ekedahl} 
the first equality follows.

For the second equality let $a\in j_*W_n\Omega^q_U$ be a local section with $\ul{p}(a)\in W_{n+1}\Omega^q_{(X,D)}$. As $\ul{p}$ is injective we have to show
that $a\in W_n\Omega^q_{(X,D)}$. By Proposition \ref{prop:strHWM}  it suffices to show
\eq{lem:KerFnFnd-mod1}{\ul{p}(W_n\Omega^q_L)\cap \Fil^p_{r} W_{n+1}\Omega^q_L=\ul{p}(\Fil^p_r W_n\Omega^q_L),}
for all $q$, $r\ge 0$ and $L\in \Phi$. For $r=0$ this follows from Lemma \ref{upreg} and Lemma \ref{lem:Filp}. 
Assume $r\ge 1$ and let $a\in W_n\Omega^q_L$ with $\ul{p}(a)\in \Fil^p_{r} W_{n+1}\Omega^q_L$. It follows from Lemma \ref{lem:char}
that there exist elements $b\in \Fil^p_rW_n\Omega^q_L$ and $c\in \Fil^p_{r-1}W_{n+1}\Omega^q_L$ such that 
\[\ul{p}(a)=c+ \ul{p}(b).\]
By induction over $r$ we find an element $e\in \Fil^p_{r-1}W_{n}\Omega^q_L$ such that $c=\ul{p}(e)$.
The  injectivity of $\ul{p}$ therefore yields $a=e+b\in \Fil^p_rW_n\Omega^q_L$.
Hence the second equality holds.
\end{proof}

\begin{prop}\label{prop:HWM-ex}
 The sequence 
\[0\lra B_n\Omega^{q+1}_{(X,D)}\lra \frac{W_{n+1}\Omega^q_{(X,D)}}{\ul{p}W_n\Omega^q_{(X,D)}}\xr{F^n} Z_n\Omega^q_{(X,D)}\lra 0\]
is a short exact sequence of $W_{n+1}\sO_X$-modules, where the first map on the left is given by $F^{n-1}d(\beta)\mapsto V(\beta)$, 
and the $W_{n+1}\sO_X$-module structure on $Z_n$ and $B_n$ is induced by $F^n:W_{n+1}\sO_X\to \sO_X$.
\end{prop}
\begin{proof}
This is a modulus version of \cite[Lemma 0.6]{Ekedahl}. After the previous work the proof is analogous:
consider the diagram (with the notation from Lemma \ref{lem:KerFnFnd-mod})
\eq{prop:HWM-ex1}{\xymatrix{
0\ar[r] & \Ker (F^n_{(X,D)})\cap \Ker(F^n_{(X,D)}d)\ar@{^(->}[d]^{\iota}\ar[r] 
       &  W_{n+1}\Omega^q_{(X,D)}\ar[r]\ar@{=}[d] 
       & \frac{W_{n+1}\Omega^q_{(X,D)}}{\ul{p}W_n\Omega^q_{(X,D)}}\ar[r]\ar[d]^{F^n} & 0\\
0\ar[r] & \Ker(F^n_{(X,D)})\ar[r]
      & W_{n+1}\Omega^q_{(X,D)}\ar[r]^{F^n}
      & Z_n\Omega^q_{(X,D)}\ar[r] &0.      
}}
Here the top row is exact by Lemma \ref{lem:KerFnFnd-mod}, the bottom row is exact by definition, and  
$\iota$ is the inclusion. By \cite[I, (3.11.3)]{IlDRW} and the inclusion $\Omega^q_{(X,D)}\subset j_*\Omega^q_U$ we have 
\[\Ker(F^n_{(X,D)})=\Ker(F^n_U)\cap W_{n+1}\Omega^q_{(X,D)}= j_*(VW_n\Omega^q_U)\cap W_{n+1}\Omega^q_{(X,D)}.\]
We have the following obvious chain of inclusions
\[B_n\Omega^{q+1}_{(X,D)}= F^n d V W_n\Omega^q_{(X,D)}
\subset F^nd (j_*(VW_n\Omega^q_U)\cap W_{n+1}\Omega^q_{(X,D)})\subset j_*(B_n\Omega^{q+1}_U)\cap \Omega^{q+1}_{(X,D)}\]
and by Lemma \ref{lem:BnZn}  we have equality everywhere. 
We obtain an exact sequence
\[0\to \Ker (F^n_{(X,D)})\cap \Ker(F^n_{(X,D)}d)\xr{\iota} \Ker(F^n_{(X,D)})\xr{F^nd} B_n\Omega^{q+1}_{(X,D)}\to 0.\]
Hence applying the Snake Lemma to  \eqref{prop:HWM-ex1} yields the exact sequence from the statement.
\end{proof}

The proof of Theorem \ref{thm:strHWM} is now complete. 
The following proposition is a finite level version of the exact sequence 
\cite[II, (1.2.2)]{IllRay} in the modulus setting.
\begin{prop}\label{prop:long-mod-seq}
The following sequence is exact for all $r\ge 1$
\[0\to W_n\Omega^q_{(X,D)}\xr{\ul{p}^r} W_{n+r}\Omega^q_{(X,D)}\xr{(F^n, F^nd)} W_r\Omega^q_{(X,D)}\oplus W_r\Omega^{q+1}_{(X,D)}
\xr{dV^n-V^n} W_{n+r}\Omega^{q+1}_{(X,D)}.\]
\end{prop}
\begin{proof}
Denote by $j:U=X\setminus D\inj X$ the open immersion.
We first consider the case $r=1$. Thus we have to show that the sequence
\[0\to W_n\Omega^q_{(X,D)}\xr{\ul{p}} W_{n+1}\Omega^q_{(X,D)}\xr{(F^n, F^nd)} \Omega^q_{(X,D)}\oplus \Omega^{q+1}_{(X,D)}
\xr{dV^n-V^n} W_{n+1}\Omega^{q+1}_{(X,D)}\]
is exact. 
The injectivity of $\ul{p}$ follows from the injectivity of $\ul{p}: W_n\Omega^q_{U} \to W_{n+1}\Omega^q_U$, see \cite[I, Proposition 3.4]{IlDRW}.
The exactness at $W_{n+1}\Omega^q_{(X,D)}$ follows from Proposition \ref{prop:HWM-ex} (use $F^ndV(\beta)=F^{n-1}d\beta$).
Let $\alpha\in \Omega^q_{(X,D)}$ and $\beta\in \Omega^{q+1}_{(X,D)}$ be local sections. It remains to show
\eq{prop:long-mod-seq1}{dV^n(\alpha)=V^n(\beta)\Longrightarrow \alpha=F^n(\gamma) \text{ and } \beta=F^nd(\gamma), 
\text{ for some }\gamma\in W_{n+1}\Omega^q_{(X,D)}.}
By \cite[I, (3.10.3)]{IlDRW} we find $\gamma_0\in W_{n+1}\Omega^q_U$ such that 
\[dV^n(\alpha)=dV^nF^n(\gamma_0)= V^n(F^nd(\gamma_0))=V^n(\beta).\]
By \cite[I, Remarques 3.21.1]{IlDRW} we find $\delta\in W_{n+2}\Omega^q_U$ and $\epsilon\in W_n\Omega^q_U$ such that 
\[\alpha= F^n(\gamma_0) +F^{n+1}(\delta)\quad \text{and}\quad \beta= F^n d(\gamma_0)+F^ndV(\epsilon).\]
Setting $\gamma_1:=\gamma_0 +F(\delta)+V(\epsilon)\in W_{n+1}\Omega^q_U$ we have 
\[\alpha= F^n(\gamma_1)\in \Omega^q_{(X,D)}\cap j_* (Z_n\Omega^q_U)\quad \text{and}\quad 
\beta=F^nd(\gamma_1)\in \Omega^q_{(X,D)}\cap j_* (B_{n+1}\Omega^q_U).\]
By Lemma \ref{lem:BnZn} we find $\gamma'\in W_{n+1}\Omega^q_{(X,D)}$ and $\gamma''\in W_{n+1}\Omega^q_{(X,D)}$
such that 
\[\alpha=F^n(\gamma') \quad \text{and} \quad \beta=F^nd(\gamma'').\]
By \cite[I, Remarques 3.21.1]{IlDRW} we find (new) $\delta$, $\epsilon$ on $U$ with 
\eq{prop:long-mod-seq2}{\gamma_1=\gamma'+V(\delta)= \gamma''+F(\epsilon).}
Hence
\[F^n(\gamma'-\gamma'')= F^{n+1}(\epsilon)\in j_*(Z_{n+1}\Omega^q_U)\cap \Omega^q_{(X,D)}\]
and
\[F^n d(\gamma''-\gamma')= F^{n-1}d (\delta)\in j_*(B_n\Omega^{q+1})\cap \Omega^{q+1}_{(X,D)}.\]
By Lemma \ref{lem:BnZn} and \cite[I, Remarques 3.21.1]{IlDRW} we find $\epsilon'\in W_{n+2}\Omega^q_{(X,D)}$, 
$\delta'\in W_n\Omega^q_{(X,D)}$, and  $\eta, \zeta\in W_{n+1}\Omega^q_U$ such that
\[F(\epsilon)= F(\epsilon')+p\zeta\quad \text{and}\quad V(\delta)=V(\delta')+p\eta.\]
Thus \eqref{prop:long-mod-seq2} yields
\[\gamma''-\gamma' +F(\epsilon')-V(\delta')=p(\eta-\zeta)=\ul{p}(R(\eta-\zeta))\in \ul{p}(j_*W_n\Omega^q_U)\cap W_{n+1}\Omega^q_{(X,D)}.\]
By Lemma \ref{lem:KerFnFnd-mod} we find $\theta\in W_n\Omega^q_{(X,D)}$ such that
\[\ul{p}\theta = \ul{p}(R(\eta-\zeta)).\]
Set
\[\gamma:=\gamma'+V(\delta')+\ul{p}\theta= \gamma''+F(\epsilon')\in W_{n+1}\Omega^q_{(X,D)}.\]
Then 
\[\alpha=F^n(\gamma_1)=F^n(\gamma'+V(\delta'))=F^n(\gamma)\quad\text{and}\quad
\beta= F^nd(\gamma_1)=F^n d(\gamma'' +F(\epsilon))=F^n d(\gamma).\]
This shows \eqref{prop:long-mod-seq1} and completes the proof in the case $r=1$.

Now assume $r \ge 2$ and consider the following diagram in which we drop for readability the subscript $(X,D)$ everywhere 
\[
\xymatrix{
        &     0\ar[d] & 0\ar[d] & 0\ar[d] & 0\ar[d]\\
0\ar[r] & W_n\Omega^q\ar@{=}[d]\ar[r]^{\ul{p}^{r-1}} &
W_{n+r-1}\Omega^q\ar[r]^-{(F^n, F^n d)}\ar[d]^{\ul{p}} &
W_{r-1}\Omega^q\oplus W_{r-1}\Omega^{q+1}\ar[r]^-{dV^n-V^n}\ar[d]^{\ul{p}}&
W_{n+r-1}\Omega^{q+1}\ar[d]^{\ul{p}}\\
0\ar[r] & W_n\Omega^q\ar[r]^{\ul{p}^r} &
W_{n+r}\Omega^q\ar[r]^-{(F^n, F^n d)}\ar[d]_{F^{n+r-1}}^{F^{n+r-1}d} &
W_r\Omega^q\oplus W_r\Omega^{q+1}\ar[r]^-{dV^n-V^n}\ar[d]_{(F^{r-1}, F^{r-1}d)}^{(F^{r-1}, F^{r-1}d)}&
W_{n+r}\Omega^{q+1}\ar[d]_{F^{n+r-1}}^{F^{n+r-1}d}\\
  &   &
\Omega^q\oplus \Omega^{q+1}\ar[r]^-{\varphi}\ar[d]^{dV^{n+r-1}-V^{n+r-1}} & 
(\Omega^q\oplus \Omega^{q+1})\oplus (\Omega^{q+1}\oplus \Omega^{q+2})\ar[r]^-{\psi}&
\Omega^{q+1}\oplus \Omega^{q+2}\\
&  & 
W_{n+r}\Omega^{q+1},
}
\]
where $\varphi$ and $\psi$ are given by
\[\varphi(a,c)=((a,0), (c,0)) \quad \text{and}\quad \psi((a,b),(c,d))=(b,-d).\]
We observe 
\begin{itemize}
\item the diagram is commutative;
\item the columns are exact by the case $r=1$;
\item the first row is exact by induction over $r$;
\item the third row is split exact.
\end{itemize}
We want to show that the middle row is exact on $(X,D)$. Clearly it is a complex and  $\ul{p}^r$ is injective by the same argument 
as in the case $r=1$. The exactness at $W_{n+r}\Omega^q_{(X,D)}$ follows from an easy diagram chase.
The exactness at $W_r\Omega^q_{(X,D)}\oplus W_r\Omega^{q+1}_{(X,D)}$ can be checked directly by a diagram chase as well 
once we observed that for $\alpha\in W_r\Omega^q_{(X,D)}$ and $\beta\in W_r\Omega^{q+1}_{(X,D)}$ with $dV^n(\alpha)=V^n(\beta)$
we have 
\[dV^{n+r-1}(F^{r-1}(\alpha))=p^{r-1}dV^n(\alpha)=p^{r-1}V^n(\beta)=V^{n+r-1}(F^{r-1}(\beta)).\]
This completes the proof.
\end{proof}

\section{Applications to Hodge-Witt cohomology with modulus}
Before we proceed to study Hodge-Witt sheaves with zeros along $D$ in section \ref{sec:strHWZ},
we  draw some consequences from the structural results in the previous section.
The main result is  Theorem \ref{lem:HWM-gamma}, which 
makes it possible to apply general results on cube invariant sheaves with transfers to Hodge-Witt sheaves with modulus,
see \ref{para:pbf-buf-gt}.

Throughout this section we assume $X\in \Sm$ pure dimensional with $\dim X=N$
and denote by $D$ an effective Cartier divisor on $X$ such that $D_{\red}$ is an SNCD.

\begin{para}\label{para:MNST-gamma}
Recall  the category $\uMCor$  from \cite[Definition 1.3.1]{KMSYI}: 
its objects are modulus pairs $(Y,E)$ (see \ref{para:mod-pairs})
and the morphisms are left  proper admissible correspondences. A {\em modulus presheaf with transfers} is an additive 
contravariant functor from $\uMCor$  to abelian groups and the category of modulus presheaves with transfers is denoted by $\uMPST$. 
If $G\in \uMPST$ and $(Y,E)$ is a modulus pair, then the assignment
\[(\text{\'etale $Y$-schemes})\ni (v:V\to Y)\mapsto G(V, v^*E),\]
defines a (Nisnevich) presheaf on $Y$ denoted by $G_{(Y,E)}$. By definition $G\in \uMPST$ is a sheaf  
if $G_{(Y,E)}$ is a Nisnevich sheaf on $Y$, for each modulus pair $(Y,E)$. 
The category of modulus sheaves with transfers is denoted by $\uMNST$. 
By \cite[Theorem 2]{KMSYI} there is a sheafification functor $\uMPST\to \uMNST$, which is an exact left adjoint to 
 the inclusion $\uMNST \to \uMPST$. The full subcategory of $\uMPST$ with objects the $(\P^1,\infty)$-invariant 
 modulus presheaves with transfers having semipurity and $M$-reciprocity is denoted by $\CI^{\tau, {\rm sp}}$, 
 see \cite[Definition 1.31]{Saito-Purity}. We set $\CI^{\tau, {\rm sp}}_{\Nis}=\CI^{\tau, {\rm sp}}\cap \uMNST$.
 By \cite[Theorem 0.4]{Saito-Purity} the sheafification restricts to $\CI^{\tau, {\rm sp}}\to \CI^{\tau, {\rm sp}}_{\Nis}$.
Note that if $F\in \RSC_{\Nis}$, then $\tF_c$, $\uomega^{\CI}F\in \CI^{\tau, {\rm sp}}_{\Nis}$, 
see \ref{para:RSC} and \ref{para:Fc}. In particular,
\[\widetilde{W_n\Omega^q}_c\in \CI^{\tau, {\rm sp}}_{\Nis},\]
see \ref{para:Fc} and Theorem \ref{thm:HW-modulus} for the notation (here $c$ is the conductor defined by 
the filtration from Definition \ref{defn:p-sat}), and 
$(\widetilde{W_n\Omega^q}_c)_{(X,D)}$ as defined above coincides with $W_n\Omega^q_{(X,D)}$ as defined in \ref{para:HW-not}.

The assignment $(Y,E)\mapsto K^M_r(Y\setminus E)$ defines a modulus sheaf with transfers denoted by $\uomega^*K^M_r$.
For $G\in \CI^{\tau, {\rm sp}}_{\Nis}$   and $r\ge 0$ we define
\[\gamma^r G:=\uHom_{\uMPST}(\uomega^*K^M_r, G).\]
By, e.g., \cite[Corollary 4.5]{BRS}, we have $\gamma^r G\in \CI^{\tau, {\rm sp}}_{\Nis}$ and 
by \cite[Theorem 6.3]{BRS} there is  a canonical isomorphism for any $s\ge r$
\eq{para:MNST-gamma1}{(\gamma^r G)_{(X,D)}\cong R^r\pi_* G_{(\P^s\times X,\pi^*D)},}
where $\pi: \P^s\times X\to X$ is the projection. 
\end{para}

\begin{lem}\label{lem:gamma-exact}
Let $ G'\to G\to G''$ be a sequence in $\CI^{\tau,{\rm sp}}_{\Nis}$.
Assume that 
\[0\to G'_{(\P^s\times X, \pi^*D)}\to G_{(\P^s\times X, \pi^*D)}\to G_{(\P^s\times X, \pi^*D)}''\to 0\]
is an exact sequence of Nisnevich sheaves on $\P^s\times X$ for all $s\ge 0$.
Then the following is an exact sequence of Nisnevich sheaves on $X$ for all $r\ge 0$
\[0\to (\gamma^r G')_{(X,D)}\to (\gamma^rG)_{(X,D)}\to (\gamma^r G'')_{(X,D)}\to 0.\]
\end{lem}
\begin{proof}
From \eqref{para:MNST-gamma1} we get a long exact sequence for any $s\ge 0$
\[0\to G'_{(X,D)}\to G_{(X,D)}\to G''_{(X,D)}\to (\gamma^1G')_{(X,D)}\to\ldots 
\to (\gamma^sG)_{(X,D)}\to (\gamma^s G'')_{(X,D)}.\]
Thus it remains to show that  
$(\gamma^rG)_{(X,D)}\to (\gamma^r G'')_{(X,D)}$  is surjective, for every $r\le s$.
By \eqref{para:MNST-gamma1} this is equivalent to the surjectivity of
\[R^r\pi_* G_{(\P^r\times X,\pi^*D)}\to R^r\pi_* G''_{(\P^r\times X,\pi^*D)}.\]
This holds as $R^{r+1}\pi_*G'_{(\P^{r}\times X, \pi^*D)}=0$, by the blow-up formula \cite[Theorem 6.3]{BRS}. 
\end{proof}

For $D=\emptyset$ the following theorem is \cite[Theorem 11.8]{BRS}.

\begin{thm}\label{lem:HWM-gamma}
Let $\eta=\dlog c_1(\sO(1))\in H^1(\P^r_X, W_n\Omega^1_{\P^r_X, \log})$, 
where $c_1(\sO(1))\in H^1(\P^r_X, \G_m)$  
is first chern class of $\sO(1)$, and denote by $\eta^r\in H^r(\P^r_X, W_n\Omega^r_{\P^r_X,\log})$ 
the $r$-fold cup product of $\eta$.  Then via \eqref{para:MNST-gamma1} 
cupping with $\eta^r$ induces an isomorphism
\[-\cup \eta^r: W_n\Omega^{q-r}_{(X,D)}\xr{\simeq} R^r\pi_* W_n\Omega^q_{(\P^r\times X,\pi^*D)}
\cong (\gamma^rW_n\Omega^q)_{(X,D)}.\]
Furthermore this isomorphism is compatible with $F$, $V$, $R$, $\ul{p}$, and $d$.
\end{thm}
\begin{proof}
The last statement follows directly from the usual compatibilities of $V$, $F$, $R$, $\ul{p}$, and $d$ with the $\dlog$-map
and the naturality of the isomorphism  \eqref{para:MNST-gamma1}.
Set
\[\sX=(X,D),\quad \text{and}\quad \sP^r=(\P^r\times X, \pi^*D).\]
For $n=1$ observe that by  \eqref{para:HW-not5} we have a Künneth formula
\[\Omega^q_{\sP^r}=\bigoplus_{a+b=q} \pi^*\Omega^a_{\sX}\otimes_{\sO_{\sP^r}}\rho^*\Omega^b_{\P^r},\]
where $\rho: \P^r\times X\to \P^r$ denotes the projection.
Hence in this case the statement holds by the projection formula and the classical result for $\P^r$.
For $q> N+r$ both sides of the first isomorphism in the statement are zero.
Set
\[G^q_n:=\widetilde{W_n\Omega^q}_c ,\quad H^q:=\widetilde{\Omega^q}_c\oplus \widetilde{\Omega^{q+1}}_c.\]
We have $G^q_n$, $H^q\in \CI^{\tau, {\rm sp}}_{\Nis}$.
Denote by $I'$ the presheaf image of $(F^n, F^nd): G^q_{n+1}\to H^q$ and 
by $J'$ the presheaf image of $dV^n-V^n: H^q\to G^{q+1}_{n+1}$.
It follows that  $I'$, $J'\in \CI^{\tau, {\rm sp}}$, c.f. \cite[Lemma 1.32]{Saito-Purity}.
Thus by \cite[Theorem 0.4]{Saito-Purity} the sheafifications $I$ and $J$ lie in $\CI^{\tau, {\rm sp}}_{\Nis}$ 
and are the sheaf-theoretic images of the morphisms $(F^n, F^nd)$ and $dV^n-V^n$, respectively.
By \Cref{prop:long-mod-seq} we have short exact sequences for all $s\ge 0$
\[0\to G^q_{n, \sP^s}\xr{\ul{p}} G^q_{n+1, \sP^s}\to I_{\sP^s}\to 0, \quad \text{and}\quad
0\to I_{\sP^s}\to H^q_{\sP^s}\to J_{\sP^s}\to 0.\]
Hence Lemma \ref{lem:gamma-exact} gives a commutative diagram with exact rows 
\[\xymatrix{
0\ar[r] &
(\gamma^rG^q_n)_{\sX}\ar[r]^{\ul{p}}&
(\gamma^rG^q_{n+1})_{\sX}\ar[rr]^{(F^n, F^n d)}& &
(\gamma^r H^q)_{\sX}\ar[r]^{dV^{n}-V^n}&
(\gamma^r J)_{\sX}\ar[d]\\
        &     &   & &    &    (\gamma^r G^{q+1}_{n+1})_{\sX}\\
0\ar[r] &
G^{q-r}_{n, \sX}\ar[r]^{\ul{p}}\ar[uu]^{\alpha_1}&
G^{q-r}_{n+1, \sX}\ar[rr]^{(F^n, F^n d)}\ar[uu]^{\alpha_2}& &
H^{q-r}_{\sX}\ar[r]^{dV^{n}-V^n}\ar[uu]^{\alpha_3}&
G^{q+1-r}_{n+1, \sX},\ar[u]^{\alpha_4} 
}\]
where the upwards pointing vertical maps are all induced by cupping with  $\eta^r$ and \eqref{para:MNST-gamma1}.
Now $\alpha_1$ is an isomorphism by induction over $n$, 
$\alpha_3$ is an isomorphism by the case $n=1$, and $\alpha_4$ is an isomorphism
by descending induction over $q$. Hence $\alpha_2$ is an isomorphism as well.
\end{proof}

\begin{para}\label{para:pbf-buf-gt}
Thanks to Theorem \ref{lem:HWM-gamma} the projective bundle formula, the blow-up formula, and the 
Gysin triangle for general $G\in \CI^{\tau, {\rm sp}}_{\Nis}$ from \cite{BRS} 
can be made more explicit for $G=\widetilde{W_n\Omega^q}_c\in \CI^{\tau, {\rm sp}}_{\Nis}$:
\begin{enumerate}[label=(\arabic*)]
\item\label{para:pbf} Let $V$  be a locally free $\sO_X$-module of rank $r+1$ and denote by
$\pi: \P(V)\to X$ the corresponding projective bundle. Then there is a canonical isomorphism
in $D(W_n\sO_X)$
\[R\pi_* W_n\Omega^q_{(\P(V), \pi^*D)}\cong \bigoplus_{j=0}^{r} W_n\Omega^{q-j}_{(X,D)}[-j].\]
The morphism is induced by cupping with $c_1(\sO(1))^j\in \CH^j(\P(V))$, see \cite[Theorem 6.3]{BRS}.
In that theorem the isomorphism takes only place in the derived category of abelian sheaves but it is direct to check
that in the case at hand it is  $W_n\sO_X$-linear. Furthermore it is compatible with $F$, $V$, $R$, $\ul{p}$, and $d$.
For $D=\emptyset$ this isomorphism is due to Gros \cite[I, Theorem 4.1.11]{Gros}.
\item\label{para:buf}
Let $i:Z\inj X$ be a closed immersion of a smooth closed subscheme of pure codimension $r$ which intersects $D$ transversally,
i.e., the scheme-theoretic intersection $Z\cap D_1\cap\ldots\cap D_j$ is either empty or smooth of codimension $j$ in $Z$,
for any number of irreducible components $D_1,\ldots, D_j$ of $D_{\red}$. Let $\rho:\tilde{X}\to X$ be the blow-up in $Z$.
Then there is a canonical isomorphism in $D(W_n\sO_X)$
\[R\rho_* W_n\Omega^q_{(\tilde{X}, \rho^*D)}\cong W_n\Omega^q_{(X,D)}
\oplus 
\bigoplus_{j=1}^{r-1} i_* W_n\Omega^{q-j}_{(Z, D_{|Z})}[-j],\]
see \cite[Corollary 7.3]{BRS}. For $D=\emptyset$, this isomorphism is due to Gros \cite[IV, Corollaire 1.1.11]{Gros}.
\item\label{para:gt}
Let the notation be as in \ref{para:buf} and denote by $E=\rho^{-1}(Z)$ the exceptional divisor.
Then there is a canonical distinguished triangle in $D(W_n\sO_X)$ called the Gysin triangle
\[i_*W_n\Omega^{q-r}_{(Z, D_{|Z})}[-r]\xr{g}  
W_n\Omega^q_{(X,D)}\xr{\rho^*} R\rho_*W_n\Omega^q_{(\tilde{X}, \rho^*D+E)}\xr{\partial} 
i_*W_n\Omega^{q-r}_{(Z, D_{|Z})}[-r+1], \]
where $g=g_{(X,D)/(Z, D_{|Z})}$ denotes the Gysin morphism, see \cite[Theorem 7.16]{BRS}.
For $D=\emptyset$, this triangle is also spelled out in \cite[Corollary 11.10(2)]{BRS}.
\end{enumerate}
\end{para}

The Gysin triangle gives the following Lefschetz type theorem.

\begin{thm}\label{thm:top-Lef}
Assume  that $X$ is additionally projective and let $H\subset X$ be a smooth hypersurface section which intersects $D$ transversally
and satisfies
\[H^j(X, \Omega^N_{(X,D)}(H))=0, \quad \text{for all } j\ge 1,\]
where $\Omega^N_{(X,D)}(H)=\Omega^N_{(X,D)}\otimes_{\sO_X}\sO_X(H)$.
Then  the Gysin map
\[H^{j-1}(H, W_n\Omega^{N-1}_{(H, D_{|H})})\lra H^j(X, W_n\Omega^N_{(X,D)})\]
is an isomorphism for $j\ge 2$ and is surjective for $j=1$.
\end{thm}
\begin{proof}
By the Gysin triangle from \ref{para:pbf-buf-gt}\ref{para:gt} we have to show the vanishing
\eq{thm:top-Lef1}{H^j(X, W_n\Omega^{N}_{(X, D+H)})=0, \quad \text{for }j\ge 1.}
As $Z_n\Omega^N_X= \Omega^N_X$, we have 
\[Z_n\Omega^N_{(X,D+H)}= \Omega^N_{(X,D+H)}= \Omega^N_{(X,D)}(H),\]
where the second equality follows from $\Omega^N_X(\log (D+H))=\Omega^N_X(D_{\red}+H)$ and 
the first equality  holds by Lemma \ref{lem:BnZn}.
Thus \eqref{thm:top-Lef1} follows by induction over $n$ from Proposition \ref{prop:HWM-ex} 
noting that $B_{n}\Omega^{N+1}_{(X,D+H)}=0$.
\end{proof}

Finally we mention the following version of Serre type vanishing.

\begin{thm}\label{thm:Serre-van}
Let $H\subset X$ be a smooth and ample divisor intersecting $D$ transversally. Then for any large enough $m$ we have 
\[H^j(X, W_n\Omega^q_{(X,D+p^{n+1} mH)})=0, \quad \text{for all } j\ge 1.\]
\end{thm}
\begin{proof}
With the notation from \ref{para:Omega-n} and \ref{para:twisted-BZ} we have 
\[\Omega^q_n(D' , p(D_{n+1}+mE))=\Omega^q_n(D', pD_{n+1})\otimes_{\sO_X}\sO_X(pmE),\]
where  $D=D' +p^{n+1}D_{n+1}$ is a \pdd, and from \eqref{para:twisted-BZ3} we obtain inductively
\[B^q_{j,n}(D', p(D_{n+1}+mE))= B^q_{j,n}(D', pD_{n+1})\otimes_{\sO_X}\sO_X(pmE),\]
for all $j=0,\ldots, n$, and similarly for $B$ replaced by $Z$. 
Thus by Serre vanishing the cohomology of the sheaves $B^q_{j,n}(D', p(D_{n+1}+mE))$ and $Z^q_{j,n}(D', p(D_{n+1}+mE))$
vanishes in positive degrees,  for all $j$ and all $q$ and for any large enough $m$.
Hence the statement follows from Proposition \ref{prop:HWM-ex} together with Proposition \ref{prop:BZnn}.
\end{proof}

\section{The structure of Hodge-Witt sheaves with zeros}\label{sec:strHWZ}
By definition  a Hodge-Witt form has modulus $(X,D)$ if it is a regular form on $X\setminus D$ with certain pole constraints along $D$. 
In this section we introduce and study a notion of forms vanishing along $D$. We will see in the next section 
that these two notions of poles and zeros match up under duality.
The main result in this section is Theorem \ref{thm:HW-zeros}, which gives precise information about the structure of 
Hodge-Witt sheaves with zeros along $D$.

\medskip

Throughout this section we assume $X\in \Sm$ and we let $D$ be an effective Cartier divisor on $X$ such that $D_{\red}$ is an SNCD.
We remind the reader of our convention $\log(D)=\log(D_{\red})$.

\begin{defn}\label{defn:HW-zeros}
Suppose $D=\sum_i D_i$ with  $D_{i,\red}$ smooth.
Define 
$$W_n\Omega^q_{(X,-D)}:=\Ker\left(W_n\Omega^q_X\lra \bigoplus_i W_n\Omega^q_{D_i}\right),$$
where we view $D_i$ as a closed (in general non-reduced) subscheme of $X$.
For $n=1$, we write 
$$\Omega^q_{(X,-D)}:=W_1\Omega^q_{(X,-D)}.$$
\end{defn}

\begin{rmk}\label{rmk:Hyodo}
It follows from \cite[Lemme 3.15.1]{Mokrane}, see also \cite[Corollary 6.28]{Nakkajima}, that
the sheaf $W_n\Omega^{q}_{(X, -D_{\red})}$ coincides with  $W_n\Omega^q_X(-\log D_{\rm red})$ defined in \cite[1.]{Hyodo}.
To our knowledge  the sheaf $W_n\Omega^q_{(X,-D)}$ was not considered before for a  non-reduced divisor.
Moreover it seems it cannot be defined or studied by the machinery introduced in \cite[6.]{Nakkajima} as 
the Cartier isomorphism from \cite[(4.2.1.3)]{DI} does not extend to this situation.
For example if $D=nD_0$ with $D_0$ smooth connected, then the inverse Cartier operator induces an isomorphism (e.g. Lemma \ref{lem:twisted-qis})
\[\Omega^q_X(\log D_0)(-nD_0)\xr{\simeq}  H^q(F_*(\Omega^\bullet(\log D_0)(-pn+(p-1)D_0)))\]
and only for $n=1$ the right hand side is isomorphic to $H^q(F_*(\Omega^\bullet_X(\log D_0)(-nD_0)))$. 
\end{rmk}

\begin{para}\label{para:HW-zeros-FVRpd}
Note that since the $D_i$ may be non-reduced many of the structure results for the de Rham-Witt complex given, e.g., in 
\cite{IlDRW} or \cite{IllRay} cannot be applied to $W_n\Omega^q_{D_i}$. 
It is however direct from the definition that $W_n\Omega^q_{(X,-D)}$
is a  $W_{n}\sO_X$-submodule of $W_n\Omega^q_X$ and that the maps $F$, $V$, $R$, $\ul{p}$, $d$ restrict to give maps
\[\xymatrix{ 
W_{n-1}\Omega^q_{(X,-D)}\ar@<1ex>[r]^-{\ul{p},\, V} & W_n\Omega^q_{(X,-D)}\ar@<1ex>[l]^-{R,\, F}\ar[r]^d & W_n\Omega^{q+1}_{(X,-D)}.
}\]    
In fact $F$, $V$, $R$, and $d$ are defined on $W_\cdot\Omega^\bullet_T$ for all $\F_p$-schemes $T$, 
by \cite[I, Th\'eor\`eme 1.3 and Th\'eor\`eme 2.17]{IlDRW}. For the definition of $\ul{p}$ note, that 
$K_{n,D_i}:=\Ker(W_n\Omega^*_X\to W_n\Omega^*_{D_i})$ is the differential graded ideal of $W_n\Omega^*_X$ generated by
$W_n(I_{D_i})$, where $I_{D_i}\subset \sO_X$ is the ideal sheaf of the closed subscheme $D_i\subset X$,
see \cite[Lemma 1.2.2]{He04}. As the restriction map $W_{n+1}(I_{D_i})\to W_n(I_{D_i})$ is surjective,
it follows directly from the definition that $\ul{p}: W_{n}\Omega^*_X\to W_{n+1}\Omega^*_X$ 
restricts to $K_{n, D_i}\to K_{n+1, D_i}$ and hence also to $W_n\Omega^*_{(X,-D)}=\cap_i K_{n, D_i}\to W_{n+1}\Omega^*_{(X,-D)}$. 
\end{para}

\begin{lem}\label{lem:n=1zeros}
Let $D=D_0+pD_1$ be a $p$-divisibility decomposition, see  \ref{para:div-not}. 
Then 
$$\Omega^q_{(X,-D)}=\Omega^q_X(\log D_0)(-D).$$
\end{lem}
\begin{proof}
Let $A$ be a smooth $k$-algebra and $t\in A$ such that $R=A/(t)$ is smooth as well.
In this case there is an isomorphism $A/(t^m)\cong R[t]/(t^m)$, for all $m\ge 1$. We find
\eq{lem:n=1zeros1}{\Omega^q_{A/(t^m)}\cong 
\begin{cases}
(R[t]/(t^{m})\otimes_R\Omega^q_R) \oplus (R[t]/(t^{m-1})\otimes_R \Omega^{q-1}_R)dt, & \text{if } (m,p)=1\\
(R[t]/(t^{m})\otimes_R\Omega^q_{R}) \oplus (R[t]/(t^{m})\otimes_R \Omega^{q-1}_R)dt, &\text{if } p|m.
\end{cases}
}
Hence $\Omega^q_{A/(t^m)}$ is a free $R$-module and the map 
$\Omega^q_{A/(t^m)}\to \Omega^q_{A/(t^m)}\otimes_R {\rm Frac}(R)$ is injective. 

As $\Omega^q_X(\log D_0)(-D)$ is locally free and clearly  maps into $\Omega^q_{(X,-D)}$, we see from the above
that it suffices to prove the equality around the generic points of $D$. In which case it follows directly from \eqref{lem:n=1zeros1}.
\end{proof}

\begin{prop}\label{prop:n=1duality} 
Suppose $X$ has pure dimension $N$. For all $n\ge 0$ the map
\eq{prop:n=1duality1}{
F^n_{X*}\Omega^q_{(X,-D)}\simeq F^n_{X*}\sH om_{\sO_X}(\Omega_{(X,D)}^{N-q},\Omega^N_X)\xra{} 
\sH om_{\sO_X}(F^n_{X*}\Omega_{(X,D)}^{N-q},\Omega^N_X),}
given by
\[\alpha\mapsto (\beta\mapsto C^n(\alpha\wedge \beta)),\]
is an isomorphism of locally free $\sO_X$-modules, where $C^n: F^n_{X*}\Omega^N_X\to \Omega^N_X$
denotes the Cartier operator and $\Omega^q_{(X,D)}$ is defined  in \ref{para:HW-not}.
\end{prop} 
\begin{proof}
Let $D=D_0+pD_1$ be a $p$-divisibility decomposition. In the case $n=0$ the isomorphism
follows from twisting the isomorphism 
$\Omega^q_X(\log D_0)\cong \sH om (\Omega^{N-q}_X(\log D_0), \Omega^N_{X}(\log D_0))$
by $\sO_X(-D)$ together with the explicit formulas given in Lemma \ref{lem:n=1zeros}, in \eqref{para:HW-not5}, and 
the isomorphism $\Omega^N_X(\log D_0)\cong \Omega^N_X(D_{0,\red})$. Applying $F^n_{X*}$ yields the first isomorphism in
\eqref{prop:n=1duality1}. The second map is an isomorphism as well. Indeed as $X$ is smooth over a perfect field of characteristic $p$
we have a natural isomorphism of $\sO_X$-modules $(F^n_X)^!\Omega^N_X\cong\Omega^N_X$ and the composition
\eq{prop:n=1duality2}{F^n_{X*}\Omega^N_X\cong F^n_{X*}(F^n_X)^!\Omega^N_X\xr{\Tr_{F^n_X}} \Omega^N_X,}
where $\Tr_{F^n_X}$ is the counit of adjunction, is equal to the $n$-fold iteration of the Cartier operator $C^n$,
e.g. \cite[II, Lemma 2.1]{Ekedahl}. Hence the second map is an isomorphism by adjunction.
\end{proof}

\begin{para}\label{para:Omega/BZ}
By \cite[I, Corollaire 3.9]{IlDRW} we have an exact sequence of $W_{n+1}\sO_X$-modules
\eq{para:Omega/BZ1}{
0\lra \frac{F_{X*}^n\Omega^q_{X}}{B_n\Omega^q_{X}}
\xra{V^n}
\Ker\left(W_{n+1}\Omega^q_{X}\xra{R} W_n\Omega^q_{X}\right)
\xra{\beta}
\frac{F_{X*}^n\Omega^{q-1}_{X}}{Z_n\Omega^{q-1}_{X}}
\lra 0,}
where the $W_{n+1}\sO_X$-module structure on the two outer $\sO_X$-modules is induced by $R^n:W_{n+1}\sO_X\to \sO_X$ and 
the map $\beta$ is given by $\beta(V^n (a)+dV^n(b))=b$.
Moreover we have an injection of $W_{n+1}\sO_X$-modules
\eq{para:Omega/BZ2}{0\to \frac{F_{X*}^n\Omega^{q-1}_{X}}{Z_n\Omega^{q-1}_{X}}\xra{dV^{n-1}} F_*W_n\Omega^q_X.}

For $n=0$ and $q\ge 0$, define
the $\sO_X$-modules 
$$(\Omega/B)^q_{0,(X,-D)}=\Omega^q_{(X,-D)}\text{\quad and \quad}(\Omega/Z)^{q-1}_{0, (X-D)}=0.$$
For $n\ge 1$, $q\ge 0$ we define the $W_{n+1}\sO_X$-modules $(\Omega/B)^q_{n,(X,-D)}$ and $(\Omega/Z)^{q-1}_{n, (X-D)}$
by requiring that the following two diagrams are cartesian in the category of $W_{n+1}\sO_X$-modules
\eq{para:Omega/BZ3}{
\xymatrix{
(\Omega/B)^q_{n, (X,-D)}\ar[r]\ar[d] & W_{n+1}\Omega^q_{(X,-D)}\ar[d]\\
\frac{F_{X*}^n\Omega^q_{X}}{B_n\Omega^q_{X}}\ar[r]^{V^n} & W_{n+1}\Omega^q_X,
}\qquad 
\xymatrix{
(\Omega/Z)^{q-1}_{n, (X,-D)}\ar[r]\ar[d] & F_*W_n\Omega^{q}_{(X,-D)}\ar[d]\\
\frac{F_{X*}^n\Omega^{q-1}_{X}}{Z_n\Omega^{q-1}_{X}}\ar[r]^{dV^{n-1}} & F_*W_n\Omega^{q}_X.
}
}
Note that we can view $(\Omega/B)^q_{n,(X,-D)}$ and $(\Omega/Z)^{q-1}_{n,(X,D)}$ as $\sO_X$-submodules of 
$F_{X*}^n\Omega^q_{X}/B_n\Omega^q_{X}$ and $F_{X*}^n\Omega^{q-1}_{X}/Z_n\Omega^q_{X}$, respectively.
In the following we will denote the top horizontal maps in the two diagrams as well by $V^n$ and $dV^{n-1}$, respectively.

Some consequences from this definition:
\begin{enumerate}[label= (\alph*)]
\item 
\label{para:Omega/BZa}
It is direct from the definition that
$$(\Omega/B)^q_{1,(X,-D)}=\{
a\in F_{X*}\Omega^q_X/\sB^q(\Omega^{\bullet}_X)
\mid 
V(a)\in W_2\Omega^q_{(X,-D)}\},$$
$$(\Omega/Z)^{q-1}_{1, (X-D)}=
\frac{F_{X*}\Omega^{q-1}_{(X,-D)}}{Z_1\Omega^{q-1}_X\cap \Omega^{q-1}_{(X,-D)}}
=\frac{F_{X*}\Omega^{q-1}_{(X,-D)}}{\sZ^{q-1}(\Omega^{\bullet }_{(X,-D)})}.$$

\item 
\label{para:Omega/BZb}
For any $n\ge 1$, the  diagram 
$$\xymatrix{
(\Omega/B)^{q-1}_{n, (X,-D)}\ar@{.>}[r]\ar[d] \ar@/^1.5pc/[rr]^{Fd\circ V^{n}}&
(\Omega/Z)^{q-1}_{n, (X,-D)}\ar[r]\ar[d] &
F_*W_n\Omega^{q}_{(X,-D)}\ar[d]
\\
\frac{F_{X*}^{n}\Omega^{q-1}_{X}}{B_{n}\Omega^{q-1}_{X}}\ar@{->>}[r]&
\frac{F_{X*}^{n}\Omega^{q-1}_{X}}{Z_n\Omega^{q-1}_{X}}\ar[r]^{dV^{n-1}} & F_*W_n\Omega^{q}_X.
}$$
shows that the surjection $\Omega^q_{X}/B_{n}\Omega^q_{X}\to \Omega^q_{X}/Z_{n}\Omega^q_{X}$ induces a well-defined map
\eq{para:Omega/BZ-hn}{(\Omega/B)^q_{n,(X,-D)}\to (\Omega/Z)^q_{n, (X,-D)}.}

\item 
\label{para:Omega/BZc}
For $n\ge 0$, 
the diagram
$$\xymatrix{
F_{X*}(\Omega/B)^q_{n, (X,-D)}\ar@{.>}[r]\ar[d]\ar@/^1.5pc/[rr]^{V\circ V^n}&
(\Omega/B)^q_{n+1, (X,-D)}\ar[r]\ar[d] & W_{n+2}\Omega^q_{(X,-D)}\ar[d]
\\
\frac{F_{X*}^{n+1}\Omega^q_{X}}{F_{X*}B_{n}\Omega^q_{X}}\ar@{->>}[r]&
\frac{F_{X*}^{n+1}\Omega^q_{X}}{B_{n+1}\Omega^q_{X}}\ar[r]^{V^{n+1}} & W_{n+2}\Omega^q_X,
}$$ 
shows that the surjection $\Omega^q_{X}/B_{n}\Omega^q_{X}\to \Omega^q_{X}/B_{n+1}\Omega^q_{X}$ induces a well-defined map
$$
F_*(\Omega/B)^{q-1}_{n,(X,-D)}\ra 
(\Omega/B)^{q-1}_{n+1,(X,-D)}.
$$
Similarly, for $n\ge 1$, 
the diagram
$$\xymatrix{
(\Omega/Z)^{q-1}_{n+1, (X,-D)}\ar@{.>}[r]\ar[d] \ar@/^1.5pc/[rr]^{F\circ dV^{n}}&
F_{X*}(\Omega/Z)^{q-1}_{n, (X,-D)}\ar[r]\ar[d] &
F^2_*W_n\Omega^{q}_{(X,-D)}\ar[d]
\\
\frac{F_{X*}^{n+1}\Omega^{q-1}_{X}}{Z_{n+1}\Omega^{q-1}_{X}}\ar@{->>}[r]&
\frac{F_{X*}^{n+1}\Omega^{q-1}_{X}}{F_{X*}Z_n\Omega^{q-1}_{X}}\ar[r]^{dV^{n-1}} & F^2_*W_n\Omega^{q}_X.
}$$
shows that the surjection $\Omega^q_{X}/Z_{n+1}\Omega^q_{X}\to \Omega^q_{X}/Z_{n}\Omega^q_{X}$ induces a well-defined map
\eq{prop:BZnn-zeros3}{(\Omega/Z)^{q-1}_{n+1,(X,-D)}\to F_*(\Omega/Z)^{q-1}_{n,(X,-D)}.}
Here we have used the fact that $F_*: (W_{n+1}\sO_X\text{-mod})\to (W_n\sO_X\text{-mod})$ has a left adjoint and hence the 
diagrams in \eqref{para:Omega/BZ3} stay cartesian after applying $F_*$.

\end{enumerate}
\end{para}

The aim of this section is to prove the following theorem:

\begin{thm}\label{thm:HW-zeros}
Let the notation be as above and assume $X$ is of pure dimension $N$.
\begin{enumerate}[label=(\arabic*)]
\item\label{thm:HW-zeros1} There is an exact sequence of $W_{n+1}\sO_X$-modules
$$0\lra (\Omega/B)^q_{n, (X,-D)} \xra{V^n}
\Ker\left(W_{n+1}\Omega^q_{(X,-D)}\xra{R} W_n\Omega^q_{(X,-D)}\right)
\xra{\beta}
(\Omega/Z)^{q-1}_{n, (X,-D)}
\lra 0,
$$
where the map $\beta$ is induced from the map $\beta$ in \eqref{para:Omega/BZ1}.
\item\label{thm:HW-zeros1.5} The restriction map $R:W_{n+1}\Omega^q_{(X,-D)}\to W_n\Omega^q_{(X,-D)}$ is surjective.
\item\label{thm:HW-zeros2}
The isomorphism \eqref{prop:n=1duality1} induces isomorphisms of locally free coherent $\sO_X$-modules
\[
(\Omega/B)^q_{n, (X,-D)} \xra{\simeq} \sH om(Z_n\Omega^{N-q}_{(X,D)},\Omega_X^N)
,\quad 
(\Omega/Z)^{q-1}_{n, (X,-D)}
\xra{\simeq}
\sH om(B_n\Omega^{N-q+1}_{(X,D)},\Omega_X^N),\]
where $Z_n\Omega^{N-q}_{(X,D)}$ and $B_n\Omega^{N-q+1}_{(X,D)}$ are defined in \eqref{thm:strHWM1}.
\end{enumerate}
\end{thm}

For $D=\emptyset$ part \ref{thm:HW-zeros1} holds by \cite[I, Corollaire 3.9]{IlDRW}, part \ref{thm:HW-zeros1.5} holds by definition of the 
de Rham-Witt complex, and  part \ref{thm:HW-zeros2} by \cite[II, Lemma 2.2.20]{Ekedahl}.
It will take the rest of this section  to prove this theorem.
Part \ref{thm:HW-zeros1} is proven in Proposition \ref{prop-ex-seq-zeros} and part \ref{thm:HW-zeros1.5} in Lemma \ref{lem:R-surj-zeros},
finally, part \ref{thm:HW-zeros2} follows from Proposition \ref{prop:BZnn-zeros} and Proposition \ref{prop:BZ-duality} together with 
Proposition \ref{prop:BZnn}.

\begin{para}\label{para:Omega-zeros}
Let $E$ be an effective Cartier divisor on $X$ such that $D_{\red}+E_\red$ is a reduced SNCD.
We will use the following notation
\eq{para:Omega-zeros1}{
\Omega^q_{n}(-D,-E):=\Omega^q_X(\log D)(-\lceil D/p^n \rceil-E),\qquad n,\,q\ge 0.}
This is a locally free coherent $\sO_X$-module. 
Let $D=D'+p^nD_n$ be a \pdd.
Similar to \ref{para:Omega-n}, we observe
\begin{enumerate}[label= (\alph*)]
\item
\label{para:Omega-zeros2}
$\Omega^q_{n}(-D',-D_n-E)\subset \Omega^q_{n}(-D,-E).$
This inclusion is strict if $D_n\neq 0$.
\item\label{para:Omega-zeros4}
There is a well-defined differential map 
$d:\Omega^q_{n}(-D, -pE)\to \Omega^{q+1}_{n}(-D,-pE)$
induced from the differential map of $\Omega^\bullet_X$.
Moreover we have an inclusion of complexes
$$\Omega^\bullet_{n-1}(-D',-pD_n-pE)\subset \Omega^\bullet_{n-1}(-D,-pE).$$
\item 
The natural map 
\[\Omega^\bullet_X(\log D)(-p \lceil D/p^n\rceil-pE)\to \Omega^\bullet_{n-1}(-D,-pE)\]
is a quasi-isomorphism  by Lemma \ref{lem:twisted-qis}. Hence the inverse Cartier operator induces an isomorphism of $\sO_X$-modules
\eq{para:Omega-zeros6}{ 
C^{-1}: \Omega^q_{n}(-D,-E)
\xr{\simeq} 
\sH^q(F_{X*}\Omega^\bullet_{n-1}(-D,-pE)).}
Moreover we obtain a commutative diagram of locally free $\sO_X$-modules
\eq{para:Omega-zeros7}{
\xymatrix{
\Omega^q_{n}(-D',-D_n-E)
\ar[r]^-{C^{-1}}_-{\simeq}
\ar@{^(->}[d]   & 
\sH^q(F_{X*}
\Omega^\bullet_{n-1}(-D', -pD_n-pE))
\ar[d]
\\
\Omega^q_{n}(-D, -E)
\ar[r]^-{C^{-1}}_-{\simeq} & 
\sH^q(F_{X*}\Omega^\bullet_{n-1}(-D, -pE)).
}
}
As a consequence the vertical map on the right is injective as well.
\end{enumerate}
\end{para}

\begin{para}\label{para:C-exs-eq-zeros}
Let  $X$, $D$, and $E$ be as in \ref{para:Omega-zeros} above.
Let $n\ge 1$. 
The inverse Cartier isomorphism defines the following exact sequence of locally free $\sO_X$-modules
\eq{para:C-ex-seq-zeros1}{
0\to \Omega^q_n(-D, -E)\xr{C^{-1}} \frac{F_{X*}\Omega^q_{n-1}(-D, -pE)}{\sB^q(F_{X*}\Omega^\bullet_{n-1}(-D, -pE))}
\to
\frac{F_{X*}\Omega^q_{n-1}(-D, -pE)}{\sZ^q(F_{X*}\Omega^\bullet_{n-1}(-D, -pE))}
\to 0.}
The fact that these modules are locally free follows via descending induction on $q$ by considering as well the exact sequence
\[0\to
\frac{F_{X*}\Omega^q_{n-1}(-D, -pE)}{\sZ^q(F_{X*}\Omega^\bullet_{n-1}(-D, -pE))}
\xr{d}
F_{X*}\Omega^{q+1}_{n-1}(-D, -pE)
\to 
\frac{F_{X*}\Omega^{q+1}_{n-1}(-D, -pE)}{\sB^{q+1}(F_{X*}\Omega^\bullet_{n-1}(-D, -pE))}
\to 0.\]

\end{para}

\begin{lem}\label{lem:dzeros}
Let $X$, $D$, and $E$ be as above.
Then for $n\ge 0$
\[\sB^{q}(\Omega^\bullet_X)\cap \Omega^{q}_n(-D,-pE)=\sB^{q}(\Omega^\bullet_n(-D,-pE)).\]
In particular 
\[\sB^{q}(\Omega^\bullet_X)\cap \Omega^{q}_{(X,-D)}=\sB^{q}(\Omega^\bullet_{(X,-D)}).\]
\end{lem}
\begin{proof}
The ``in particular'' part follows from the first statement and Lemma \ref{lem:n=1zeros}.
For the  first statement we only have to show this ``$\subset$'' inclusion.
As the $\sO_X$-module $\sB^{q}(F_{X*}\Omega^\bullet_n(-D,-pE))$ is locally free by \ref{para:C-exs-eq-zeros}
it suffices to show the inclusion around each generic point of $D$ or $E$,
c.f. the argument in the proof of \eqref{lem:BnZn1}. Since the question is moreover local in the Nisnevich topology,
we may assume that $X=\Spec R[t]$, with $R$ a smooth $k$-algebra, and $D$ or $E$ is equal to $\Div(t^r)$, $r\ge 0$.
Thus $\Omega^q_{R[t]}=\Omega^q_{R}[t]\oplus \Omega^{q-1}_{R}[t]dt$,
where for an $R$-module $M$ we denote by $M[t]$ the free $R$-module $\oplus_{j\ge 0} M t^j$.
In this situation  $\Omega^q_n(-D,-pE)$ has one of the following forms
\begin{enumerate}[label=(\alph*)]
\item\label{lem:dzeros1} 
$t^{s}\Omega^q_{R}[t]\oplus t^{s-1}\Omega^{q-1}_{R}[t]dt$, with $s=\lceil r/p^n\rceil$, 
when $D=\Div(t^r)$;
\item\label{lem:dzeros2} 
$t^{pr}\Omega^q_{R[t]}$, when $E=\Div(t^r)$.
\end{enumerate}
We consider the case \ref{lem:dzeros1}. Let 
$a=a_0+ \sum_{j\ge 1} (a_j t^j+ b_j t^{j-1}dt)\in \Omega^{q-1}_{R[t]}$, with $a_j\in \Omega^{q-1}_R$ and  
$b_j\in \Omega^{q-2}_R$, and assume $da\in t^{s}\Omega^q_{R}[t]\oplus t^{s-1}\Omega^{q-1}_{R}[t]dt$.
Set 
\[a':= a_0 +
\sum_{\stackrel{1\le j\le s-1}{p|j}} 
(a_jt^j+ b_jt^{j-1}dt)+ 
      \sum_{\stackrel{1\le j\le s-1}{p\nmid j}} \tfrac{(-1)^{q-2}}{j} d(b_j t^j) \in \Omega^{q-1}_{R[t]}.\]
Then it is direct to check that 
\[a-a'\in t^{s}\Omega^{q-1}_R[t] \oplus t^{s-1}\Omega^{q-2}_{R}[t]dt \quad \text{and}\quad da=d(a-a').\]
In case \ref{lem:dzeros2}, take $a\in \Omega^{q-1}_{R[t]}$ with $da= t^{pr}b$, with $b\in \Omega^{q}_{R[t]}$.
Applying the Cartier operator yields $0=t^rC(b)$. Thus $C(b)=0$, thus $b=dc$, for some 
$c\in \Omega^{q-1}_{R[t]}$. Hence $da=d(t^{pr}c)$.

Both cases together imply the statement.
\end{proof}

\begin{lem}
\label{lem:pHWM-zeros}
Let $D=D'+p^{n+1}D_{n+1}$ be a \pdd. 
For $\alpha\in \Omega^q_X$ we have
\[\ul{p}^n(\alpha)\in  W_{n+1}\Omega^q_{(X,-D)} \Longleftrightarrow  \alpha\in \Omega^q_{n}(-D',-pD_{n+1}).\]
Moreover, the restriction map 
\[R^j: W_{n+j+1}\Omega^q_{(X,-D)}\cap \ul{p}^n(W_{j+1}\Omega^q_X)\to W_{n+1}\Omega^q_{(X,-D)}\cap\ul{p}^n(\Omega^q_X)\]
is surjective. 
\end{lem}
\begin{proof}
We first show that for any germ $\alpha\in \Omega^q_n(-D', -pD_{n+1})_x$, with  $x\in D$ a {\em closed} point,
we find a lift $\tilde{\alpha}\in W_{j+1}\Omega^q_{X,x}$ such that $\ul{p}^n(\tilde{\alpha})\in W_{n+j+1}\Omega^q_{(X,-D),x}$.
Taking $j=0$ yields this ``$\Leftarrow$" direction and after the other direction will be proven below, it will also imply the 
``Moreover"-part.
Let $A=\sO_{X,x}$. We find an \'etale map $k[t_1,\ldots, t_d]\to A$, such that 
$t_1,\ldots,t_d$ form a regular sequence of parameters for $A$ and on $\Spec A$ we have
\[D'=\Div(t_1^{m_1}\cdots t_r^{m_r}) \quad\text{and}\quad D_{n+1}=\Div(t_{r+1}^{m_{r+1}}\cdots t_s^{m_s}),\]
with $0\le r\le s\le d$ and $p^{n+1}\nmid m_i$, for $i=1,\ldots,r$ (with the convention that $r=0$ means that $D'=\emptyset$ and $r=s$ means that 
$D_{n+1}=\emptyset$). 
Thus 
\[D=\Div(t_1^{m_1}\cdots t_r^{m_r}\cdot t_{r+1}^{p^{n+1}m_{r+1}}\cdots t_s^{p^{n+1}m_s}).\]
A basis of the free $A$-module 
$\Omega^q_n(-D',-pD_{n+1})_x=
\Omega_X^q(\log D')(-\lceil D'/p^n\rceil-pD_{n+1})_x$ 
is given by 
\[e_{I,J} = 
t_1^{\lceil\frac{m_1}{p^n}\rceil}\cdots 
t_r^{\lceil\frac{m_r}{p^n}\rceil} 
t_{r+1}^{pm_{r+1}}\cdots t_{s}^{pm_{s}}\cdot
\dlog \ul t_I\dlog (\ul{1+t}_J),\]
where $I$ and $J$ run through the tuples $I=(1\le i_1<\ldots<i_{q_1}\le r)$ and $J=(r+1\le j_1<\ldots< j_{q_2}\le d)$ with $q_1+q_2=q$
and where we use the notations 
\[ \ul{t}_I= \{t_{i_1}, \ldots, t_{i_{q_1}}\}\in K^M_{q_1}(A[\tfrac{1}{t_{i_1}\cdots t_{i_{q_1}}}])
\quad \text{and}\quad \ul{1+ t}_J= \{ 1+t_{j_1},\ldots, 1+ t_{j_{q_2}}\}\in K^M_{q_2}(A).\]
Set
\[\tilde{e}_{I,J}:= [t_1]^{\lceil\frac{m_1}{p^n}\rceil}\cdots [t_r]^{\lceil\frac{m_r}{p^n}\rceil} [t_{r+1}]^{p m_{r+1}}\cdots 
[t_{s}]^{pm_{s}}\cdot \dlog \ul t_I\dlog (\ul{1+t}_J)\in W_{j+1}\Omega^q_A.\]
Note that $\ul p^{n}: \iota_*W_{j+1}\Omega^q_X \ra W_{n+j+1}\Omega^q_X$ is $W_{n+j+1}\sO_X$-linear, where $\iota:W_{j+1}X\ra W_{n+j+1}X$ 
denotes the map of schemes induced by the restriction $R^j:W_{n+j+1}\sO_X\ra W_{j+1}\sO_X$.
Since  $W_{n+j+1}\Omega^q_{(X,-D),x}$ is a $W_{n+j+1}A$-submodule of $W_{n+j+1}\Omega^q_{X,x}$, it suffices to show that
$\ul{p}^n(\tilde{e}_{I,J})\in W_{n+j+1}\Omega^q_{(X,-D),x}$. Indeed, if $\alpha=\sum_{I,J} a_{I,J} e_{I,J}\in \Omega^{q}_n(-D',-pD_{n+1})_x$, 
$a_{I,J}\in A$, is an arbitrary element, then we can choose any lifts $\tilde{a}_{I,J}\in W_{j+1}(A)$ of $a_{I,J}$, 
and the element $\tilde{\alpha}:=\sum_{I,J}\tilde{a}_{I,J} \tilde{e}_{I,J}$ is a lift of $\alpha$ and satisfies 
$\ul{p}^n(\tilde{\alpha})\in W_{n+j+1}\Omega^q_{(X,-D),x}$.

Thus it remains to show that the element
\[\ul p^{n}(\tilde{e}_{I,J})= p^n [t_1]^{\lceil\frac{m_1}{p^n}\rceil}\cdots [t_r]^{\lceil\frac{m_r}{p^n}\rceil} [t_{r+1}]^{p m_{r+1}}\cdots 
[t_{s}]^{pm_{s}}\cdot \dlog \ul t_I\dlog (\ul{1+t}_J)\in W_{n+j+1}\Omega^q_A\]
lies in $W_{n+j+1}\Omega^q_{(X,-D),x}$.
To this end we have to show that $\ul p^{n}(\tilde{e}_{I,J})$ vanishes in $W_{n+j+1}\Omega^q_{A/t_i^{m}}$,
where $m=m_i$, if $1\le i\le r$,  and $m=p^{n+1}m_i$, if $r+1\le i\le s$.
We  consider the following four cases:
\begin{itemize}
\item  If $r+1\le i\le s$, then
\[\ul p^{n}(\tilde{e}_{I,J})=p^n [t_i^{pm_i}]\alpha = V^n([t_i]^{p^{n+1}m_i}) \cdot \alpha,\quad \text{for some }\alpha \in W_{n+j+1}\Omega^q_A.\]
\item If $1\le i\le r$ and $i\not \in I$, then $p^n\lceil\frac{m_i}{p^n}\rceil= m_i+e$, for some $e\ge 0$, and thus
\[\ul p^{n}(\tilde{e}_{I,J})=p^n [t_i]^{\lceil\frac{m_i}{p^n}\rceil}\alpha  = V^n([t_i]^{m_i+e})\cdot \alpha,\quad 
\text{for some }\alpha \in W_{n+j+1}\Omega^q_A.\]
\item If $1\le i\le r$, $i\in I$, and $p^n\nmid m_i$, then we have $p^n\lceil m_i/p^n\rceil=m_i+e$, for some $e\ge 1$, and thus
\[\ul p^{n}(\tilde{e}_{I,J})= p^n [t_i]^{\lceil\frac{m_i}{p^n}\rceil}\dlog t_i \cdot \beta = V^n([t_i]^{m_i} [t_i]^e\dlog t_i)\cdot \beta,\quad 
\text{for some }\beta \in W_{n+j+1}\Omega^{q-1}_A.\]
\item If $1\le  i \le r$, $i\in I$, and $p^n | m_i$, then $\mu= m_i/p^n$ is prime to $p$ and we find
\[\ul p^{n}(\tilde{e}_{I,J})= p^n [t_i]^{\mu}\dlog t_i \cdot \beta= \tfrac{1}{\mu} dV^n([t_i]^{m_i})\cdot \beta, \quad 
\text{for some }\beta\in W_{n+j+1}\Omega^{q-1}_A.\]
\end{itemize}
Therefore the vanishing of $\ul p^{n}(\tilde{e}_{I,J})$ in $W_{n+j+1}\Omega^q_{A/t_i^m}$ holds in all four cases.

\medskip

Next we prove this ``$\Rightarrow$" direction. Since $\Omega^q_n(-D', -pD_{n+1})$ is locally free on a smooth scheme,
it suffices to show the statement around all generic points of $D$. Since the statement is furthermore Nisnevich local 
we may assume $X=\Spec R[t]$, with $R$ a smooth $k$-algebra, and $D=\Div(t^m)$.
Let $R_{n+1}$ be a smooth lift of $R$ over $W_{n+1}(k)$. By \cite[III, (1.5)]{IllRay} and \cite[Proposition 8.4]{BER}
there is a unique  injective map
\[\tilde{F}^{n+1}: W_{n+1}\Omega^q_{R[t]}\to \Omega^q_{R_{n+1}[t]}/\sB^q(\Omega^\bullet_{R_{n+1}[t]}),\] 
which makes the following diagram commutative
\[\xymatrix{
W_{n+2}\Omega^q_{R_{n+1}[t]}\ar[r]^{F^{n+1}}\ar[d] & \Omega^q_{R_{n+1}[t]}\ar[d] \\
W_{n+1}\Omega^q_{R[t]}\ar[r]^{\tilde{F}^{n+1}} & \Omega^q_{R_{n+1}[t]}/\sB^q,
}\]
where we use that the de Rham-Witt complex exists for $R_{n+1}[t]$ as well, see \cite{HeMa04}.
By \cite[Lemma 1.2.2]{He04} 
\[W_{n+1}\Omega^q_{(R[t], t^m)}:=\Ker(W_{n+1}\Omega^q_{R[t]}\to W_{n+1}\Omega^q_{R[t]/t^m})\]
is the degree $q$ part of the differential graded ideal 
generated by $W_{n+1}(t^mR[t])$, i.e.,
\[W_{n+1}\Omega^q_{(R[t], t^m)}= W_{n+1}(t^m R[t])\cdot W_{n+1}\Omega^q_{R[t]}+ 
d(W_{n+1}(t^m R[t]))\cdot W_{n+1}\Omega^{q-1}_{R[t]}.\]
Thus $\tilde{F}^{n+1}(W_{n+1}\Omega^q_{(R[t], t^m)})$ is contained in the image of the following  group in $\Omega^q_{R_{n+1}[t]}/\sB^q$
\[t^{pm}\Omega^q_{R_{n+1}[t]} + m t^{pm}\dlog(t)\Omega^{q-1}_{R_{n+1}},\]
where the second summand contains  $\tilde{F}^{n+1}(dV^n(t^mR)W_{n+1}\Omega^{q-1}_R)$. 
More precisely, set
\[\Omega^q_{e}:= t^e(\Omega^q_{R_{n+1}}\oplus \Omega^{q-1}_{R_{n+1}}\dlog t), \quad \text{if } p^{n+1}\nmid e\ge 1,\]
and 
\[\Omega^q_{0,0}:= \Omega^q_{R_{n+1}}, \quad 
\Omega^q_{e,0}:= t^e\Omega^q_{R_{n+1}}, \quad 
\Omega^q_{e,1}:= t^e\Omega^{q-1}_{R_{n+1}}\dlog t, \quad \text{if } p^{n+1}\mid e\ge 1.\]
Then we obtain a direct sum decomposition of complexes
\[\Omega^{\bullet}_{R_{n+1}[t]} = \bigoplus_{e\ge 1,\, p^{n+1}\nmid e} \Omega^{\bullet}_e \oplus 
\bigoplus_{e\ge 0, \,p^{n+1}\mid e} \Omega^{\bullet}_{e,0} \oplus \bigoplus_{e\ge 1,\, p^{n+1}\mid e} \Omega^{\bullet}_{e,1}\]
and hence also a direct sum decomposition
\[\frac{\Omega^q_{R_{n+1}[t]}}{\sB^q(\Omega^{\bullet}_{R_{n+1}[t]})} = 
\bigoplus_{e\ge 1,\, p^{n+1}\nmid e} \frac{\Omega^q_e}{\sB^q(\Omega^{\bullet}_e)} \oplus 
\bigoplus_{e\ge 0, \,p^{n+1}\mid e} \frac{\Omega^q_{e,0}}{\sB^q(\Omega^{\bullet}_{e,0})} \oplus 
\bigoplus_{e\ge 1,\, p^{n+1}\mid e} \frac{\Omega^{q}_{e,1}}{\sB^q(\Omega^{\bullet}_{e,1})}.\]
With this notation the discussion above shows that $\tilde{F}^{n+1}(W_{n+1}\Omega^q_{(R[t], t^m)})$ is contained
in
\eq{lem:pHWM-zeros1}{
\bigoplus_{e\ge pm,\, p^{n+1}\nmid e} \frac{\Omega^q_e}{\sB^q(\Omega^{\bullet}_e)} \oplus 
\bigoplus_{e\ge pm, \,p^{n+1}\mid e} \frac{\Omega^q_{e,0}}{\sB^q(\Omega^{\bullet}_{e,0})} \oplus 
\bigoplus_{e\ge pm+1,\, p^{n+1}\mid e} \frac{\Omega^{q}_{e,1}}{\sB^q(\Omega^{\bullet}_{e,1})}
\oplus \frac{m\Omega^q_{pm,1}}{\sB^q(\Omega^\bullet_{pm,1})\cap m\Omega^q_{pm,1}},
}
where the last summand only occurs if $p^{n+1}|pm$ and $p^{n+1}\nmid m$.

We show that 
\[\ul{p}^n:\Omega^q_{R[t]}/\Omega^q_n(-D',-pD_{n+1})\to W_{n+1}\Omega^q_{R[t]}/W_{n+1}\Omega^q_{(R[t], t^m)}\]
is injective. Note that in the above situation we have  two cases, where $s=\lceil m/p^n\rceil$,
\[\Omega^q_n(-D',-pD_{n+1})=
\begin{cases}
t^s(\Omega^q_{R}[t]+ \Omega^q_R[t]\dlog t)    & \text{if } p^{n+1}\nmid m\\
t^s\Omega^q_{R[t]} & \text{if } p^{n+1}\mid m.
\end{cases}\]
Thus any element in the quotient $\Omega^q_{R[t]}/\Omega^q_n(-D',-pD_{n+1})$ can be written in the form
\[\alpha = a_0+ \sum_{j\ge 1}^{s-1} t^j (a_j + b_j \dlog t)+ b_s t^s \dlog t,\]
for $a_i\in \Omega^q_R$ and $b_j\in \Omega^{q-1}_R$, where $b_s=0$, if $p^{n+1}\nmid m$. 
We obtain 
\[\tilde{F}^{n+1}(\ul{p}^n\alpha)\equiv  \sum_{j=0}^{s-1} p^nt^{jp^{n+1}}\tilde{F}^{n+1}(\tilde{a}_j) + 
\sum_{j=1}^s p^nt^{jp^{n+1}}\tilde{F}^{n+1}(\tilde{b}_j)\dlog t \quad\text{mod }\sB^q,\]
where $\tilde{a}_i$ and $\tilde{b}_j$ are lifts of $a_i$ and $b_j$ to $W_{n+1}\Omega^q_R$ and $W_{n+1}\Omega^{q-1}_R$, respectively.
As $s-1<m/p^n$ we obtain 
\[\tilde{F}^{n+1}(\ul{p}^n\alpha)\in 
\bigoplus_{\stackrel{0\le e< pm}{p^{n+1}\mid e}}\frac{\Omega^q_{e,0}}{\sB^q(\Omega^\bullet_{e,0})}\oplus
\bigoplus_{\stackrel{1\le e< pm}{p^{n+1}\mid e}}\frac{\Omega^q_{e,1}}{\sB^q(\Omega^\bullet_{e,0})}
\oplus \frac{\Omega^q_{pm,1}}{\sB^q(\Omega^q_{pm,1})},\]
where the last summand only occurs if $p^{n+1}\mid m$. 
Now assume $\ul{p}^n\alpha\in W_{n+1}\Omega^q_{(R[t], t^m)}$. 
Then $\tilde{F}^{n+1}(\ul{p}^n\alpha)$ lies in \eqref{lem:pHWM-zeros1} as well. This implies
$\tilde{F}^{n+1}(\ul{p}^n\alpha)=0$ and hence by the injectivity of $\tilde{F}^{n+1}\ul{p}^n$, also
that $\alpha=0$. This completes the proof.
\end{proof}

\begin{lem}\label{lem:VdV-zeros-claim}
Let $a\in \Omega^{q}_X$ be a local section satisfying $dV^{n-1}(a)\in W_n\Omega^{q+1}_{(X,-D)}$.
Then there exist local sections $a_j\in W_{j+1}\Omega^q_X$, $j=0,\ldots, n-1$,
and $b\in W_{n+1}\Omega^q_X$ such that $p^j a_j=\ul{p}^j R^j(a_j)\in W_{j+1}\Omega^q_{(X,-D)}$ and 
\[a= \sum_{j=0}^{n-1} F^j(a_j) + F^{n}(b).\]
\end{lem}
\begin{proof}
We show, for all $s=-1,\ldots, n-1$ we find $a_j\in W_{j+1}\Omega^q_X$, $j=0,\ldots, s$, as in the statement
such that 
\[a\equiv \sum_{j=0}^s F^j(a_j) \quad {\rm mod}\quad F^{s+1}W_{s+2}\Omega^q_X.\]
There is nothing to show for $s=-1$. Assume $0\le s\le n-1$ and the statement is true for $s-1$, i.e.,
we find $a_j$, $j=0,\ldots, s-1$, as in the claim and $b\in W_{s+1}\Omega^q_X$ such that 
\[a= \sum_{j=0}^{s-1} F^j(a_j)+ F^s(b) \quad \text{in } \Omega^q_X.\]
Applying $dV^{s}$ yields 
\[\ul{p}^s dR^s(b) = dV^s(a)-\sum_{j=0}^{s-1} dV^{s-j}(\ul{p}^j R^j(a_j))\in W_{s+1}\Omega^{q+1}_{(X,-D)},\]
where we use $dV^s(a)=F^{n-1-s}dV^{n-1}(a)\in W_{s+1}\Omega^{q+1}_{(X,-D)}$.
Hence by Lemma  \ref{lem:dzeros} and Lemma \ref{lem:pHWM-zeros} 
we find an element $a_s\in W_{s+1}\Omega^q_{X}$ with 
$\ul{p}^sR^s(a_s)\in W_{s+1}\Omega^q_{(X,-D)}$ such that $dV^s(F^s(a_s))=dV^s(F^s(b))$.
As the kernel of $dV^s$ is $F^{s+1}W_{s+2}\Omega^{q-1}_X$ the statement follows.
\end{proof}

\begin{lem}\label{lem:VdV-zeros}
The following two equalities of subsheaves of $W_{n+1}\Omega^q_X$ hold
\[W_{n+1}\Omega^{q}_{(X,-D)}\cap V^n(\Omega^q_X)=\sum_{j=0}^n V^j\left(W_{n+1-j}\Omega^q_{(X,-D)} \cap \ul{p}^{n-j}(\Omega^q_X)\right)\]
and 
\[W_{n+1}\Omega^{q}_{(X,-D)}\cap dV^n(\Omega^{q-1}_X)=\sum_{j=0}^n dV^j\left(W_{n+1-j}\Omega^{q-1}_{(X,-D)} \cap \ul{p}^{n-j}(\Omega^{q-1}_X)\right).\]
\end{lem}
\begin{proof}
We show the first equality. This ``$\supset$" inclusion is obvious. We show the other inclusion.
Let $a\in \Omega^{q}_X$ be a local section with $V^n(a)\in W_{n+1}\Omega^q_{(X,-D)}$.
As $dV^{n-1}(a)=FdV^n(a)\in W_{n}\Omega^{q+1}_{(X,-D)}$,  
we find elements $a_j\in W_{j+1}\Omega^q_X$ and an element $b\in W_{n+1}\Omega^q_X$ as in Lemma \ref{lem:VdV-zeros-claim}.
Applying $V^n$ to the equality in that lemma yields
\[\ul{p}^n R^n(b)=p^n b= \sum_{j=0}^{n-1}V^{n-j}(\ul{p}^jR^j(a_j))- V^n(a)\in W_{n+1}\Omega^q_{(X,-D)},\]
which proves the other inclusion in the first equality.

For the second equality we take $a\in \Omega^{q-1}_X$ such that $dV^n(a)\in W_{n+1}\Omega^q_{(X,-D)}$.
Thus Lemma \ref{lem:VdV-zeros-claim} yields the existence 
of $a_j\in W_{j+1}\Omega^{q-1}_X$, $j= 0,\ldots, n$, and $b\in W_{n+2}\Omega^{q-1}_X$ such that
$\ul{p}^j R^j(a_j)\in W_{j+1}\Omega^{q-1}_{(X,-D)}$ and 
\[a=\sum_{j=0}^n F^j(a_j)+F^{n+1}(b).\]
As $dV^n(F^{n+1}(b))=p^{n+1}Fd(b)=0$ we find 
\[dV^n(a)= \sum_{j=0}^n dV^{n-j}(\ul{p^j}R^j(a_j)),\]
which proves the second equality.
\end{proof}

\begin{prop}\label{prop-ex-seq-zeros}
There is a short exact sequence of $W_{n+1}\sO_X$-modules
\eq{prop-ex-seq-zeros0}{0\to (\Omega/B)^q_{n, (X,-D)}\xr{V^n} {\rm gr}^n_{(X,-D)}\xr{\beta} (\Omega/Z)^{q-1}_{n, (X,-D)}\to 0,}
where ${\rm gr}^n_{(X,-D)}:=\Ker(R: W_{n+1}\Omega^q_{(X,-D)}\to W_n\Omega^q_{(X,-D)})$  and the map 
$\beta$ is induced by $V^n(a)+dV^n(b)\mapsto b$.    
Furthermore the following diagram of sheaves of abelian groups is cartesian
\eq{prop-ex-seq-zeros00}{\xymatrix{
dV^n(\Omega^{q-1}_{X})\cap W_{n+1}\Omega^q_{(X,-D)}\ar[d]\ar[r]^-F &
W_n\Omega^q_{(X,-D)}\ar[d]\\
dV^n(\Omega^{q-1}_X)\ar[r]^-F&
W_n\Omega^q_X.
}}
\end{prop}
\begin{proof}
Consider the following diagram of $W_{n+1}\sO_X$-modules
\eq{prop-ex-seq-zeros1}{
\xymatrix{
(\Omega/B)^q_{n, (X,-D)}\ar[r]\ar[d] &
{\rm gr}^n_{(X,-D)}\ar[d]\ar@{.>}[r]\ar@/^1.5pc/[rr]^F &
(\Omega/Z)^{q-1}_{n,(X-D)}\ar[d]\ar[r]&
F_*W_n\Omega^q_{(X,-D)}\ar[d]\\
\frac{F^n_{X*}\Omega^q_X}{B_n\Omega^q_X}\ar[r]^{V^n}&
{\rm gr}^n_X\ar[r]^{\beta}\ar@/_1.5pc/[rr]_F&
\frac{F^n_{X*}\Omega^{q-1}_X}{Z_n\Omega^{q-1}_X}\ar[r]^-{dV^{n-1}}&
F_*W_n\Omega^q_X.
}}
The  two outer squares  are cartesian in the category of $W_{n+1}\sO_X$-modules by definition, and also in the category
of sheaves of abelian groups as the forgetful functor $(W_{n+1}\sO_X\text{-mod})\to (\text{abelian-sheaves})$ has a  left adjoint. 
Clearly the square with $F$ as horizontal maps is commutative
and we have $dV^{n-1}\circ \beta=F$. Hence by the definition of $(\Omega/Z)^{q-1}_{n, (X,-D)}$ as a pullback in \eqref{para:Omega/BZ3}
there is a unique dotted arrow which makes the middle square commutative.  
Thus all the maps in \eqref{prop-ex-seq-zeros0} are defined. The injectivity of $V^n$ holds by definition, the exactness in the middle of
\eqref{prop-ex-seq-zeros0} follows from the exactness of \eqref{para:Omega/BZ1}
and the fact that the left square is cartesian.
We claim that the diagram of sheaves of abelian groups 
\eq{prop-ex-seq-zeros05}{\xymatrix{
dV^n(\Omega^{q-1}_{X})\cap W_{n+1}\Omega^q_{(X,-D)}\ar[d]\ar[r]^-{\beta_0} &
(\Omega/Z)^{q-1}_{n,(X-D)}\ar[d]\\
dV^n(\Omega^{q-1}_X)\ar[r]^-{\beta_0}&
\frac{\Omega^{q-1}_X}{Z_n\Omega^{q-1}_X}
}}
is cartesian, where $\beta_0$ is induced by $\beta$, i.e., $\beta_0(dV^n(a))=a$.
If the claim is true, then the top $\beta_0$ is surjective as the bottom one is and hence also the 
dotted arrow in \eqref{prop-ex-seq-zeros1} is surjective. Moreover the square \eqref{prop-ex-seq-zeros00} is cartesian, 
as it is a composition of two cartesian squares.

We prove the claim. Let $a\in \Omega^{q-1}_X$ be a local section which mod $Z_n\Omega^{q-1}_X$ lies in the image
$\Im((\Omega/Z)^{q-1}_{n,(X,-D)}\hra \Omega^{q-1}_X/Z_n\Omega^{q-1}_X)$.
By \eqref{prop-ex-seq-zeros1} we have $dV^{n-1}(a)\in W_n\Omega^q_{(X,-D)}$.
As $Z_n\Omega^q_X=F^n(W_{n+1}\Omega^q_X)$, Lemma \ref{lem:VdV-zeros-claim} yields the 
existence of local sections $a_j\in W_{j+1}\Omega^q_X$, $j=0,\ldots, n-1$, such that $p^j a_j=\ul{p}^j R^j(a_j)\in W_{j+1}\Omega^q_{(X,-D)}$ and 
\[a=\sum_{j=0}^{n-1} F^j(a_j) \quad \text{in } \Omega^q_X/Z_n\Omega^q_X.\]
Moreover 
\[dV^n(\sum_{j=0}^{n-1} F^j(a_j))= \sum_{j=0}^{n-1}d V^{n-j}(\ul{p}^jR^j(a_j))\in W_{n+1}\Omega^q_{(X,-D)}\]
and it maps to $a$ under $\beta_0$. Hence  \eqref{prop-ex-seq-zeros05} is cartesian.
\end{proof}

\begin{lem}\label{lem:R-surj-zeros}
The restriction map $R:W_{n+1}\Omega^q_{(X,-D)}\to W_n\Omega^q_{(X,-D)}$ is surjective.    
\end{lem}
\begin{proof}
We show 
\eq{lem:R-surj-zeros1}{R^r:W_{n+r}\Omega^q_{(X,-D)} \ra W_n \Omega^q_{(X,-D)}\quad \text{is surjective, for all } r\ge 1.}
For $n=1$ this holds by the ``Moreover"-part of Lemma \ref{lem:pHWM-zeros} (take $(j,n)$ there as $(r,0)$ here). 
Assume $n\ge 2$. Let $w\in W_n\Omega^q_{(X,-D)}$ be a local section.
Set $w_{n-1}:=R(w)\in W_{n-1}\Omega^q_{(X,-D)}$.
By induction there exists an $w_{n+r}\in W_{n+r}\Omega^q_{(X,-D)}$  with $R^{r+1}(w_{n+r})=w_{n-1}$.
Set $w_n:=R^r(w_{n+r})\in W_{n}\Omega^q_{(X,-D)}$.
Since both $w$ and $w_n$ are lifts of $w_{n-1}$, we have $R(w-w_n)=0$.
Thus there exist elements $a\in\Omega^q_X$ and $b\in \Omega^{q-1}_X$ with 
\[V^{n-1}(a)+dV^{n-1}(b)=w-w_n\in W_n\Omega^q_{(X,-D)}.\]
As $F(V^{n-1}(a))=0$, the cartesian diagram \eqref{prop-ex-seq-zeros00} yields 
$dV^{n-1}(b)\in  W_n\Omega^q_{(X,-D)}$ and hence also $V^{n-1}(a)\in W_n\Omega^q_{(X,-D)}$.
By Lemma \ref{lem:VdV-zeros} we find elements $a_j\in \Omega^q_X$ and $b_j\in \Omega^{q-1}_X$ such that 
\[V^{n-1}(a)= \sum_{j=0}^{n-1}V^j(\ul{p}^{n-1-j}a_j)\quad \text{and}\quad \ul{p}^{n-1-j}a_j\in W_{n-j}\Omega^q_{(X,-D)}, \text{ for all }j,\]
and 
\[dV^{n-1}(b)=\sum_{j=0}^{n-1} dV^j(\ul{p}^{n-1-j}b_j)\quad \text{and}\quad \ul{p}^{n-1-j}b_j\in W_{n-j}\Omega^{q-1}_{(X,-D)}, \text{ for all }j.\]
By the ``Moreover"-part of  Lemma \ref{lem:pHWM-zeros} we find sections $\alpha_j\in W_{n-j+1}\Omega^q_{(X,-D)}$ and 
$\beta_j\in W_{n-j+1}\Omega^{q-1}_{(X,-D)}$ with $R(\alpha_j)=\ul{p}^{n-j-1}a_j$ and $R(\beta_j)=\ul{p}^{n-j-1} b_j$, 
for all $j=0,\ldots, n-1$. Altogether we find that the element
\[R^{r-1}(w_{n+r})+ \sum_{j=0}^{n-1}V^j(\alpha_j) +\sum_{j=0}^{n-1}V^j(\beta_j)\]
lies in $W_{n+1}\Omega^q_{(X,-D)}$ and lifts $w$. This completes the proof.
\end{proof}

The proof of part \ref{thm:HW-zeros1} and \ref{thm:HW-zeros1.5} of Theorem \ref{thm:HW-zeros} is now complete.
To prove the last part of that theorem we need to study certain variants of $(\Omega/Z)^q$ and $(\Omega/B)^q$
which are defined using iterates of the twisted inverse Cartier operator from \ref{para:Omega-zeros}, similar to what we did in 
section \ref{sec:strHWM}.

\begin{para}\label{para:twisted-BZ-zeros}
Let  $X$, $D$, and $E$ be as in \ref{para:Omega-zeros} above. Let $n\ge 0$ and let 
\[D=D_0+pD_1+\ldots+ p^n D_n\]
be a \pdd. 
In \ref{para:twisted-BZ} we have defined divisors
$\underline D_j$ and $\underline D^j$ 
such that 
$D=\ul{D}_{n-j}+ p^{n-j}\ul{D}^j$
is a \pdd, for any $0\le j\le n$.\footnote{Recall that $\ul{D}^j$ is divisible by $p$ and that the actual \pdd \, is
$D=\ul{D}_{n-j}+ p^{n-j+1}(\frac{1}{p}\ul{D}^j)$.}

For $1\le j\le n$, 
define $\cO_X$-modules
\[\Omega^q_{j,n,-}:=F^j_{X*}\Omega^q_{n-j}(-\ul{D}_{n-j}, -\ul{D}^j-p^jE),\quad \sB^q_{j,n,-}:=\sB^q(\Omega^\bullet_{j,n,-}), \quad 
\sZ^q_{j,n,-}:=\sZ^q(\Omega^\bullet_{j,n,-}),\]
we furthermore set $\sB^q_{0,n,-}:=0$.
Define $\sO_X$-modules
$(\Omega/B)^q_{j,n}(-D,-E)$, 
$(\Omega/Z)^q_{j,n}(-D,-E)$
and natural maps 
\eq{para:twisted-BZ-zeros0}{
\Omega^q_{j,n,-}/\sB^q_{j,n,-}
\lra (\Omega/B)^q_{j,n}(-D,-E) 
\lra (\Omega/Z)^q_{j,n}(-D,-E) 
}
by setting 
\[(\Omega/B)^q_{0,n}(-D,-E):= \Omega^q_n(-D, -E)\surj 0=:
(\Omega/Z)^q_{0,n}(-D,-E),\]
and for $j\ge 1$ recursively by the condition that the  two squares on the left in the following diagram are pushout squares 
\eq{para:twisted-BZ-zeros1}{
\xymatrix@C=1.2em{
0\ar[r] & 
{\scriptstyle F^{j-1}_{X*}\Omega^q_{n-j+1}(-\ul{D}_{n-j}, -D_{n-j+1}-\ul{D}^{j-1}-p^{j-1}E)}
\ar[r]^-{C^{-1}}\ar[d]&
\Omega^q_{j,n,-}/\sB^q_{j,n,-}
\ar[d]\ar[r]&
\Omega^q_{j,n,-}/\sZ^q_{j,n,-}
\ar@{=}[d]\ar[r] &0
\\
0\ar[r] & (\Omega/B)^q_{j-1,n}(-D,-E)
\ar[r]\ar[d] &
(\Omega/B)^q_{j,n}(-D,-E)
\ar[r]\ar[d]&
\Omega^q_{j,n,-}/\sZ^q_{j,n,-}
\ar@{=}[d]\ar[r] &0
\\
0\ar[r] & 
(\Omega/Z)^q_{j-1,n}(-D,-E)
\ar[r] &
(\Omega/Z)^q_{j,n}(-D,-E)
\ar[r]&
\Omega^q_{j,n,-}/\sZ^q_{j,n,-}
\ar[r] &0,
}}
where the top horizontal exact sequence is induced by \eqref{para:C-ex-seq-zeros1} and 
the top vertical map on the left is the composition of  the natural map
\[F^{j-1}_{X*}\Omega^q_{n-j+1}(-\ul{D}_{n-j}, -D_{n-j+1}-\ul{D}^{j-1}-p^{j-1}E) \hra 
\Omega^q_{j-1,n,-}\lra \frac{\Omega^q_{j-1,n,-}}{\sB^q_{j-1,n,-}}\]
with the first map of \eqref{para:twisted-BZ-zeros0} (with $j$ replaced by $j-1$).

Some consequences of this definition:
\begin{enumerate}[label= (\alph*)]
    \item
    \label{para:twisted-BZ-zerosa2}
    By definition we have for $n\ge 1$
     \[(\Omega/Z)^q_{1,n}(-D,-E)= \Omega^q_{1,n,-}/\sZ^q_{1,n,-}.\]  
    Moreover we have
    \[(\Omega/B)^q_{1,n}(-D,-E)=
    \frac{
    \sZ^q(F_{X*}\Omega^\bullet_{n-1}(-D,-pE))
    +\Omega^q_{1,n,-}}{\sB^q(F_{X*}\Omega^\bullet_{n-1}(-D,-pE))}.\]
     For the latter equality note that the right hand side (RHS) fits into a short exact sequence
    $$0\to \Omega^q_n(-D,-E)
    \xr{\alpha} 
    \text{RHS}
    \xr{\beta}
    \Omega^q_{1,n,-}/\sZ^q_{1,n,-}
    \to 0,$$
    where $\alpha$ is given by $C^{-1}$ from \eqref{para:Omega-zeros6}
     composed with the natural inclusion  and the map $\beta$ is defined by
     $\beta(a+b)=b$, for $a\in\sZ^q(F_{X*}\Omega^q_{n-1}(-D,-pE))$, and $b\in \Omega^q_{1,n,-}$.
     It is direct to check that $\beta$ is well-defined, that the sequence is exact, and that this implies the above equality.
    \item
    \label{para:twisted-BZ-zerosb} 
    The sheaves $(\Omega/B)^q_{j,n}(-D,-E)$, $(\Omega/Z)^q_{j,n}(-D,-E)$ are locally free coherent $\sO_X$-modules. 
    For $j=0$ this follows from the definition 
    and for $j\ge 1$ by induction from \eqref{para:twisted-BZ-zeros1} and 
    the fact that the $\sO_X$-modules $\Omega^q_{j,n,-}/\sZ^q_{j,n,-}$ are locally free, see \ref{para:C-exs-eq-zeros}.
    Their restrictions to $U=X\setminus (D+E)$ are quotients of $\Omega^q_U$.
    The natural maps 
    \eq{para:twisted-BZ-zerosb1}{(\Omega/B)^q_{j,n}(-D,-E)\surj (\Omega/Z)^q_{j,n}(-D,-E)}
    are surjective.   Moreover the $j$-fold iteration of the  Cartier operator induces an isomorphism
    \[ (Z/B)^q_{j,n}(-D,-E):=\Ker \left((\Omega/B)^q_{j,n}(-D,-E)\to (\Omega/Z)^q_{j,n}(-D,-E)\right)\xra[\simeq]{C^{j}}
    \Omega^q_n(-D,-E).\]
    For $j=0$ this holds by definition, and for $j\ge 1$ it follows by induction from applying the Snake Lemma to the 
    two lower rows in \eqref{para:twisted-BZ-zeros1}.
    \item
    \label{para:twisted-BZ-zerosc} 
    We have natural maps  of  $\sO_X$-modules 
   \eq{para:twisted-BZ-zerosc1}{F_{X*}(\Omega/B)^q_{j-1,n-1}(-\ul{D}_{n-1}, -\ul{D}^1-pE)\to (\Omega/B)^q_{j,n}(-D,-E)}
   and surjections of $\sO_X$-modules
    \eq{para:twisted-BZ-zerosc2}{(\Omega/Z)^q_{j,n}(-D,-E)\surj F_{X*}(\Omega/Z)^q_{j-1,n-1}(-\ul{D}_{n-1}, -\ul{D}^1-pE).}
In the case $j=1$ this follows from the definition
and for $j\ge 2$ it follows by induction from the fact 
the top sequence in  \eqref{para:twisted-BZ-zeros1} constructed for $D$, $E$, $n$, $j$ is  equal to $F_{X*}$ 
applied to the sequence constructed for $\ul{D}_{n-1}$, $\ul{D}^1+pE$, $n-1$, $j-1$.
By definition the $(j-1)$-fold iteration of the Cartier operator induces an isomorphism
\[ C^{j-1}:\Ker\eqref{para:twisted-BZ-zerosc2}    \xra{\simeq}   (\Omega/Z)^q_{1,n}(-D,-E).\]
\end{enumerate}
\end{para}

\begin{rmk}
Note that the first map in \eqref{para:twisted-BZ-zeros0} is in general not surjective, which under the duality from 
Proposition \ref{prop:BZ-duality} relates to the fact  that 
\[Z^q_{j,n}(D,E)\subset F^j_{X*}\Omega^q_{n-j}(\ul{D}_{n-j}, \ul{D}^j+ p^j E)\]
is not a subbundle, see Remark \ref{rmk:BZ-loc-free}.
\end{rmk}

\begin{lem}\label{lem:BZ-reg}
In the situation of \ref{para:twisted-BZ-zeros} we have natural inclusions
\[(\Omega/B)^q_{j,n}(-D,-E)\subset F^j_{X*}\Omega^q_X/B_j\Omega^q_X, \qquad 
(\Omega/Z)^q_{j,n}(-D,-E)\subset F^j_{X*}\Omega^q_X/Z_j\Omega^q_X.\]
\end{lem}
\begin{proof}
The statements hold by definition for $j=0$. The general statement follows by induction over $j$. We explain the case $(\Omega/Z)$
and a similar argument also works for $(\Omega/B)$. 
We have a morphism of exact sequences (we drop the Frobenius twists)
\[\xymatrix{
0\ar[r] & 
(\Omega/Z)^q_{j-1,n}(-D,-E)
\ar[r] \ar[d]&
(\Omega/Z)^q_{j,n}(-D,-E)
\ar[r]\ar@{.>}[d]&
\Omega^q_{j,n,-}/\sZ^q_{j,n,-}
\ar[r]\ar[d] &0\\
0\ar[r] & 
\Omega^q_X/Z_{j-1}\Omega^q_X
\ar[r]^{C^{-1}} &
\Omega^q_X/Z_j\Omega^q_X
\ar[r]&
\Omega^q_X/Z_1\Omega^q_X
\ar[r] &0,
}\]
where left vertical map is the inclusion which we have by induction, the right vertical map is the natural inclusion, and 
the middle map exists by the definition of $(\Omega/Z)^q_{j,n}(-D,-E)$ as a pushout and the fact that there is a natural map
$\Omega^q_{j,n,-}/\sB^q_{j,n,-}\to \Omega^q_X/Z_j\Omega^q_X$. This yields the statement.
\end{proof}

\begin{prop}\label{prop:BZ-duality}
Let $X$ be  of pure dimension $N$. The isomorphism (cf. Proposition \ref{prop:n=1duality})
\[F^j_{X*}\Omega^q_X\to \sH om(F^j_{X*}\Omega^{N-q}_X, \Omega^N_X), \quad 
\alpha\mapsto \left(\beta\mapsto C^j(\alpha\wedge \beta)\right),\]
induces isomorphisms
\eq{prop:BZ-duality1}{(\Omega/B)^q_{j,n}(-D,-E)\cong \sH om(Z^{N-q}_{j,n}(D,E), \Omega^N_X)}
and
\eq{prop:BZ-duality2}{(\Omega/Z)^q_{j,n}(-D,-E)\cong \sH om(B^{N-q}_{j,n}(D,E), \Omega^N_X),}
where $Z^{N-q}_{j,n}(D,E)$ and $B^{N-q}_{j,n}(D,E)$ are defined in \ref{para:twisted-BZ}.
\end{prop}
First we record the following lemma whose proof is similar to the one of Proposition \ref{prop:n=1duality}.
\begin{lem}\label{lem:n=1duality'}
The map
\[F^j_*\Omega^q_n(-D,-E)\xr{\simeq} \sH om(F^j_*\Omega^{N-q}_n(D,E), \Omega^N_X), 
\quad \alpha\mapsto (\beta\mapsto C^j(\alpha\wedge\beta)),\]
is an isomorphism for all $j$, $n\ge 0$.
\end{lem}

\begin{proof}[Proof of Proposition \ref{prop:BZ-duality}]
This is a generalization of \cite[II, Lemma 2.2.20]{Ekedahl}.
For $j=0$ both sides of \eqref{prop:BZ-duality2} are zero and   \eqref{prop:BZ-duality1} holds by Lemma  \ref{lem:n=1duality'}.
Set 
\[\sD:= \sH om(-,\Omega^N_X): (\text{loc. free }\sO_X\text{-mod})^{\rm op}\to (\text{loc. free }\sO_X\text{-mod}).\] 
Note that $\sD$ is dualizing, in the sense that the natural map $\sE\to \sD(\sD(\sE))$ is an isomorphism for any locally free $\sO_X$-module $\sE$.
Hence $\sD$ maps a pushout square of locally free $\sO_X$-modules to a pullback square and vice versa.
Thus it remains to show that $\sD$ applied to the bottom exact sequence in \eqref{para:twisted-BZ3} 
(with $N-q$ instead of $q$) is isomorphic to the  top sequence in \eqref{para:twisted-BZ-zeros1}, for all $j\ge 1$.
Then the result follows by induction over $j$ directly from the definitions.
The looked for duality between the two exact sequences is standard, e.g. \cite[Lemma 1.7]{Milne76}. For convenience of the reader we give the proof
in this twisted situation.
In the following we assume $j\ge 1$. 
By Lemma \ref{lem:n=1duality'} the map 
$$
b^q_j: \Omega^q_{j,n,-}  \lra \sD(\Omega^{N-q}_{j,n}),\qquad 
b^q_j(\alpha)(\beta)=C^j(\alpha\wedge \beta) \quad (\alpha\in \Omega^q_{j,n,-},\, 
\beta\in \Omega^{N-q}_{j,n}).
$$
is an isomorphism. If $\alpha\in \sB^q_{j,n,-}$ and $\beta\in \sZ^{N-q}_{j,n}$,
then $\alpha=d\alpha'$ and  $d\beta=0$ and hence  
\[b^q_j(\alpha)(\beta)= C^j(d\alpha'\wedge \beta)=C^j(d(\alpha'\wedge \beta))=0.\]
Similarly, the above vanishing holds if $\alpha\in \sZ^q_{j,n,-}$  and $\beta\in \sB^{N-q}_{j,n}$.
Therefore the map $b^q_j$ induces  well-defined maps 
$$c^q_j: 
\Omega^q_{j,n,-}/\sB_{j,n,-}^q
\lra 
\sD(\sZ^{N-q}_{j,n})
\quad 
\text{and}
\quad
a^{q}_j:
\Omega^q_{j,n,-}/\sZ_{j,n,-}^q
\lra 
\sD(\sB_{j,n}^{N-q}).
$$
Consider the following two diagrams,
\eq{lem:BZ-duality3}{
\xymatrix{
0\ar[r]&
\Omega^{q-1}_{j,n,-}/\sZ^{q-1}_{j,n,-}
\ar[r]^-{d}\ar[d]_{a^{q-1}_j}& 
\Omega^q_{j,n,-}
\ar[r]\ar[d]_{b^{q}_j}^{\simeq}&
\Omega^q_{j,n,-}/\sB_{j,n,-}^q
\ar[r]\ar[d]_{c^{q}_j}&0
\\
0\ar[r]&
\sD(\sB_{j,n}^{N-q+1})
\ar[r]^-{(-1)^q d^\vee}&
\sD(\Omega^{N-q}_{j,n})
\ar[r]&
\sD(\sZ^{N-q}_{j,n})
\ar[r]&
0
}}
and
\eq{lem:BZ-loc-free-zeros4}{
\xymatrix@C=1em{
0\ar[r]&
{\scriptstyle 
F^{j-1}_{X*}\Omega^q_{n-j+1}(-\ul{D}_{n-j}, -D_{n-j+1}-\ul{D}^{j-1}-p^{j-1}E)}
\ar[r]^-{C^{-1}}\ar[d]_{e^q}^{\simeq}&
\Omega^q_{j,n,-}/\sB_{j,n,-}^q
\ar[r]\ar[d]_{c^q_j}&
\Omega^q_{j,n,-}/\sZ^q_{j,n,-}
\ar[r]\ar[d]_{a^q_j}&0
\\
0\ar[r]&
{\sD\left(\scriptstyle F^{j-1}_{X*}\Omega^{N-q}_{n-j+1}(\ul{D}_{n-j}, D_{n-j+1}+\ul{D}^{j-1}+p^{j-1}E)\right)}
\ar[r]^-{C^\vee}&
\sD(\sZ^{N-q}_{j,n})
\ar[r]&
\sD(\sB_{j,n}^{N-q})
\ar[r]&
0.
}}
where $e^q$ is the isomorphism from Lemma  \ref{lem:n=1duality'}, which is induced by $\alpha\mapsto C^{j-1}(\alpha\wedge \beta)$.
The two diagrams clearly commute and the rows are exact. In fact the top sequence of \eqref{lem:BZ-loc-free-zeros4}
is the top sequence of \eqref{para:twisted-BZ-zeros1} and the bottom sequence of \eqref{lem:BZ-loc-free-zeros4}
is $\sD$ applied to the bottom sequence of \eqref{para:twisted-BZ3}.  Thus by the above it is only left to show that 
$a^q_j$ and $c^q_j$ are isomorphisms for all $q\ge 0$ and $j\ge 1$. But this follows by descending induction over $q$ from 
\eqref{lem:BZ-duality3} together with \eqref{lem:BZ-loc-free-zeros4}. This completes the proof of the proposition.
\end{proof}

\begin{prop}\label{prop:BZnn-zeros}
Let $D=D'+p^{n+1}D_{n+1}$ be a \pdd.
Consider the $\sO_X$-modules $(\Omega/Z)^{q-1}_{n,(X,-D)}$ and $(\Omega/B)^q_{n,(X,-D)}$ defined in \eqref{para:Omega/BZ3}.
There are natural  isomorphisms of $\sO_X$-modules for $n\ge 0$
\eq{prop:BZnn-zeros1}{(\Omega/Z)^{q-1}_{n,n}(-D', -pD_{n+1})\xr{\simeq}(\Omega/Z)^{q-1}_{n,(X,-D)}.}
and
\eq{prop:BZnn-zeros2}{
(\Omega/B)^q_{n,n}(-D', -pD_{n+1})\xr{\simeq} (\Omega/B)^q_{n,(X,-D)} .}
In particular the sheaves on the right hand side are locally free coherent $\sO_X$-modules.
\end{prop}

\begin{proof}
The last statement follows from the two isomorphisms and \ref{para:twisted-BZ-zeros}\ref{para:twisted-BZ-zerosb}.
We first prove \eqref{prop:BZnn-zeros1}.
In case $n=0$ both sides are zero by definition. 
The case $n=1$ follows from \ref{para:Omega/BZ}\ref{para:Omega/BZa} and \ref{para:twisted-BZ-zeros}\ref{para:twisted-BZ-zerosa2}.

Let $n\ge 2$.
Let $D=D_0+pD_1+\ldots+ p^{n+1}D_{n+1}$ be a \pdd.
Assume that by induction we have a natural isomorphism
$$(\Omega/Z)^{q-1}_{n-1,n-1}(-\ul{D}_{n-1}, -pD_n-p^2D_{n+1})\xr{a_{n-1}} (\Omega/Z)^{q-1}_{n-1, (X,-D)}.$$
\begin{claim}\label{prop:BZnn-zeros-claim1}
There exits a natural morphism $a_n$ which makes the following diagram commutative
\eq{prop:BZnn-zeros-claim10}{\xymatrix{
(\Omega/Z)^{q-1}_{n,n}(-\ul{D}_n,-pD_{n+1})\ar[d]^{a_n}\ar@{->>}[r]^-{\eqref{para:twisted-BZ-zerosc2}}&
F_{*}(\Omega/Z)^{q-1}_{n-1,n-1}(-\ul{D}_{n-1}, -pD_n-p^2D_{n+1})\ar[d]^{a_{n-1}}\\
(\Omega/Z)^{q-1}_{n, (X,-D)}\ar[r]^-{\eqref{prop:BZnn-zeros3}} &
F_*(\Omega/Z)^{q-1}_{n-1, (X,-D)}.
}}
\end{claim}
We prove the claim. The following diagram of sheaves of abelian groups clearly commutes
\[\xymatrix{
(\Omega/Z)^{q-1}_{n,n}(-\ul{D}_n,-pD_{n+1})\ar[rr]^-{a_{n-1}\circ \eqref{para:twisted-BZ-zerosc2}}\ar[d]& &
(\Omega/Z)^{q-1}_{n-1, (X,-D)}\ar[r]^{dV^{n-2}}&
W_{n-1}\Omega^q_{(X,-D)}\ar[d]
\\
\Omega^{q-1}_X/Z_n\Omega^{q-1}_X\ar[rr]^{dV^{n-1}}& &
dV^{n-1}(\Omega^q_X)\ar[r]^F&
W_{n-1}\Omega^q_X,
}\]
where the left vertical map is the natural inclusion from Lemma \ref{lem:BZ-reg}.
As \eqref{prop-ex-seq-zeros00} is cartesian we obtain a morphism
$(\Omega/Z)^{q-1}_{n,n}(-\ul{D}_n,-pD_{n+1})\to dV^{n-1}(\Omega^{q-1}_X)\cap W_n\Omega^q_{(X,-D)}$,
composing with the inclusion $dV^{n-1}(\Omega^{q-1}_X)\cap W_n\Omega^q_{(X,-D)}\subset W_n\Omega^q_{(X,-D)}$
yields the top horizontal morphism in the following commutative diagram of  $W_{n+1}\sO_X$-modules 
\[\xymatrix{
(\Omega/Z)^{q-1}_{n,n}(-\ul{D}_n,-pD_{n+1})\ar[r]^-{dV^{n-1}}\ar[d]&
F_*W_n\Omega^q_{(X,-D)}\ar[d]\\
F^n_{X*}\Omega^{q-1}_X/Z_n\Omega^{q-1}_X\ar[r]^-{dV^{n-1}}&
F_*W_n\Omega^q_X.
}\]
Now the  map $a_n$ in the claim exists by definition of $(\Omega/Z)^{q-1}_{n,(X,-D)}$ in \eqref{para:Omega/BZ3}.
The diagram \eqref{prop:BZnn-zeros-claim10} commutes by construction, or maybe easier, the fact that both lines of the diagram
map injectively into $F^n_{X*}\Omega^{q-1}_X/Z_n\to F^n_{X*}\Omega^{q-1}_X/ F_{X*}Z_{n-1}$.

It remains to show that $a_n$ is an isomorphism. As $a_{n-1}$ is an isomorphism by induction and the map \eqref{para:twisted-BZ-zerosc2}
is surjective it remains to show that $a_n$ induces an isomorphism between the kernels of the horizontal maps in \eqref{prop:BZnn-zeros-claim10}.
By construction 
\[\Ker \eqref{prop:BZnn-zeros3}= \{a\in Z_{n-1}\Omega^{q-1}_X/Z_n\Omega^{q-1}_X\mid dV^{n-1}(a)\in W_n\Omega^q_{(X,-D)}\}.\]
We note
\begin{itemize}
\item $Z_{n-1}\Omega^{q-1}_X=F^{n-1}W_n\Omega^{q-1}_X$;
\item the $(n-1)$-fold iterated Cartier operator $C^{(n-1)}$ maps $F^{n-1}(b)$ to $R^{n-1}(b)$, for $b\in W_n\Omega^{q-1}_X$;
\item $dV^{n-1}(F^{n-1}(b))=\ul{p}^{n-1}d R^{n-1}(b)$.
\end{itemize}
Hence $C^{(n-1)}$ induces an isomorphism
\[C^{(n-1)}: \Ker \eqref{prop:BZnn-zeros3}\xr{\simeq} \{b\in \Omega^{q-1}_X/Z_1\Omega^{q-1}_X\mid  \ul{p}^{n-1}d b\in W_n\Omega^q_{(X,-D)}\}.\]
By Lemma \ref{lem:pHWM-zeros} and Lemma \ref{lem:dzeros} we can rewrite this as
\[C^{(n-1)}: \Ker \eqref{prop:BZnn-zeros3}\xr{\simeq} 
\frac{\Omega^q_{n-1}(-\ul{D}_{n-1}, -pD_n-p^2D_{n+1})}{\sZ^{q-1}(\Omega^\bullet_{n-1}(-\ul{D}_{n-1}, -pD_n-p^2D_{n+1}))}.\]
This together with \ref{para:twisted-BZ-zeros} \ref{para:twisted-BZ-zerosa2} and \ref{para:twisted-BZ-zerosc} implies that $a_n$ induces an isomorphism
on the kernels of the horizontal maps in \eqref{prop:BZnn-zeros-claim10}. Thus we have proven the isomorphism \eqref{prop:BZnn-zeros1} in general.

\medskip

Next \eqref{prop:BZnn-zeros2}.
First note that for  $a\in W_j\Omega^{q-1}_X$,  we have $V^j(F^{j-1}da)=p^jdV(a)=0$.
As $B_j\Omega^q_X=F^{j-1}d(W_j\Omega^{q-1}_X)$ we obtain a well-defined morphism of $W_{n+1}\sO_X$-modules
\[\ul{p}^{n-j}V^j=V^j\ul{p}^{n-j}: F^j_{X*}\Omega^q_X/B_j\Omega^q_X\to W_{n+1}\Omega^q_X.\]
\begin{claim}\label{prop:BZnn-zeros-claim2}
For  $j=0,\ldots, n$ the above map induces well-defined morphisms of $W_{n+1}\sO_X$-modules
\[\ul{p}^{n-j}V^j: (\Omega/B)^q_{j,n}(-D', -pD_{n+1})\to W_{n+1}\Omega^q_{(X,-D)}.\]
\end{claim}
We prove the claim. For $j=0$, this follows from Lemma \ref{lem:pHWM-zeros} and the definition of $(\Omega/B)^q_{0,n}$.
Let $j\ge 1$ and consider the following diagram
\[\xymatrix{
F^{j-1}_{X*}\Omega^{q}_{n-j+1}(- \ul{(D')}_{n-j}, -D_{n-j+1}- \ul{(D')}^{j-1}-p^j D_{n+1})\ar[r]^-{C^{-1}}\ar[d]&
\Omega^q_{j,n,-}/\sB^{q}_{j,n,-}\ar[d]^{V^j\ul{p}^{n-j}} \\
(\Omega/B)^q_{j-1,n}(-D',-pD_{n+1})\ar[r]^-{\ul{p}^{n-j+1}V^{j-1}} &
W_{n+1}\Omega^q_{(X,-D)},
}\]
where
\[\Omega^q_{j,n,-}= F^j_{X*}\Omega^q_{n-j}(-\ul{(D')}_{n-j}, -\ul{(D')}^j-p^{j+1}D_{n+1})\quad \text{and}\quad \sB^q_{j,n,-}=d\Omega^{q-1}_{j,n,-}.\]
Here the bottom horizontal map exists by induction and the existence of the right vertical map follows from Lemma \ref{lem:pHWM-zeros}, which implies
$\ul{p}^{n-j}(\Omega^q_{j,n,-})\subset W_{j+1}\Omega^q_{(X,-D)}$. It is direct to check that the above diagram commutes.
Thus by the definition of $(\Omega/B)^q_{j,n}(-D',-pD_{n+1})$ as a pushout we obtain the map from the statement of the claim.

Thus we obtain a commutative diagram 
\[\xymatrix{
(\Omega/B)^q_{n,n}(-D',-pD_{n+1})\ar[r]^-{V^n}\ar[d] &
 W_{n+1}\Omega^q_{(X,-D)}\ar[d]\\
 F^n_{X*}\Omega^q_X/B_n\Omega^q_X\ar[r]^-{V^n}&
 W_{n+1}\Omega^q_X,
}\]
where the top map is induced from Claim \ref{prop:BZnn-zeros-claim2} with $j=n$. 
The definition of $(\Omega/B)^q_{n, (X,-D)}$ as a pullback in \eqref{para:Omega/BZ3} yields therefore a natural map
\[b_n: (\Omega/B)^q_{n,n}(-D',-pD_{n+1})\to (\Omega/B)^q_{n, (X,-D)}. \]
Altogether we obtain the following diagram of $\sO_X$-modules
\eq{prop:BZnn-zeros4}{\xymatrix{
(\Omega/B)^q_{n,n}(-D',-pD_{n+1})\ar[r]^{\eqref{para:twisted-BZ-zerosb1}}\ar[d]^{b_n}&
(\Omega/Z)^q_{n,n}(-D', -pD_{n+1})\ar[d]^{a_n}\\
(\Omega/B)^q_{n, (X,-D)}\ar[r]^{\eqref{para:Omega/BZ-hn}}&
(\Omega/Z)^q_{n, (X,-D)}.
}}
The diagram commutes as both lines of the diagram
map injectively into $F^n_{X*}\Omega^{q-1}_X/B_n\to F^n_{X*}\Omega^q_X/Z_n$.
As \eqref{para:twisted-BZ-zerosb1} is surjective and $a_n$ is an isomorphism as shown above, it remains
to show that $b_n$ induces an isomorphism between the kernels of the horizontal maps.
By construction we have 
\[\Ker(\eqref{para:Omega/BZ-hn})=\{a\in Z_n\Omega^q_X/B_n\Omega^q_X\mid V^n(a)\in W_{n+1}\Omega^q_{(X,-D)}\}.\]
We note
\begin{itemize}
\item $Z_n\Omega^q_X/B_n\Omega^q_X=C^{-n}(\Omega^q_X)$;
\item $V^n(C^{-n}(b))=\ul{p}^n b.$
\end{itemize}
Thus the $n$-fold iterated Cartier operator $C^n$ induces an isomorphism
\[C^n: \Ker(\eqref{para:Omega/BZ-hn})\xr{\simeq} \{b\in \Omega^q_X\mid \ul{p}^n b\in W_{n+1}\Omega^q_{(X,-D)}\}.\]
By Lemma \ref{lem:pHWM-zeros} the right hand side is isomorphic to $\Omega_n^q(-D', -pD_{n+1})$ and thus it follows from 
\ref{para:twisted-BZ-zeros}\ref{para:twisted-BZ-zerosb}, that $b_n$ induces an ismorphism between the kernels of the horizontal maps
of the diagram \eqref{prop:BZnn-zeros4}. This completes the proof.
\end{proof}

The proof of Theorem \ref{thm:HW-zeros} is now complete.

\begin{remark}
\label{rmk:pHWM-VdV-zeros}
Using \Cref{prop:BZnn-zeros}, we can give an explicit description for $(\Omega/B)^q_{n,(X,-D)}$: 
$$(\Omega/B)^q_{n,(X,-D)}
=\frac{\sum_{j=0}^n C^{-j}(F^{n-j}_{X*}\Omega^q_{j}(-\ul D_j, -\ul D^{n+1-j}))}{\sB^q_{n}(\Omega^\bullet_n(-D',-pD_{n+1}))}.
$$
Here
$\sB^q_{n}(\Omega^\bullet_n(-D',-pD_{n+1}))$ is defined inductively so that
$\sB^q_{0}(\Omega^\bullet_n(-D',-pD_{n+1}))=0,$ and
$$
C^{-1}:\sB^q_{j}(\Omega^\bullet_n(-D',-pD_{n+1}))\stackrel{\simeq}{\lra}
\frac{\sB^q_{j+1}(\Omega^\bullet_n(-D',-pD_{n+1}))}{\sB^q_{1}(\Omega^\bullet_n(-D',-pD_{n+1}))},$$
for every $j\in [0,n-1]$.
\end{remark}

\section{Duality for Hodge-Witt cohomology with modulus: finite level}\label{sec:duality}
In this section we prove the duality between Hodge-Witt sheaves with modulus and their counterpart  with zeros
and draw some first consequences. This is one of the main results of the paper, see Theorem \ref{thm:duality-mod}.

Throughout this section we assume $X\in \Sm$ is of pure dimension $\dim X=N$
and we let $D$ be an effective Cartier divisor on $X$ such that $D_{\red}$ is an SNCD.

\begin{lem}\label{lem:mult}
Set $U=X\setminus D$.
The multiplication map
$$W_n\Omega^q_U\times W_n\Omega^{N-q}_U\ra W_n\Omega^N_U$$
induces a well-defined map
\eq{lem:mult0}{W_n\Omega^q_{(X,-D)}\times W_n\Omega^{N-q}_{(X,D)}\ra W_n\Omega^N_X}
of Nisnevich sheaves on $X$.
\end{lem}
\begin{proof}
Let $\eta\in D^{(0)}$ be a geometric point. Set  $\sO_L:=\sO_{X, \eta}^h$, it is a henselian discrete valuation ring 
with maximal ideal $\fm_L$ and fraction field $L$. 
Let $r\ge 1$ be the multiplicity of $D$ at $\eta$ and set
\[W_n\Omega^q_{(\sO_L, \fm_L^r)}:=(W_n\Omega^q_{(X,-D),\eta})^h=\Ker (W_n\Omega^{q}_{\sO_L}\to W_n\Omega^{q}_{\sO_L/\fm_L^r}).\]
In view of the definition of $W_n\Omega^{N-q}_{(X,D)}$ in \ref{para:HW-not} and as 
$W_n\Omega^N_X$ is a successive extension of locally free $\sO_X$-modules (see \cite[I, Corollaire 3.9]{IlDRW}), it suffices to show that multiplication induces a well-defined map
\[W_n\Omega^q_{(\sO_L, \fm_L^r)}\times \Fil^p_rW_n\Omega^{N-q}_L\to W_n\Omega^N_{\sO_L},\]
with $\Fil^p_rW_n\Omega^{N-q}_L$ as in Definition \ref{defn:p-sat}.
By the equality $\ul{p}^s(\Fil_rW_{n-s}\Omega^{N-q}_L)=p^s\Fil_{rp^s}W_{n}\Omega^{N-q}_L$ 
(see Lemma \ref{lem:fil-FVR}) and the formulas 
\[\alpha \cdot\ul{p}^s(\beta)= \ul{p^s}(R^s(\alpha)\cdot \beta)\quad \text{and}\quad  
\alpha\cdot d(\gamma)= (-1)^q (d(\alpha\cdot \gamma)+ - d(\alpha)\cdot \gamma), \]
for $\alpha\in W_n\Omega^q_{(\sO_L,\fm_L^r)}$, $\beta\in \Fil_{r} W_{n-s}\Omega^{N-q}_L$, 
and $\gamma \in \fil_rW_n\Omega^{N-q-1}_L$,
it suffices to show that multiplication induces a well-defined map
\[W_n\Omega^q_{(\sO_L,\fm_L^r)}\times \fil_rW_n\Omega^{N-q}_L \to W_n\Omega^N_{\sO_L}.\]
By \cite[Lemma 1.2.2]{He04} $W_n\Omega^q_{(\sO_L,\fm_L^r)}$ is the degree $q$-part of the differential graded ideal in $W_n\Omega^*_{\sO_L}$ generated by $W_n(\fm_L^r)$. Thus any element  $\alpha\in W_n\Omega^q_{(\sO_L,\fm_L^r)}$ is a sum of elements 
\[ V^i([a])\alpha' \quad \text{and}\quad  dV^i([a])\alpha'', \quad 
\text{with } v_L(a)\ge r,\, \alpha'\in W_n\Omega^q_{\sO_L},\, \alpha''\in W_n\Omega^{q-1}_{\sO_L}, 
\]
where $v_L: L\to \Z\cup\{\infty\}$ denotes the discrete valuation on $L$.
Let
\[m=\min\{v_p(r), n\}.\]
By Definition \ref{defn:fil} any element in $\fil_rW_n\Omega^{N-q}_L$ is a sum of elements
\begin{enumerate}[label=(\arabic*)]
\item\label{lem:mult1} $\beta=V^j([b])\dlog u$, with $p^{n-j-1}v_L(b)\ge -r+1$, $u\in K^M_{N-q}(L)$, $j\in \{0,\ldots, n-1\}$;
\item\label{lem:mult2} $\gamma=V^j([c])\dlog v$, with $p^{n-j-1}v_L(c)\ge -r$, $v\in K^M_{N-q}(\sO_L)$, $j\in \{n-m,\ldots, n-1\}$.
\end{enumerate}
Thus it suffices to show that for $a\in \sO_L$ with $v_L(a)\ge r$ and $\beta$ and $\gamma$ as in \ref{lem:mult1} and \ref{lem:mult2}, respectively,
we have 
\[V^i([a])\beta,\quad  V^i([a])\gamma \in W_n\Omega^{N-q}_{\sO_L},  \quad 
dV^i([a])\beta, \quad dV^i([a])\gamma\in W_n\Omega^{N-q+1}_{\sO_L}.\]
We consider the various cases separately.

\begin{description}
\item[1st case] $V^i([a])\beta$.
We have 
\[v_L(a^{p^j}b^{p^i})\ge p^j r+ p^i \tfrac{-r+1}{p^{n-j-1}}\ge p^j\ge 1\Longrightarrow 
V^i([a])\beta= V^{i+j}([a^{p^j}b^{p^i}]\dlog u)\in W_n\Omega^N_{\sO_L}.\]

\item[2nd case] $V^i([a])\gamma$.
We have 
\[v_L(a^{p^j}c^{p^i})\ge p^jr+ p^i\tfrac{-r}{p^{n-j-1}}\ge 0\Longrightarrow 
V^i([a])\gamma=V^{i+j}([a^{p^j}c^{p^i}]\dlog v)\in W_n\Omega^N_{\sO_L}.\]

\item[3rd case] $dV^i([a])\beta$.

\begin{description}
    \item[1st subcase] $n-1\ge j\ge i\ge 0$.
We have 
\[v_L(ba^{p^{j-i}})\ge \tfrac{-r+1}{p^{n-j-1}}+p^{j-i}r>0\Longrightarrow 
dV^i([a])\beta= V^j([b][a]^{p^{j-i}}\dlog \{a, u\})\in W_n\Omega^N_{\sO_L}.\]
\item[2nd subcase] $n-1\ge i\ge j\ge 0$.
We have 
\mlnl{v_L(ab^{p^{i-j}})\ge r+ p^{i-j}\tfrac{-r+1}{p^{n-j-1}}\ge 1\Longrightarrow \\
dV^i([a])\beta= V^j\left(dV^{i-j}([ab^{p^{i-j}}]\dlog u)- V^{i-j}([a b^{p^{i-j}}]\dlog\{b,u\})\right)\in W_n\Omega^N_{\sO_L}.}
\end{description}
\item[4th case] $dV^i([a])\gamma$.
\begin{description}
\item[1st subcase] $n-1\ge j\ge n-m$ and $\min\{n-2, j\}\ge i$.
We have
\[v_L(ca^{p^{j-i}})\ge \tfrac{-r}{p^{n-j-1}} + p^{j-i}r>0\Longrightarrow dV^i([a])\gamma = V^j([ca^{p^{j-i}}]\dlog\{a,v\})\in W_n\Omega^N_{\sO_L}.\]
\item[2nd subcase] $n-1=j=i\ge n-m$.
We have 
\[dV^{n-1}([a])\gamma=V^{n-1}([ca]\dlog a\dlog v)\]
and $v_{L}(ca)\ge 0$. If $v_L(ca)\ge 1$, then $dV^{n-1}([a])\gamma$ is clearly regular.
If $v_L(ca)=0$, we have $v_L(a)=r$ and $v_p(r)\ge m\ge 1$. Thus in this case $\dlog a\in \Omega^1_{\sO_L}$ and 
hence also $dV^{n-1}([a])\gamma\in W_n\Omega^N_{\sO_L}$.
\item[3rd subcase] $n-1>i\ge j\ge n-m$.
We have
\mlnl{v_L(a c^{p^{i-j}})\ge r-\tfrac{rp^{i-j}}{p^{n-j-1}}>0 \Longrightarrow \\
dV^i([a])\gamma= V^j\left(dV^{i-j}([ac^{p^{i-j}}]\dlog v)- V^{i-j}([ac^{p^{i-j}}]\dlog\{c,v\})\right)\in W_n\Omega^N_{\sO_L}.}
\item[4th subcase] $n-1=i\ge j\ge n-m$.
We have
\[dV^{n-1}([a])\gamma= p^jdV^{n-1}([ac^{p^{n-1-j}}]\dlog v)-V^{n-1}([a c^{p^{n-j-1}}]\dlog c\dlog v)\]
and $v_L(ac^{p^{n-j-1}})\ge 0$.  If $v_L(ac^{p^{n-j-1}})\ge 1$, then $dV^{n-1}([a])\gamma$ is clearly regular.
If $v_L(ac^{p^{n-j-1}})=0$, then $p^{n-j-1}v_L(c)=-r$.
Since $m\le v_p(r)= n-j-1+ v_p(v_L(c))$, we see that $p|v_L(c)$. Thus $\dlog c\in \Omega^1_{\sO_L}$
and hence also $dV^{n-1}([a])\gamma\in W_n\Omega^N_{\sO_L}$.
\end{description}
\end{description}
This completes the proof of the lemma.
\end{proof}

\begin{para}\label{para:W-duality}
Consider the scheme $W_nX=(X,W_n\sO_X)$ and denote by $\pi_n: W_nX\to \Spec W_n(k)$ the structure map, it is separated and of finite-type.
Recall that  $W_n\Omega^q_{(X,D)}$  and $W_n\Omega^q_{(X,-D)}$ can be viewed as  coherent sheaves on $W_nX$, see
Proposition \ref{prop:strHWM} and Definition \ref{defn:HW-zeros}.

By \cite[I, Theorem 4.1]{Ekedahl} there is a canonical isomorphism
\eq{para:W-duality1}{W_n\Omega^N_X[N]\xr{\simeq}\pi_n^!W_n(k),}
where $\pi_n^!$ denotes the twisted inverse image from Grothendieck duality, see \cite{Ha66}, \cite{Conrad}, also \cite[4.8 and 4.10]{Lipman}.
Hence
\[\sD_{X,n}:=R\sH om_{W_n\sO_X}(-, W_n\Omega^N_X): D^b_c(W_n\sO_X)\to D^b_c(W_n\sO_X)\]
is a dualizing functor in the sense that the canonical map $\id_{D^b_c(W_n\sO_X)}\to \sD_{X,n}\circ \sD_{X,n}$ is an isomorphism, see, e.g., \cite[V, \S 10]{Ha66}. 
If $f:X\to Y$ is a proper morphism between smooth $k$-schemes and $Y$ is of pure dimension $e$, 
then  $f_n: W_nX\to W_nY$ is proper as well and Grothendieck duality yields an isomorphism
\eq{para:W-duality2}{Rf_{n*}(\sD_{X,n}(-))\cong \sD_{Y, n}(Rf_{n*}(-)) [-r],}
where $r=N-e$ is the relative dimension of $f$.
\end{para}

The following is one of the main results of the paper. For $D=\emptyset$ it is due to Ekedahl, see \cite[II]{Ekedahl},
for $D$ is reduced, it is  \cite[Theorem 5.3(1)]{Nakkajima}, cf. \cite[(3.3.1)]{Hyodo2}.
\begin{thm}\label{thm:duality-mod}
The multiplication map \eqref{lem:mult0} induces  isomorphisms, for all $q$, $n$,
\eq{thm:duality-mod1}{W_n\Omega^{N-q}_{(X,D)}\xr{\simeq} \sD_{X,n}(W_n\Omega^q_{(X,-D)})}
and 
\eq{thm:duality-mod2}{W_n\Omega^q_{(X,-D)}\xr{\simeq}\sD_{X,n}(W_n\Omega^{N-q}_{(X,D)}).}
\end{thm}
\begin{proof}
As $\sD_{X,n}$ is dualizing we obtain the isomorphism \eqref{thm:duality-mod2}  from  \eqref{thm:duality-mod1} by applying $\sD_{X,n}$.
In view of Theorem \ref{thm:strHWM} and Theorem \ref{thm:HW-zeros} the proof of \eqref{thm:duality-mod1}
is along the lines of Ekedahl's proof. We explain the strategy and the required modifications in the following. The case $n=1$ holds by Proposition \ref{prop:n=1duality}.
For $m\in \{0,\ldots, n\}$ denote by $i_m: W_mX\inj W_{n+1}X$ the closed immersion induced by the restriction
$R^{n+1-m}: W_{n+1}\sO_X\to W_m\sO_X$.
By Theorem \ref{thm:HW-zeros}, the sheaf $W_{n+1}\Omega^{q}_{(X,-D)}$ is a successive extension of $W_{n+1}\sO_X$-modules of the from
$i_{1*}M$, with $M$ a coherent locally free $\sO_X$-module. 
Thus 
\eq{eq:ext0}{
\sE xt^i_{W_{n+1}\sO_X}(W_{n+1}\Omega^q_{(X,-D)}, W_{n+1}\Omega^N_X)=0, \text{for all $i\ge 1$,}
}
by the same proof as in \cite[II, Lemma 2.2.7]{Ekedahl}.
Moreover we get that the restriction map 
\eq{thm:duality-mod3}{
\sD_{n+1,X}(W_n\Omega^q_{(X,-D)})\to j_*\sD_{n+1,U}(W_{n+1}\Omega^q_U)}
along the open immersion $j:U=X\setminus D\inj X$ is injective. Set
\[\gr^n_{(X,-D)}:=\Ker(R: W_{n+1}\Omega^q_{(X,-D)}\to i_{n*}W_n\Omega^q_{(X,-D)})\]
and
\[\gr^n_{1,(X,D)}:=\Coker(\ul{p}: i_{n*}W_n\Omega^{N-q}_{(X,D)}\to W_{n+1}\Omega^{N-q}_{(X,D)}).\]
Consider  the following  diagram 
\[\xymatrix{
0\ar[r]& 
i_{n*}W_n\Omega^{N-q}_{(X,D)}\ar[r]^{\ul{p}}\ar[d] &
W_{n+1}\Omega^{N-q}_{(X,D)}\ar[r]\ar[d] &
\gr^n_{1, (X,D)}\ar[r]\ar[d] &
0\\
0\ar[r]&
i_{n*}\sD_{X,n}(W_n\Omega^q_{(X,-D)})\ar[r]&
\sD_{X,{n+1}}(W_{n+1}\Omega^q_{(X,-D)})\ar[r]&
\sD_{X, n+1}(\gr^n_{(X,-D)})\ar[r]&
0,
}\]
where the vertical maps are induced by the multiplication map \eqref{lem:mult0} and the horizontal map on the bottom left side
is induced by the duality isomorphism
$i_{n*}\circ \sD_{X,n}=\sD_{X,n+1}\circ i_{n*}$ 
composed with the dual of the restriction $R$. 
As the restriction map \eqref{thm:duality-mod3} is injective it follows from \cite[II, (2.2.8), (2.2.9)]{Ekedahl} that the diagram commutes.
(Note that it is a diagram of sheaves by the vanishing of the Ext-groups shown above.)
As the lines are exact it remains by induction to show that the right vertical arrow is an isomorphism.
In view of the exact sequences in Theorem \ref{thm:strHWM} and Theorem \ref{thm:HW-zeros}\ref{thm:HW-zeros1}
the same argument as in \cite[II, Lemma 2.2.17]{Ekedahl} (see also the explanation around \eqref{prop:n=1duality2})
reduces us to show that the maps
\eq{thm:duality-mod4}{ Z_n\Omega^{N-q}_{(X,D)} \to \sH om((\Omega/B)^q_{n, (X,-D)},\Omega_X^N), \quad \text{and} \quad
B_n\Omega^{N-q+1}_{(X,D)} \to \sH om((\Omega/Z)^{q-1}_{n, (X,-D)},\Omega_X^N),}
induced by $\alpha\mapsto (\beta\mapsto C^n(\alpha\wedge\beta))$, 
are isomorphisms. These maps are exactly $\sD_{X,1}$ applied to the isomorphisms 
in part \ref{thm:HW-zeros2} of Theorem \ref{thm:HW-zeros} and hence are isomorphisms as well.
This completes the  proof.
\end{proof}

\begin{cor}\label{cor:pbf-buf-zeros}
The projective bundle formula and the blow-up formula from \ref{para:pbf-buf-gt}, \ref{para:pbf} and \ref{para:buf}
also hold for $D$ replaced by $-D$.
\end{cor}

Grothendieck duality \eqref{para:W-duality2} yields:

\begin{cor}\label{cor:duality-mod}
Assume additionally that $X$ is proper over $k$. Then there is a canoncial isomorphism of finite $W_n(k)$-modules 
\[H^j(X, W_n\Omega^q_{(X,-D)})\cong \Hom_{W_n(k)}\left(H^{N-j}(X, W_n\Omega^{N-q}_{(X,D)}), W_n(k)\right).\]
\end{cor}

\begin{cor}\label{cor:covMCor}
Assume $X$ is proper and let $(Y,E)$  be another modulus pair with $Y\in \Sm$ proper of pure dimension $N_Y$ 
and $E_{\red}$ an SNCD and  denote by $\uMCor((X,D), (Y,E))$ the group of left proper admissble correspondences, 
see \cite[Definition 1.3.1]{KMSYI}.
Then for any correspondence $\alpha\in \uMCor((X,D), (Y,E))$ we have a natural map
\[\alpha_* : H^N(X, W_n\Omega^{N-q}_{(X,-D)})\to H^{N_{Y}}(Y, W_n\Omega^{N_Y-q}_{(Y,-E)})\]
such that:
\begin{enumerate}
    \item If $(Z,F)$  is a modulus pair with $Z\in \Sm$ proper of pure dimension $N_Z$ and $F_{\red}$ an SNCD, and 
    $\beta\in \uMCor((Y,E), (Z,F))$ then 
    \[\beta_*\alpha_*=(\beta\circ \alpha)_*: 
    H^N(X, W_n\Omega^{N-q}_{(X,-D)})\to H^{N_{Z}}(Y, W_n\Omega^{N_Z-q}_{(Z,-F)}).\]
    \item If $\alpha$ is induced by the graph of a morphism $f:X\to Y$ satisfying $D\ge f^*E$, then $\alpha_*=f_*$,
    the dual of $f^*$.
\end{enumerate}
In particular, if $f:X\to Y$ is a morphism which induces an isomorphism $X\setminus f^{-1}(E)\xr{\simeq} Y\setminus E$,
then $f_*$ induces an isomorphism $f_*: H^N(X, W_n\Omega^q_{(X, f^*E)})\xr{\simeq}  H^N(Y, W_n\Omega^q_{(Y,E)})$.
\end{cor}
\begin{proof}
Use Corollary \ref{cor:duality-mod} and the corresponding properties of $H^0(X, W_n\Omega^q_{(X,D)})$.
\end{proof}

Recall that following  an idea of Deligne \cite[Appendix]{Ha66} Hartshorne defines in \cite[\S 2]{Ha72} 
cohomology with compact support for coherent sheaves on schemes which are separated and of finite type over a field $k$.
In fact the same approach works for schemes of finite type over an artinian ring, as it is only used 
that the cohomology of a coherent sheaf on a proper scheme has finite length over the base ring, to guarantee 
the exactness of certain limits. Let $A$ be an artinian ring,  $V\to \Spec A$ be a separated morphism of finite type,
and $G$ a coherent sheaf on $V$, then the corresponding cohomology with compact support is denoted by 
$H^j_c(V, G)$.

\begin{cor}\label{cor:HW-compactsupp}
With the above notation and assuming that $X$ is proper there are canonical isomorphisms
\[H^j_c(W_n(U), W_n\Omega^q_U)\cong \varprojlim_r H^{j}(X, W_n\Omega^q_{(X,-rD)})\cong 
\Hom_{W_n(k)}\left(H^{N-j}(U, W_n\Omega^{N-q}_U), W_n(k)\right).\]
\end{cor}
\begin{proof}
In view of the definition of $H^j_c(V,G)$ and Theorem \ref{thm:duality-mod} we only have to check that 
the pro-systems of coherent $W_n\sO_X$-modules $\{ W_n\Omega^q_{(X,-rD)}\}_r$ and 
$\{W_n(I^r)\cdot W_n\Omega^q_X\}_r$ are isomorphic, where $I=\Ker(\sO_X\to \sO_D)$.
To this end, note that by definition we have an inclusion $W_n(I^r)\cdot W_n\Omega^q_X\subset W_n\Omega^q_{(X,-rD)}$,
for all $r\ge 1$. Furthermore
we have an inclusion 
\eq{cor:HW-compactsupp1}{W_n\Omega^q_{(X, -p^{n+1}rD)}\subset W_n(I^{p r})\cdot W_n\Omega^q_X.} 
To check this  we may assume $X=\Spec A$ and $rD=\Div(f)$.
It follows from the recursive definition in  \eqref{para:twisted-BZ-zeros1} that 
$(\Omega/B)^q_{j,n}(0, -prD)$ is a quotient of $f^{p^{j+1}r}\Omega^q_A$ and similarly with 
$(\Omega/Z)^{q-1}_{j,n}(0, -prD)$. Thus ${\rm gr}^n_{(X,-D)}$ is 
by Propositions \ref{prop-ex-seq-zeros} and  \ref{prop:BZnn-zeros} generated by elements 
\[V^n([f^{p^{n+1}r}]\alpha)= [f]^{pr}V^n(\alpha)\quad \text{and}\quad
  dV^n([f^{p^{n+1}r}]\beta)= [f]^{pr}dV^n(\beta),\]
where $\alpha\in \Omega^q_A$, $\beta\in \Omega^{q-1}_A$. Thus \eqref{cor:HW-compactsupp1} holds by induction.
\end{proof}

The following generalizes Ekedahl's isomorphism \eqref{para:W-duality1} to the case of thickenings of smooth $k$-schemes.
\begin{cor}\label{cor:twisted-inverse-image-Fil}
Let $H\subset X$ be smooth closed subscheme of codimension 1. Denote by $H_r\subset X$ the closed subscheme defined by the ideal
sheaf $\sO_X(-rH)$, for $r\ge 1$. Denote by $f:W_n(H_r)\to \Spec W_n(k)$ the structure map.
Then there is an isomorphism in $D(W_n\sO_{H_r})$
\[f^!W_n(k)\cong \frac{W_n\Omega^N_{(X,rH)}}{W_n\Omega^N_X}[N-1].\]
\end{cor}
\begin{proof}
Denote by $i: H_r\inj X$ the closed immersion. Ekedahl's isomorphism \eqref{para:W-duality1} and Grothendieck duality yield
an isomorphism 
\[\sD_{X,n}(i_*W_n\sO_{H_r})[N]\cong i_*f^!W_n(k).\]
Now apply $\sD_{X,n}$ to the exact sequence
\[0\to W_n\sO_{(X,-rH)}\to W_n\sO_X\to i_*W_n\sO_{H_r}\to 0\]
and use Theorem \ref{thm:duality-mod} to obtain the isomorphism from the statement.
\end{proof}

\section{Duality for Hodge-Witt cohomology with modulus: infinite level}
Throughout this section  $X\in \Sm$ has pure dimension $\dim X=N$ and $D$ is an effective Cartier divisor on $X$ 
such that $D_{\red}$ is an SNCD. We denote by 
\eq{i-restr}{i: W_{n-1}X\inj W_nX}
the closed immersion induced by the restriction $R: W_n\sO_X\to W_{n-1}\sO_X$ and by 
\eq{sigma-Frob}{\sigma: W_nX\to W_nX}
the finite morphism induced by $W_n(F_X^*): W_n\sO_X\to W_n\sO_X$, where $F_X: X\to X$ denotes the absolute Frobenius on $X$.
With this  notation, e.g.,  the Verschiebung $V: (\sigma i)_* W_{n-1}\Omega^q_X\to W_n\Omega^q_X$ and the differential
$d: (\sigma^n)_*W_n\Omega^q_X\to (\sigma^n)_*W_n\Omega^{q+1}_X$ are $W_n\sO_X$-linear morphisms.

\begin{para}\label{para-dir-inv-dRW}
We recall some definitions from \cite[III, Definition 2.1]{Ekedahl}.
An {\em inverse de Rham-Witt system} on $X$ consists of a family of quasi-coherent graded $W_n\sO_X$-modules $M_n$, $n\ge 1$, together with 
morphisms of graded $W_n\sO_X$-modules
\[R: M_{n}\to i_*M_{n-1}, \, F: M_n\to (\sigma i)_*M_{n-1},\, V: (\sigma i)_*M_{n-1}\to M_n,\,  d:(\sigma^n)_*M_n\to (\sigma^n)_*M_n(1),\]
where $M_0=0$ and for a graded module $M$ we denote by $M(1)$ the graded module with $M(1)^q=M^{q+1}$, satisfying the following identities
\[RF=FR,\quad RV=VR,\quad Rd=dR, \quad FV=p,\quad VF=p, \quad d^2=0,\quad FdV=d.\]
A {\em direct de Rham-Witt system}  consists of a family of quasi-coherent graded $W_n\sO_X$-modules $M_n$, $n\ge 1$, together with maps $F$, $V$, $d$ as above, 
and (instead of the maps $R$)  morphisms of graded $W_n\sO_X$-modules
\[\ul{p}: i_*M_{n-1}\to M_n,\]
such that $\ul{p}$, $F$, $V$, and $d$ satisfy the same identities as above with $R$ replaced by $\ul{p}$.
We call a direct (resp. inverse) de Rham-Witt system {\em coherent} if each $M_n$ is a coherent $W_n\sO_X$-module.
We denote by 
\[\idrw_X\quad \text{and} \quad \ddrw_X\]
the categories of inverse- and direct de Rham-Witt systems, respectively, with morphisms defined in the obvious way,
and by $\idrw_{X,c}$ and $\ddrw_{X,c}$ their full subcategories of coherent objects.
These are abelian categories and the forgetful functors to the product category $\prod_n (\text{graded } W_n\sO_X\text{-modules})$ sending
$( (M_n)_n, F,V,d, R \text{(resp. $\ul{p}$)})$ to $(M_n)_n$ is exact and conservative. 
For $M\in \ddrw$ we denote by $M(j)$ the direct de Rham-Witt system with shifted grading $M_n(j)^q=M_n^{q+j}$, and where $R_{M(j)}=R_M$, $F_{M(j)}=F_M$, $V_{M(j)}=V_M$, 
and $d_{M(j)}=(-1)^j d_M$. 

Assume $X$ is proper. Then the global section functor $\Gamma(X,-)$ derives to
\[R\Gamma(X,-): D^b(\idrw_{X,c})\to D^b(\idrw_{k,c}),\quad D^b(\ddrw_{X,c})\to D^b(\ddrw_{k,c}),\]
where $\idrw_{k,c}:=\idrw_{\Spec k, c}$ and similarly with $\ddrw$, 
see \cite[III, (2.5.2)]{Ekedahl}.
\end{para}

\begin{example}\label{exa-dir-invDRW}
We have
 \[ W_{\infty\bullet}\Omega^*_{(X,\pm D)}:=((W_n\Omega^*_{(X,\pm D)})_n,\, \ul{p}, \,F, \,V,\, d) 
 \in \ddrw_{X,c},\]
and
 \[W_{\bullet}\Omega^*_{(X,\pm D)}:=((W_n\Omega^*_{(X,\pm D)})_n, \,R,\,F,\, V,\, d)\in\idrw_{X,c}.\]
\end{example}

\begin{para}\label{para-dualizing-Witt-sys}
Recall from \cite[III, Example 2.2.1]{Ekedahl} that the  top Witt forms $W_n\Omega^N_X$, $n\ge 1$, together with the maps $\ul{p}$, $V$, and 
the $n$-fold Cartier operators $C^n: (\sigma^n)_*W_n\Omega^N_X\to W_n\Omega^N_X$ form a  dualizing system in the sense of 
\cite[III, Definition 2.2]{Ekedahl}. Hence we obtain a well-defined functor
\[(\ddrw_{X,c})^o\to \idrw_{X}, \quad M=(M_n)_n \mapsto (\uHom_{W_n\sO_X}(M_n, W_n\Omega^N_X))_n.\]
We only remark that the maps  $R$, $F$, $V$, and $d$ on $(\uHom_{W_n\sO_X}(M_n, W_n\Omega^N_X))_n$
are induced by $\ul{p}_M^*$, $V_M^*$, $F_M^*$, and $d_M^*$, respectively, and refer to \cite[p. 205 -- 206]{Ekedahl} for the details, see also \cite[1.6.6]{CR12}.

Let $E_{n,X}$ be the Cousin complex of $W_n\Omega^N_X$. By \cite[III, Example 2.2]{Ekedahl}
$E=(E_{n,X})_n$ has the structure of a complex of dualizing systems (and is a Witt residual complex in the sense of \cite[Definition 1.8.3]{CR12}).
We obtain a functor between the derived categories
\[\sD_{X,\bullet}: D^b(\ddrw_{X,c})^o\to D^b(\idrw_X), \quad M\mapsto \sD_{X,\bullet}(M),\]
where 
\[(\sD_{X,\bullet}(M))_n= \uHom_{W_n\sO_X}(M_n, E_{n,X})=\sD_{X,n}(M_n) \quad \text{in } D^b_c(W_n\sO_X\text{-mod}),\]
with $\sD_{X,n}$ as in \ref{para:W-duality} and where the degree $q$-part of $\sD_{X,n}(M_n)$
is $\sD_{X,n}(M_n^{-q})$.
\end{para}

\begin{cor}\label{cor:duality-W-systems}
The multiplication map induces  isomorphisms  in $D^b(\idrw_{X,c})$
\[W_\bullet\Omega^*_{(X,-D)}\cong \sD_{X,\bullet}(W_{\infty\bullet}\Omega^*_{(X,D)})(-N)
\quad \text{and}\quad
W_\bullet\Omega^*_{(X,D)}\cong \sD_{X,\bullet}(W_{\infty\bullet}\Omega^*_{(X,-D)})(-N).\]
\end{cor}
\begin{proof}
 By Theorem \ref{thm:duality-mod} we have the isomorphisms of graded $W_n\sO_X$-modules
 \[W_n\Omega^*_{(X,\pm D)}\cong \sH om_{W_n\sO_X}(W_{n}\Omega^*_{(X,\mp D)}, W_n\Omega^N_X)(-N),\]
 for all $n$.   Thus it suffices to show that these isomorphisms are compatible with the maps 
 $R,F,V,d$ defined on both sides. As before this can be checked after restriction to $U=X\setminus D$
 and hence follows from \cite[III, Proposition 2.4]{Ekedahl}
\end{proof}

\begin{para}\label{para:CDR}
Recall from \cite[I, (1.1)]{IllRay} that the Cartier-Dieudonn\'e-Raynaud ring $\mathcal{R}$ 
is the graded (non-commutative) $W=W(k)$-algebra,
$\sR=\sR^0\oplus \sR^1$ generated by symbols  $F$, $V\in \sR^0$ and $d\in \sR^1$ subject to the relations 
\[ \sigma(a) F= F a,\, aV=V \sigma(a), \, da=ad,\, FV=p=VF,\, dd=0,\, FdV=d,\]
for $a\in W$, where $\sigma:W\to W$ denotes the Frobenius lift.

We obtain a functor $\lim: \idrw_{X,c}\to {\rm Sh}(X,\sR)$, the category of sheaves of 
$\sR$-modules on $X$, where the limit is taken along the restriction map $R$.
This functor derives to
\[R\lim: D^b(\idrw_{X,c})\to D^b(X,\sR),\]
see \cite[III, (2.6.2)]{Ekedahl}. For $M\in {\rm Sh}(X,\sR)$, the shift $M(j)$ is defined by
$M(j)^q=M^{j+q}$ and the action of $F,V,d$ on $M(j)$ is induced by the action of $F,V, (-1)^j d$ on $M$.
\end{para}

\begin{prop}\label{prop:modulus-limit}
The following equalities hold in $D^b(X,\sR)$
\[R\lim W_\bullet\Omega^*_{(X,-D)}=W\Omega^*_{(X,-D)}
:=\Ker\left(W\Omega^*_X\to \bigoplus_i W\Omega^*_{D_i}\right),\]
where we write $D=\sum_i D_i$ with $D_{i,\red}$ smooth, and 
\[R\lim W_\bullet\Omega^*_{(X,D)}=W\Omega^*_X(\log D),\]
the log de Rham-Witt complex, see, e.g., \cite[Proposition-Definition 3.10]{Matsuue}.
In particular, 
\[R\lim W_\bullet\Omega^*_{(X,D)}=R\lim W_\bullet\Omega^*_{(X,D_\red)}.\]
\end{prop}
\begin{proof}
The statement for $W_\bullet\Omega^*_{(X,-D)}$ follows directly 
from the surjectivity of $R$, see Theorem \ref{thm:HW-zeros}\ref{thm:HW-zeros1.5}, and 
the left exactness of $\lim$. 
By \eqref{rmk:HW-modulus2} we have a natural inclusion,
\eq{prop:modulus-limit1}{W_\bullet\Omega^*_X(\log D)\inj W_\bullet\Omega^*_{(X,D)}.}
Let $M>0$ be a positive integer such that $p^M$ is strictly larger than all the multiplicities of $D$.
We claim
\[R^M(W_{n+M}\Omega^*_{(X,D)})=W_n\Omega^*_X(\log D),\quad \text{for all }n\ge 1.\]
By the above we only have to prove this ``$\subset$" inclusion.
It suffices to prove this after pullback to henselian dvf's  $L$ of geometric type.
Thus by Theorem \ref{thm:HW-modulus} it suffices to show for $r<p^M$ and $s=0,\ldots, n-1$,
\eq{prop:modulus-limit2}{R^M(p^s\fil_{rp^s}W_{n+M}\Omega^*_L)\subset W_n\Omega^*_{\sO_L}(\log).}
Under our assumptions on $M$, the inequality 
\[\frac{-rp^s}{p^{n+M-i-1}}\le -1  \]
is only possible if $i>n-s-1$. But in this case 
$p^sV^i([a])=0$ in $W_n(L)$, for $a\in L$. 
In view of the definition of $\Fil^p_rW_n\Omega^q_L$ in Definition \ref{defn:p-sat} this yields
the inclusion \eqref{prop:modulus-limit2}. 
It follows that the inclusion \eqref{prop:modulus-limit1} is an isomorphism of pro-objects and hence
\[R\lim W_\bullet\Omega^*_{(X,D)}= R\lim W_\bullet\Omega_X^*(\log D)
=W\Omega_X^*(\log D),\]
where the last equality follows from the surjectivity of  the restriction map $R$ 
on $W_\bullet\Omega^*_X(\log D)$. 
\end{proof}

\begin{thm}\label{thm:limit}
Assume additionally that $X$ is proper and set
\[W_\infty\Omega^*_{(X,D)}:=\colim \, W_{\infty\bullet}\Omega^*_{(X,D)},\]
see Example \ref{exa-dir-invDRW}.
Then we have isomorphisms in $D^b(\sR)$
\[R\Gamma(X, W\Omega^*_{(X,-D)})\cong  R\Hom_{W}(R\Gamma(X, W_\infty\Omega^*_{(X,D)}), K/W )(-N)[-N],\]
where  $K=W[1/p]$.
In particular we have isomorphisms of $\cR^0$-modules, for $i,q\ge 0$,
\[H^i(X, W\Omega^q_{(X,-D)})\cong \Hom_W(H^{N-i}(X, W_\infty\Omega^{N-q}_{(X,D)}), K/W).\]    
\end{thm}
\begin{proof}
Applying $R\lim\circ R\Gamma(X,-)$ to the first isomorphism in Corollary \ref{cor:duality-W-systems} 
gives
\begin{align*}
R\Gamma(X, W\Omega^*_{(X,-D)}) 
&\cong
R\Gamma(X, R\lim W_\bullet\Omega^*_{(X,-D)})
\quad\text{(by \Cref{prop:modulus-limit})}
\\
&\cong R\lim R\Hom_{W_n}(R\Gamma(X, W_n\Omega^*_{(X,D)}), W_n)(-N)[-N]
\quad\text{(by \eqref{para:W-duality2})}
\\
& \cong R\lim R\Hom_{W}(R\Gamma(X, W_n\Omega^*_{(X,D)}), K/W)(-N)[-N],
\end{align*}
where the transition maps in the last limit are given by precomposition with $\ul{p}$.
Here the third isomorphism holds by definition of the restriction maps on $\sD_{X,\bullet}(M)$ in \ref{para-dualizing-Witt-sys}, 
see \cite[III, (2.3.2)]{Ekedahl}.
Let $\sU$ be an affine covering of $X$ and denote by $C_n=C^\bullet(\sU, W_n\Omega^*_{(X,D)})$ the \v{C}ech complex
(it is a complex of graded $W_n[d]$-modules).
Then $(C_n)_n$ is a representative 
of $R\Gamma(W_{\infty\bullet}\Omega_{(X,D)}^*)$ in $D^b(\ddrw_{\Spec k, c})$
and $\colim_{\ul{p}}\, C_n$ is a representative of $R\Gamma(W_\infty\Omega^*_{(X,D)})$.
As $\ul{p}: W_n\Omega^*_{(X,-D)}\to W_{n+1}\Omega^*_{(X,-D)}$ is injective,  the induced map $\ul{p}: C_n\to C_{n+1}$
is term-wise injective. 
As $K/W$ is an injective $W$-module, $\ul{p}^*: \Hom(C_{n+1}, K/W)\to \Hom(C_n, K/W)$ is term-wise surjective and thus
(e.g. \cite[\href{https://stacks.math.columbia.edu/tag/07KW}{Tag 07KW(5)}]{stacks-project})
\begin{align*}
R\lim R\Hom_{W}(R\Gamma(X, W_n\Omega^*_{(X,D)}), K/W) & = \lim \Hom_{W}(C_n,K/W)\\
                                                     & =\Hom_W(\colim\, C_n, K/W),
\end{align*}
which implies the statement.
\end{proof}

\begin{para}
We denote by $W[F]$ the free $W$-algebra 
with generator $F$ subject to the relation $F\cdot a=\sigma(a)\cdot F$, for all $a\in W$, where $\sigma: W\to W$ is the Frobenius. 
Let
\[R\Gamma_{\rm crys}((X,-D)/W)\in D^b(W[F])\]
be the simple  complex associated to  the double complex $R\Gamma(X, W\Omega^*_{(X,-D)})$, where $F$ 
acts via the absolute Frobenius $F_X^*$, cf. \cite[III, (5.3), (5.4)]{Ekedahl}.
Furthermore denote by
\[\sF:  R\Hom_{W}(R\Gamma(X, W_\infty\Omega^{q}_{(X,D)}), K/W )\to  R\Hom_{W}(R\Gamma(X, W_\infty\Omega^{q}_{(X,D)}), K/W )\]
the morphism, which sends a map $\mu$ to the composition
\[R\Gamma(X, W_\infty\Omega^{q}_{(X,D)})\xr{p^{N-q}\cdot V} R\Gamma(X, W_\infty\Omega^{q}_{(X,D)})\xr{\mu} K/W\xr{\sigma} K/W.\]
Representing $R\Gamma(X, W_\infty\Omega^{q}_{(X,D)})$ by a \v{C}ech complex $C^{\bullet, q}_{\infty}$ 
we find that $(\Hom_W(C^{\bullet, q}_\infty, K/W), \sF)$ is a complex of $W[F]$-modules.
As 
\[d: \Hom_W(C^{\bullet, q}_\infty, K/W)\to \Hom_W(C^{\bullet, q-1}_\infty, K/W)\]
is given by (see \cite[1.6.6]{CR12})
\[\mu\mapsto (-1)^{q+1}\mu\circ d\]
it follows from the relation $Vd=pdV$ that 
$(\Hom_W(C^{\bullet,\bullet}_\infty, K/W), \sF)$ is a double complex of $W[F]$-modules.
We denote the image of the associated simple complex in $D^b(W[F])$ by
\[R\Hom_W(R\Gamma_{\rm crys}((X,D)/W_\infty), K/W).\]
\end{para}

\begin{cor}\label{cor:limit-crys}
Assume $X$ is proper.
Multiplication in $W_\bullet\Omega^*_X$ induces an isomorphism in $D^b(W[F])$
\[R\Gamma_{\rm crys}((X,-D)/W)\cong R\Hom_W(R\Gamma_{\rm crys}((X,D)/W_\infty), K/W)[-2N].\]  
\end{cor}
\begin{proof}
The isomorphism follows by applying the  functor  which sends a double complex to its associated simple complex from 
the isomorphism in Theorem \ref{thm:limit}. 
For the compatibility with $F_X^*$ on the left and $\sF$ on the right we  observe:
\begin{itemize}
  \item for $\alpha \in W_{n+m}\Omega^q_{(X,-D)}$ and $\beta\in W_n\Omega^{N-q}_{(X,D)}$ ($m\ge 1$) we have 
           \[V(R^m(F_X^*(\alpha))\beta)= V(p^qF(R^{m-1}\alpha)\beta)=R^{m-1}(\alpha)\cdot p^qV(\beta);\]
  \item under the identification $W_n(k) \cong \frac{1}{p^n}W/W\subset K/W$ the inverse map to $V:W_n(k)\to W_{n+1}(k)$
  is induced by $\sigma: K/W\to K/W$.
\end{itemize}
Now the statement follows from this and the definition of $\sF$.
\end{proof}

\begin{para}\label{para:coh-R-mod}
Set $\sR_n=\sR/(V^n\sR+dV^n\sR)$. It is a right $\sR$-module and a left $W_n[d]$-module.
By \cite[(2.3)]{Illusie-Ekedahl}  $\sR_n\otimes_{\sR}$  derives to 
\[\sR_n\otimes^L_{\sR}: D^b(\sR)\to D^b(W_n[d]).\]
Recall that a complex $M\in D^b(\sR)$ (or an $\sR$-module $M$) is called {\em coherent} (in the sense of Illusie-Raynaud-Ekedahl)
if the following two conditions are satisfied
\begin{enumerate}
\item $H^i(\sR_n\otimes^L_{\sR} M)$ is a finitely generated $W_n$-module for all $i$, and 
\item $M$ is complete, i.e., $M=R\lim_n (\sR_n\otimes^L_{\sR} M)$.
\end{enumerate}
By \cite[p. 190]{Ekedahl} (see also \cite[2.4]{Illusie-Ekedahl}) a complex $M$ is coherent if and only if $H^i(M)$ is a coherent $\sR$-module
in the sense of Illusie-Raynaud \cite[I, D\'efinition 3.9]{IllRay}. The coherent subcomplexes form a triangulated subcategory
of $D^b(\sR)$, which is denoted by $D^b_c(\sR)$.
Ekedahl shows in \cite[IV, Proposition 1.1]{Ekedahl} that there is a well-defined functor
\eq{para:coh-R-mod1}{\sD: D^b_c(\sR)^{\rm op}\to D^b_c(\sR), \quad 
M\mapsto \sD(M)=R\lim R\Hom_{W_n}(\sR_n\otimes^L_{\sR} M, W_n),}
which is dualizing in the sense that it satisfies $\sD\circ\sD=\id$. Moreover for any $M\in D^b_c(\sR)$
\eq{para:coh-R-mod2}{\sR_n\otimes^L_{\sR} \sD(M)=R\Hom_{W_n}(\sR_n\otimes^L_{\sR} M, W_n)\quad \text{in } D^b(W_n[d]).}
\end{para}

In case $D=\emptyset$ the following corollary is \cite[II, Th\'eor\`eme 2.2]{IllRay}, \cite[III, Theorem 2.9]{Ekedahl},
and \cite[II, (1.4.3)]{IllRay}. The last part on the compatibility with 
$\sR_n\otimes_{\sR}^L$ follows also from \cite[Theorem 6.24]{Nakkajima}, but the proof is different,
see also \cite[Lemme 1.3.3]{Mokrane}.
\begin{cor}\label{thm:Dred-coh}
Assume $X$ is proper. We have 
\[R\Gamma(X, W\Omega^*_{(X, \pm D_{\rm red})})  \in D^b_c(\sR),\]
and 
\eq{thm:Dred-coh0}{R\Gamma(X, W\Omega^*_{(X, -D_{\rm red})})\cong\sD(R\Gamma(X, W\Omega^*_{(X, D_{\rm red})}))(-N)[-N],}
where $\sD$ is the dualizing functor \eqref{para:coh-R-mod1}. 
Moreover,
\[\sR_n\otimes^L_{\sR} R\Gamma(X, W\Omega^*_{(X, \pm D_{\rm red})})= R\Gamma(X, W_n\Omega^*_{(X, \pm D_{\rm red})}).\]
\end{cor}
\begin{proof}
Write $D_{\red}=\sum_{j=1}^r D_j$, with $D_j$ smooth and connected. 
Set 
\[D^{(s)}=\coprod_{(j_1<\ldots< j_s)} D_{j_1}\cap \ldots\cap D_{j_s}.\]
By \cite[9.]{Matsuue} (see also \cite[Proposition 9.3]{Nakkajima})  and \eqref{rmk:HW-modulus2} 
there is an increasing filtration (the weight filtration) 
$\{P_s\}_{s=0,\ldots, r}$ of $\sR$-submodules of $W\Omega^*_{(X,D_{\rm red})}$ 
with $P_0=W\Omega^*_X$ and $P_r=W\Omega^*_{(X, D_{\rm red})}$, which fits into exact sequences of $\sR$-modules
\eq{thm:Dred-coh1}{0\to P_{s-1}\to P_s\xr{\rho} W\Omega^*_{D^{(s)}}(-s)\to 0.}
As it is not stated explicitly in {\em loc. cit.} we remark that $\rho$ is a morphism of $\sR$-modules.
Indeed,  the iterated residue map $\rho$ is induced by the inverse of the isomorphism 
\[W\Omega^*_{D^{(s)}}(-s)\to P_s/P_{s-1},\]
which sends a local form $\alpha\in W\Omega^{q-s}_{D_{j_1}\cap \ldots\cap D_{j_s}}$ to $ \dlog \{t_{j_1},\ldots, t_{j_s}\}\cdot \tilde{\alpha}$,
where $\tilde{\alpha}\in W\Omega^{q-s}_X$ is a lift of $\alpha$ and $t_j$ is a local equation for $D_j$.
This is clearly a morphism of $\sR$-modules (see \ref{para:CDR} for the ``$(-s)$" shift of an $\sR$-module).
As $R\Gamma(D^{(s)}, W\Omega^*_{D^{(s)}})\in D^b_c(\sR)$, for all $s=0,\ldots, r$, by \cite[II, Th\'eor\`eme 2.2]{IllRay}, we find 
$R\Gamma(X, W\Omega^*_{(X, D_{\rm red})})\in D^b_c(\sR)$.
Similarly it follows from \cite[II, (1.4.3)]{IllRay} (applied to $W\Omega^*_{D^{(s)}}(-s)$) that we have 
\[\sR_n\otimes^L_{\sR} R\Gamma(X, W\Omega^*_{(X, D_{\rm red})})= R\Gamma(X, W_n\Omega^*_{(X, D_{\rm red})}).\]
Hence applying $R\lim\circ R\Gamma$ to \eqref{thm:duality-mod2} (for $D_{\rm red}$)
yields \eqref{thm:Dred-coh0} in view of the discussion  in \ref{para-dir-inv-dRW}, \ref{para-dualizing-Witt-sys}, \ref{para:CDR}.
The remaining statements for $-D_{\red}$ follow from \eqref{para:coh-R-mod1} and \eqref{para:coh-R-mod2}.
\end{proof}

\begin{para}\label{para:conseq-coh}
We can now use the results on coherent $\sR$-modules  from  \cite{IllRay} and \cite{Ekedahl} to get many consequences, which
however only work in the case $D$ is reduced or we work up to bounded torsion\footnote{We say that a certain statement 
in an additive category $C$ is true up to bounded torsion, 
if it is true in the localized category $C_{\Q}$ which has the same objects but the Hom's are tensored with $\Q$.}. 
Here is a sample of corollaries, where we assume that $X$ is proper additionally to our standing assumptions:
\begin{enumerate}[label=(\arabic*)]
\item\label{para:conseq-coh1} $R\Gamma(X, W\Omega^*_{(X,D_{\rm red})})$ is a complex of  $\sR$-modules of level $N$, 
i.e., it is represented by a complex of $\sR$-modules
which vanishes in degrees outside $[0,N]$ and on which $F$ is invertible in degree $N$. 
(That  $F$ is invertible in degree $N$ in the case at hand, follows by induction from the 
exact sequences \eqref{thm:Dred-coh1} and  \cite[I, Proposition 3.7]{IlDRW}.) 
Therefore by \cite[III, 5.]{Ekedahl} taking the simple complex of the duality isomorphism  
in Theorem \ref{thm:Dred-coh} we get an isomorphism  of $F$-crystals of level $N$
\[R\Gamma(X, W\Omega^\bullet_{(X, -D_{\rm red})})\cong R\Hom_{W}(R\Gamma(X, W\Omega^\bullet_{(X, D_{\rm red})}), W)[-2N].\]
\item\label{para:conseq-coh2} We have 
\[R\Gamma(X, W\Omega^\bullet_{(X, D_{\rm red})})=R\Gamma(X, W\Omega^\bullet_X(\log D))\cong R\Gamma_{\text{log-crys}}((X,D)/W),\]
where the equality holds by \eqref{rmk:HW-modulus2} and the isomorphism, e.g., by \cite[Theorem 7.2]{Matsuue}.
It is a perfect complex of $W$-modules by, e.g.,  \cite[II, Th\'eor\`eme 2.7]{IlDRW} and \eqref{thm:Dred-coh1}.
Thus duality yields an isomorphism
\[R\Gamma(X, W\Omega^\bullet_{(X, -D_{\rm red})})\cong R\Hom_W(R\Gamma_{\text{log-crys}}((X,D)/W), W)[-2N].\]
Hence $R\Gamma(X, W\Omega^\bullet_{(X, -D_{\rm red})})$ is a perfect complex as well.
\item\label{para:conseq-coh3} As the cokernel of the natural inclusion $W\Omega^*_{(X, -D)}\to W\Omega^*_{(X, -D_{\rm red})}$ is annihilated by 
some fixed $p$-power (take $p^r$ greater than all the multiplicities in $D$) it follows that
\[R\Gamma(X, W\Omega^*_{(X, -D)})\to R\Gamma(X, W\Omega^*_{(X, -D_{\rm red})})\]
is an isomorphism up to bounded $p$-primary torsion.
As the complex on the right is a coherent $\sR$-module, it follows that
the cohomology groups $H^i(X, W\Omega^q_{(X,-D)})$ are finitely generated $W$-modules up to bounded torsion, and 
$R\Gamma(X, W\Omega^\bullet_{(X, -D)})$ is a perfect complex of $W$-modules up to bounded (but possibly infinitely generated) torsion. 
\item\label{para:conseq-coh4} As observed in \cite[Remark 5.4(1)]{Nakkajima} it follows from Shiho's comparison 
of log-crystalline cohomology for $(X,D)$
with rigid cohomology, see \cite[Corollary 2.4.13 and Theorem 3.1.1]{Shiho}, 
and Berthelot's Poincar\'e duality for rigid cohomology, see \cite[Th\'eor\`eme 2.4]{Berthelot-PD-KunnethRig}, that the above
yields an isomorphism
\[R\Gamma(X, W\Omega^\bullet_{(X,-D)})\otimes_W K= R\Gamma_{c, {\rm rig}}(U),\]
where $U=X\setminus D$ and $K=\Frac(W)$. Thus for any effective Cartier divisor $D$ with SNC support equal to $X\setminus U$ 
the complex $R\Gamma(X, W\Omega^\bullet_{(X,-D)})$ is  perfect up to bounded torsion and  is an integral model 
for the compactly supported rigid cohomology of $U$. 
\item\label{para:conseq-coh5} By \cite[Corollary 2.5.4]{Illusie-Ekedahl} and the above the slope spectral sequence
\[E^{i,j}_1= H^i(X, W\Omega^j_{(X,-D)})\Longrightarrow H^*(X, W\Omega^\bullet_{(X,-D)})\]
degenerates up to bounded torsion and thus by the same argument
\footnote{This argument uses that $V$ acts topologically nilpotent on  $H^i(X, W\Omega^{\le j}_{(X, -D_{\rm red})})$ in the $p$-adic topology, 
which follows from \cite[II, Proposition 2.10]{IlDRW} and the resolution 
$W\Omega^{\le j}_{(X, -D_{\rm red})}\to W\Omega^{\le j}_{D^\bullet}$, see \cite[Lemme 3.15.1]{Mokrane} and \cite[Corollary 6.28]{Nakkajima}}
as in \cite[III, Corollary 3.4]{Bloch-DRW} 
(see  also \cite[II, Corollary 3.5]{IlDRW}) we obtain the isomorphism
\[H^{n}_{c, {\rm rig}}(U)_{[j,j+1[}\cong H^i(X, W\Omega^{n-j}_{(X,-D)})\otimes_W K,\]
where the left hand side denotes the part of the compactly supported rigid cohomology of $U$ on which the Frobenius acts with slope $\lambda$ with
$j\le \lambda <j+1$. Note that the above isomorphism also appears in the proof of \cite[Theorem 5.9]{Nakkajima} (there for $D$ reduced). 
In particular $H^i(X, W\sO_{(X,-D)})\otimes_W K$ is equal to the compactly supported Witt vector cohomology of $U$ defined in 
\cite{BBE}.
\item\label{para:conseq-coh6} Assume $U$ is affine (e.g. $D$ is ample) then, the same argument as in \cite[Corollary 1.2]{BBE} yields
\eq{para:conseq-coh6.1}{H^{i}(X, W\Omega^{j}_{(X,-D)})\otimes_W K=0,\quad \text{for all } i+j<\dim X.}
As rigid cohomology of an affine $k$-scheme vanishes above the dimension a similar argument using the slope decomposition of $H^*(X, W\Omega^\bullet_X(\log D))$
yields the vanishing
\eq{para:conseq-coh6.2}{H^{i}(X, W\Omega^{j}_{(X,D)})\otimes_W K=0,\quad \text{for all } i+j>\dim X,}
these two vanishing results can be viewed as a kind of Kodaira-Akizuki-Nakano vanishing. See also  \cite[Theorem 1.1]{Tanaka}, where a different kind
of Kodaira vanishing for Witt vector cohomology is proven. The above two  vanishing results depend only on $D_{\rm red}$ and 
can be found a least implicitly in the literature and are probably not so surprising. 
What might be more intriguing is that by \eqref{para:conseq-coh6.1} and \ref{para:conseq-coh3}
we know that $H^{i}(X, W\Omega^{j}_{(X,-D)})$, $i+j<\dim X$, is annihilated by some fixed $p$-power, say $p^M$. 
Hence by  Theorem \ref{thm:limit} 
\[p^M\cdot \colim_{\ul{p}}\, H^i(X, W_\bullet\Omega^j_{(X,D)})=0, \quad \text{for all } i+j>\dim X.\]
Note that $D$ does not need to be  reduced here and that $M$ might depend on $D$.
\end{enumerate}
\end{para}

\begin{rmk}\label{rmk:not-coh}
\begin{enumerate}[label=(\arabic*)]
\item\label{rmk:not-coh1}
We remark that the complex $\sR_n\otimes_{\sR}^LR\Gamma(X, W\Omega^*_{(X,-D)})$ will for non-reduced $D$ in general not be
isomorphic to $R\Gamma(X, W_n\Omega^*_{(X,-D)})$. 
For example take $X$ smooth projective with a smooth effective and ample divisor $H$ such that $H^i(X,W\sO_X)$ does not vanish for some $i>0$.
Choose $r\ge 1$ such that $H^i(X, \Omega^j_{(X,-rH)})=0$ for all $i<\dim X=N$ and for all $j$. (This is always possible by Serre vanishing and duality.)
We claim that in this situation 
$$\sR_n\otimes_{\sR}^L R\Gamma(X, W\Omega^*_{(X,-rH)})\neq R\Gamma(X, W_n\Omega^*_{(X,-rH)})
\text{\quad for all $n$.}$$
Indeed, else it would follow from Ekedahl's Nakayama Lemma \cite[I, Proposition 1.1, ii)]{EkedahlII}
that we have $H^i(X, W_n\Omega^j_{(X,-rH)})=0$, for all $i<N$ and all $j$, $n$. 
By Corollary \ref{cor:duality-mod} we get  $H^i(X, W_n\Omega^j_{(X,rH)})=0$ for all $i>0$ and all $j$, $n$.
In particular 
\[H^i(X, W\sO_X)=H^i(X, R\lim_n W_n\sO_{(X,rH)})=\lim_n H^i(X, W_n\sO_{(X,rH)})=0, \quad \text{for all } i,\]
where the first equality holds by Proposition \ref{prop:modulus-limit}, and the second follows from the fact that the 
$H^i(X, W_n\sO_{(X,rH)})$ are $W$-modules of finite length, for all $i$.
This contradicts our choice of $X$. 
\item\label{rmk:not-coh2}
Item \ref{rmk:not-coh1} in particular implies that in general 
$$\sR_n\otimes^L_{\sR} W\Omega^*_{(X,-D)}\neq W_n\Omega^*_{(X,-D)}.$$
That Nakkajima's method \cite[Theorem 6.24]{Nakkajima} cannot be applied is realted to 
Remark \ref{rmk:Hyodo}. We analyze the tor-terms in case $D_{\red}$ is smooth:
define the quotient $Q$ by the exact sequence
\[0\to W\Omega^*_{(X,-D)}\to W\Omega^*_{(X, -D_{\rm red})}\to Q\to 0.\]
Definition \ref{defn:HW-zeros} together with the Snake Lemma give the exact sequence
\[0\to Q\to W\Omega^*_D\to W\Omega^*_{D_{\rm red}}\to 0.\]
We have 
\[\sR_n\otimes_{\sR}^L W\Omega^*_{D_{\rm red}}= W_n\Omega^*_{D_{\rm red}} \quad \text{and} \quad 
\sR_n\otimes_{\sR}^L W\Omega^*_{(X,-D_{\rm red})}= W_n\Omega^*_{(X,-D_{\rm red})},\]
where the first equality holds by \cite[II, Th\'eor\`eme (1.2)]{IllRay} and the second equality by \cite[Theorem 6.24]{Nakkajima}.
Moreover 
\[\Tor_0^{\sR}(\sR_n, W\Omega^*_D)= W_n\Omega^*_D,\]
by, e.g., \cite[Lemma 3.2.4]{HeMa03}.
Altogether we get an exact sequence
\[0\to \Tor_1^{\sR}(\sR_n, W\Omega^*_D)\to \Tor_0^{\sR}(\sR_n, W\Omega^*_{(X,-D)})\to W_n\Omega^*_{(X,-D_{\red})}\to W_n\Omega^*_D
\to W_n\Omega^*_{D_{\red}}\to 0\]
and equalities 
\[\Tor_1^{\sR}(\sR_n, W\Omega^*_{(X,-D)})=\Tor_2^{\sR}(\sR_n, W\Omega^*_D), \quad 
\Tor_{j}^{\sR}(\sR_n, W\Omega^*_{(X,-D)})=0,\quad j\neq 0,1.\]
Here the vanishing follows from \cite[I, Corollaire (3.3)]{IllRay}, which can be also used to express
$\Tor_j^{\sR}(\sR_n, W\Omega^*_D)$, $j=1$, $2$,  more explicitly.  
\end{enumerate}
\end{rmk}

\section{Milne-Kato duality with modulus}
In this section we generalize some of the duality results from Milne \cite{Milne76} and 
\cite{Kato-dualityI}, \cite{Kato-dualityII} to the modulus setup. The main result is  Theorem \ref{thm:Zpnqfrp-duality}.

We continue to assume that $X$ is smooth of pure dimension $N$ and $D$ is an effective Cartier divisor on $X$ such that 
$D_{\red}$ is an SNCD. We denote by $j: U=X\setminus D\inj X$ the open immersion of the complement.

\begin{para}\label{para:et-poles-zeros}
We note that the Nisnevich sheaves $W_n\Omega^q_{(X,\pm D)}$ defined by
$V\mapsto \Gamma(V, W_n\Omega^q_{(V, \pm D_{|V})})$
are in fact  \'etale sheaves of coherent $W_n\sO_X$-modules.
Indeed, if $u:V\to X$ is \'etale, then so is $u_n: W_nV\to W_nX$ and it follows from 
the \'etale base change for the de Rham-Witt complex (see \cite[I, Proposition 1.14]{IlDRW})
and the Definition \ref{defn:HW-zeros} (in case of $-D$), and the Propositions \ref{prop:BZnn} and \ref{prop:HWM-ex} (in case of $+D$) that we have
\[\Gamma(V, W_n\Omega^q_{(V, \pm D_{|V})})=\Gamma(V, u_n^*W_n\Omega^q_{(X,\pm D)}).\]
In particular we have 
\eq{para:et-poles-zeros1}{R\e_*W_n\Omega^q_{(X,\pm D)}=W_n\Omega^q_{(X,\pm D)},}
where $\e:X_{\et}\to X_{\Nis}$ is the change-of-sites morphism.
\end{para}

\begin{para}\label{para_Witt-Cartier}
Recall that the Cartier operator on the de Rham-Witt complex is the morphism on $X_{\et}$
\eq{Cartier}{C: F(W_{n+1}\Omega^q_X)\to W_n\Omega^q_X, \quad F(a)\mapsto R(a).}
The so called inverse Cartier operator is given by
\eq{invCartier}{C^{-1}: W_n\Omega^q_X\to W_n\Omega^q_X/dV^{n-1}\Omega^{q-1}_X, \quad a\mapsto F(\tilde{a}),\quad \text{where } R(\tilde{a})=a.}
Note that these are morphisms of $W_n\sO_X$-modules when 
$F(W_{n+1}\Omega^q_X)$ is viewed as a submodule of $F_*W_n\Omega^q_X$ and $W_n\Omega^q_X/dV^{n-1}\Omega^{q-1}_X$ as a quotient-module of 
$F_*W_n\Omega^q_X$. Furthermore, we remark that  $W_n\Omega^N=F(W_{n+1}\Omega^N_X)$ and that the Cartier operator
coincides in this case with the  composition
\ml{topCartier}{C: F_{X*}W_n\Omega^N_X\cong F_{X*}\pi_n^!W_n(k)[-N]\cong F_{X*}\pi_n^! F_k^!W_n(k)[-N]\\
\cong F_{X*}F_X^!\pi_n^!W_n(k)[-N] \xr{\tr_{F_X}} \pi_n^!W_n(k)[-N]=W_n\Omega^N_X,}
where $F_X: W_n(X)\to W_n(X)$ denotes the morphism induced by the absolute Frobenius (similarly for $F_k: W_n(\Spec k)\to W_n(\Spec k)$),
${\tr_{F_X}}: F_{X*}F_X^!\to \id$ denotes the counit of adjunction, and the other notation is as in \ref{para:W-duality},
see \cite[II, Lemma 2.1]{Ekedahl}\footnote{In {\em loc. cit.} the statement is for $C^n$, the $n$-fold iterate of the Cartier operator, 
but the same proof works here.}.
\end{para}

\begin{lem}\label{lem:CC-}
Set
\[(FW_{n+1}\Omega^q)_{(X,D)}:= j_*F(W_{n+1}\Omega^{q}_U)\cap W_{n}\Omega^q_{(X,D)},\]
and 
\[ (W_n\Omega^q/dV^{n-1})_{(X,-D)}:= \frac{W_n\Omega^q_{(X,-D)}}{dV^{n-1}\Omega^{q-1}_X\cap W_n\Omega^q_{(X,-D)}}.\]
These extend naturally to coherent $W_n\sO_X$-modules on $X_\et$.
In particular
\[R\e_*(FW_{n+1}\Omega^q)_{(X,D)}=(FW_{n+1}\Omega^q)_{(X,D)},\]
and
\[R\e_*(W_n\Omega^q/dV^{n-1})_{(X,-D)}=(W_n\Omega^q/dV^{n-1})_{(X,-D)},\]
where $\e:X_{\et}\to X_{\Nis}$ is the change-of-sites morphism.
Furthermore, on $X_\et$ the Cartier operator on $j_*F(W_{n+1}\Omega^q_U)$ restricts to 
\eq{lem:CC-1}{C: (FW_{n+1}\Omega^q)_{(X,D)}\to W_n\Omega^q_{(X,D)}}
and  the inverse Cartier operator $C^{-1}$ on $W_n\Omega^q_X$ restricts to 
\eq{lem:CC-2}{C^{-1}: W_n\Omega^q_{(X,-D)}\to (W_n\Omega^q/dV^{n-1})_{(X,-D)}.}
\end{lem}
\begin{proof}
Since $FW_{n+1}\Omega^q_{U}=\Ker(F^{n-1}d:W_n\Omega^q_{U}\ra \Omega^q_U)$, by \cite[I, (3.21.1.1)]{IlDRW}, the sheaf
\eq{lem:CC-3}{(FW_{n+1}\Omega^q)_{(X,D)}=\Ker(F^{n-1}d: F_*W_n\Omega^q_{(X,D)}\ra F^n_*\Omega^q_{(X,D)})}
is a kernel between coherent (\'etale) sheaves and hence is coherent. Similarly, the sheaf
$(W_n\Omega^q/dV^{n-1}\Omega^{q-1})_{(X,-D)}$ is a quotient of the coherent \'etale sheaf $F_*W_n\Omega^q_{(X,-D)}$ and hence is coherent.
The existence of \eqref{lem:CC-2} follows from the surjectivity of $R$ onto $W_n\Omega^q_{(X,-D)}$, see Lemma \ref{lem:R-surj-zeros}.
For the pole side, note that $C$ defines a morphism $C: F(W_{n+1}\Omega^q)\to W_n\Omega^q$ of Nisnevich sheaves with transfers 
and hence induces a well-defined map
$C:\uomega^{\CI}F(W_{n+1}\Omega^q)\to \uomega^{\CI}W_n\Omega^q$. 
For $q\ge 1$ we have $\uomega^{\CI}F(W_{n+1}\Omega^q)_{(X,D)}=(FW_{n+1}\Omega^q)_{(X,D)}$ by definition (see also \ref{para:RSC})
and Theorem \ref{thm:HW-modulus}, which yields \eqref{lem:CC-1}. 
By the definition in \eqref{para:HW-not2}, the case for $q=0$ is reduced to show the following:
Let $L$ be a henselian dvf and let $a\in W_{n+1}(L)$ such that $F(a)\in \fil^p_rW_n(L)$, for $r\ge 2$, 
then $R(a)\in \fil^p_rW_n(L)$ (see \ref{defn:fil} and \ref{defn:p-sat} for the notation).
By definition we find $b_s\in \fil_{rp^s}W_n(L)$, $s=0,\ldots, n-1$, such that
\[F(a)=b_0+\sum_{s\ge 1} p^s b_s= b_0+\sum_{s\ge 1} F(p^{s-1}V(b_s)).\]
Hence there exists an $a_0\in W_{n+1}(L)$ with $F(a_0)=b_0\in \fil_r W_n(L)\subset \fillog_r W_n(L)$.
It follows directly from the definition that this implies $a_0\in \fillog_rW_{n+1}(L)$ and hence
$R(a_0)\in \fillog_{\lfloor r/p\rfloor}W_n(L)\subset \fillog_{r-1}W_n(L)\subset \fil_rW_n(L)$, see \eqref{para:fillog1.5}.
Moreover, for $s\ge 1$, we have  $V(b_s)\in \fil_{rp^s}W_{n+1}(L)$ and hence $R(V(b_s))\in \fil_{rp^{s-1}}W_n(L)$.
Thus 
\[R(a)=R(a_0)+\sum_{s\ge 1} p^{s-1} R(V(b_s))\in \fil^p_rW_n(L).\] 
This completes the proof of \eqref{lem:CC-1}.
\end{proof}

\begin{defn}\label{defn:Zpnq}
We define  the following two complexes of abelian sheaves on $X_{\et}$
\[\Z/p^n(q)_{(X, D)}:= \left((FW_{n+1}\Omega^q)_{(X,D)}\xr{1-C} W_n\Omega^q_{(X,D)}\right)[-q]\]
and 
\[\Z/p^n(q)_{(X,-D)}:=\left(W_n\Omega^q_{(X,-D)}\xr{C^{-1}-1} (W_n\Omega^q/dV^{n-1})_{(X,-D)}\right)[-q].\]
Both complexes sit in degree $[q, q+1]$.
\end{defn}

\begin{para}\label{para:Zpnq}
We make some comments on the complexes defined above:
\begin{enumerate}[label=(\arabic*)]
\item\label{para:Zpnq1} If $D=\emptyset$, then the two complexes are isomorphic to $W_n\Omega^q_{X,\log}[-q]$ in 
the derived category of abelian sheaves on $X_\et$ (see, e.g., \cite[Lemma 4.1.5]{Kato-dualityI}), 
which is isomorphic to the motivic complex $\Z/p^n(q)_X$, by \cite[Theorem 8.3]{GL}.
\item\label{para:Zpnq2} Let $\e:X_\et \to X_{\Nis}$ be the change-of-sites map then by Lemma \ref{lem:CC-}
\[R\e_*\Z/p^n(q)_{(X,D)}= \uomega^{\CI}\left((FW_{n+1}\Omega^q)\xr{1-C} W_n\Omega^q\right)_{(X,D)}[-q], \quad q\ge 1,\]
for $q=0$ the complex is still a two-term complex of cube invariant Nisnevich sheaves with transfers.
\item \label{para:Zpnq3}As the open immersion $j: U\inj X$ is affine we have in $D(X_\et)$ 
\[\colim_r\, \Z/p^n(q)_{(X,rD)}= Rj_*\Z/p^n(q)_U.\]
\item \label{para:Zpnq4} Complexes similar to $\Z/p^n(q)_{(X,-D)}$ were already considered in \cite{JSZ}, \cite{Morrow}, and \cite{Gupta-Krishna}.
More precisely, in \cite[(3.1.3)]{JSZ} the analog complex with $W_n\Omega^q_{(X,-D)}$ replaced by
$W_n\Omega^q_{X}(\log D)\otimes_{W_n\sO_X}W_n\sO_X(-D)$ is used, and in \cite[below (7.12)]{Gupta-Krishna} (see also \cite[2.4]{Morrow}) the sheaf 
$\Ker(W_n\Omega^q_{X}\to W_n\Omega^q_D)$ is considered instead. The pro-systems over $(X,-rD)$, $r\ge 1$,  
defined by the respective complexes are isomorphic
but for fixed $D$ the complexes differ. We do not know whether they are quasi-isomorphic, but we expect them not to be.
\item \label{para:Zpnq5} For $q=0$ a complex similar to $\Z/p^n(0)_{(X,D)}$ is considered in \cite[(3-2)]{KeS-Lefschetz}, 
where the filtration $\fillog_rW_n(L)$ was used  instead of $\fil^p_rW_n(L)$ which is used here, 
also there the sequence with $F=C^{-1}$ is considered instead of the one with $C$ considered here.
\end{enumerate}
\end{para}

\begin{lem}\label{lem:Zpnq-}
Set
\[W_n\Omega^q_{(X,-D),\log}:=\sH^q(\Z/p^n(q)_{(X,-D)}).\]
Then the natural map
\[W_n\Omega^q_{(X,-D),\log}[-q]\xr{\simeq} \Z/p^n(q)_{(X,-D)}\]
is an isomorphism in the derived category $D(X_\et)$.
\end{lem}
\begin{proof}
We have to show that $C^{-1}-1: W_n\Omega^q_{(X,-D)}\to (W_n\Omega^q/dV^{n-1})_{(X,-D)}$ is surjective in the \'etale topology.
For $n=1$ this holds by the  following claim, which is a version of the surjectivity statement in \cite[Theorem 1.2.1]{JSZ}:
\begin{claim}\label{cl:Cinv-1surjn=1}
Let $A$ be an SNCD, and $B$ be an effective Cartier divisor with $(A+B)_\red$  an SNCD. 
Then
$$\Omega^q_X(\log A)(-B)\xra{C^{-1}-1}
\Omega^q_X(\log A)(-B)/(d\Omega^{q-1}_X\cap \Omega^q_X(\log A)(-B))$$
is surjective on $X_\et$.
\end{claim}
Indeed, let $R$ be the strict henselization of a local ring at a closed point of $X$ and choose a system $t_1,\dots , t_d$ 
of regular parameters of $R$ such that 
$$A=\Div(t_1\dots t_e), \quad
B=\Div(t_1^{r_1}\dots t_e^{r_e}t_{e+1}^{r_{e+1}}\dots t_f^{r_{f}})$$
where $e,f$ are integers such that $0\le e\le f\le d$, and $r_1,\dots r_f$ are non-negative integers.
Let  $I^q$ be the set of strictly increasing functions $\{1,\dots, q\}\ra \{1,\dots, d\}.$
Let $s\in I^q$ and let $\theta$ be the biggest integer in the image of $s$ such that $s(\theta)\le e$.
Denote
$$\omega_s=\dlog t_{s(1)}\dots \dlog t_{s(\theta)} \dlog(1+t_{s(\theta+1)})\dots \dlog (1+t_{s(q)}). $$
Then the $\omega_s$, $s\in I^q$, form an $R$-basis of $\Omega^q_R(\log A)$. Moreover,
$\Omega^q_R(\log A)(-B)=
(\pi)\cdot \Omega^q_{R}(\log A)$ with $\pi =t_1^{r_1}\dots t_e^{r_e}t_{e+1}^{r_{e+1}}\dots t_f^{r_{f}}$.
Since $\omega_s$ is invariant under $C^{-1}$, it suffices to show that for any
$a\in \pi R$, there exists an element $b\in \pi R$ such that $a=b^p-b$.
Such $b$ exists by \cite[Claim 1.2.2]{JSZ}. This shows the claim.

Now we do induction on $n$ and assume $n\ge 2$. As $R$ is surjective by Lemma \ref{lem:R-surj-zeros}, it suffices to show that
$\Ker(R:W_n\Omega^q_{(X,-D)}\ra W_{n-1}\Omega^q_{(X,-D)})$ lies in  the image of $C^{-1}-1$. 
Every element of this kernel is of the form $V^{n-1}(\alpha)+dV^{n-1}(\beta)$ with $\alpha\in \Omega^q_X$ and $\beta\in \Omega^{q-1}_X$.
As $FV^{n-1}(\alpha)=0$, \eqref{prop-ex-seq-zeros00} yields $dV^{n-1}(\beta)\in W_n\Omega^q_{(X,-D)}$.
Hence $V^{n-1}(\alpha)\in W_n\Omega^q_{(X,-D)}$ as well and  it suffices to show $V^{n-1}(\alpha)$ lies in the image of $C^{-1}-1$.
Combining \Cref{lem:VdV-zeros} with \Cref{lem:pHWM-zeros}, we see that $V^{n-1}(\alpha)$ can be written as a sum 
\[V^{n-1}(\alpha)=\ul{p}^{n-1}(\alpha_0)+\sum_{j=1}^{n-1} V^j(\alpha_j)\]
with $\alpha_0\in \Omega^{q}_{n-1}(-D',-pD_n)$, where $D=D'+p^n D_n$ is a \pdd, and $\alpha_j\in W_{n-j}\Omega^q_{(X,-D)}$, for $j\ge 1$. 
By induction $\alpha_j$ lies in the image of $C^{-1}-1$, for $j\ge 1$, and $\alpha_0$ lies in the image  by Claim \ref{cl:Cinv-1surjn=1},
hence so does $V^{n-1}(\alpha)$.
\end{proof}

\begin{rmk}\label{rmk:Zpnq+}
\begin{enumerate}
\item In \cite[Definition 1.1.1]{JSZ} a sheaf $W_n\Omega^q_{(X|D),\log}$ is defined and 
in \cite[(5.3)]{Gupta-Krishna} a sheaf $W_n\Omega^q_{(X,D),\log}$ is defined. These are related to 
the complexes mentioned in \ref{para:Zpnq} \ref{para:Zpnq4} in the same way as $W_n\Omega^q_{(X,-D)}$ is defined in
Lemma \ref{lem:Zpnq-}. We don't know what is the precise relation between these three sheaves, but the pro-systems for  $(X,rD)$, $r\ge 1$,
are isomorpic.
 \item A similar argument as for Lemma \ref{lem:Zpnq-}  shows that in $D(X_{\et})$ we have an isomorphism 
\[\sH^q(\Z/p^n(q)_{(X,D_{\rm red})})[-q]\cong \Z/p^n(q)_{(X, D_{\red})}.\]
(This would correspond to the case $B=0$ in Claim \ref{cl:Cinv-1surjn=1}.)
However, even in the \'etale toplogy the sheaf $\sH^{q+1}(\Z/p^n(q)_{(X,D)})$ does in general not vanish if $D$ is not reduced.
The main problem when trying to adapt the proof of Claim \ref{cl:Cinv-1surjn=1} is that for $a\in \frac{1}{\pi} \cdot A$ 
the Artin-Schreier covering defined by $b^p-b=a$ will (wildly) ramify along the vanishing locus of $\pi$ and is  in particular not \'etale over $A$. 
Also it is direct to show that we have 
\eq{rmk:Zpnq+1}{\sH^q(\Z/p^n(q)_{(X,D)})=j_*W_n\Omega^q_U,}
which only depends on the support of $D$. Thus all the information about the multiplicities of the components of $D$ is stored in the \'etale sheaf
$\sH^{q+1}(\Z/p^n(q)_{(X,D)})$, or the  Nisnevich sheaf  $R^{q+1}\e_*(\Z/p^n(q)_{(X,D)})$.
The study of this sheaf was initiated by Kato in \cite{Kato-Swan} as becomes apparent from  the following Lemma.

\end{enumerate}
\end{rmk}

\begin{lem}\label{lem:relKato}
Let $\eta\in X$ be a generic point of $D$ and set $L=\Frac(\sO_{X,\eta}^h)$. Denote by $r\ge 1$ the multiplicity of $D$ at $\eta$.
Set $H^{q+1}(\Spec L_\et, \Z/p^n(q)):=H^{q+1}_n(L)$. Then
\[\Im(\fil_rW_n\Omega^q_L\to H^{q+1}_n(L))= R^{q+1}\e_*(\Z/p^n(q)_{(X,D)})_{\eta}^h,\]
see Definition \ref{defn:fil} for notation.
In particular, if $\fil^{\rm Kato}_r H^{q+1}_n(L)$ denotes the filtration defined in \cite[Definition 2.1]{Kato-Swan}, then 
\[\fil^{\rm Kato}_{r-1} H^{q+1}_n(L)\subset R^{q+1}\e_*(\Z/p^n(q)_{(X,D)})_{\eta}^h\subset \fil^{\rm Kato}_r H^{q+1}_n(L).\]
\end{lem}
\begin{proof}
By definition 
\eq{lem:relKato1}{R^{q+1}\e_*(\Z/p^n(q)_{(X,D)})_{\eta}^h= \frac{\Fil^p_rW_n\Omega^q_{L}}{(1-C)(FW_{n+1}\Omega^q_L)\cap \Fil^p_rW_n\Omega^q_L}.}
For $a\in \Fil_{p^sr}W_n\Omega^q_L$ and $s\ge 1$ we have 
\[p^sa= F^sV^s(a)\equiv V^s(R^s(a)) \quad \text{mod } (1-C)(FW_{n+1}\Omega^q_L).\]
It follows from the definition 
(see \eqref{para:fillog1.5})
that we have  $R^s(\Fil_{p^sr}W_n\Omega^q_L)\subset \Fil_r W_{n-s}\Omega^q_L$, hence
$\Fil_rW_n\Omega^q_{L}$ surjects onto the quotient \eqref{lem:relKato1}. Moreover for $b\in \fil_rW_n\Omega^{q-1}_L$ we have 
\[db=FdV(b)\equiv dVR(b)\quad \text{mod } (1-C)(FW_{n+1}\Omega^q_L).\]
Iterating this we see that $db$ vanishes in the quotient, which yields the first statement. 
By \cite[Theorem (3.2)]{Kato-Swan} we have a surjection
\[\fillog_rW_n\Omega^q_L\surj \fil^{\rm Kato}_r H^{q+1}_n(L),\]
for all $r\ge 0$, which gives the second statement.
\end{proof}

\begin{prop}\label{prop:dualityF-mod-dVn}
There is an isomorphism in $D(W_n\sO_{X_\et})$
\[(FW_{n+1}\Omega^q)_{(X,D)}\xr{\simeq} \sD_{X,n}\left((W_n\Omega^{N-q}/dV^{n-1})_{(X,-D)}\right),\]
induced by $F(\alpha)\mapsto (\beta\mapsto C(F(\alpha)\beta))$, where
 $\sD_{X,n}=R\cHom_{W_n{\sO_X}}(-, W_n\Omega^N_X)$ denotes the dualizing functor considered  in \ref{para:W-duality} 
(and extended to the \'etale site).
\end{prop}
\begin{proof}
Consider the following two short exact sequences of $W_n\cO_X$-modules
\eq{eq:Fndses}{0\ra 
(FW_{n+1}\Omega^q)_{(X,D)}
\ra
F_{X*}W_n\Omega^q_{(X,D)}
\xra{F^{n-1}d} 
i_*B_n\Omega^{q+1}_{(X,D)}
\ra 0,
}
\eq{eq:dVnses}{
0\ra 
i_*(\Omega/Z)^{N-q-1}_{n,(X,-D)}
\xra{dV^{n-1}}
F_{X*}W_n\Omega^{N-q}_{(X,-D)}
\ra 
(W_n\Omega^{N-q}/dV^{n-1})_{(X,-D)}
\ra 0,}
see  \eqref{lem:CC-3}, \eqref{thm:strHWM1}, and \eqref{para:Omega/BZ3}, where $i:X\inj W_n X$ denotes the closed immersion induced by $R^{n-1}$.
We have an isomorphism
\[a: F_{X_*}W_n\Omega^q_{(X,D)}\xr{\simeq}F_{X*}\sD_{X,n}(W_n\Omega^{N-q}_{(X,-D)})\cong \sD_{X,n}(F_{X*}W_n\Omega^{N-q}_{(X,-D)}),\]
which by Theorem \ref{thm:duality-mod} and \eqref{topCartier} is induced by $\alpha\mapsto (\beta\mapsto C(\alpha \beta))$.
Moreover, there is an isomorphism 
\[b: i_*B_n\Omega^{q+1}_{(X,D)}\xr{\simeq} i_*D_{X,1}((\Omega/Z)^{N-q-1}_{n, (X,-D)})\cong D_{X,n}(i_*(\Omega/Z)^{N-q-1}_{n, (X,-D)}),\]
which is a composition of the isomorphism \eqref{thm:duality-mod4} with the duality isomorphism and which is induced by
\[ i_*B_n\Omega^{q+1}_{(X,D)}\to \cHom_{W_n\sO_X}(i_*(\Omega/Z)^{N-q-1}_{n, (X,-D)}, W_n\Omega^N_X), \quad
\alpha\mapsto (\gamma\mapsto \ul{p}^{n-1}C^n(\alpha \gamma)),\]
see \cite[II, (2.2.5)]{Ekedahl} for the fact that the second isomorphism in the defintion of $b$ is induced by 
$\ul{p}^{n-1}:i_*\Omega^N_X\to W_n\Omega^N_X$.
Finally the following diagram  commutes (note that the top row is concentrated in degree zero)
\[\xymatrix{
\sD_{X,n}(F_{X*}W_n\Omega^{N-q}_{(X,-D)})\ar[rr]^{(dV^{n-1})^\vee} & &
\sD_{X,n}(i_*(\Omega/Z)^{N-q-1}_{n, (X,-D)})\\
F_{X_*}W_n\Omega^q_{(X,D)}\ar[u]^{a}\ar[rr]^{F^{n-1}d}& &
i_*B_n\Omega^{q+1}_{(X,D)}\ar[u]^{(-1)^{q-1} b}.
}\]
Indeed, the commutativity follows directly from the relations $\ul{p}^{n-1}C^{n-1}=V^{n-1}$ and $C\circ dV^{n-1}=0$, 
which follow from the definition of the Cartier operator.
Taking everything together we obtain an isomorphism $\eqref{eq:Fndses}\xr{\simeq} \sD_{X,n}(\eqref{eq:dVnses})$,
which yields  the isomorphism of the statement.
\end{proof}

\begin{para}\label{para:frp}
To extend Milne's classical duality result for \'etale cohomology with mod $p^n$ coefficients 
(see \cite[Theorem 5.2]{Milne76} and \cite[Theorem 1.11]{Milne86}) to the modulus setup we will use
Kato's general dualizing formalism from \cite{Kato-dualityI, Kato-dualityII}, 
which we recall in the following.
A morphism of $\F_p$-schemes $T\to X$ is relatively perfect if and only if 
the relative Frobenius $F_{T/X}: T\to T\times_{X, F_X} X$ is an isomorphism.
If $T\to X$ is additionally flat, then so is $W_n(T)\to W_n(X)$ and the following diagrams are cartesian, 
\eq{para:frp1}{
\xymatrix{
W_n(T)\ar[r]^{W_n(F_T)}\ar[d] &
W_n(T)\ar[d]\\
W_n(X)\ar[r]^{W_n(F_X)} &
W_n(X),}
\quad
\xymatrix{
W_{n}(T)\ar[r]\ar[d] &
W_{n+1}(T)\ar[d]\\
W_n(X)\ar[r] &
W_{n+1}(X),
}
}
where the horizontal maps in the diagram on the right hand side are the closed immersions 
induced by the restriction $R$,
see \cite[Lemma 2]{Kato80} and \cite[0, Proposition 1.5.8]{IlDRW}. 

Let $X_{\FRP}$ be the site with underlying category, the category of 
all flat and relatively perfect $X$-schemes, equipped with the \'etale topology.
In particular, the underlying category of $X_{\FRP}$ contains all \'etale $X$-schemes and we get a fully faithful functor 
of the  underlying catgeories
\[u: X_{\et}\to X_{\FRP},\]
which is continuous, cocontinuous, and commutes with fiber products. 
The induced restriction functor on the topoi $\widetilde{X_{\FRP}}\to \widetilde{X_{\et}}$, $S\mapsto u^*S=S_{|X_{\et}}$,
therefore is exact and has an exact left adjoint, see \cite[I, Proposition 5.4, 4) and III, Proposition 2.6]{SGA4-1}. 
Hence for any complex $L\in D(X_{\FRP},\Z/p^n)$ in the derived category of $\Z/p^n$-sheaves on $X_{\FRP}$ we get a canonical isomorphism  
\eq{para:frp1.5}{R\Gamma(X_{\FRP}, L)\cong R\Gamma(X_{\et}, L_{|X_{\et}})}
in the derived category of $\Z/p^n$-modules.

For any $n\ge 1$, the association
$T\mapsto \Gamma(T, W_n\cO_T)$ defines a sheaf on $X_{\FRP}$, we denote it by $W_n\cO_X^{\FRP}$.
Now let $\tau=\Zar$, $\Nis$ or $\et$.
For any $W_n\cO_X$-module $M$ on $X_\tau$, the association 
$$(g:T\ra X)\mapsto \Gamma(T,(W_ng)^*M)$$
defines a presheaf of $W_n\cO_X^{\FRP}$-modules on $X_{\FRP}$, we denote its sheafification by $M^{\FRP}$.
When $M$ is (quasi-)coherent as a $W_n\cO_X$-module on $X_\tau$, we say that $M^{FRP}$ is \textit{(quasi-)coherent} as a $W_n\cO_X^{\FRP}$-module on $X_{\FRP}$.
The functor $M\mapsto M^{\FRP}$ from the category of $W_n\cO_X$-modules on the site $X_\tau$ to the category of $W_n\cO_X^{\FRP}$-modules is exact.
In particular, the association 
$$(g:T\ra X)\mapsto \Gamma(T, (W_ng)^*W_n\Omega^q_{(X,\pm D)})$$
defines a sheaf of $W_n\cO_X^{\FRP}$-modules on $X_{\FRP}$, 
which we denote by $W_n\Omega^{q,\FRP}_{(X,\pm D)}$. 
The maps $F,V,R,\ul p, d, C, C^{-1}$ on $W_\bullet\Omega^*_X$ 
extend naturally to maps on $W_\bullet\Omega^{*,\FRP}_{(X,\pm D)}$; 
this follows  from the cartesian diagrams \eqref{para:frp1} and the flatness of 
$W_n(g): W_n(T)\to W_n(X)$, see also \ref{para:ZpnqFRP} below where we spell 
out how to extend $F$ and $C$.
In particular, when $D=\emptyset$, we write
$W_n\Omega^{\bullet,\FRP}_{X}:=W_n\Omega^{\bullet,\FRP}_{(X,\emptyset)}$. 
Since $(W_ng)^*W_n\Omega^{\bullet}_{X}\simeq W_n\Omega^q_T$ for any 
flat and relatively perfect morphism $g:T\ra X$, 
the sheaf $W_n\Omega^{\bullet,\FRP}_{X}$ agrees with the presheaf given by 
$T\mapsto \Gamma(T,W_n\Omega^q_T)$, cf. the discussion in \cite[\S 4]{Kato-dualityI} or the proof of \cite[Proposition 1.7]{LZ}.

Denote by $D_0(X_{\FRP},\Z/p^n)$ the full triangulated subcategory of $D(X_{\FRP}, \Z/p^n)$ generated by coherent $\cO_X^{\FRP}$-modules 
regarded as complexes of $\Z/p^n$-modules in degree $0$.
Note that $W_n\Omega^{q, \FRP}_{(X,\pm D)}\in D_0(X_{\FRP}, \Z/p^n)$, 
by the Theorems \ref{thm:strHWM} and  \ref{thm:HW-zeros}. Set
$$\mathbb{D}_{X,n}(-):=R\sH om_{D(X_{\FRP},\Z/p^n)}(-, W_n\Omega^{N,\FRP}_{X,\log}),$$
where $W_n\Omega^{N,\FRP}_{X,\log}$ is defined by the short exact sequence
\eq{para:frp2}{0\to W_n\Omega^{N,\FRP}_{X,\log}\to W_n\Omega^{N,\FRP}_{X}\xr{C-1} W_n\Omega^{N,\FRP}_{X}\to 0.}
This functor agrees with Kato's functor defined in \cite[\S 6]{Kato-dualityI} 
 and differs from the one defined in  \cite[\S 3]{Kato-dualityII} by a shift $[N]$.
Kato proves the following results:
\begin{enumerate} 
\item 
If $M\in D^b_c(X_{\Zar},W_n\cO_X)$, then 
the natural map 
$W_n\Omega^{N,\FRP}_X\ra W_n\Omega^{N,\FRP}_{X,\log}[1]$
from \eqref{para:frp2} induces an isomorphism in $D_0(X_{\FRP},\Z/p^n)$
\eq{eq:KatoDualityCoh}{
\cD_{n,X}(M)^{\FRP}
\simeq 
\D_{n,X}(M^{\FRP})[1],
}
where $\cD_{n,X}=R\cHom_{W_n\sO_X}(-, W_n\Omega^N_X)$ as in \ref{para:W-duality}, see \cite[Proposition 6.1]{Kato-dualityI}.
\item 
The functor $\D_{X,n}$ restricts to an involutive endo-functor of $D_0(X_{\FRP},\Z/p^n)$, i.e.,
\eq{eq:KatoDualitySelfDual}{\D_{X,n}\circ \D_{X,n}\simeq \id:
D_0(X_{\FRP},\Z/p^n)\ra D_0(X_{\FRP},\Z/p^n),
}
see \cite[Theorem (0.1) and Proposition (3.4)]{Kato-dualityII}.
\item 
Let $f:X\to Y$ be a proper morphism of relative dimension $r$ between two pure dimensional 
smooth $k$-schemes. There is a natural isomorphism
\eq{eq:KatoDuality}{ 
Rf_*\D_{n,X}( M)
\simeq
\D_{n, Y}( Rf_*M)[-r]
}
in $D_0(Y_{\FRP},\Z/p^n)$,
for any $M\in D_0(X_{\FRP},\Z/p^n)$, see \cite[Theorem (0.2)]{Kato-dualityII}.
\end{enumerate}

The following Lemma is implicitly in \cite{Kato-dualityI} and \cite{Milne76}.

\begin{lem}\label{lem:Kato-duality-frp-et}
Assume $k$ is a finite field and   $X$ is proper over $k$.
Let $M\in D_0(X_{\FRP}, \Z/p^n)$. There is a canonical isomorphism 
\[R\Gamma(X_{\FRP}, \D_{X,n}(M))\cong R\Hom_{D(\Z/p^n)}\left(R\Gamma(X_{\et}, M_{|X_{\et}}), \Z/p^n\right)[-N-1].\]
\end{lem}
\begin{proof}
We start by giving a canonical map from left to right. As $k$ is finite, we can view $X$ also as 
a smooth proper scheme of relative dimension $N$ over $S_0:=\Spec \F_p$.
Note that \eqref{para:frp2} for $S_0$ instead of $X$ yields in view of \eqref{para:frp1.5} a decomposition
\eq{lem:Kato-duality-frp-et1}{R\Gamma(S_{0,\FRP}, W_n\Omega^{0,\FRP}_{S_0,\log})
= W_n(\F_p)[0] \oplus W_n(\F_p)[-1].}
Let $\pi: X\to S_0=\Spec \F_p$ be the structure map.
We define the map from the statement to be the composition
\begin{align}
R\Gamma(X_{\FRP}, \D_{X,n}(M)) & \xr{\simeq} R\Gamma(S_{0, \FRP}, \D_{S_0,n}(R\pi_*M))[-N]
           \label{lem:Kato-duality-frp-et2}\\
                              & \to R\Hom_{D(\Z/p^n)}(R\Gamma(S_{0,\FRP}, R\pi_*M), R\Gamma(S_{0,\FRP}, 
                              W_n\Omega^{0,\FRP}_{S_0,\log}))[-N]\notag\\
                              &\cong R\Hom_{D(\Z/p^n)}(R\Gamma(X_{\et}, M_{|X_{\et}}), W_n(\F_p)[0] 
                              \oplus W_n(\F_p)[-1])[-N]\notag\\
                              &\to R\Hom_{D(\Z/p^n)}(R\Gamma(X_{\et}, M_{|X_{\et}}), \Z/p^n)[-N-1],\notag   
\end{align}
where the first isomorphism is induced from the duality isomorphism \eqref{eq:KatoDuality}, the second map  is the natural transformation
$R\Gamma(T,-)\circ R\cHom(-,-)\to R\Hom(R\Gamma(T,-), R\Gamma(T,-))$, 
the third isomorphism  is \eqref{para:frp1.5} and \eqref{lem:Kato-duality-frp-et1}, 
and the last map is the projection to the second summand plus the identification $W_n(\F_p)=\Z/p^n$.
We note that this morphism defines for $M$ variable an exact functor of triangulated categories 
$D_0(X_{\FRP}, \Z/p^n)\to D(\Z/p^n)$. Thus it suffices to show that it is an isomorphism
for $M=L^{\FRP}$, with $L\in D^b_c(X_{\Zar}, W_n\sO_X)$. For such $L$ we have  isomorphisms
\begin{align}
R\Gamma(X_{\FRP}, \D_{X,n}(L^{\FRP})) & \cong R\Gamma(X_{\FRP}, \cD_{X,n}(L)^{\FRP})[-1]   
\label{lem:Kato-duality-frp-et3}\\
 & \cong R\Gamma(X_{\et}, (\cD_{X,n}(L)^{\FRP})_{|X_{\et}})[-1]\notag\\
 &\cong R\Gamma(X_{\Zar}, \cD_{X,n}(L))[-1]\notag\\
 &\cong R\Hom_{W_n(\F_p)}(R\Gamma(X_{\Zar}, L), W_n(\F_p))[-N-1]\notag\\
 &\cong R\Hom_{W_n(\F_p)}(R\Gamma(X_{\et}, (L^{\FRP})_{|X_{\et}}), \Z/p^n)[-N-1],\notag
\end{align}
where  the first isomorphism is \eqref{eq:KatoDualityCoh}, the second is \eqref{para:frp1.5},
the third isomorphism holds since $\cD_{X,n}(L)\in D^b_c(X_{\Zar},W_n\sO_X)$,
the fourth isomorphism  is Ekedahl-Grothendieck duality and the last isomorphism is clear.
It follows from the construction of the duality isomorphism \eqref{eq:KatoDuality} and 
the isomorphism \eqref{eq:KatoDualityCoh} that the two compositions \eqref{lem:Kato-duality-frp-et2}
and \eqref{lem:Kato-duality-frp-et3} agree for $M=L^{\FRP}$.
Indeed by \cite[(3.2.2), (3.2.3)]{Kato-dualityII}
we have a commutative diagram
\[\xymatrix{
R\pi_*W_n\Omega^{N,\FRP}_X[N]\ar[r]^-{\Tr}\ar[d] &
W_n\Omega^{0,\FRP}_{S_0}\ar[d]\\
R\pi_*W_n\Omega^{N,\FRP}_{X,\log}[N+1]\ar[r]^-{\Tr} & W_n\Omega^{0,\FRP}_{S_0,\log}[1],
}\]
where the horizontal maps are the trace maps defining  the respective duality isomorphism and 
the vertical maps are induced by \eqref{para:frp2}. With this commutativity and the formula $R\pi_*(L^{\FRP})=(R\pi_*(L))^{\FRP}$ 
(see \cite[Lemma 2.6]{Kato-dualityII}) it is straightforward to
check that the above compositions agree. This completes the proof.
\end{proof}

\end{para}
\begin{para}\label{para:ZpnqFRP}
We extend the complexes from Definition \ref{defn:Zpnq} to the site $X_{\FRP}$ by setting
\[\Z/p^n(q)_{(X, D)}^{\FRP}:= \left((FW_{n+1}\Omega^q)_{(X,D)}^{\FRP}\xr{1-C} W_n\Omega^{q,\FRP}_{(X,D)}\right)[-q],\]
and 
\[\Z/p^n(q)_{(X,-D)}^{\FRP}:=\left(
W_n\Omega^{q,\FRP}_{(X,-D)}
\xr{C^{-1}-1} 
(W_n\Omega^q/dV^{n-1})_{(X,-D)}^{\FRP}\right)[-q].\]    
We spell out the extension of $1-C$ to $X_{\FRP}$.
Let $g: T\to X$ be relatively perfect and flat.
Then the Frobenius on $W_{\bullet}\Omega_{(X,\pm D)}^{q,\FRP}$ is induced by the composition
\[W_{n+1}(g)^* W_{n+1}\Omega^q_{(X,\pm D)}\xr{W_{n+1}(g)^*(F)} W_{n+1}(g)^* F_*W_n\Omega^q_{(X,\pm D)}\cong 
F_* W_{n}(g)^*W_n\Omega^q_{(X,\pm D)},\]
which on elements is given by 
\[W_{n+1}\sO_T\otimes_{W_{n+1}\sO_X} W_{n+1}\Omega^q_{(X,\pm D)}\ni a\otimes \alpha \mapsto 
F(a)\otimes F(\alpha)\in W_n\sO_T\otimes W_n\Omega^q_{(X,\pm D)}.\]
Here the isomorphism on the right hand side comes from the cartesian diagrams \eqref{para:frp1}
and the flatness of $W_{n+1}(g)$. Hence 
\[\Gamma(T,(FW_{n+1}\Omega^q)_{(X,D)}^{\FRP})=\Gamma(T, F(W_{n+1}(g)^*W_{n+1}\Omega^q_{(X,D)})).\]
The map $(FW_{n+1}\Omega^q)_{(X,D)}^{\FRP}\xr{1-C} W_n\Omega^{q,\FRP}_{(X,D)}$ is therefore induced by
\[F(a)\otimes F(\alpha)\mapsto F(a)\otimes F(\alpha) - R(a)\otimes R(\alpha).\]
Similarly for the map $C^{-1}-1$ in the definition of $\Z/p^n(q)_{(X,-D)}^{\FRP}$.
\end{para}

\begin{lem}\label{lem:pairing-triangle}
Let $\alpha :A^0\to A^1$ and $\beta: B^0\to B^1$ be two morphisms of sheaves of abelian groups  on some site $X$.
Let $C^\bullet=(C^0\xr{d} C^1)$ be a two term complex of  sheaves of abelian groups on $X$ sitting in degree $[0,1]$.
Assume we are given pairings
\[\langle-,-\rangle^0_{0,0}: A^0\otimes_{\Z} B^0\to C^0, \quad \langle-,-\rangle^1_{1,0}: A^1\otimes_{\Z} B^0\to C^1,\quad
\langle-,-\rangle^1_{0,1}: A^0\otimes_{\Z} B^1\to C^1,\]
such that 
\[\langle\alpha(-),-\rangle^1_{1,0}+\langle-,\beta(-)\rangle^1_{0,1}= d\circ \langle-,-\rangle^0_{0,0}: A^0\otimes_{\Z} B^0\to C^1.\]
Then we obtain the following morphism between distinguished triangles in the homotpoy category of complexes of abelian sheaves on $X$
\footnote{We use the usual sign conventions, see, e.g., \cite[\href{https://stacks.math.columbia.edu/tag/0FNG}{Tag 0FNG}]{stacks-project}.}
\[
\xymatrix{
A^0[0]\ar[r]^{\alpha}\ar[d]^{\pi'} &
A^1[0]\ar[r]\ar[d]^{\pi''}&
\cone(\alpha)\ar[d]^{\pi}\ar[r]&
A^0[1]\ar[d]^{\pi'[1]} \\
\cHom^\bullet(B^1, C^\bullet)[1]\ar[r]^{\beta^*[1]}&
\cHom^\bullet(B^0, C^\bullet)[1]\ar[r]&
\cone(\beta^*[1])\ar[r]&
\cHom^\bullet(B^1, C^\bullet)[2],
}\]
where the non-labled horizontal maps are the natural morphisms in the respective cone sequence
and the vertical maps are defined as follows:
the morphism $\pi$ is in degree -1 given by
\[A^0\to \cHom(B^0,C^0)\oplus \cHom(B^1, C^1), \quad 
x\mapsto \left(-\langle x,-\rangle^0_{0,0}, -\langle x,-\rangle^1_{0,1}\right),\]
and in degree 0  by
\[A^1\to \cHom(B^0, C^1), \quad y\mapsto \langle y,-\rangle^1_{1,0},\]
and $\pi'$ and $\pi''$ are (in degree 0) given by 
\[\pi'(x)=-\langle x,-\rangle^1_{0,1}\quad \text{and} \quad \pi''(y)=\langle y ,-\rangle^1_{1,0}.\]
Furthermore, there is an isomorphism of complexes
\eq{lem:pairing-triangle1}{\cone(\beta^*[1])\xr{\simeq} \cHom^\bullet(\cone(\beta)[-1], C^\bullet)[1],}
which is the identity in the degrees -2 and 0 and which in degree -1 is given by
$\id\oplus (-\id)\in \End(\cHom(B^0, C^0)\oplus \cHom(B^1,C^1))$.
\end{lem}
\begin{proof}
This is all direct to check. We  remark that the square on the left commutes only up to homotopy, 
more precisely if
$h: A^0\to (\cHom(B^0, C^\bullet)[1])^{-1}= \cHom(B^0, C^0)$ is defined by $h(x)=-\langle x , - \rangle^0_{0,0}$, then we have 
$\pi''\circ \alpha- (\beta^*)\circ \pi'=(d[1])\circ h$. 
\end{proof}

\begin{para}\label{para:pairingZpnqFRP}
We define the following pairings on $X_{\FRP}$:
\begin{align*}
\langle - ,- \rangle^0_{0,0}: &(FW_{n+1}\Omega^q)_{(X,D)}^{\FRP}
\otimes_{\Z} 
W_n\Omega^{N-q,\FRP}_{(X,-D)}
\to
W_n\Omega^{N,\FRP}_X,   & &
\alpha\otimes\beta\mapsto \alpha\beta,\\
\langle - ,- \rangle^1_{1,0}:&
W_n\Omega^{q,\FRP}_{(X,D)}
\otimes_{\Z} 
W_n\Omega^{N-q,\FRP}_{(X,-D)}
\to W_n\Omega^{N,\FRP}_X, & &\alpha\otimes \beta\mapsto \alpha\beta,\\
\langle - ,- \rangle^1_{0,1}:&
(FW_{n+1}\Omega^q)_{(X,D)}^{\FRP}
\otimes_{\Z} 
(W_n\Omega^q/dV^{n-1})_{(X,-D)}^{\FRP}
\to
W_n\Omega^{N,\FRP}_X, & &
\alpha\otimes\beta\mapsto C(\alpha\beta).
\end{align*}
These satisfy
\[\langle (C-1)(x) ,y \rangle^1_{1,0}+ \langle x ,(1-C^{-1}(y)) \rangle^1_{0,1}= 
(C-1)\langle x ,y \rangle^0_{0,0},\]
for all local sections $x\in (FW_{n+1}\Omega^q)_{(X,D)}^{\FRP}$ and $y\in W_n\Omega^{N-q,\FRP}_{(X,-D)}$.
Note that
\[\Z/p^n(q)^{\FRP}_{(X,D)}[q]=
\cone\left((FW_{n+1}\Omega^q)^{\FRP}_{(X,D)}\xr{C-1} W_n\Omega^{q,\FRP}_{(X,D)}\right)[-1],\]
\[\Z/p^n(N-q)^{\FRP}_{(X,-D)}[N-q]=
\cone\left( W_n\Omega^{q,\FRP}_{(X,-D)}\xr{1-C^{-1}} (W_n\Omega^q/dV^{n-1})_{(X,-D)}^{\FRP}\right)[-1].\]
Hence we obtain a natural morphism in $D_0(X_{\FRP}, \Z/p^n)$
\begin{align}
\Z/p^n(q)^{\FRP}_{(X,D)}[q] & =\cone(C-1)[-1]
\label{para:pairingZpnqFRP1}\\ 
& \lra R\cHom\left(\cone(1-C^{-1})[-1], (W_n\Omega^{N,\FRP}_{X}\xr{C-1} W_n\Omega^{N,\FRP}_{X})\right)
\notag\\
& \cong \D_{n,X}\left(\Z/p^n(N-q)_{(X,-D)}^{\FRP}[N-q]\right),\notag
\end{align}
where the morphism in the second line is induced by the composition of $\pi[-1]$ 
from Lemma \ref{lem:pairing-triangle} with the isomorphism \eqref{lem:pairing-triangle1}
and the isomorphism in the third line is induced by the exact sequence \eqref{para:frp2}.
\end{para}

In view of Kato's duality isomorphism \eqref{eq:KatoDuality} the following 
Theorem generalizes Milne's duality, see \cite[Theorem 5.2]{Milne76} (for smooth projective surfaces), 
\cite[Theorem 1.11]{Milne86} (for  general smooth proper $k$-schemes), and \cite[Theorem 4.3]{Kato-dualityI}, 
to the case of non-empty modulus.
\begin{thm}\label{thm:Zpnqfrp-duality}
\label{thm:FRPDuality}
The morphism \eqref{para:pairingZpnqFRP1} induces isomorphisms in $D_0(X_{\FRP},\Z/p^n)$
$$\Z/p^n(q)^{\FRP}_{(X,D)}[q]
\xr{\simeq}
\D_{n,X}\left(\Z/p^n(N-q)_{(X,-D)}^{\FRP}[N-q]\right),$$
$$\Z/p^n(N-q)_{(X,-D)}^{\FRP}[N-q]
\xr{\simeq}
\D_{n,X}\left(\Z/p^n(q)^{\FRP}_{(X,D)}[q]\right).$$
\end{thm}
\begin{proof}
By \eqref{eq:KatoDualitySelfDual} it suffices to show that \eqref{para:pairingZpnqFRP1} is an isomorphim.
By Lemma \ref{lem:pairing-triangle} we have a morphism of distinguished triangles 
$$\resizebox{\displaywidth}{!}{
\xymatrix{
W_n\Omega^{q,\FRP}_{(X,D)}[-1]
\ar[d]^{\pi''[-1]}
\ar[r]&
\Z/p^n(q)_{(X, D)}^{\FRP}[q]
\ar[d]
\ar[r]&
(FW_{n+1}\Omega^q)_{(X,D)}^{\FRP}
\ar[d]^{\pi'}
\ar[r]&
\\
\D(W_n\Omega^{N-q,\FRP}_{(X,-D)})
\ar[r]&
\D\left(\Z/p^n(N-q)_{(X,-D)}^{\FRP}[N-q]\right)
\ar[r]&
\D((W_n\Omega^{N-q}/dV^{n-1})_{(X,-D)}^{\FRP})[1] 
\ar[r]&,
}
}$$
where we write $\D$ instead of $\D_{n,X}$.
By definition the map $\pi''[-1]$ factors as
\begin{align*} 
W_n\Omega^{q,\FRP}_{(X,D)}[-1]
&
\xra[\simeq]{\text{Thm. }\ref{thm:duality-mod}}
\cD_{n,X}(W_n\Omega^{N-q}_{(X,-D)})^{\FRP}[-1]
\\
&
\xra[\simeq]{\eqref{eq:KatoDualityCoh}}
\D_{n,X}(W_n\Omega^{N-q,\FRP}_{(X,-D)})
\end{align*}
and hence is an isomorphism.
Similarly, the morphism $\pi'$ factors as  
\begin{align*} 
(FW_{n+1}\Omega^q)_{(X,D)}^{\FRP} 
&
\xra[\simeq]{\text{Prop. \ref{prop:dualityF-mod-dVn}}}
\cD_{n,X}((W_n\Omega^{N-q}/dV^{n-1})_{(X,-D)})^{\FRP}
\\
&
\xra[\simeq]{\eqref{eq:KatoDualityCoh}}
\D_{n,X}((W_n\Omega^{N-q}/dV^{n-1})^{\FRP}_{(X,-D)})[1]
\end{align*}
and hence is an isomorphism as well. 
This completes the proof of the theorem.
\end{proof}

\begin{cor}\label{cor:Zpnq-duality-finite}
Assume additionally that $k$ is a finite field and $X$ is proper over $k$.
Then the morphism \eqref{para:pairingZpnqFRP1} induces an isomorphism in $D(\Z/p^n)$
\[R\Gamma(X_{\et}, \Z/p^n(q)_{(X,D)})
\cong R\Hom_{\Z/p^n}\left(R\Gamma(X_{\et}, \Z/p^n(N-q)_{(X,-D)}), \Z/p^n\right)[-2N-1].\]
In particular we obtain isomorphisms of finite groups for all $i$
\[H^{i+q}(X_{\et}, \Z/p^n(q)_{(X,D)})\cong 
\Hom_{\Z/p^n}(H^{2N-i-q+1}(X_{\et}, \Z/p^n(N-q)_{(X,-D)}), \Z/p^n).\]
\end{cor}
\begin{proof}
As $(\Z/p^n(q)^{\FRP}_{(X,\pm D)})_{|X_{\et}}=\Z/p^n(q)_{(X,\pm D)}$, this follows from 
Theorem \ref{thm:Zpnqfrp-duality}, Lemma \ref{lem:Kato-duality-frp-et}, and \eqref{para:frp1.5}.
\end{proof}

\begin{rmk}\label{rmk:fin-duality}
By Lemma  \ref{lem:Zpnq-} the isomorphism in Corollary \ref{cor:Zpnq-duality-finite} also induces an isomorphism
\[H^{N-i+1}(X_{\et}, W_n\Omega^{N-q}_{(X,-D),\log})\cong H^{i+q}(X_{\et}, \Z/p^n(q)_{(X,D)})^{\vee},\]
where the $(-)^\vee$ on the right hand side denotes the  Pontryagin dual. 
Taking the limit and using \ref{para:Zpnq}\ref{para:Zpnq3} yields the isomorphism
\[\lim_r H^{N-i+1}(X_{\et}, W_n\Omega^{N-q}_{(X,-rD),\log}) \cong  H^i(U_{\et}, W_n\Omega^q_{U,\log})^\vee,\]
which gives back  \cite[Theorem 2]{JSZ}.
\end{rmk}

\section{Weights zero and one}
In this section we consider the complex $\Z/p^n(q)_{(X,D)}$ for $q=0, 1$ and and apply 
the finite duality from the previous section. 
The main results are Corollaries \ref{cor:duality-pi} and \ref{cor:Zpn0-zeros}, and Theorem \ref{thm:duality-Br-surf}.
The assumptions on $X$ and $D$ are the same as in the previous sections.

\subsection{Weight zero}

\begin{lem}\label{lem:Zpn0-poles}
Assume $X$ is irreducible.
For all $i$ we have natural maps
\eq{para:Zpn0-poles1}{H^i(X_{\et},\Z/p^n)\to H^i(X_{\et}, \Z/p^n(0)_{(X,D)})\to H^i(U_{\et}, \Z/p^n)}
inducing an isomorphism 
\[\colim_r H^i(X_{\et}, \Z/p^n(0)_{(X,r D)})\cong H^i(U_{\et}, \Z/p^n).\]
Moreover, the maps \eqref{para:Zpn0-poles1} are injective for $i=1$ and identify with the inclusions
\[H^0(X_{\Nis},R^1\e_*(\Z/p^n)_{X_{\et}})\to H^0(X_{\Nis}, R^1\e_*\Z/p^n(0)_{(X,D)})\to 
H^0(X_{\Nis}, j_*R^1\e_*(\Z/p^n)_{U_{\et}}),\]
where $\epsilon: X_{\et}\to X_{\Nis}$ is the change of sites map.
\end{lem}
\begin{proof}
The natural maps in \eqref{para:Zpn0-poles1}  are induced from the usual 
Artin-Schreier-Witt sequence on $X_{\et}$ and $U_{\et}$
and the isomorphism  comes from \ref{para:Zpnq}\ref{para:Zpnq3}. 
For the injectivity note that 
\[R^0\epsilon_*\Z/p^n(0)_{(X,D)}=(\Z/p^n)_{X_{\Nis}}.\]
Hence $H^i(X_{\Nis},R^0\epsilon_*\Z/p^n(0)_{X,D})=0= H^i(U_{\Nis}, R^0\epsilon_*\Z/p^n_{|U_\et})$, for all $i\neq 0$,
which yields the identification of \eqref{para:Zpn0-poles1} for $i=1$ with the maps in the second part of the statement.
Using the usual Artin-Schreier-Witt sequence and the definition of $\Z/p^n(0)_{(X,D)}$
(see Lemma \ref{lem:CC-} and Definition \ref{defn:Zpnq}) we see that this latter maps are injective.
\end{proof}

\begin{lem}\label{lem:Zpn0-zeros}
There is distinguished triangle in $D(X_{\et}, \Z/p^n)$
\[\Z/p^n(0)_{(X,-D)}\to (\Z/p^n)_{X_{\et}}\to (\Z/p^n)_{D_{\et}}\xr{+1},\]
where the second map is the natural restriction map.
In particular, the natural map 
\[\Z/p^n(0)_{(X,-D)}\to \Z/p^n(0)_{(X,-D_{\rm red})}\] 
is an isomorphism.
\end{lem}
\begin{proof}
By the Artin-Schreier-Witt sequence the complex of \'etale sheaves
\eq{lem:Zpn0-zeros1}{\Z/p^n(0)_Y:=\left(W_n\sO_Y\xr{F-1} W_n\sO_Y\right)}
sitting in degrees $[0,1]$ is quasi-isomorphic to the constant sheaf $(\Z/p^n)_{Y_\et}$, for any $\F_p$-scheme $Y$. 
Thus the first statement follows from the exact sequence of complexes
\[0\to \Z/p^n(0)_{(X,-D)}\to \Z/p^n(0)_X\to \Z/p^n(0)_D\to 0.\]
For the exactness on the left note that a local section $a\in W_n\sO_X$ maps to zero in $W_n\sO_D$ if and only if 
it maps to zero in all $W_n\sO_{D_i}$, where $D=\sum_i D_i$, with $D_{i,\red}\in \Sm$. 
The last statement follows from the equivalence of topoi  ${\rm Sh}(D_{\et})\cong {\rm Sh}(D_{\red, \et})$,
see \cite[VIII, Th\'eor\`eme 1.1]{SGA4-2}.
\end{proof}

\begin{para}\label{para:pi1-pole-zero}
 Set 
 \[\pi_1^{\rm ab}(X,D)^p:= \lim_n \Hom (H^1(X_{\et}, \Z/p^n(0)_{(X,D)}), \Q/\Z),\]
 where the transition maps in the limit are induced by $\ul{p}: \Z/p^n(0)_{(X,D)}\to \Z/p^{n+1}(0)_{(X,D)}$.
 It follows from Lemma \ref{lem:Zpn0-poles}, Lemma \ref{lem:relKato}, and \cite[Proposition 2.5]{KeS-Lefschetz} that this group
 is equal to the maximal pro-$p$ quotient of the group $\pi_1^{\rm ab}(X,D)$ considered in
 \cite[(2-4)]{KeS-Lefschetz} (see also \cite{KeS}). It classifies  finite \'etale covers
 of $U$ which have an abelian $p$-group as Galois group and whose ramification is bounded by $D$.

Moreover, set 
\[\pi_1^{\rm ab}(X,-D)^p:=\lim_p \Hom(H^1(X_{\et}, \Z/p^n(0)_{(X,-D)}), \Q/\Z).\]
By Lemma \ref{lem:Zpn0-zeros}
\[\pi_1^{\rm ab}(X,-D)^p=\pi_1^{\rm ab}(X,-D_{\red})^p\]
and if $X$ is proper and $|D|$ connected  we have an exact sequence  of profinite groups
\eq{para:pi1-pole-zero1}{\pi_1^{\rm ab}(D)^p\to \pi_1^{\rm ab}(X)^p\to \pi_1^{\rm ab}(X,-D)^p\to 0.}
Note that $\pi_1^{\rm ab}(X, -D)$ classifies finite \'etale  covers of $X$, which have an abelian $p$-group 
as Galois group and split completely when restricted to $D$. 
\end{para}

Now Corollary \ref{cor:Zpnq-duality-finite} directly gives:

\begin{cor}\label{cor:duality-pi}
Assume $k$ is finite and $X$ is proper.
Then 
\[\pi_1^{\rm ab}(X,D)^p= \lim_n H^{2N}(X_{\et}, \Z/p^n(N)_{(X,-D)})=H^N(X_{\et}, W\Omega^N_{(X,-D), \log}),\]
where $W\Omega^N_{(X,-D), \log}=\lim_n W_n\Omega^N_{(X,-D), \log}$, and
\[\pi_1^{\rm ab}(X, -D)^p= \lim_n H^{2N}(X_{\et}, \Z/p^n(N)_{(X,D)}).\]
\end{cor}

\begin{rmk}
Taking the limit over $\{rD\}_r$ in the first equality above gives back the isomrophism
from \cite[below Theorem 2]{JSZ}.
\end{rmk}

\begin{cor}\label{cor:Zpn0-zeros}
Let $k$ be a finite field, let $X$ be proper  and let $H\subset X$ be a smooth ample divisor intersecting $D$ 
transversally. Then the natural map
\[H^i(X_{\et}, \Z/p^n(0)_{(X,-D)})\to H^i(H_{\et}, \Z/p^n(0)_{(H, -D_{|H})})\]
\begin{itemize}
    \item is injective, if $N>i$, and 
    \item is bijective, if $N> i+1$.
\end{itemize}
\end{cor}
\begin{proof}
With the notation from \eqref{lem:Zpn0-zeros1} we define for a closed immersion of $\F_p$-schemes
$Z\subset Y$ the complex $\Z/p^n(0)_{(Y,-Z)}$ of sheaves on $Y_{\et}$ by the exact sequence
\[0\to \Z/p^r(0)_{(Y,-Z)}\to \Z/p^n(0)_Y\to \Z/p^n(0)_Z\to 0.\]
By (the proof of) Lemma \ref{lem:Zpn0-zeros} the complex $\Z/p^n(0)_{(Y,-Z)}$ coincides with the one defined in 
Definition \ref{defn:Zpnq} in case $Y$ is smooth and $Z_{\red}$ is an SNCD. 
Furthermore, up to quasi-isomorphism, these complexes only depend on $Y_{\red}$ and $Z_{\red}$.
For $r\ge 1$ we denote in the following by $D+rH$ also the closed subscheme of $X$ 
defined by the the effective Cartier divisor.
It follows directly from the definition and an application of the Snake Lemma, 
that we have a short exact sequence of complexes on $X_{\et}$
\eq{cor:Zpn0-zeros1}{0\to \Z/p^n(0)_{(X,-D-rH)}\to \Z/p^n(0)_{(X,-D)}\to \Z/p^n(0)_{(D+rH,-D)}\to 0.}
There is a natural morphism of complexes 
\eq{cor:Zpn0-zeros2}{\Z/p^n(0)_{(D+rH,-D)}\to \Z/p^n(0)_{(rH,-D_{|rH})},}
which is an isomorphism. Indeed the map is induced by the natural surjection $\sO_{D+rH}(-D)\to \sO_{rH}(-D_{|rH})$
and it suffices  to check that this latter map is an isomorphism. This is a local question and we can assume
that $X=\Spec A$ with a  factorial local noetherian ring $A$ 
and $D$ and $H$ are given by elements $d$ and $h$ in $A$,
which have no prime divisors in common (as $H$ intersects $D$ transversally). In this case
the global sections of $\sO_{D+rH}(-D)$ are
\[\Gamma(\sO_{D+rH}(-D))= dA/dh^rA= dA/(dA\cap h^rA)= (dA+h^rA)/ h^rA=\Gamma(\sO_{rH}(-D_{|rH})).\]
This shows that \eqref{cor:Zpn0-zeros2} is an isomorphism. Thus  \eqref{cor:Zpn0-zeros1} yields an exact sequence
\mlnl{H^i(X_{\et}, \Z/p^n(0)_{(X,-D-rH)})\to H^i(X_{\et}, \Z/p^n(0)_{(X,-D)})\to H^i(H_{\et}, \Z/p^n(0)_{(rH, -D_{|rH})})
\\
\to H^{i+1}(X_{\et}, \Z/p^n(0)_{(X,-D-rH)}).}
As the third term only depends on $(rH)_{\red}=H$ (see Lemma \ref{lem:Zpn0-zeros}) it suffices to show
for the injectivity statement that the term on the left hand side vanishes for $r\gg 0$ and $N>i$, 
and for the bijectivity that additionally the term on the right hand side vanishes for $r\gg 0$ and $N>i+1$.
By  Corollary \ref{cor:Zpnq-duality-finite} the statement follows from  the vanishing
\eq{cor:Zpn0-zeros3}{ H^{2N-i+1}(X_{\et}, \Z/p^n(N)_{(X,D+rH)})=0, \quad \text{if } N>i \text{ and } r\gg 0.}
By definition of $\Z/p^n(N)_{(X, D+rH)}$ this last cohomology group fits into an exact sequence
\[H^{N-i}(X, W_n\Omega^N_{(X, D+rH)})\to H^{2N-i+1}(X_{\et}, \Z/p^n(N)_{(X,D+rH)})\to H^{N-i+1}(X, W_n\Omega^N_{(X, D+rH)}).\]
Thus the vanishing  \eqref{cor:Zpn0-zeros3} follows from the Serre-type vanishing Theorem \ref{thm:Serre-van}.
\end{proof}

\begin{rmk}\label{rmk:Lefschetz}
\begin{enumerate}
    \item The Lefschetz type statement that 
    \[\pi_1^{\rm ab}(H, -D_{|H})^p\lra \pi_1^{\rm ab}(X, -D)^p\]
    is surjective, if $N\ge 2$, and bijective, if $N\ge 3$, which is implied by Corollary \ref{cor:Zpn0-zeros},
    of course follows also directly from  the general Lefschetz Theorem for the \'etale fundamental group and
    the sequence \eqref{para:pi1-pole-zero1}. 
    \item In \cite[Theorem 1.1]{KeS-Lefschetz} (see also \cite{KeS-Lefschetz-Err}), the authors prove a
    Lefschetz Theorem for $\pi_1^{\rm ab}(X,D)$. Note that we cannot get a proof of this statement directly from
    the duality statements which we developed here so far as it would require to deal with complexes like
    ''$\Z/p^n(0)_{(X, H-D)}$'', i.e., we  need to allow divisors which are neither effective nor anti-effective.
    Note also that in the Lefschetz statement for $\pi_1^{\rm ab}(X,D)$ from \cite{KeS-Lefschetz} 
    one needs $H$ to be {\em sufficiently} ample, whereas 
    in Corollary \ref{cor:Zpn0-zeros} any effective ample $H$ works. 
\end{enumerate}    
\end{rmk}

\subsection{Weight one}

\begin{para}\label{para:wt1}
We will use the notation 
\eq{para:wt1-0}{\sO_{(X,-D)}^\times:=\Ker(\sO_X^\times\to\sO_D^\times)\quad \text{and} \quad
(\sO^\times/p^n)_{(X,-D)}:=\Ker(\sO^\times_X/p^n\to \sO_{D}^\times/p^n).}
Note that the sheaf on the left is in the literature often denoted by variants of $\sO_{(X,D)}^\times$, 
but in order to have the notation  consistent  with the other parts of the paper, we put here the minus in front of the $D$.
\end{para}

\begin{prop}\label{prop:wt1}
There is a  quasi-isomorphism induced by the $\dlog$-map
\[(\sO^\times/p^n)_{(X,-D)}\xr{\simeq} \Z/p^n(1)_{(X,-D)}[1],\]
which gives rise to a distinguished triangle
\eq{prop:wt1-1}{\sO^\times_{(X,-\lceil D/p^n\rceil)}\xr{F^n} \sO^\times_{{(X,-D)}} 
\lra \Z/p^n(1)_{(X,-D)}[1]\xr{+1}. }
\end{prop}

\begin{proof}
By definition and the surjectivity of the restriction map $\sO_X^\times\to \sO^\times_D$ we have
\[(\sO^\times/p^n)_{(X,-D)}=\frac{\sO^\times_{(X,-D)}}{(\sO^\times_X)^{p^n}\cap \sO_{(X,-D)}^\times}=
\coker(F^n : \sO^\times_{(X,-\lceil D/p^n\rceil)}\to \sO^\times_{(X,-D)}).\]
By \cite[1.4, Lemme 2]{CTSS} together with  \cite[I, Proposition 3.23.2]{IlDRW} we have an exact sequence
\eq{prop:wt1-1.5}{0\to \sO^\times_X/(\sO^\times_X)^{p^n}\xr{\dlog} W_n\Omega^1_X\xr{C^{-1}-1} W_n\Omega^1_X/dV^{n-1}\sO_X.}
In view of Lemma \ref{lem:Zpnq-} it remains to show that the natural map 
\eq{prop:wt1-2}{\dlog: \sO^\times_{(X,-D)}\to W_n\Omega^1_{(X,-D)}\cap W_n\Omega^1_{X,\log},} 
induced by restriction of $\dlog:\sO_{X}^\times\to W_n\Omega^1_{X,\log}$, is surjective.

We do induction on $n$. For $n=1$, this essentially follows from \cite[Theorem 1.2.1]{JSZ}. In fact, 
by \textit{loc. cit.} the map
$$\dlog: \cO_{(X,-D)}^\times \ra \Omega_X^1(\log D)(-D)$$
is surjective. As this map factors via the natural inclusion
$\Omega_X^1(\log E)(-D)\inj \Omega_X^1(\log D)(-D)$, for any reduced effective Cartier divisor $E$, satisfying 
$D_{0,\red}\subset E\subset D_{\red}$, where $D=D_0+pD_1$ is a \pdd, we get the surjectivity 
\eq{eq:dlogsurjzerosE}{\dlog: \cO_{(X,-D)}^\times \surj \Omega_X^1(\log E)(-D)\cap \Omega^1_{X,\log}.
}
As $\Omega_X^1(\log D_0)(-D)=\Omega^1_{(X,-D)}$, by Lemma \ref{lem:n=1zeros}, we in particular obtain 
\eqref{prop:wt1-2} for $n=1$

The induction step is adapted from the proof of \cite[Theorem 2.3.1]{JSZ}.
Denote 
$$M_n:=\mathrm{Image}\left(\dlog:
\sO^\times_{(X,-D)} \ra 
W_n\Omega^1_{(X,-D)}\cap W_n\Omega^1_{X,\log}
\right).$$
The restriction map $R$ clearly induces a surjection $R:M_n\ra M_{n-1}$.
Consider the following commutative diagram
$$\xymatrix{
&&
M_n\ar@{->>}[r]^R\ar@{^(->}[d]&
M_{n-1}\ar@{^(->}[d]
\\
0\ar[r]&
\Omega^1_{X,\log}\ar[r]^{\ul p^{n-1}}&
W_n\Omega^1_{X,\log}\ar[r]^R&
W_{n-1}\Omega^1_{X,\log}\ar[r]&
0.
}$$
The second row is exact by \cite[1.4, Lemme 3]{CTSS}.
Take $x\in W_n\Omega^1_{(X,-D)}\cap W_n\Omega^1_{X,\log}$. Then $R(x)\in W_{n-1}\Omega^1_{(X,-D)}\cap W_{n-1}\Omega^1_{X,\log}$, and hence by induction hypothesis, $R(x)\in M_{n-1}$.
Since $R:M_n\ra M_{n-1}$ is surjective, there exists an element $y\in M_n$ such that $R(y)=R(x)$. 
By the exactness of the second row, $x-y=\ul p^{n-1}(z)$, for some $z\in \Omega^1_{X,\log}$.
In particular, $\ul p^{n-1}(z)\in W_n\Omega^1_{(X,-D)}$ and  thus \Cref{lem:pHWM-zeros} yields
$z\in \Omega^1_{n-1}(-D',-pD_n)$, where $D=D'+p^nD_n$ is a \pdd.
By  the surjectivity of \eqref{eq:dlogsurjzerosE} we find
\[z=\dlog u, \quad \text{for some } u\in \sO^\times_{(X, -\lceil D'/p^{n-1}\rceil -pD_n)}.\]
Hence 
$\ul p^{n-1}(z) =\dlog u^{p^{n-1}} \in M_n$ and therefore $x=y+\ul p^{n-1}(z)\in M_n$.
This finishes the induction step.
\end{proof}

\begin{rmk}
Let $\Z_{(X|D)}(r)$ be the relative motivic cycle complex from \cite{BiSa}. 
By \cite[Theorem 4.1]{BiSa} there is a quasi-isomorphism  of complexes of \'etale sheaves
\[\Z_{(X|D)_{\et}}(1)[1]\simeq \sO_{(X,-D)}^\times.\]
And hence Proposition \ref{prop:wt1} gives rise to a distinguished triangle of complexes
of \'etale sheaves 
\[\Z_{(X|\lceil D/p^n \rceil)_{\et}}(1)\xr{\sV^n} \Z_{(X|D)_{\et}}(1)\to \Z/p^n(1)_{(X,-D)}\xr{+1},\]
where $\sV^n$ is induced by pullback along 
$\id_X\times F^n :X\times_{\F_p} (\P^1_{\F_p})^r\to X\times_{\F_p} (\P^1_{\F_p})^r$, cf. \cite[Lemma 5.7]{KP-dga}.
Thus the notation $\Z/p^n(1)_{(X,-D)}$ might be a bit misleading as we rather have
\[\Z/p^n(1)_{(X,-D)}=\, ``\Z_{(X|D)}(1)/\sV^n \Z_{(X|\lceil D/p^n \rceil)}(1)",\]
which does not agree with $\Z_{(X|D)}(1)/p^n\Z_{(X|D)}(1)$. (It does, if $D=\emptyset$.)
It is interesting to note that the difference between the $p$-adic and the $V$-adic filtraton, which
we have on the de Rham-Witt complex (or just on the Witt vectors of a non-perfect field), seems to be visible
also on the motivic complex with modulus.
\end{rmk}

\begin{para}\label{para:Brauer-pz}
We define the Brauer group of $X$ with ramification bounded by $D$, by
\[\Br(X,D):= H^0(U, R^2\e_* \Q/\Z(1)'_{U})\oplus H^0(X, R^2\e_*\Q_p/\Z_p(1)_{(X,D)}),\]
where $U=X\setminus D$, $\e:X_\et\to X_{\Nis}$ is the change of sites map, and 
\[\Q/\Z(1)'_{U}= \colim_{n'}\, \mu_{n',U}, \quad  
\Q_p/\Z_p(1)_{(X,D)}=\colim_{\ul{p}}\, \Z/p^n(1)_{(X,D)}\]
with the colimit on the left over all $n'$ which are prime to $p$.

We define the Brauer group with zeros along $D$ by 
\[\Br(X,-D):= H^2(X_\et, \sO_{(X,-D)}^\times).\]

We remark:
\begin{enumerate}
\item There are natural inclusions $\Br(X)\subset \Br(X,D)\subset \Br(U)$ and by \ref{para:Zpnq}\ref{para:Zpnq3}
\[\colim_r \Br(X,rD)=\Br(U).\]
\item We have an exact sequence
\eq{para:Brauer-pz1}{\Pic(X)\to \Pic(D)\to \Br(X,-D)\to \Br(X)\to \Br(D),}
where for a singular scheme $Z$ we use the definition $\Br(Z)=H^2(Z_{\et},\G_m)$.
\item Assume $p^n$ is bigger than all the multiplicities of $D$. 
Then the $n$-power Frobenius induces by \eqref{prop:wt1-1} induces a morphism $F^n$ fitting 
in the following commutative diagram
\eq{para:Brauer-pz3}{\xymatrix{
\Br(X,-D)\ar[rr]^{p^n}\ar[rd] &  &\Br(X,-D)\\
 & \Br(X, -D_{\red}).\ar[ru]_{F^n}
}}
\end{enumerate}
\end{para}

\begin{remark}\label{rmk:Brauer-pz}
As is well-known $\Br(X)$ and $\Br(U)$ are torsion groups and hence so is $\Br(X,D)$.
But note that $\Br(X,-D)$ is in general not a torsion group as becomes apparent from the exact sequence 
\eqref{para:Brauer-pz1}. For example if $X=\P^2$ and $D=E$ is an elliptic curve we have 
$E(k)\subset \Br(\P^2, -E)$. Even if $k$ is a finite field the subgroup $\coker(\Pic(X)\to \Pic(D))$ might not be torsion due 
to a contribution by $\Pic_{D/k}(k)/\Pic^0_{D/k}(k)$. 

Now assume that $k$ is a finite field, $p\neq 2$, and $X$ is a smooth projective surface.
In this case the Tate conjecture for divisors is equivalent to the finiteness of $\Br(X)$
which is also equivalent to the finiteness of its the $p$-primary torsion subgroup
$\Br(X)[p^\infty]$, see \cite[Theorem 4.1]{Milne-ConjArtin}.
We note:
\begin{enumerate}[label=(\arabic*)]
    \item\label{rmk:Brauer-pz1} $\Br(X)$ is finite if and only if $\Br(X,-D)[p^\infty]$ is finite.
    \item\label{rmk:Brauer-pz2} 
    If the reduced irreducible components of $D$ are linearly independent in ${\rm Num}(X)\otimes \Q$
    (divisors modulo numerical equivalence) and $\Br(X)$ is finite, then $\Br(X,-D)$ is finite as well.
\end{enumerate}
Indeed, first note that $\Br(D)=0$ as follows from the Brauer-Hasse-Noether Theorem, 
see also \cite[Remarques (2.5), b)]{Grothendieck-BrauerIII}, together with \cite[3.6.6]{CTS}. 
Thus the finiteness of $\Br(X,-D)[p^\infty]$ implies the finiteness of $\Br(X)$. 
On the other hand, $\Pic_{D/k}^0$ is a smooth algebraic group scheme (e.g. \cite[\S 9.2]{BLR}) 
and hence $\Pic_{D/k}^0(k)$ is a finite group. Moreover
$(\Pic_{D/k}/\Pic^0_{D/k})(k)$ is a finitely generated group and hence so is 
$\Pic(D)/\Pic_{D/k}^0(k)$. It follows that  $\Im(\Pic(D)\to \Br(X,-D))\cap \Br(X,-D)[p^\infty]$ is finite, which yields the
``if" direction in \ref{rmk:Brauer-pz1}. For \ref{rmk:Brauer-pz2} it remains to observe that under the assumptions made there
the index of $\Pic(X)\to \Pic(D)/\Pic_{D/k}^0(k)$ is finite.
\end{remark}

For $D=\emptyset$ the following duality statement is the $p$-part of Milne's duality for the Brauer group of a smooth 
projective surface over a finite field, see \cite[Theorem 2.4]{Milne-ConjArtin}.

\begin{thm}\label{thm:duality-Br-surf}
Assume $k$  is a finite field and  $X$ is a smooth proper surface (i.e. $N=2$). 
Then there is a canonical isomorphism of profinite groups
\[ \frac{\Br(X,-D)[p^\infty]}{(\Br(X, -D)[p^\infty])_{\rm div}}\xr{\simeq} 
\Hom\left(\frac{\Br(X,D)[p^\infty]}{(\Br(X,D)[p^\infty])_{\rm div}}, \Q/\Z\right),\]
where the index ``div" refers  to the divisible part, and $M[p^\infty]$ denotes the $p$-primary torsion in $M$.
\end{thm}
\begin{proof}
The proof follows the same strategy as the one in \cite[Theorem 2.4]{Milne-ConjArtin}.
First of all note that the right hand side of the isomorphism in the statement is a profinite group by
Pontryagin duality, hence it suffices to prove the isomorphism of abstract groups.
We start by showing that we have an exact sequence
\eq{thm:duality-Br-surf1}{ 0\to \Pic(U)/p^n \to H^2(X_{\et}, \Z/p^n(1)_{(X,D)})\to H^0(X,R^2\e_*(\Z/p^n(1)_{(X,D)}))\to 0.}
To this end we first claim that we have 
\[R^1\e_*(\Z/p^n(1)_{(X,D)})=j_*(\sO_{U_\Nis}^\times/p^n)= Rj_*(\sO_{U_\Nis}^\times/p^n).\]
The first equality follows directly from the definition, see also \eqref{rmk:Zpnq+1}, and 
the exact sequence \eqref{prop:wt1-1.5}. 
The second equality follows from the vanishing $R^i j_*\sO^\times_U=0$, $i\ge 1$, which relies on the smoothness of $X$.
Therefore we get a distinguished triangle
\[R j_*(\sO^\times_{U_\Nis}/p^n)[-1]\to R\e_*(\Z/p^n(1)_{(X,D)})\to R^2\e_*(\Z/p^n(1)_{(X,D)})[-2]\xr{+1}.\]
As $H^r(U_{\Nis}, \sO_U^\times)=0$, for $r\ge 2$, we get the exact sequence \eqref{thm:duality-Br-surf1}.
Taking the colimit over $n$ of \eqref{thm:duality-Br-surf1} gives an exact sequence
\eq{thm:duality-Br-surf2}{0\to \Pic(U)\otimes_{\Z} \Q_p/\Z_p\to H^2(X_{\et}, \Q_p/\Z_p(1))\to \Br(X,D)[p^\infty]\to 0.}

Now we consider the negative divisor.
The distinguished triangle \eqref{prop:wt1-1} yields an exact sequence
\ml{thm:duality-Br-surf3}{0\to \frac{\Br(X,-D)}{F^n\Br(X, -\lceil D/p^n\rceil)}\to H^3(X_{\et}, \Z/p^n(1)_{(X,-D)})\\
\to \Ker (F^n: H^3(X_{\et},\sO^\times_{(X,-\lceil D/p^n\rceil)})\to H^3(X_{\et},\sO^\times_{(X,-D)}))\to 0.}
Assume $n$ is large enough so that $\lceil D/p^n\rceil=D_\red$. 
Note that the quotient sheaf $\sQ=\sO_{(X,-D_\red)}^\times/ \sO^\times_{(X,-D)}$ 
is a successive extension of coherent $\sO_D$-modules.
As $D$ is a curve we get $H^i(X_{\et}, \sQ)=0$, $i\ge 2$, and hence an isomorphism
\[H^3(X_{\et}, \sO^\times_{(X,-D)})\xr{\simeq} H^3(X_{\et}, \sO^\times_{(X,-D_{\red})}).\]
Thus the term on the right  of \eqref{thm:duality-Br-surf3} is equal to the $p^n$-torsion 
in $H^3(X_{\et}, \sO^\times_{(X,-D_\red)})$.
Furthermore for fixed $n_0$ with $\lceil D/p^{n_0} \rceil=D_{\red}$ we have by \eqref{para:Brauer-pz3} the following inclusions
\[p^n\Br(X,-D)\subset F^n(\Br(X,-D_{\red}))\subset F^{n_0}(p^{n-n_0}\Br(X,-D_\red))
\subset p^{n-n_0}\Br(X,-D), \]
for $n\ge n_0$.
Thus $\{F^n\Br(X,-D_{\red})\}_{n\ge n_0}$ and $\{p^n\Br(X,-D)\}_n$ define the same topology on $\Br(X,-D)$. 
As \eqref{thm:duality-Br-surf3} is an exact sequence of finite groups 
taking the limit over $n$ yields  an exact sequence
\[0\to \lim_n \frac{\Br(X,-D)}{p^n\Br(X, -D)}\to \lim_n H^3(X_{\et}, \Z/p^n(1)_{(X,-D)})\to 
T_pH^3(X_{\et}, \sO^\times_{(X,-D_\red)})\to 
0, \]
where $T_pM$ denotes the Tate module of $M$. 
As $T_pH^3(X_{\et}, \sO^\times_{(X,-D_\red)})$ is torsion-free and 
$\Pic(U)\otimes_{\Z} \Q/\Z$ is divisible  the statement follows from Corollary \ref{cor:Zpnq-duality-finite}
by the same argument as in \cite[Theorem 2.4]{Milne-ConjArtin}.
\end{proof}

\begin{rmk}\label{rmk:BrA}
Tate has shown that if $X$ is an abelian surface over a finite field, 
then $\Br(X)[\ell^\infty]$ ($\ell\neq p$) is finite, 
see \cite[Corollary]{Tate}. By \cite[Theorem 4.1]{Milne-ConjArtin} 
we also have the finiteness  of $\Br(X)[p^\infty]$.
Hence by Remark \ref{rmk:Brauer-pz} and Theorem \ref{thm:duality-Br-surf} we find that 
$\Br(X, D)[p^\infty]$ is the direct sum of a divisible group with a finite group, 
for all effective Cartier divisors
$D$ on $X$ with $D_{\red}$ a simple normal crossing divisor.    
\end{rmk}


\end{document}